\documentclass[a4paper,10pt]{article}
\pdfoutput=1
\usepackage{graphics}
\usepackage{graphicx}
\usepackage{float}

\usepackage[english]{babel}
\usepackage{latexsym}
\usepackage{amssymb}
\usepackage{amsmath,amsthm}
\usepackage{amsfonts}

\usepackage[utf8]{inputenc}
\usepackage[T1]{fontenc}
\usepackage{tikz,tikz-3dplot}
\usetikzlibrary{fadings}
\tikzfading[name=fade out,
            inner color=transparent!0,
            outer color=transparent!100]
\tikzfading[name=fade right,
            left color=transparent!0,
            right color=transparent!100]

\usepackage{bm}
\usepackage{eso-pic}
\usepackage{color} 
\usepackage{xcolor-material}

\usepackage{booktabs}
\usepackage{caption}
\usepackage[a4paper,top=2cm,bottom=2cm,left=2.6cm,right=2.6cm,bindingoffset=5mm]{geometry}

\theoremstyle{plain}
\newtheorem{thm}{Theorem}[section]
\newtheorem{lem}{Lemma}[section]
\newtheorem{coro}{Corollary}[section]
\newtheorem{propos}{Proposition}[section]

\theoremstyle{definition}
\newtheorem{defn}{Definition}[section]
\newtheorem{exe}{Example}[section]
\newtheorem{rmq}{Remark}[section]

\newcommand{\pf}{{\bf Proof.~}} 
\newcommand{\thmref}[1]{Theorem~\ref{#1}}
\newcommand{\lemref}[1]{Lemma~\ref{#1}}
\newcommand{\cororef}[1]{Corollary~\ref{#1}}

\newcommand{\defnref}[1]{Definition~\ref{#1}}
\newcommand{\rmqref}[1]{Remark~\ref{#1}}
\newcommand{\exeref}[1]{Example~\ref{#1}}

\def\scat#1#2%
 {\hbox{\vrule\vbox to#2pc{\hrule\hbox to#1pc{\hfil} 
  \vfil\hrule}\vrule}}

\def\findimlogo
 {\scat{.25}{.5}}

\def\findim
 {\hfill\findimlogo\smallskip\par}

\def\ellerovesciata#1#2%
 {\vbox{\hbox to#1pc{\vbox to#2pc{\vfil} 
  \hfil\vrule}\hrule}}

\def\ctr
 {\ellerovesciata{.5}{.5}}

\def\Supp
 {\mathop{\rm Supp\,}\nolimits}

\def\lin
 {\mathop{\rm lin}\nolimits}

\def\relint
 {\mathop{\rm relint\,}\nolimits}
 
\def\relbd
 {\mathop{\rm relbd\,}\nolimits}

\def\conv
 {\mathop{\rm conv\,}\nolimits}

\def\rk
 {\mathop{\rm rank\,}\nolimits}

\def\tr
 {\mathop{\rm tr\,}\nolimits}
 
\def\V
 {\mathop{\rm V\,}\nolimits}

\def\re
 {\mathop{\rm Re\,}\nolimits}

\def\im
 {\mathop{\rm Im\,}\nolimits}

\def\aff
 {\mathop{\rm aff}\nolimits}

\def\vol
 {\mathop{\rm vol\,}\nolimits}

\def\Hess
 {\mathop{\rm Hess}\nolimits}

\def\Subd
 {\mathop{\rm Subd}\nolimits}

\def\Arg
 {\mathop{\rm Arg}\nolimits}

 \def\sec
 {\mathop{\textsl{II}}\nolimits}

 \def\Bigcirc
 {\mathop{\raisebox{-0.5mm}{\scalebox{1.7}{$\bigcirc$}}}}

\newcommand{\de}{\mathrm d}
 
\newcommand{\dc}{{\mathrm d}^{\mathrm c}}
 
\newcommand{\ddc}{{\mathrm {dd}}^{\mathrm c}}

\def\hfl#1#2{\smash{\mathop{\hbox to 12mm{\rightarrowfill}}
 \limits^{\scriptstyle#1}_{\scriptstyle #2}}}

\def\Hfl#1#2{\smash{\mathop{\hbox to 15mm{\rightarrowfill}}
 \limits^{\scriptstyle#1}_{\scriptstyle #2}}}

\def\deq{\mathrel{\hbox{\raise .09ex\hbox{$:$}$=$}}}

\def\eqd{\mathrel{\hbox{$=$\raise .09ex\hbox{$:$}}}}

\begin{document}
\title{\Large \bf Kazarnovski{\v\i} mixed pseudovolume.
}
\author{\scshape{James Silipo}
}
\maketitle
\begin{abstract}
The present paper is a detailed and almost self-contained introduction to Kazarnovski{\v{\i}} mixed pseudovolume. It provides new formulas and new results about Alexandroff-Fenchel type inequalities for Kazarnovski{\v{\i}} mixed pseudovolume.
\end{abstract}

\tableofcontents

%

\newpage

\section{Introduction}

The present article provides a detailed introduction to \emph{Kazarnovski{\v\i} mixed pseudovolume} as a generalization to $\mathbb C^n$ of Minkowski mixed volume in $\mathbb R^n$. 

Given a convex body (i.e. a compact convex subset) $A\subset\mathbb C^n$, the support function of $A$ is the convex, positively $1$-homogeneous function $h_A$ defined, for $z\in \mathbb C^n$, by $h_A(z)=\max_{u\in A}\re\langle z,u\rangle$, where $\langle\,,\rangle$ is the standard hermitian product on $\mathbb C^n$. Suppose $h_A\in{\mathcal C}^2(\mathbb C^n\setminus\{0\})$ and let $\dc=i(\bar\partial-\partial)$, let also $B_{2n}$ denote the unit full-dimensional ball about the origin of $\mathbb C^n$ and $\varkappa_{2n}$ its Lebesgue measure, then the \emph{Kazarnovski{\v\i} pseudovolume} of $A$ is the non-negative real number $P_n(A)$ defined as 
\begin{equation}\label{a}
P_n(A)=\dfrac{1}{n!\varkappa_n}\int_{B_{2n}}(\ddc h_A)^{\wedge n}\,.
\end{equation}
Accordingly, if $A_1,\ldots,A_n$ are convex bodies whose support functions belong to ${\mathcal C}^2(\mathbb C^n\setminus\{0\})$, their \emph{mixed pseudovolume} is the non-negative real number $Q_n(A_1,\ldots,A_n)$ defined as
\begin{equation}\label{b}
Q_n(A_1,\ldots,A_n)=\dfrac{1}{n!\varkappa_n}\int_{B_{2n}}\ddc h_{A_1}\wedge\ldots\wedge\ddc h_{A_n}\,.
\end{equation}
The mixed version polarizes the unmixed one, i.e. $Q_n(A,\ldots,A)=P_n(A)$, so $Q_n$ is multilinear with respect to positive scalar multiplication and Minkowski addition. 

The notion of pseudovolume can be defined without the regularity assumption on support functions, it requires a regularization argument on positive currents that, for a polytope $\Gamma$, yields a remarkable combinatorial formula:
\begin{equation}\label{c}
P_n(\Gamma)
\quad=
\sum_{\Delta\in{\mathcal B_{\rm ed}}(\Gamma,n)}
\varrho(\Delta)\vol_n(\Delta)\psi_\Gamma(\Delta)\,,
\end{equation}
where ${\mathcal B_{\rm ed}}(\Gamma,n)$ is the set of \emph{equidimensional} $n$-faces of $\Gamma$, i.e. those faces $\Delta$ of $\Gamma$ spanning a real affine subspace $\aff_\mathbb R \Delta$ and a complex one $\aff_\mathbb C \Delta$ which satisfy the dimensional equalities $\dim_\mathbb R\aff_\mathbb R \Delta=\dim_\mathbb C\aff_\mathbb C \Delta=n$; $\varrho(\Delta)$ is a weight ascribed to every equidimensional face and equal to the jacobian determinant of a projection (cf. section~\ref{distorzione}), $\vol_n(\Delta)$ is the $n$-dimensional volume of $\Delta$ and $\psi_\Gamma(\Delta)$ is the outer angle of $\Gamma$ at $\Delta$ (cf. section~\ref{preliminari}).

For convex bodies included in $\mathbb R^n$, Kazarnovski{\v\i} mixed pseudovolume reduces to Minkowski mixed volume, which yields an integral formula for the latter similar to that involving superforms obtained by Larson~\cite{Lar}. 

Kazarnovski{\v\i} psedovolume can be even defined for a non convex compact set $A\subset\mathbb C^n$ with interior points and sufficiently regular boundary as
\begin{equation}\label{d}
P_n(A)=\dfrac{1}{n!\varkappa_n}
\int_{\partial A}
\alpha_{\partial A}\wedge(\de \alpha_{\partial A})^{\wedge (n-1)}\,,
\end{equation}
where $\alpha_{\partial A}$ is the differential 1-form whose value on the point $\zeta\in\partial A$ equals 
the projection on the characteristic direction of the tangent space $T_\zeta \partial A$ (cf. section~\ref{pS}), or equivalently as 
\begin{equation}\label{e}
P_n(A)=\dfrac{2^{n-1}}{n\varkappa_n}
\int_{\partial A}
{\mathcal E}_{\partial A}\upsilon_{\partial A}\,,
\end{equation}
where ${\mathcal E}_{\partial A}$ denotes the determinant of the Levi form of $\partial A$ and $\upsilon_{\partial A}$ is the standard volume form on the hypersurface $\partial A$.

The resulting object is very interesting: as shown by Alesker~\cite{Ale}, it yields a valuation on convex bodies, i.e. 
\begin{equation}\label{f}
P_n(A_1\cup A_2)+P_n(A_1\cap A_n)=P_n(A_1)+P_n(A_2)\,,
\end{equation}
whenever $A_1,A_2$ and $A_1\cup A_2$ are convex. This valuation has some remarkable properties, namely: continuity with respect to Hausdorff metric and invariance for translations or unitary transformations. $P_n$ is not monotonically increasing unless $n=1$, so it cannot be extended to a measure on $\mathbb C^n$. However, just like ordinary volume, the pseudovolume of a polytope is positive if and only if the polytope spans an $n$-dimensional complex affine subspace. 

\subsection{Statement of main results}
The paper contains the proofs of well known and new results about Kazarnovski{\v\i} mixed pseudovolume. The proofs of well known ones have been included because they are hardly available elsewhere. Among the new results, some are easy but very useful.

\begin{itemize}
\item Well known results:
\begin{itemize} 
\item the equivalence of the above mentioned definitions of pseudovolume (cf. \thmref{thmlevi}, \cororef{kaza}, \thmref{azzok1}, \thmref{azzok2}, \rmqref{vcorrk});
\item the vanishing condition for mixed pseudovolume of polytopes (cf. \cororef{ndeg});
\item the non monotonicity of $P_n$ for $n>1$ (cf. \exeref{nonmon});
\item the valuation property, for every integer $0\leqslant k\leqslant n$, of the mapping $A\mapsto (\ddc h_A)^k$ on the space ${\mathcal K}(\mathbb C^n)$ of compact convex subsets of $\mathbb C^n$ with values in that of positive currents (cf. \cororef{valk});
\end{itemize}
\item New results:
\begin{enumerate}
\item For every $A_1,\ldots,A_k\in{\mathcal K}(\mathbb C^n)$, (cf. \cororef{strano1})
\begin{align}
Q_n(A_1,\ldots,A_k,\underbrace{B_{2n},\ldots,B_{2n}}_{(n-k) \text{ times}})
=
\dfrac{1}{2^{n-k}n!\varkappa_n}
\int_{B_{2n}}
\left(
\bigwedge_{\ell=1}^k\ddc h_{A_\ell}
\right)
\wedge
(\ddc h^2_{B_{2n}})^{\wedge(n-k)}\,;
\end{align}

\item Alexandroff-Fenchel inequality for Kazarnovski{\v\i} mixed pseudovolume is generally false (cf. \rmqref{nokaf}), however the following weaker version holds true:
\begin{equation}
Q_n(A,\underbrace{B_{2n},\ldots,B_{2n}}_{(n-1) \text{ times}})^2\geq \left(\dfrac{2n-2}{2n-1}\right)Q_n(A,A,\underbrace{B_{2n},\ldots,B_{2n}}_{(n-2) \text{ times}})P_n(B_{2n})\,,
\end{equation}
for every $A\in{\mathcal K}(\mathbb C^n)$;
\item for every integer $0\leqslant k\leqslant n$ and every $\Gamma_1,\ldots,\Gamma_k\in{\mathcal P}(\mathbb C^n)$\begin{equation}
\ddc h_{\Gamma_1}\wedge\ldots\wedge\ddc h_{\Gamma_k}=k!\sum_{\Delta\in{\mathcal B}_{\rm ed}(\Gamma,k)}
\varrho(\Delta)V_k(\Delta_1,\ldots,\Delta_k)[K_\Delta]\wedge \upsilon _{E_\Delta^\prime}\,,
\end{equation}
where $\Gamma=\sum_{\ell=1}^k\Gamma_\ell$, and for every $\Delta\in{\mathcal B}_{\rm ed}(\Gamma,k)$, $\Delta=\sum_{\ell=1}^k\Delta_\ell$ is the unique sum decomposition, $V_k(\Delta_1,\ldots,\Delta_k)$ is the $k$-dimensional Minkowski mixed volume on the linear space $E_\Delta$ parallel to $\Delta$, $[K_\Delta]$ is the integration current on the dual cone to $\Delta$ and $\upsilon_{E_\Delta^\prime}$ is a suitably chosen volume form on $E_\Delta^\prime=E_\Delta^\perp\cap(E_\Delta+ i E_\Delta)$;
\item if $A_1,\ldots,A_n\in{\mathcal K}(\mathbb C^n)$ have ${\mathcal C}^2(\mathbb C^n\setminus\{0\})$ support functions with complex hessian matrices $\Hess_\mathbb C h_{A_1},\ldots,\Hess_\mathbb C h_{A_n}$, then (cf. (\ref{kmix-D-f}))
\begin{equation}
Q_n(A_1,\ldots,A_n)=\dfrac{4^n}{\varkappa_n}\int_{B_{2n}} D_n(\Hess_\mathbb C h_{A_1},\ldots,\Hess_\mathbb C h_{A_n})\,\upsilon_{2n}\,,
\end{equation}
where $D_n$ is the mixed discriminant and $\upsilon_{2n}$ the standard volume form on $\mathbb C^n$;
\item if $A_1,\ldots,A_n\in{\mathcal K}(\mathbb C^n)$ are such that $h_{A_1}\in{\mathcal C}^1(\mathbb C^n\setminus\{0\})$ and $h_{A_\ell}\in{\mathcal C}^2(\mathbb C^n\setminus\{0\})$ for $2\leqslant \ell\leqslant n$, then (cf. \lemref{dj})
\begin{equation}
Q_n(A_1,A_2,\ldots,A_n)
=
\dfrac{4^{n-1}}{\varkappa_n}
\int_{\partial B_{2n}}
\iota^*_{\partial B_{2n}}
D_n(N_{A_1},\Hess_\mathbb C h_{A_2},\ldots,\Hess_\mathbb C h_{A_n})\upsilon_{\partial B_{2n}}
\,,
\end{equation}
where $i_{\partial B_{2n}}:\partial B_{2n}\to\mathbb C^n$ is the inclusion, $\upsilon_{\partial B_{2n}}$ is the standard volume form on $\partial B_{2n}$ and $N_{A_1}$ is the matrix defined by (\ref{matriceN});
\item if $A_1,\ldots,A_n\in{\mathcal K}(\mathbb C^n)$ have ${\mathcal C}^2(\mathbb C^n\setminus\{0\})$ support functions, then (cf. \lemref{djj})
\begin{equation}
Q_n(A_1,A_2,\ldots,A_n)
=
\dfrac{4^{n-1}}{\varkappa_n}
\int_{\partial B_{2n}}
\iota^*_{\partial B_{2n}}
D_n(Y_{A_1},\Hess_\mathbb C h_{A_2},\ldots,\Hess_\mathbb C h_{A_n})\upsilon_{\partial B_{2n}}
\,,
\end{equation}
where $Y_{A_1}$ is the matrix defined by (\ref{matriceY});
\item if $A_1,\ldots,A_n\in{\mathcal K}(\mathbb C^n)$,  $h_{A_1}\in{\mathcal C}^1(\mathbb R^n\setminus\{0\})$, $h_{A_\ell}\in{\mathcal C}^2(\mathbb R^n\setminus\{0\})$, for $2\leqslant \ell\leqslant n$, and $\Xi(A_2,\ldots,A_n)(u)$ is the matrix defined on page \pageref{matriceXi}, then (cf. \cororef{corobuono}) the condition $\re\langle\nabla h_{A_1}(u),\Xi(A_2,\ldots,A_n)(u)\, u\rangle=0$ for almost all $u\in\partial B_{2n}$, implies 
\begin{equation}
Q_n(A_1,\ldots,A_n)
=
\dfrac{4^{n-1}}{\varkappa_n}
\int_{\partial B_{2n}}
h_{A_1}
D_n(I_n,\Hess_\mathbb C h_{A_2},\ldots,\Hess_\mathbb C h_{A_n})
\upsilon_{\partial B_{2n}}\,;
\end{equation}
\item if $A_1,\ldots,A_n\subset\mathbb R^n$ are convex bodies such that $h_{A_1}\in{\mathcal C}^1(\mathbb R^n\setminus\{0\})$ and $h_{A_\ell}\in{\mathcal C}^2(\mathbb R^n\setminus\{0\})$ for $2\leqslant \ell\leqslant n$, then (cf. (\ref{nuova}))
\begin{equation}
V_n(A_1,\ldots,A_n)
=
\int_{\partial B_{2n}}
h_{A_1}
D_n(I_n,\Hess_\mathbb R h_{A_2},\ldots,\Hess_\mathbb R h_{A_n})\upsilon_{\partial B_{2n}},
\end{equation}
where $V_n$ is the $n$-dimensional Minkowski mixed volume,  $I_n$ is the $n$-by-$n$ identity matrix and $\Hess_\mathbb R h_{A_2},\ldots,\Hess_\mathbb R h_{A_n}$ are the real hessian matrices of the respective support functions;
\item let $A\in{\mathcal K}(\mathbb C^n)$ span a real affine subspace $\aff_\mathbb R A$ of dimension $2n-1$ and let let $E_A$ be the $\mathbb R$-linear subspace parallel to $\aff_\mathbb R A$. If there exists a real hyperplane $H$ such that $E_A\cap H=E_A\cap iE_A$ then (cf. \lemref{luglio2019})
$$
P_n(A)\geqslant\max\{P_n(A\cap H^+), P_n(A\cap H^-)\}.
$$
\end{enumerate}
\end{itemize}

\subsection{Origin of the notion}
Although this notion is interesting in integral, tropical and complex analytic geometry, the available literature on the subject is still relatively scarce.  
The theory of mixed pseudovolume has been inspired by the so-called BKK theorem about the number of solutions to a generic system $P\subset\mathbb C[z_1^\pm,\ldots,z_n^\pm]$ of $n$ Laurent polynomials due to Bernstein~\cite{Be}, Kushnirenko~\cite{Ku} and Khovanski~\cite{Kho1}. Recall that, if  $p(z)=\sum_{\lambda\in\Lambda_p}c_\lambda z^\lambda=\sum_{\lambda\in\Lambda_p}c_\lambda z_1^{\lambda_1}\cdots z_n^{\lambda_n}$ is a Laurent polynomial, the \emph{Newton polytope} $\Delta_p$ of $p$ is the convex hull  in $\mathbb R^n$ of its (finite) support of summation $\Lambda_p\subset\mathbb Z^n$. According to the BKK theorem, the number of solutions in the torus $(\mathbb C^*)^n$ to a generic system of $n$ Laurent polynomials $P$ equals the Minkowski mixed volume of the $n$ polytopes $\Delta_p$, as $p\in P$. 

Kazarnovski{\v\i}~\cite{Ka1}, by investigating the asymptotic distribution of the solutions to a finite generic system $F$ of exponential sums, (i.e. entire functions of the form $f(z)=\sum_{\lambda\in\Lambda_f}c_\lambda e^{\langle \cdot,\lambda\rangle}$, where $\Lambda_f\subset\mathbb C^n$ is a finite subset called the \emph{spectrum} of $f$, $c_\lambda\in\mathbb C^*$ and $\langle\,,\rangle$ is the standard hermitian product of $\mathbb C^n$), found that such an asymptotic distribution is controlled by the geometry of the convex hulls $\Delta_f\subset\mathbb C^n$ of the spectra $\Lambda_f$, $f\in F$. In this exponential setting, the mixed pseudovolume plays the role of mixed volume in the polynomial one and reduces to it when the spectra are included in $\mathbb R^n$, in particular the BKK theorem appears as a special instance of Kazarnovsk{\v\i}i's theorem when the spectra are included in $\mathbb Z^n$.

The link between the asymptotic distribution of zeros and the geometry of the exponents comes from King's formula for the integration current associated to a complete intersection of analytic hypersurfaces. According to~\cite{GK}, such a current is nothing but the wedge product of the integration currents associated to the involved hypersurfaces. By Lelong-Poincare formula~\cite{Le}, the integration current of a single analytic hypersurface $\{f=0\}$ equals, up to a normalization constant, the current $\ddc \log \vert f\vert$. Loosely speaking, the exponential shape of $f$ and the presence of a logarithm in Lelong-Poincaré formula explain why the asymptotic distribution of the zeros of the system is controlled by the geometry of its exponents. More precisely, let $\#F=n$, $F=\{f_1,\ldots,f_n\}$ and, for $1\leq \ell\leq n$, let $\Lambda_\ell=\Lambda_{f_\ell}$, $\Delta_\ell=\Delta_{f_\ell}$ and $h_\ell=h_{\Delta_\ell}$. If the zero set $V(F)=\{z\in\mathbb C^n\mid f_\ell(z)=0,\;1\leq\ell\leq n\}$ is a complete intersection, then the current $\bigwedge_{\ell=1}^n t^{-1}\ddc \log \vert f_\ell(tz)\vert$ weakly converges to the current $\bigwedge_{\ell=1}^n \ddc h_\ell(z)$, as $t\to +\infty$, and it turns out (non trivially) that the mixed pseudovolume $Q_n(\Delta_1,\ldots,\Delta_n)$ affects the asymptotic distribution of $V(F)$. 
The requirement that $V(F)$ is a complete intersection is crucial for the preceding method to work. If one fixes both $\# F$ and the number of points in the spectra, this requirement depends on both the coefficients and the exponents of the exponential sums, it is generically satisfied and can be seen as a sort of \emph{non-singularity} of $V(F)$. The non-singularity conditions assumed in~\cite{Ka1}, though generalizing those considered in the polynomial case, are merely sufficient conditions which are far from being necessary. 

Though mixed pseudovolume is just a tool in the theory of exponential sums, it seems to deserve at least the same attention as the theory it originated from. 
The notion of pseudovolume first appeared in the 1981 note~\cite{Ka1} and, as early as its introduction, it was known to the Russian mathematical community as witnessed by Shabat, who mentioned pseudovolume in~\cite{Sha} p.~199.
In 1984 Kazarnovski{\v\i} gave a more detailed construction of pseudovolume in a second paper~\cite{Ka2}, however this second article does not provide proofs for all the theorems stated in the first one. The subsequent papers~\cite{Ka3,Ka4,Ka5,Ka6,Ka7,Ka8,Ka9} deal with interesting analytic and combinatorial investigations in which, however, pseudovolume does not play a prominent role. 
Since 1984 Kazarnovski{\v\i} pseudovolume seemed to have been forgotten until 2003, when the notion got new attention in integral geometry~\cite{Ale}. Then the papers~\cite{Ka10,Ka11} by Kazarnovski{\v\i} briefly provided further details about pseudovolume and its vanishing condition in the polytopal case. In his most recent paper~\cite{Ka12}, Kazarnovski{\v\i} builds an intersection theory for exponential sums in which the intersection index of $n$ exponential hypersurfaces is equal to the mixed pseudovolume of their Newton polyhedra.

In 2001 Yger drew my attention towards the geometry of exponential sums and I started my PhD under his supervision by conducting a thorough study of Kazarnovski{\v\i}'s seminal paper \cite{Ka1} and looking for fully detailed proofs of the main results included therein. The present paper is a substantial expansion of the first chapter of my thesis~\cite{Sil}, its intentionally elementary style aims at making the reading as accessible as possible, the exposition is sometimes redundant: some results, especially in low dimension, can be deduced by more general ones, nevertheless I inserted them to ease the understanding or to emphasize the use of particular methods which are valid just in special situations. I also tried to give a self-contained presentation with the following two notable exceptions: several complex variables and the theory of currents, at least at an introductory level. The reader should be familiar with these topics and is possibly referred to~\cite{Dem},~\cite{Ran} or \cite{LeGr}.

\subsection{Structure of the paper}
The paper is organized as follows. Section~\ref{preliminari} collects some notions and facts about convexity which will be used in the sequel. Section~\ref{folume} presents the construction of a class of valuations on convex bodies of $\mathbb R^n$ generalizing ordinary intrinsic volumes and mixed volume. This construction, briefly outlined in~\cite{Ka1} and further studied in~\cite{Wan}, admits Kazarnovski{\v\i} pseudovolume as a special instance. Section~\ref{distorzione} is devoted to the description of the weight function $\varrho$ involved in~\eqref{c}.
This section does not go beyond elementary linear algebra, however it provides a useful machinery used throughout the article. Section~\ref{corpi-complessi} describes some typical features of convex bodies of $\mathbb C^n$ thanks to the information from the preceding two sections. Section~\ref{pS} deals with the notion of pseudovolume on the class of compact subsets of $\mathbb C^n$ with non empty interior and sufficiently regular boundary, it contains the proof of the equality of \eqref{d} and~\eqref{e} as well as that of~\eqref{a} and~\eqref{d} for strictly convex bodies with sufficiently regular boundary. 
Section~\ref{ksez} treats the convex case in full generality, it includes a proof of~\eqref{c} (following a suggestion by Kazarnovski{\v\i} himself) that does not seem to have appeared yet.
Section~\ref{chk} is dedicated to a closer study of the current $(\ddc h_A)^{\wedge k}$ in the polytopal case. Such a current will be represented as a linear combination of some special currents obtained as products of an integration current and a volume form (cf. \eqref{lambdak}). A similar representation is briefly described in~\cite{Ka5} and~\cite{Ka9}. By virtue of such representation it is shown that $(\ddc h_A)^{\wedge k}$ enjoys a valuation property  (already mentioned in~\cite{Ka6} and ~\cite{Ka7}) implying that revealed by Alesker for $P_n$. 
Section~\ref{kmd} clarifies the relations between mixed pseudovolume and mixed discriminants. 
Section~\ref{esempi} provides several examples, including those showing that $P_n$ is neither rational on lattice polytopes nor orthogonally invariant. 
Section~\ref{u} collects some useful notions, facts and computations only occasionally used in the paper. 
Finally section~\ref{open} presents some open problems and partial answers about non vanishing, monotonicity and Alexandroff-Fenchel type inequality for $Q_n$. The article ends with a list of non standard notation employed throughout the paper.

\section{Preliminaries on convex bodies in $\mathbb R^n$}\label{preliminari}

Let $n\in\mathbb N^*$ and give $\mathbb R^n$ the euclidean structure induced by the standard scalar product $(\,,)$\label{scalarr}.
If $A\subseteq\mathbb R^n$ is a non-empty subset, 
let $\chi_A$\label{chi} be its \emph{indicator function}, i.e. the function defined by $\chi_A(x)=1$ if $x\in A$ and $\chi_A(x)=0$ if $x\notin A$. Let also $\aff_\mathbb R A$\label{affr} denote the real affine subspace of $\mathbb R^n$ spanned by $A$, whereas let $E_A$\label{giacituraA} be the $\mathbb R$-linear subspace of $\mathbb R^n$ parallel to $\aff_\mathbb R A$. The set of interior points and the set of boundary points in the relative topology of $A$, as a subset of $\aff_{\mathbb R} A$, are referred to as the \textit{relative interior} and the \textit{relative boundary} of $A$, they are respectively noted $\relint A$\label{relintA} and $\relbd A$\label{relbdA}.
For any $0\leq k\leq n$, let ${\mathcal G}_k(\mathbb R^n)$\label{grassk} be the Grassman manifold of the (real) linear $k$-dimensional subspaces of $\mathbb R^n$ endowed with the usual topology, let also ${\mathcal G}(\mathbb R^n)$\label{wgrass} denote the whole grassmannian, i.e. the topological sum of the ${\mathcal G}_k(\mathbb R^n)$'s.  For every $E\in{\mathcal G}(\mathbb R^n)$, by $E^\perp$\label{perp} we will denote the orthogonal complement of $E$. A linear subspace $E\in{\mathcal G}(\mathbb R^n)$ is \textit{full-dimensional} if its dimension equals $n$.

A \textit{convex body} is a compact convex subset of $\mathbb R^n$. The set of convex bodies of $\mathbb R^n$ is noted ${\mathcal K}(\mathbb R^n)$\label{konvexr}, whereas ${\mathcal C}(\mathbb R^n)$\label{compactr} denotes that of nonempty compact subsets of $\mathbb R^n$. 
A \textit{polytope} is the convex hull of a finite, possibly empty,\footnote{By convention, one sets $\conv(\varnothing)=\{0\}$.} subset of $\mathbb R^n$. A \textit{lattice polytope} is the convex hull of a finite, possibly empty, subset of $\mathbb Z^n$. Every polytope can be realized as a (bounded) polyhedral set, i.e. a finite intersection of half-spaces.
The polytopes of $\mathbb R^n$ form a set denoted ${\mathcal P}(\mathbb R^n)$\label{polytopesr}. For a non-empty subset $A\subseteq \mathbb R^n$, the notions of dimension and full-dimensionality introduced for $\mathbb R$-linear subspaces of $\mathbb R^n$, can be simply extended by just invoking the subspace $E_A$. In particular, this can be done for convex bodies.  If $A\in{\mathcal K}(\mathbb R^n)$, set $d_A=\dim_\mathbb R A=\dim_\mathbb R E_A$\label{rdimA}. If $d\in\mathbb N^*$, a $d$-polytope $\Gamma$ is a polytope for which $d_\Gamma=d$. For $d=0,1,2,3,4$, a $d$-polytope is respectively referred to as a \textit{point, segment, polygon, polyhedron, polychoron}. If $v,w\in\mathbb R$, the segment with $v$ and $w$ as endpoints will be denoted $[v,w]$.\label{segmento} A \emph{simplex} is the convex hull of affinely independent elements of $R^n$, a $d$-simplex is a $d$-dimensional simplex, it is the convex hull of $d+1$ affinely independent points.

If $A\in{\mathcal K}(\mathbb R^n)$, let us denote $h_A$ the support function of $A$, i.e.
$
h_A(z)=\sup_{v\in A}(z,v)\,.
$\label{supportf} For example, if $B_{n}$\label{fullnball} is the full-dimensional unit ball in $\mathbb R^n$, then $h_{B_{n}}(z)=\Vert z\Vert$. 
The support function of a convex body is a positively $1$-homogeneous convex function. Remark that if $x=x_1+x_2\in \mathbb R^n=E_A\oplus E^\perp_A$, with $x_1\in E_A$ and $x_2\in E^\perp_A$ then $h_A(x)=h_A(x_1)$, so that $h_A$ actually depends on $d_A$ variables. 
Given $v\in\mathbb R^n\setminus\{0\}$, the corresponding \textit{supporting hyperplane} $H_A(v)$ for the convex body $A$ is given by
$
H_A(v)=\{z\in\mathbb C^n\mid ( z, v)=h_A(v)\}\,.
$\label{supportH}
Observe that, if $v=0$, the definition of $H_A(v)$ still makes sense, nevertheless it is no longer a hyperplane because $H_A(0)=\mathbb R^n$.

A \textit{face} $\Delta$ of $A$, (in symbols $\Delta\preccurlyeq A$\label{facce} or $A\succcurlyeq\Delta$), is a convex subset such that $u,v\in A$ and $(u+v)/2\in \Delta$ imply $u,v\in \Delta$. This notion applies also to the case of unbounded polyhedral sets. Every convex body admits two trivial faces, namely $\varnothing$ and $A$; these faces are referred to as the \textit{improper faces} of $A$ in contrast to the \textit{proper} ones which are non-empty and different from the convex body $A$ itself. For any proper face $\Delta$ of $A$ we will write $\Delta\prec A$ or $A\succ\Delta$.

An \textit{exposed face} $\Delta$ of $A$ is any intersection of the type $\Delta=A\cap H_A(v)$, for some $v\in\mathbb R^n$. By taking $v=0$, it follows that the non-empty improper face is an exposed face. The empty face is considered an exposed face too. For polytopes the distinction between faces and exposed faces is not necessary, since the two notions coincide. 
\begin{figure}[ht]
\centering
\begin{tikzpicture}[scale=0.4]
\fill[blue!20!white]  (4.5,4) arc (270:90:4)  -- (8.5,12) arc (90:-90:4) -- (4.5,4);
\draw[very thick,orange] (4.5,4) arc (270:90:4)  -- (8.5,12) arc (90:-90:4) -- (4.5,4);
\draw[ultra thick] (4.5,12) -- (8.5,12);
\fill (4.5,12) circle (3pt);
\fill (8.5,12) circle (3pt);
\end{tikzpicture}
\caption{The end-points of the black segment are faces but not exposed ones.}
\label{exp-faces}
\end{figure}
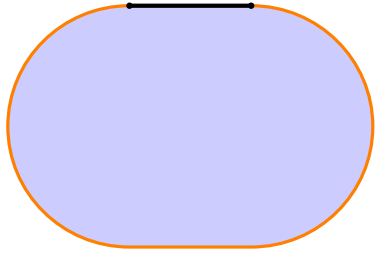

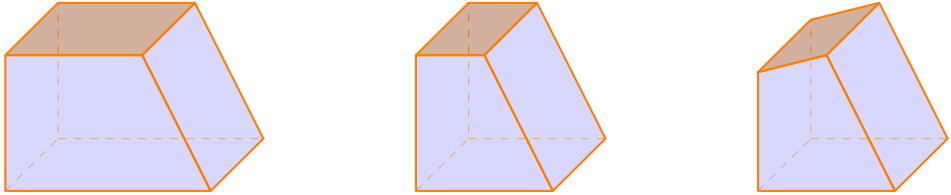
\begin{figure}[ht]
\begin{center}
\begin{tikzpicture}[scale=1]
\fill[fill=blue!20!white, semitransparent,line cap=round, line join=round] (-10,2,-2)--(-8,2,-2)--(-7,0,-2)--(-10,0,-2)--cycle;
\fill[fill=blue!20!white, semitransparent,line cap=round, line join=round] (-10,0,-2)--(-7,0,-2)--(-7,0,0)--(-10,0,0)--cycle;
\fill[fill=blue!20!white,semitransparent,line cap=round, line join=round] (-10,2,-2)--(-10,2,0)--(-10,0,0)--(-10,0,-2)--cycle;
\draw[orange,dashed] (-10,0,0)--(-10,0,-2)--(-10,2,-2);
\draw[orange,dashed] (-10,0,-2)--(-7,0,-2);
\fill[blue!20!white, semitransparent,line cap=round, line join=round] (-7,0,-2)--(-7,0,0)--(-8,2,0)--(-8,2,-2)--cycle;
\fill[brown, semitransparent,line cap=round, line join=round] (-10,2,-2)--(-8,2,-2)--(-8,2,0)--(-10,2,0)--cycle;
\fill[blue!20!white, semitransparent,line cap=round, line join=round] (-10,2,0)--(-8,2,0)--(-7,0,0)--(-10,0,0)--cycle;
\draw[orange,thick,line cap=round, line join=round] (-7,0,-2)--(-7,0,0)--(-8,2,0)--(-8,2,-2)--cycle;
\draw[orange,thick,line cap=round, line join=round] (-10,2,-2)--(-8,2,-2)--(-8,2,0)--(-10,2,0)--cycle;
\draw[orange,thick, line cap=round, line join=round] (-10,2,0)--(-8,2,0)--(-7,0,0)--(-10,0,0)--cycle;
\fill[fill=blue!20!white, semitransparent,line cap=round, line join=round] (-4,2,-2)--(-3,2,-2)--(-2,0,-2)--(-4,0,-2)--cycle;
\fill[fill=blue!20!white, semitransparent,line cap=round, line join=round] (-4,0,-2)--(-2,0,-2)--(-2,0,0)--(-4,0,0)--cycle;
\fill[fill=blue!20!white,semitransparent,line cap=round, line join=round] (-4,2,-2)--(-4,2,0)--(-4,0,0)--(-4,0,-2)--cycle;
\draw[orange,dashed] (-4,0,0)--(-4,0,-2)--(-4,2,-2);
\draw[orange,dashed] (-4,0,-2)--(-2,0,-2);
\fill[blue!20!white, semitransparent,line cap=round, line join=round] (-2,0,-2)--(-2,0,0)--(-3,2,0)--(-3,2,-2)--cycle;
\fill[brown, semitransparent,line cap=round, line join=round] (-4,2,-2)--(-3,2,-2)--(-3,2,0)--(-4,2,0)--cycle;
\fill[blue!20!white, semitransparent,line cap=round, line join=round] (-4,2,0)--(-3,2,0)--(-2,0,0)--(-4,0,0)--cycle;
\draw[orange,thick,line cap=round, line join=round] (-2,0,-2)--(-2,0,0)--(-3,2,0)--(-3,2,-2)--cycle;
\draw[orange,thick,line cap=round, line join=round] (-4,2,-2)--(-3,2,-2)--(-3,2,0)--(-4,2,0)--cycle;
\draw[orange,thick, line cap=round, line join=round] (-4,2,0)--(-3,2,0)--(-2,0,0)--(-4,0,0)--cycle;
\fill[blue!20!white, semitransparent,line cap=round, line join=round] (1,0,-2)--(1,0,0)--(1,1.75,0)--(1,1.75,-2)--cycle;
\fill[blue!20!white, semitransparent,line cap=round, line join=round] (1,0,-2)--(3,0,-2)--(2,2,-2)--(1,1.75,-2)--cycle;
\fill[blue!20!white, semitransparent,line cap=round, line join=round] (3,0,-2)--(1,0,-2)--(1,0,0)--(3,0,0)--cycle;
\draw[orange,dashed] (1,0,0)--(1,0,-2)--(1,1.75,-2);
\draw[orange,dashed] (1,0,-2)--(3,0,-2);
\fill[brown, semitransparent,line cap=round, line join=round] (1,1.75,-2)--(2,2,-2)--(2,2,0)--(1,1.75,0)--cycle;
\fill[blue!20!white, semitransparent,line cap=round, line join=round] (1,0,0)--(3,0,0)--(2,2,0)--(1,1.75,0)--cycle;
\fill[blue!20!white, semitransparent,line cap=round, line join=round] (3,0,-2)--(3,0,0)--(2,2,0)--(2,2,-2)--cycle;
\draw[orange,thick,line cap=round, line join=round] (1,1.75,-2)--(2,2,-2)--(2,2,0)--(1,1.75,0)--cycle;
\draw[orange,thick,line cap=round, line join=round] (1,0,0)--(3,0,0)--(2,2,0)--(1,1.75,0)--cycle;
\draw[orange,thick,line cap=round, line join=round] (3,0,-2)--(3,0,0)--(2,2,0)--(2,2,-2)--cycle;
\end{tikzpicture}
\end{center}
\caption{The three polytopes are combinatorially isomorphic but only the first two are strongly combinatorially isomorphic because the brown face of the third polytope is not parallel to the brown ones of the first two.}
\label{isom}
\end{figure}

If $\Gamma\in{\mathcal P}(\mathbb R^n)$, any non-empty face of $\Gamma$ is itself a polytope.
If $\Delta_1\preccurlyeq\Delta_2\preccurlyeq\Gamma$, then $\Delta_1\preccurlyeq\Gamma$ and $E_{\Delta_1}\subseteq E_{\Delta_2}\subseteq E_\Gamma$. 
The \textit{boundary complex} ${\mathcal B}(\Gamma)$\label{bc} of a polytope $\Gamma$ is the set of its non-empty faces. Two polytopes $\Gamma_1$ and $\Gamma_2$ are \emph{combinatorially isomorphic} if there exists a bijective and inclusion-preserving map from ${\mathcal B}(\Gamma_1)$ to ${\mathcal B}(\Gamma_2)$. The polytopes $\Gamma_1$ and $\Gamma_2$ are \emph{strongly combinatorially isomorphic} if there exists a combinatorial isomorphism $\psi:{\mathcal B}(\Gamma_1)\to{\mathcal B}(\Gamma_2)$ such that for any $\Delta\in {\mathcal B}(\Gamma_1)$ the affine subspace $\aff_\mathbb R \Delta$ is parallel to the affine subspace $\aff_\mathbb R \psi(\Delta)$. Two triangles are always combinatorially isomorphic but they are strongly combinatorially isomorphic if and only if the sides of the first triangle are respectively parallel to those of the second one. This shows that the preceding notions of isomorphism are different; Figure~\ref{isom} provides a three-dimensional example.

For every $0\leq k\leq \dim_\mathbb R\Gamma$, ${\mathcal B}(\Gamma,k)$\label{bck} will denote the subsets of ${\mathcal B}(\Gamma)$ consisting of $k$-faces (i.e. $k$-dimensional faces), the cardinal number of the set ${\mathcal B}(\Gamma,k)$ is denoted $f_k(\Gamma)$. A \textit{facet} is a face belonging to ${\mathcal B}(\Gamma,d_\Gamma-1)$, whereas a \textit{ridge} is a face belonging to ${\mathcal B}(\Gamma,d_\Gamma-2)$, a $1$-face is called an \textit{edge}, a $0$-face is called a \textit{vertex}. 
A polytope is \emph{simplicial} if all its proper faces are simplices. A $d$-polytope is \emph{simple} if each of its vertices belongs to $d-1$ facets.
The \textit{face vector} of a polytope $\Gamma\in{\mathcal P}(\mathbb R^n)$ is the vector 
\begin{equation*}\label{facevector}
f(\Gamma)=(f_0(\Gamma),\ldots,f_{d_\Gamma-1}(\Gamma))\,.
\end{equation*}
Setting $f_{-1}(\Gamma)=f_{d_\Gamma}(\Gamma)=1$ and $f_k(\Gamma)=0$ for $k<-1$ or $k>d_\Gamma$, it is well known that the natural numbers $f_k(\Gamma)$ cannot be arbitrary, indeed they  satisfy the Euler's relation
\begin{equation*}
\sum_{k=-1}^{d_\Gamma}(-1)^kf_k(\Gamma)=0\,,
\end{equation*}
as well as the inequalities
\begin{equation*}
\binom{d_\Gamma+1}{k+1}\leq f_k(\Gamma)\leq\binom{f_0(\Gamma)}{k+1} \,.
\end{equation*}

An \emph{oriented} $d$-polytope $\Gamma\in{\mathcal P}(\mathbb R^n)$ is a polytope such that  $E_\Gamma$ is an oriented linear subspace. To any complete flag of faces of $\Gamma$
\begin{equation*}
\Delta_0\prec\Delta_1\prec\ldots\prec\Delta_{d-2}\prec\Delta_{d-1}\prec \Delta_d=\Gamma
\end{equation*} 
starting with a vertex and ending with $\Gamma$, we can associate an orthonormal basis of $E_\Gamma$ made of outer normal unit vectors. 
In fact, for any vertex $\Delta_0$ of $\Gamma$, let $\Delta_1$ be an edge of $\Gamma$ admitting $\Delta_0$ as a vertex, then we have $E_{\Delta_0}=\{0\}$ and the trivial orthogonal decomposition \begin{equation*}
E_{\Delta_1}=E_{\Delta_0}\oplus (E_{\Delta_0}^\perp\cap E_{\Delta_1})=\{0\}\oplus(\mathbb R^n\cap E_{\Delta_1})=\{0\}\oplus E_{\Delta_1}\,.
\end{equation*} 
By induction suppose we have defined a sequence $\Delta_0\prec\Delta_1\prec\ldots\prec\Delta_{k-1}$ of faces of $\Gamma$ such that $\dim_\mathbb R E_{\Delta_\ell}=\ell$, for every $\ell\in\{0,1,\ldots,k-1\}$. Then, for every $k$-face $\Delta_k\preccurlyeq\Gamma$ admitting $\Delta_{k-1}$ as a facet, there is an orthogonal decomposition 
\begin{equation*}
E_{\Delta_k}=E_{\Delta_{k-1}}\oplus (E_{\Delta_{k-1}}^\perp\cap E_{\Delta_{k}})\,,
\end{equation*}
where the second direct summand is $1$-dimensional. The process ends after $d$ steps and yields the orthogonal decomposition
\begin{equation*}
E_\Gamma=\bigoplus_{k=1}^{d} (E_{\Delta_{k-1}}^\perp\cap E_{\Delta_k})\,,
\end{equation*}
where each summand is $1$-dimensional. 
For every $k\in\{1,\ldots,d\}$, the affine subspace $\aff_\mathbb R\Delta_{k-1}$ relatively bounds two relatively open half-spaces of $\aff_\mathbb R \Delta_k$ one of which does not intersect $\Delta_k$. The \emph{outer unit normal vector}\label{ounv} to the face $\Delta_{k-1}$ is the the unit vector $u_{\Delta_{k-1},\Delta_k}\in E_{\Delta_{k-1}}^\perp\cap E_{\Delta_k}$ pointing towards the relatively open half-space of $\aff_\mathbb R \Delta_k$ (relatively bounded by $\aff_\mathbb R\Delta_{k-1}$) that does not intersect $\Delta_k$.

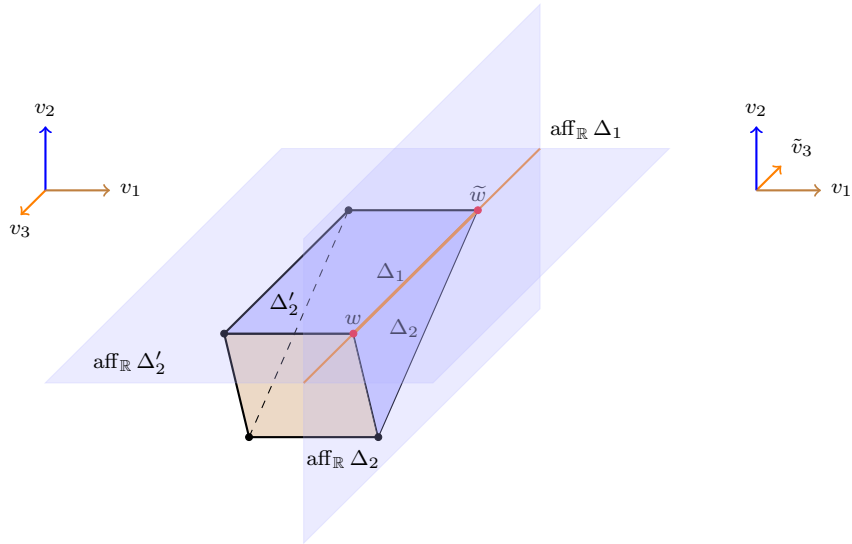
\begin{figure}[ht]
\begin{center}
\begin{tikzpicture}[scale=0.85]
\draw[->,thick,orange] (-2,0,0)-- (-2,0,1) node[below,black]{\footnotesize$v_3$};
\draw[->,thick,blue] (-2,0,0)-- (-2,1,0) node[above,black]{\footnotesize$v_2$};
\draw[->,thick,brown] (-2,0,0)-- (-1,0,0) node[right,black]{\footnotesize$v_1$};
\draw[->,thick,orange] (9,0,0)-- (9,0,-1) node[above right,black]{\footnotesize$\tilde v_3$};
\draw[->,thick,blue] (9,0,0)-- (9,1,0) node[above,black]{\footnotesize$v_2$};
\draw[->,thick,brown] (9,0,0)-- (10,0,0) node[right,black]{\footnotesize$v_1$};
\filldraw[draw=black,fill=brown!30!white,thick] (5,2,11)--(7,2,11)--(7,0,10)--(5,0,10)--cycle;
\filldraw[draw=black,fill=blue!20!white,thick] (5,2,11)--(7,2,11)--(7,2,6)--(5,2,6)--cycle;
\draw[dashed] (5,0,10)--(5,2,6);
\filldraw[blue!20!white] (7,2,6)--(7,2,11)--(7,0,10)--cycle;
\draw (7,2,6)--(7,2,11)--(7,0,10)--cycle;
\draw[thick, orange] (7,2,13)--(7,2,3.5) node[above right, black]{\footnotesize$\aff_\mathbb R\Delta_1$};
\draw[very thick, orange] (7,2,6)--(7,2,11);
\draw (7,2,8.5) node[left]{\footnotesize$\Delta_1$};
\draw (7.2,1.125,8.5) node[left]{\footnotesize$\Delta_2$};
\filldraw (5,2,11) circle(1.5pt);
\filldraw (5,2,6) circle(1.5pt);
\filldraw[red] (7,2,11) circle(1.5pt) node[above, black]{\footnotesize$w$};
\filldraw[red] (7,2,6) circle(1.5pt) node[above, black]{\footnotesize$\widetilde w$};
\filldraw (5,0,10) circle(1.5pt);
\filldraw (7,0,10) circle(1.5pt);
\filldraw[blue!30!white, nearly transparent] (7,4.25,13)--(7,-0.5,13)-- (7,-0.5,3.5)--(7,4.25,3.5)--cycle;
\filldraw[blue!30!white, nearly transparent] (3,2,13)--(9,2,13)-- (9,2,3.5)--(3,2,3.5)--cycle;
\draw (5.3,-1.5,7) node{\footnotesize$\aff_\mathbb R \Delta_2$};
\draw (2,0,7) node{\footnotesize$\aff_\mathbb R \Delta_2^\prime$};
\draw (5.35,1.5,8.5) node[left]{\footnotesize$\Delta_2^\prime$};
\end{tikzpicture}
\caption{Induced orientations.}
\label{basis}
\end{center}
\end{figure}

For the sake of notation, let us set $v_{d-k}=u_{\Delta_{k-1},\Delta_k}$, for every $k\in\{0,\ldots,d-1\}$. This construction provides an orthonormal basis $(v_1,\ldots,v_d)$ of $E_\Gamma$ depending on the choice of the flag. If $\Gamma$ is oriented, it is possible to choose a flag such that the obtained basis conforms to the orientation of $\Gamma$. Indeed, if the flag 
$\Delta_0\prec\Delta_1\prec\ldots\prec\Delta_{d-2}\prec\Delta_{d-1}\prec \Delta_d=\Gamma$ produces the wrong orientation, it is enough to replace the vertex $\Delta_0$ with the other vertex of $\Delta_1$. Of course when $d<n$, such a basis can be completed to a positively (resp. negative) basis of $\mathbb R^n$. A facet $\Delta$ of an oriented $d$-polytope $\Gamma$ is said to be \emph{oriented by the unit outer} (resp. \emph{inner}) \emph{normal vector} $u_{\Delta,\Gamma}$ (resp. $-u_{\Delta,\Gamma}$) if $E_\Delta$ is given a basis $(v_2,\ldots,v_d)$ such that $(u_{\Delta,\Gamma},v_2,\ldots,v_d)$ conforms to the chosen orientation of $E_\Gamma$. By induction, the choice of an orientation for $\Gamma$ implies an orientation of each of its faces. Observe, however, that a same face can get different orientations from different facets. For example, a ridge gets opposite orientations from the two facets containing it. 
Figure \ref{basis} depicts the situation for a positively oriented $3$-polytope $\Gamma\subset\mathbb R^3$. 
The basis $(v_1,v_2,v_3)$ corresponds to the flag $\Delta_0\prec\Delta_1\prec\Delta_2\prec\Gamma$. On the left, the case $\Delta_0=w$. On the right, the case $\Delta_0=\widetilde w$. The only difference between the two basis is the third vector. The first basis is positively oriented and $(v_2,v_3)$ gives the facet $\Delta_2$ the outer normal vector orientation, whereas $(v_2,\tilde v_3)$ gives $\Delta_2$ the opposite one. Moreover, the facets $\Delta_2$ and $\Delta_2^\prime$ give the common edge $\Delta_1$ opposite orientations.

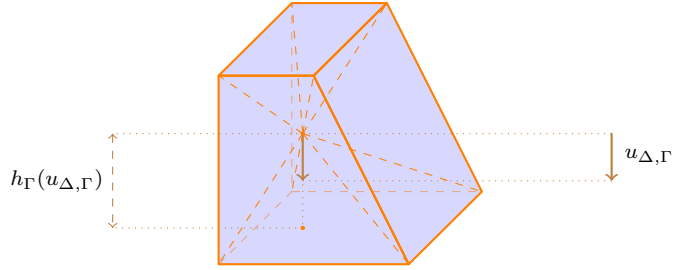
\begin{figure}[ht]
\begin{center}
\begin{tikzpicture}[scale=1.25]
\fill[fill=blue!20!white, semitransparent,line cap=round, line join=round] (-4,2,-2)--(-3,2,-2)--(-2,0,-2)--(-4,0,-2)--cycle;
\fill[fill=blue!20!white, semitransparent,line cap=round, line join=round] (-4,0,-2)--(-2,0,-2)--(-2,0,0)--(-4,0,0)--cycle;
\fill[fill=blue!20!white,semitransparent,line cap=round, line join=round] (-4,2,-2)--(-4,2,0)--(-4,0,0)--(-4,0,-2)--cycle;
\draw[orange,dashed] (-4,0,0)--(-4,0,-2)--(-4,2,-2);
\draw[orange,dashed] (-4,0,-2)--(-2,0,-2);
\fill[blue!20!white, semitransparent,line cap=round, line join=round] (-2,0,-2)--(-2,0,0)--(-3,2,0)--(-3,2,-2)--cycle;
\fill[blue!20!white, semitransparent,line cap=round, line join=round] (-4,2,-2)--(-3,2,-2)--(-3,2,0)--(-4,2,0)--cycle;
\fill[blue!20!white, semitransparent,line cap=round, line join=round] (-4,2,0)--(-3,2,0)--(-2,0,0)--(-4,0,0)--cycle;
\draw[orange,thick,line cap=round, line join=round] (-2,0,-2)--(-2,0,0)--(-3,2,0)--(-3,2,-2)--cycle;
\draw[orange,thick,line cap=round, line join=round] (-4,2,-2)--(-3,2,-2)--(-3,2,0)--(-4,2,0)--cycle;
\draw[orange,thick, line cap=round, line join=round] (-4,2,0)--(-3,2,0)--(-2,0,0)--(-4,0,0)--cycle;
\draw[orange,dashed] (-3.5,1,-1)--(-4,2,-2);
\draw[orange,dashed] (-3.5,1,-1)--(-4,2,0);
\draw[orange,dashed] (-3.5,1,-1)--(-3,2,-2);
\draw[orange,dashed] (-3.5,1,-1)--(-3,2,0);
\draw[orange,dashed] (-3.5,1,-1)--(-2,0,-2);
\draw[orange,dashed] (-3.5,1,-1)--(-2,0,0);
\draw[orange,dashed] (-3.5,1,-1)--(-4,0,-2);
\draw[orange,dashed] (-3.5,1,-1)--(-4,0,0);
\draw[brown,dotted] (-3.5,1,-1)--(-3.5,0,-1);
\draw[brown,thick,->,line cap=round,line join=round] (-3.5,1,-1)--(-3.5,0.5,-1);
\draw[brown,dotted] (-3.5,0,-1)--(-5.5,0,-1);
\draw (-5.5,0.5,-1) node[left]{\footnotesize$h_\Gamma(u_{\Delta,\Gamma})$};
\draw[brown,dotted] (-0.25,1,-1)--(-5.5,1,-1);
\draw[brown,dotted] (-0.25,0.5,-1)--(-3.5,0.5,-1);
\draw[brown,thick,->,line cap=round] (-0.25,1,-1)--(-0.25,0.5,-1);
\draw (0.5,0.75,-1) node[left]{\footnotesize$u_{\Delta,\Gamma}$};
\draw[brown,dashed,<->] (-5.5,1,-1)--(-5.5,0,-1);
\filldraw[orange](-3.5,0,-1) circle(0.5pt);
\filldraw[orange](-3.5,1,-1) circle(0.5pt);

\end{tikzpicture}
\caption{A polytope, containing the origin, decomposed into pyramids based onto its facets.}
\label{fig-piramidi}
\end{center}
\end{figure}

The $d$-volume of a $d$-polytope is related to the $(d-1)$-volume of its facets. Indeed, let $\Gamma\in{\mathcal P}(\mathbb R^n)$ and, for every $\Delta\in{\mathcal B}(\Gamma,d-1)$, let $u_{\Delta,\Gamma}\in E_\Delta^\perp\cap E_\Gamma$ be the outer unit normal vector to the facet $\Delta$. Then the following equalities are worth-noting:
\begin{equation}\label{piramidi}
\sum_{\Delta\in{\mathcal B}(\Gamma,d-1)} \vol_{d-1}(\Delta)u_{\Delta,\Gamma}=0\,,
\end{equation}
\begin{equation}\label{vol-piramidi}
\vol_{d}(\Gamma)=\frac{1}{d}\sum_{\Delta\in{\mathcal B}(\Gamma,d-1)} h_\Gamma(u_{\Delta,\Gamma})\vol_{d-1}(\Delta)\,.
\end{equation}
Figure~\ref{fig-piramidi} provides a geometric interpretation of~\eqref{piramidi}.

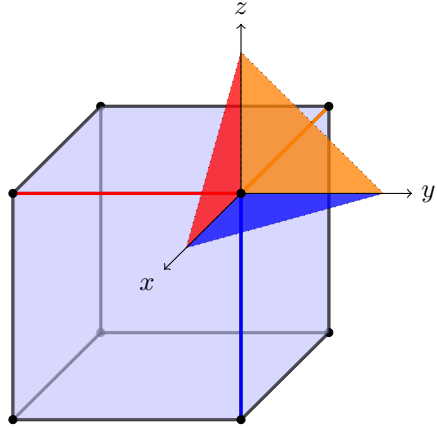
\begin{figure}[ht]
\begin{center}
\begin{tikzpicture}[scale=1.5,path fading=fade right]
\filldraw (0,0,0) circle(1pt);
\filldraw (2,0,0) circle(1pt);
\filldraw (0,2,0) circle(1pt);
\filldraw[fill=blue!20!white, draw=black,very thick,semitransparent] (0,0,0)--(0,0,2)--(2,0,2)--(2,0,0)--cycle;
\filldraw[fill=blue!20!white, draw=black,very thick,semitransparent] (2,0,0)--(2,2,0)--(0,2,0)--(0,0,0)--cycle;
\filldraw[fill=blue!20!white, draw=black,very thick,semitransparent] (0,0,0)--(0,0,2)--(0,2,2)--(0,2,0)--cycle;
\filldraw[fill=blue!20!white, draw=black,very thick,semitransparent] (0,0,2)--(2,0,2)--(2,2,2)--(0,2,2)--cycle;
\filldraw[fill=blue!20!white, draw=black,very thick,semitransparent] (2,2,2)--(2,0,2)--(2,0,0)--(2,2,0)--cycle;
\filldraw[fill=blue!20!white, draw=black,very thick,semitransparent] (0,2,0)--(0,2,2)--(2,2,2)--(2,2,0)--cycle;
\draw[red,very thick,line cap=round, line join=round] (0,2,2)--(2,2,2);
\draw[blue,very thick,line cap=round, line join=round] (2,0,2)--(2,2,2);
\draw[orange,very thick,line cap=round, line join=round] (2,2,0)--(2,2,2);
\filldraw (2,2,0) circle(1pt);
\filldraw (0,0,2) circle(1pt);
\filldraw (2,0,2) circle(1pt);
\filldraw (0,2,2) circle(1pt);
\filldraw (2,2,2) circle(1pt);
\filldraw[fill=blue, draw=black,dotted,nearly opaque] (2,2,2)--(3.25,2,2)--(2,2,3.25)--cycle;
\filldraw[fill=red, draw=black,dotted,nearly opaque] (2,2,2)--(2,3.25,2)--(2,2,3.25)--cycle;
\filldraw[fill=orange, draw=black,dotted,nearly opaque,path fading,fading transform={rotate=45}] (2,2,2)--(3.25,2,2)--(2,3.25,2)--cycle;
\filldraw (2,2,2) circle(1pt);
\draw[->] (2,2,2)--(3.5,2,2) node[right]{$y$};
\draw[->] (2,2,2)--(2,3.5,2) node[above]{$z$};
\draw[->] (2,2,2)--(2,2,3.75) node[below left]{$x$};
\end{tikzpicture}
\caption{A cube and parts of the dual cones to the edges passing through the origin.}
\label{angoli-esterni}
\end{center}
\end{figure}

For any non-empty face $\Delta$ of a polytope $\Gamma\in{\mathcal P}(\mathbb R^n)$, there is a \textit{dual cone} $K_{\Delta,\Gamma}$\label{dualcone} (or simply $K_\Delta$, if no confusion may arise) defined as
$
K_{\Delta,\Gamma}=\{ v\in\mathbb R^n\mid \Delta=\Gamma\cap H_\Gamma(v)\}
$.
For the improper face $\Gamma$, one has $K_{\Gamma,\Gamma}=E_\Gamma^\perp$ and this is the only  dual cone that is closed. If $\Delta$ is a proper $k$-face, the subset $K_{\Delta,\Gamma}$ is a non-empty convex polyhedral cone in $E_\Delta^\perp$ (i.e. an unbounded intersection of finitely many open half-spaces passing through the origin) and a relatively open $(2n-k)$-submanifold of $\mathbb R^n$. 
Figure~\ref{angoli-esterni} shows parts of the dual cones to some edges of a cube.

\begin{figure}[ht]
\begin{center}
\begin{tikzpicture}[scale=0.85]
\draw[->] (0,0,0)--(9,0,0) node[right]{$E_\Delta^\perp$};
\draw (4,3,0) node{$E_\Lambda^\perp$};
\draw[->] (0,0,0)--(0,3,0) node[above]{$E_\Lambda^\perp\cap E_\Delta$};
\draw[->] (0,0,0)--(0,0,10) node[below left]{$E_\Lambda$};
\filldraw[draw=black,fill=brown!20!white,thick] (5,2,11)--(7,2,11)--(7,0,10)--(5,0,10)--cycle;
\filldraw[draw=black,fill=brown!20!white,thick] (5,2,11)--(7,2,11)--(7,2,6)--(5,2,6)--cycle;
\draw[dashed] (5,0,10)--(5,2,6);
\filldraw[blue!20!white] (7,2,6)--(7,2,11)--(7,0,10)--cycle;
\draw (7,2,6)--(7,2,11)--(7,0,10)--cycle;
\draw[thick, orange] (7,2,12)--(7,2,-1) node[above right, black]{$\aff_\mathbb R\Lambda$};
\draw[very thick, orange] (7,2,6)--(7,2,11);
\draw (7,2,8.5) node[left]{$\Lambda$};
\draw (7,1.125,8.5) node[left]{$\Delta$};
\draw[->,thick] (0,0,0)--(0,2,0) node[left]{$h_\Delta(u_{\Lambda,\Delta})u_{\Lambda,\Delta}$};
\draw[->,thick] (0,0,0)--(0,1,0) node[left]{$u_{\Lambda,\Delta}$};
\draw[dotted] (7,0,0)--(7,2,0);
\draw[dotted] (0,2,0)--(7,2,0);
\filldraw (5,2,11) circle(1.5pt);
\filldraw (5,2,6) circle(1.5pt);
\filldraw (7,2,11) circle(1.5pt);
\filldraw (7,2,6) circle(1.5pt);
\filldraw (5,0,10) circle(1.5pt);
\filldraw (7,0,10) circle(1.5pt);
\filldraw[blue!20!white, nearly transparent] (7,3,12)--(7,-0.5,12)-- (7,-0.5,-1)--(7,3,-1)--cycle;
\draw (8,1.25,7) node{$\aff_\mathbb R \Delta$};
\draw (-2,1,2) node{$E_\Delta$};
\draw[->,thick] (0,0,0)--(7,0,0) node[below]{$p_\Delta$};
\draw[thick,->] (0,0,0)--(7,2,0) node[above]{$p_\Lambda$};
\filldraw (7,0,0) circle(1.5pt);
\filldraw (7,2,0) circle(1.5pt);
\draw (5.25,1.5,8.5) node[left]{$\Gamma$};
\end{tikzpicture}
\caption{$\Gamma$ is the light brown polytope with the blue face $\Delta$ and the orange edge $\Lambda$. The blue, transparent plane is $\aff_\mathbb R \Delta$, the orange line is $\aff_\mathbb R\Lambda$. }
\label{decomp}
\end{center}
\end{figure}
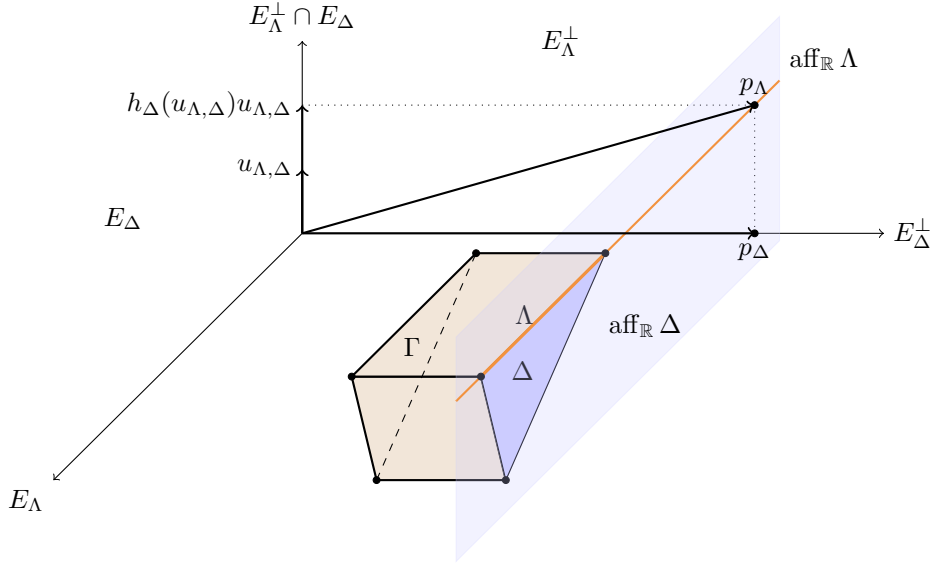

On the dual cone $K_{\Delta,\Gamma}$ the support function $h_\Gamma$ is linear and, for every $z\in K_{\Delta,\Gamma}$, $h_\Gamma(z)=h_\Delta(z)=(z,v)$, where $v$ is any point of $\Delta$. In fact, the chosen point $v\in\Delta$ may be represented as $v=v_1+v_2$, with $v_1\in E_\Delta$ and $v_2\in E^\perp_\Delta$, then, on $K_{\Delta,\Gamma}$, one has $h_\Gamma(z)=h_\Delta(z)=(z,v)=(z,v_2)$. For any non empty face $\Delta$ of $\Gamma$, let $p_\Delta$ be the only point in the intersection $\aff_\mathbb R \Delta\cap E_\Delta^\perp$, then $h_\Gamma=(\cdot\,,p_\Delta)$ on the whole $K_{\Delta,\Gamma}$. As a consequence, if $\Lambda$ is a facet of $\Delta\preccurlyeq\Gamma$ and $u_{\Lambda,\Delta}\in E^\perp_\Lambda\cap E_\Delta$ is its outer unit normal vector, there is an orthogonal decomposition: 
\begin{equation}\label{decomposition}
p_\Lambda=p_\Delta+h_\Delta(u_{\Lambda,\Delta})u_{\Lambda,\Delta}\,.
\end{equation} 
Figure~\ref{decomp} shows that the points $p_\Lambda$ and $p_\Delta$ need not belong to $\Lambda$ and $\Delta$, respectively. 

If $\Delta\preccurlyeq \Gamma$ is a $k$-dimensional face, the topological boundary of the manifold $K_{\Delta,\Gamma}$ equals its relative boundary and is given by the disjoint union of the cones which are dual to the faces of $\Gamma$ admitting $\Delta$ as a proper face, or by the union of the topological closures of the cones which are dual to the faces of $\Gamma$ admitting $\Delta$ as a facet, i.e.
\begin{equation*}
\partial K_{\Delta,\Gamma}
=
\relbd K_{\Delta,\Gamma}
=
\bigcup_{\substack{\Delta^\prime\in{\mathcal B}(\Gamma)\\ \Delta\prec\Delta^\prime}} K_{\Delta^\prime,\Gamma}
=
\bigcup_{\substack{\Delta^\prime\in{\mathcal B}(\Gamma,k+1)\\ \Delta\prec\Delta^\prime}} \overline{K}_{\Delta^\prime,\Gamma}
\,.
\end{equation*}
Observe that, if $\Delta\prec\Gamma$, the closure of the convex polyhedral cone $K_{\Delta,\Gamma}$ is spanned by the unit outer normal vectors $u_{\Delta,\Delta^\prime}$, as $\Delta^\prime$ runs in the set of those faces of $\Gamma$ admitting $\Delta$ as a facet.

For any integer $0\leq k\leq \dim\Gamma$, the \textit{$k$-star} $\Sigma_{k,\Gamma}$\label{kstar} (or simply $\Sigma_k$, if no confusion may arise) is the disjoint union of the cones which are dual to the $k$-dimensional faces of $\Gamma$. Notice that $\Sigma_k$ is a $(2n-k)$-dimensional (disconnected) manifold and
\begin{equation*}
\mathbb R^n=\bigcup_{0\leq k\leq d_\Gamma}\Sigma_k=\bigcup_{\Delta\in{\mathcal B}(\Gamma)} K_{\Delta,\Gamma}=\bigcup_{v\in{\mathcal B}(\Gamma,0)} \overline{K}_{v,\Gamma}\,,
\end{equation*}
where, for any vertex $v\in\Gamma$, $\overline{K}_{v,\Gamma}$ is the topological closure of $K_{v,\Gamma}$.
It follows that $h_\Gamma$ is a continuous piecewise linear function on the whole $\mathbb R^n$ with the following representation
\begin{equation*}
h_\Gamma(z)=\sum_{\Delta\in{\mathcal B}(\Gamma)}(p_\Delta, z)\chi_{K_{\Delta,\Gamma}}(z)\,.
\end{equation*}
By continuity, on $K_{\Lambda,\Gamma}$, one has $h_\Lambda=h_\Delta$, for any $\Delta\preccurlyeq \Gamma$ such that $\Lambda\preccurlyeq\Delta$. 

If $\Delta\in{\mathcal B}(\Gamma,k)$, the \textit{outer angle} of $\Delta$ with respect to $\Gamma$ is the number 
\begin{equation}\label{outerangle}
\psi_\Gamma(\Delta)=\varkappa^{-1}_{n-k}\vol_{n-k}(K_{\Delta,\Gamma}\cap B_{n})\,,
\end{equation} 
where, for any $0\leq\ell\leq n$, $\vol_\ell$\label{vollell} is the $\ell$-dimensional Lebesgue measure in $\mathbb R^\ell$, whereas $\varkappa_\ell$\label{varkappaell} is the $\ell$-dimensional Lebesgue measure of the $\ell$-dimensional closed unit ball $B_\ell$ about the origin of $\mathbb R^\ell$. The subscript $\Gamma$ in the notation $\psi_\Gamma(\Delta)$ will usually be dropped, unless it is necessary to avoid confusion.\footnote{The outer angle $\psi_\Gamma(\Delta)$ is also known as the \emph{external angle} and it is sometimes noted $\gamma(\Delta,\Gamma)$.}

Recall that, if $\varGamma$ denotes Euler's gamma function, then 
\begin{equation*}
\varkappa_\ell=\frac{{\pi^{(\ell/2)}}}{{\varGamma(1+(\ell/2))}}\,,
\end{equation*} 
in particular $\varkappa_0=1$, $\varkappa_1=2$ and $\varkappa_{2\ell}={\pi^\ell}/{\ell!}$.
By Fubini's theorem together with the properties of Euler's gamma and beta functions, it is possible to show that, for a $d$-polytope $\Gamma$, outer angles can be computed in the smaller space $E_\Gamma$, i.e. $\psi_\Gamma(\Delta)=\varkappa_{d-k}^{-1}\vol_{d-k}(K_{\Delta,\Gamma}\cap E_\Gamma\cap B_{n})$.
 Observe that $\psi_\Gamma(\Gamma)=1$ and, for any facet $\Delta\prec\Gamma$, $\psi_\Gamma(\Delta)=1/2$. 
 
The \emph{normal fan} $\Sigma(\Gamma)$\label{nfan} of $\Gamma$ is the union of the $k$-stars $\Sigma_{k,\Gamma}$, as $k$ runs in the set $\{0,\ldots,\dim_\mathbb R\Gamma\}$. Two polytopes $\Gamma_1$ and $\Gamma_2$ are strongly combinatorially isomorphic if and only if $\Sigma(\Gamma_1)=\Sigma(\Gamma_2)$.
 
For any compact subset $A\subset \mathbb R^n$ and any $\lambda\in\mathbb R$, one defines the \textit{multiplication of $A$ by $\lambda$} as the subset $\lambda A=\{\lambda z\in\mathbb R^n\mid z\in A\}$. Of course, if $A$ is a convex body and $\lambda\neq0$, $h_{\lambda A}=\lambda h_A$ and when $A$ is a polytope, $\lambda A$ is also a polytope. For every polytope $\Gamma$, any face of the polytope $\lambda \Gamma$ has the form $\lambda\Delta$, for a unique face $\Delta$ of $\Gamma$, so that
\begin{equation*}
\aff_{\mathbb R}\lambda\Delta=\lambda\aff_{\mathbb R}\Delta\,,
\qquad
E_{\lambda \Delta}=E_\Delta\,,
\qquad
E_{\lambda\Delta}^\perp=E_\Delta^\perp\,,
\qquad
K_{\lambda\Delta,\lambda \Gamma}=K_{\Delta,\Gamma}\,.
\end{equation*}

The convex bodies $A_1$ and $A_2$ are \textit{(positively) homothetic} if $A_1=t+\lambda A_2$, for some $t\in\mathbb R^n$ and $\lambda\in\mathbb R_{\geq 0}$, i.e. if one of them is a translate of a multiple of the other one or is reduced to a single point.

The \textit{Minkowski sum} $A=\sum_{\ell\in I}A_\ell$ of a finite family of subsets $\{A_\ell\}_{\ell\in I}$ is either the zero set $\{0\}$ in case $I=\varnothing$, or the subset $A=\{\sum_{\ell\in I} v_\ell\mid v_\ell\in A_\ell\,,\; \ell\in I\}$ in case $I$ is a non-empty finite set. The spaces ${\mathcal C}(\mathbb R^n)$ and ${\mathcal K}(\mathbb R^n)$ are stable with respect to the Minkowski addition. Of course, if $I\neq\varnothing$ and the $\{A_\ell\}_{\ell\in I}$ are convex bodies, then $h_A=\sum_{\ell\in I} h_{A_\ell}$ and when the summands are polytopes, the sum is itself a polytope. In the latter case, for any face $\Delta$ of the Minkowski sum $\Gamma=\sum_{\ell\in I}\Gamma_\ell$ there exists a unique sequence of faces $\Delta_\ell\preccurlyeq \Gamma_\ell$, $\ell\in I$ for which $\Delta=\sum_{\ell\in I}\Delta_\ell$, such uniquely determined faces are the \emph{(Minkowski) summands} of $\Delta$. As a consequence, one gets the equalities
\begin{equation}\label{utili}
\aff_{\mathbb R}\Delta=\sum_{\ell\in I} \aff_{\mathbb R}{\Delta_\ell}\,,
\qquad
E_\Delta=\sum_{\ell\in I} E_{\Delta_\ell}\,,
\end{equation}
\begin{equation}\label{utili2}
E_\Delta^\perp=\bigcap_{\ell\in I} E_{\Delta_\ell}^\perp\,,
\qquad
K_{\Delta,\Gamma}=\bigcap_{\ell\in I} K_{\Delta_\ell,\Gamma_\ell}\,.
\end{equation}
If $\# I\leq n$, an \emph{edge-sum face} $\Delta$ of a Minkowski sum of polytopes $\Gamma=\sum_{\ell\in I}\Gamma_\ell$ is a face whose summands $\Delta_\ell$, respectively, include segments $\Lambda_\ell$, $\ell\in I$, spanning linearly independent directions.

Notice that the Minkowski sum of finitely many strongly combinatorially isomorphic polytopes yields a polytope that is strongly combinatorially isomorphic to each of its summands.

When $\lambda\in\mathbb N^*$, one has $\lambda A=\sum_{\ell=1}^\lambda A$, for any convex body $A$. If $A\subset \mathbb R^n$ and $\varepsilon>0$, the convex body $(A)_\varepsilon=A+\varepsilon B_n$\label{eintorno} is called the \textit{$\varepsilon$-neighborhood} of $A$. It is just the Minkowski sum of $A$ and $\varepsilon$ times the full-dimensional closed unit ball $B_n$. In particular $(B_n)_\varepsilon$ is the full-dimensional closed ball around the origin of radius $1+\varepsilon$. Remark that the Minkowski sum of a convex body $A$ with a full-dimensional ball has the effect of smoothing the boundary of $A$. Observe also that, for every $\varepsilon>0$, the $\varepsilon$-neighborhood of any subset $A$ contains $\partial A$, in particular, for any polytope $\Gamma$, one has
\begin{equation*}
\Sigma_{k,\Gamma}\subset(\Sigma_{k,\Gamma})_\varepsilon\subseteq\left(\overline{\Sigma_{1,\Gamma}}\right)_\varepsilon\,,
\end{equation*}
for any $1\leq k\leq\dim\Gamma$. The class of convex subsets of $\mathbb R^n$ which are $\varepsilon$-neighborhood of polytopes, for some $\varepsilon>0$, will be noted ${\mathcal P}_\circ(\mathbb R^n)$\label{pinfinitor}. 

\begin{exe}\label{nsimplesso} The standard $n$-simplex.

\noindent
{\rm The standard $n$-simplex is the polytope $\Delta_n$ defined as
\begin{equation*}
\Delta_n=\left\{x\in [0,1]^{n+1}\bigg\vert\sum_{\ell=1}^{n+1}x_\ell=1\right\}\,.
\end{equation*}
For any integer $0\leq k\leq n-1$, the $k$-th component of the face vector $f(\Delta_{n})$ is given by
\begin{equation*}
f_{k}(\Delta_{n})=\binom{n+1}{k+1}\,,
\end{equation*}
and, for any face $\Delta\in{\mathcal B}(\Delta_n,k)$, one has $\vol_k(\Delta)=\sqrt{k+1}/k$.
\findim}

\end{exe}

\begin{exe}\label{ncubo} The standard $n$-cube.

\noindent
{\rm
Let $I_{n}=[-1,1]^{n}\subset\mathbb R^{n}$. For any integer $0\leq k\leq n-1$, the $k$-th component of the face vector $f(I_{n})$ is given by
\begin{equation*}
f_{k}(I_{n})=2^{n-k}\binom{n}{k}\,,
\end{equation*}
and, for any face $\Delta\in{\mathcal B}(I_n,k)$, one has $\vol_k(\Delta)=2^k$ and $\psi(\Delta)=2^{k-n}$.\findim
}
\end{exe}

\begin{exe}\label{ntetra} The standard $n$-crosspolytope.

\noindent
{\rm
Let $\Theta_{n}\subset\mathbb R^{n}$ be the convex hull of the $2^{n}$ points whose coordinates are all zero except one of them that runs in the set $\{1,-1\}$. Equivalently $\Theta_{n}=\{x\in\mathbb R^n\mid\sum_{1\leq\ell\leq n} \vert x_\ell\vert\leq 1\}$. The $k$-th component of the face vector $f(\Theta_{n})$ is given by
\begin{equation*}
f_{k}(\Theta_{n})=2^{k+1}{n\choose k+1}\,.
\end{equation*}
The unit $n$-crosspolytope can be constructed recursively as the convex hull of the union $\Theta_1\cup\Theta_{n-1}$, where $\Theta_{1}=[-1,1]\subset \mathbb R\times \{0\}^{n-1}\subset\mathbb R^{n}$ and $\Theta_{n-1}\subset \mathbb R^{n-1}\times\{0\}\subset\mathbb R^{n}$ are respectively a unit $1$-crosspolytope (i.e. a segment) and a unit $(n-1)$-crosspolytope orthogonal to each other in $\mathbb R^{n}$. It thus follows that $\vol_n(\Theta_n)=2^n/n!$, whereas $\vol_k(\Delta)=\sqrt{k+1}/k!$, for any $\Delta\in{\mathcal B}(\Theta_n,k)$, with $1\leq k<n$.
As $\Theta_n$ is a centrally symmetric polytope, the outer angle to any $k$-face depends only on the dimension of the face. If $e_\ell$, $1\leq \ell\leq n$ is the usual basis of $\mathbb R^n$, consider, for $1\leq j\leq n$, the vectors $u_j=e_{n-j+1}+\sum_{\ell=1}^{n-j} e_\ell$ and $v_j=-e_{n-j+1}+\sum_{\ell=1}^{n-j} e_\ell$ together with the corresponding supporting hyperplanes to $\Theta_n$. Each such vector spans the dual cone to a facet of $\Theta_n$ and different vectors span different facets. For example, when $n=4$, the outer angle of any $2$-face of $\Theta_4$ is $1/6$. In fact, if $\Delta$ denotes the $2$-face obtained by intersecting $\Theta_4$ with the supporting hyperplanes corresponding to $u_1$ and $v_1$, one realizes that  $\cos\widehat{u_1v_1}=(u_1,v_1)/\Vert u_1\Vert\cdot\Vert v_1\Vert=1/2$, i.e. $\widehat{u_1v_1}=\pi/3$, whence
\begin{equation*}
\psi(\Delta)=\frac{\vol_2(K_\Delta\cap B_4)}{\varkappa_2}=\frac{(1/2)(\pi/3)}{\pi}=\frac{1}{6}\,,
\end{equation*}
as asserted.\findim

}
\end{exe}
By a piecewise ${\mathcal C}^2$-hypersurface in $\mathbb R^n$ we mean a closed ${\mathcal C}^1$-hypersurface $X$ with the following properties: $\partial X=\varnothing$; $X$ admits a subset $X_c$ such that $X\setminus X_c$ is a finite disjoint union of ${\mathcal C}^2$-hypersurfaces; the $(n-1)$-dimensional Hausdorff measure of $X_c$ is zero.

A connected compact subset of $\mathbb R^n$ is said to be piecewise $2$\textit{-regular} if it has non-empty interior and its boundary is a piecewise ${\mathcal C}^2$-hypersurface. The space of piecewise $2$-regular subsets will be denoted ${\mathcal S}^{2}(\mathbb R^n)$\label{skreg}. 

Any piecewise $2$-regular subset $A$ admits a \textit{(global) defining function} $\rho_A\in{\mathcal C}^1(\mathbb R^n,\mathbb R)$\label{deffunct} such that $\rho_A\equiv0$ on $\partial A$, $\rho_A>0$ on $\mathbb R^n\setminus A$, $\rho_A<0$ on $\relint A$ and $\de\rho_A\neq0$ on $\partial A$. 
We set ${\mathcal K}^2(\mathbb R^n)={\mathcal S}^2(\mathbb R^n)\cap{\mathcal K}(\mathbb R^n)$\label{kreg2}. If $A\in{\mathcal K}^2(\mathbb R^n)$, then any of its defining functions has to be convex. Observe that ${\mathcal P}_\circ(\mathbb R^n)\subset{\mathcal K}^2(\mathbb R^n)$, indeed the boundary of an $\varepsilon$-neighborhood of a polytope is a ${\mathcal C}^1$-hypersurface but only a piecewise ${\mathcal C}^2$-hypersurface. If a convex body admits a single supporting hyperplane at each point of its boundary, then it has non-empty interior and, for convexity reasons, this means that its boundary is a ${\mathcal C}^1$-hypersurface. 

Given a fixed $\Gamma\in{\mathcal P}(\mathbb R^n)$ admitting the origin as an interior point, a global defining function of $(\Gamma)_\varepsilon$ can be described in simple combinatorial terms. Observer first of all that the subsets $\Delta+\overline{K_\Delta}$, as $\Delta$ runs in the set of non-empty faces of $\Gamma$, comprise a finite closed cover of $\mathbb R^n$. For every non empty face $\Delta\in{\mathcal B}(\Gamma)$, let $q_\Delta:\mathbb R^n\to E_\Delta^\perp$\label{q-Delta} be the orthogonal projection and let $p_\Delta=\aff_\mathbb R \Delta\cap E_\Delta^\perp$. If $\Lambda$ is a facet of $\Delta$, then, for every $z\in\Lambda+\overline{K_\Delta}$, there is an orthogonal decomposition similar to (\ref{decomposition}), namely
\begin{align}\label{decomposition2}
q_\Lambda(z)=q_\Delta(z)+h_\Delta(u_{\Lambda,\Delta})u_{\Lambda,\Delta}.
\end{align}
Figure~\ref{decomp2} depicts the situation
\begin{figure}[ht]
\begin{center}
\begin{tikzpicture}[scale=0.85]
\draw[->] (0,0,0)--(12,0,0) node[right]{$E_\Delta^\perp$};
\draw (4,3,0) node{$E_\Lambda^\perp$};
\draw[->] (0,0,0)--(0,3,0) node[above]{$E_\Lambda^\perp\cap E_\Delta$};
\draw[->] (0,0,0)--(0,0,10) node[below left]{$E_\Lambda$};
\filldraw[draw=black,fill=brown!20!white,thick] (5,2,11)--(7,2,11)--(7,0,10)--(5,0,10)--cycle;
\filldraw[draw=black,fill=brown!20!white,thick] (5,2,11)--(7,2,11)--(7,2,6)--(5,2,6)--cycle;
\draw[dashed] (5,0,10)--(5,2,6);
\filldraw[blue!20!white] (7,2,6)--(7,2,11)--(7,0,10)--cycle;
\draw (7,2,6)--(7,2,11)--(7,0,10)--cycle;
\draw[thick, orange] (8.5,2,0)--(8.5,2,8) node[below, black]{$z$};
\fill (8.5,2,8) circle(1.5pt);
\draw[very thick, orange] (7,2,6)--(7,2,11);
\draw (7,2,9.5) node[left]{$\Lambda$};
\draw (7,0.9,8.5) node[left]{$\Delta$};
\draw[->,thick] (0,0,0)--(0,2,0) node[left]{$h_\Delta(u_{\Lambda,\Delta})u_{\Lambda,\Delta}$};
\draw[->,thick] (0,0,0)--(0,1,0) node[left]{$u_{\Lambda,\Delta}$};
\draw[dotted] (8.5,0,0)--(8.5,2,0);
\draw[dotted] (0,2,0)--(8.5,2,0);
\filldraw (5,2,11) circle(1.5pt);
\filldraw (5,2,6) circle(1.5pt);
\filldraw (7,2,11) circle(1.5pt);
\filldraw (7,2,6) circle(1.5pt);
\filldraw (5,0,10) circle(1.5pt);
\filldraw (7,0,10) circle(1.5pt);
\filldraw[blue!50!orange, nearly transparent] (7,2,6)--(13,2,6)-- (13,2,11)--(7,2,11)--cycle;
\draw (11,1.25,7) node{$\Lambda+\overline{K_\Delta}$};
\draw (-2,1,2) node{$E_\Delta$};
\draw[->,thick] (0,0,0)--(8.5,0,0) node[above right]{$q_\Delta(z)$};
\draw[thick,->] (0,0,0)--(8.5,2,0) node[above right]{$q_\Lambda(z)$};
\filldraw (8.5,0,0) circle(1.5pt);
\filldraw (8.5,2,0) circle(1.5pt);
\draw (5.25,1.5,8.5) node[left]{$\Gamma$};
\draw[blue!50!orange,->,thick] (0,0,0)--(8.5,2,8);
\draw[orange,->,thick] (0,0,0)--(7,2,8);
\draw[blue,->,thick] (0,0,0)--(1.5,0,0);
\draw[blue,->,thick] (7,2,8)--(8.5,2,8);
\draw[dotted] (8.5,2,8)--(8.5,0,8)--(8.5,0,0);
\draw[dotted] (8.5,0,8)--(0,0,8);
\end{tikzpicture}
\caption{$\Gamma$ is the light brown polytope with the blue face $\Delta$ and the orange edge $\Lambda$. The violet, transparent half strip is $\Lambda+\overline{K_\Delta}$. }
\label{decomp2}
\end{center}
\end{figure}

For $z\in \mathbb R^n$, set 
\begin{equation}\label{funzione-def}
\rho_{(\Gamma)_\varepsilon}(z)=-\varepsilon^2+\sum_{\varnothing\neq\Delta\in{\mathcal B}(\Gamma)}\Vert q_\Delta(z)-p_\Delta\Vert^2\chi_{\Delta+\overline{K_\Delta}}(z).
\end{equation}
Observe first of all that, for every $z\in \Gamma$, $q_\Gamma(z)=p_\Gamma$, so $\rho_{(\Gamma)_\varepsilon}(z)=-\varepsilon^2$. 
In order to prove that $\rho_{(\Gamma)_\varepsilon}$ is well defined, let $z\in\mathbb R^n$ belong to $(\Delta_1+\overline{K_{\Delta_1}})\cap (\Delta_2+\overline{K_{\Delta_2}})$, for some non empty faces $\Delta_1$ and $\Delta_2$. If $\Lambda=\Delta_1\cap\Delta_2$, then there are two possibilities: $\Lambda$ is strictly included in both $\Delta_1$ and $\Delta_2$, or $\Lambda$ coincides with one of them, say $\Delta_1$. In the first case $z\in \Lambda$ so   $q_\Lambda(z)=p_\Lambda$, $q_{\Delta_1}(z)=p_{\Delta_1}$ and  $q_{\Delta_2}(z)=p_{\Delta_2}$. In the second case $z\in (\Lambda+\overline{K_{\Delta_2}})\setminus\Lambda$ and now $q_\Lambda(z)\neq p_\Lambda$, however, by (\ref{decomposition}) and (\ref{decomposition2}), $q_\Lambda(z)=q_{\Delta_2}(z)+h_{\Delta_2}(u_{\Lambda,\Delta_2})u_{\Lambda,\Delta_2}$ and $p_\Lambda=p_{\Delta_2}+h_{\Delta_2}(u_{\Lambda,\Delta_2})u_{\Lambda,\Delta_2}$, so
$
q_\Lambda(z)-p_\Lambda
=
q_{\Delta_2}(z)-p_{\Delta_2}
$.
In any case, the preceding discussion shows that the restrictions of $\rho_{(\Gamma)_\varepsilon}$ to the intersection of two overlapping closed subsets of the form $\Delta+\overline{K_\Delta}$ coincide, whence $\rho_{(\Gamma)_\varepsilon}$ is a well defined continuous function on the whole $\mathbb R^n$. Its convexity and differentiability follow at once from the definition so $\rho_{(\Gamma)_\varepsilon}\in{\mathcal C}^1(\mathbb R^n,\mathbb R)$ but it is only piecewisely smooth since, for every pair of distinct faces $\Delta_1$ and $\Delta_2$, it does not admit second derivatives on $(\Delta_1+\overline{K_{\Delta_1}})\cap (\Delta_2+\overline{K_{\Delta_2}})$ .

For any given $A\in{\mathcal K}(\mathbb R^n)$ recall that the \textit{subdifferential} $\Subd h_A(u)$\label{subd} of $h_A$ at a point $u\in\mathbb R^n$ is defined as
\begin{equation*}
\Subd h_A(u)=\{v\in\mathbb R^n\mid ( v,z-u)\leq h_A(z)-h_A(u)\,,\;\forall\,z\in\mathbb R^n\}\,. 
\end{equation*}
It is a non-empty set as long as $A$ is such, moreover it equals the exposed face $A\cap H_A(u)$. Its elements are called the \textit{subgradients} of $h_A$ at $u$ and $h_A$ is differentiable at $u$ if and only if $\#\Subd h_A(u)=1$. For convexity reasons, the condition $\#\Subd h_A(u)=1$ implies that $h_A$ is ${\mathcal C}^1$ at $u$. Anyway, $h_A$ will never be differentiable at the origin, unless $A$ is reduced to a single point.

A convex body is \textit{strictly convex} if its boundary contains no line segment. In particular such a body has non-empty interior. The support function of a convex body $A$ is differentiable on $\mathbb R^n\setminus\{0\}$ if and only if $h_A\in{\mathcal C}^1(\mathbb R^n\setminus\{0\},\mathbb R)$, or equivalently, if and only if $A$ is strictly convex, (cf.~\cite{Sch}). The class of strictly convex bodies of $\mathbb R^n$ will be denoted ${\mathcal K}_1(\mathbb R^n)$.\label{strict} 
Observe that ${\mathcal K}_1(\mathbb R^n)\neq{\mathcal K}^2(\mathbb R^n)$, indeed the convex bodies from the class ${\mathcal P}_\circ(\mathbb R^n)$ belong to ${\mathcal K}^2(\mathbb R^n)$ but they not strictly convex, whereas the planar convex body depicted in Fig.~\ref{stretto} is strictly convex but it does not belong to ${\mathcal K}^2(\mathbb R^2)$ as its boundary is not ${\mathcal C^1}$. Although ${\mathcal K}_1(\mathbb R^n)\neq {\mathcal K}^2(\mathbb R^n)$, taking $\varepsilon$-neighborhood is a convenient way to make ${\mathcal C^1}$ boundaries piecewise ${\mathcal C^2}$.
\begin{figure}[ht]
\centering
\begin{tikzpicture}[scale=0.4]
\filldraw[very thick,fill=blue!20!white, draw=orange] (4.5,4) cos (10.5,0) sin (4.5,-4) arc (270:90:4);
\filldraw  (10.5,0) circle(2pt);
\end{tikzpicture}
\caption{A strictly convex body which is not in ${\mathcal K^2}(\mathbb R^2)$ because of the black point.}
\label{stretto}
\end{figure}

If $\ell\in\mathbb N^*\cup\{\infty\}$, ${\mathcal K}_\ell(\mathbb R^n)$\label{ell-strict} will denote the class of convex bodies of $\mathbb R^n$ whose support function belongs to ${\mathcal C}^\ell(\mathbb R^n\setminus\{0\})$, we also set ${\mathcal K}^2_\ell(\mathbb R^n)={\mathcal K}^2(\mathbb R^n)\cap{\mathcal K}_\ell(\mathbb R^n)$\label{2-ell-strict}.

The \textit{Hausdorff metric} on ${\mathcal C}(\mathbb R^n)$ is the mapping
${\mathfrak h}:{\mathcal C}(\mathbb R^n)\times{\mathcal C}(\mathbb R^n)\rightarrow\mathbb R_{\geq0}
$\label{hausm}
defined, for any $A_1,A_2\in{\mathcal C}(\mathbb R^n)$, as ${\mathfrak h}(A_1,A_2)=\inf\{\varepsilon\in\mathbb R_{\geq 0}\mid A_1\subseteq(A_2)_\varepsilon\;{\rm and}\; A_2\subseteq(A_1)_\varepsilon\}$.
The Hausdorff metric gives ${\mathcal C}(\mathbb R^n)$ the structure of a locally compact metric space in which ${\mathcal K}(\mathbb R^n)$ is a closed subspace. The subspace ${\mathcal P}(\mathbb R^n)$ (and hence ${\mathcal P}^\infty(\mathbb R^n)$) is dense in ${\mathcal K}(\mathbb R^n)$. Moreover, any $m$-tuple of convex bodies $A_1,\ldots,A_m\in {\mathcal K}(\mathbb R^n)$ can be respectively approximated by an $m$-tuple of sequences of polytopes $\Gamma_{1,k},\ldots,\Gamma_{m,k}\in{\mathcal P}(\mathbb R^n)$ such that, for every fixed $k\in \mathbb N$, the polytopes $\Gamma_{1,k},\ldots,\Gamma_{m,k}$ are strongly combinatorially isomorphic to one-another.

The convergence of a sequence $A_m\to A$ in ${\mathcal K}(\mathbb R^n)$ is equivalent to the uniform convergence $h_{A_m}\to h_A$ on compacta of $\mathbb R^n$ of the corresponding sequence of support functions. By homogeneity, $h_{A_m}\to h_A$ uniformly on compacta of $\mathbb R^n$ if and only if $h_{A_m}\to h_A$ uniformly on the single compactum $\partial B_n$. Moreover, (for convexity reasons) $h_{A_m}\to h_A$ uniformly on $\partial B_n$ if and only if $h_{A_m}\to h_A$ converges point-wise on $\partial B_n$, (cf.~\cite{Sch}). By a regularization argument, one can show that the subset ${\mathcal K}_2(\mathbb R^n)\subset {\mathcal K}(\mathbb R^n)$ is dense in ${\mathcal K}(\mathbb R^n)$ for the Hausdorff metric. Indeed, for any convex body $A$ and every $\varepsilon\in\mathbb R_{\geq 0}$, there exist an \textit{$\varepsilon$-regularization}\label{ereg} of $A$, i.e. a convex body $R_\varepsilon A\in{\mathcal K}_\infty(\mathbb R^n)\subset{\mathcal K}_2(\mathbb R^n)$ such that $\lim_{\varepsilon\to 0}R_\varepsilon A=A$. The support function of $R_\varepsilon A$ is given by
\begin{equation*}
h_{R_\varepsilon A}(x)
=
\int_{\mathbb R^n} h_A(x+u \Vert x\Vert )\varphi_\varepsilon(\Vert u\Vert)\,\upsilon_n 
\,,
\end{equation*}
where $\upsilon_n=\de u_1\wedge\ldots\wedge\de u_n$\label{vformn} is the usual standard volume form on $\mathbb R^n$ and 
$\varphi_\varepsilon\in{\mathcal C}^\infty(\mathbb R,\mathbb R)$\label{ecutoff} is a non-negative function with (compact) support included in the interval $[\varepsilon/2,\varepsilon]$ such that 
$
\int_{\mathbb R^n} \varphi_\varepsilon(\Vert u\Vert)\,\upsilon_n=1\,.
$
A convenient choice for $\varphi_\varepsilon$ is given by the function
\begin{equation*}
\varphi_\varepsilon(x)
=
m_\varepsilon^{-1}\chi_{[\varepsilon/2,\varepsilon]}(x)\exp\left[\frac{\varepsilon^2}{(4x-3\varepsilon)^2-\varepsilon^2}\right]
\end{equation*}
where
\begin{equation*}
m_\varepsilon=\int_{\mathbb R^n}
\chi_{[\varepsilon/2,\varepsilon]}(\Vert u\Vert)\exp\left[\frac{\varepsilon^2}{(4\Vert u\Vert-3\varepsilon)^2-\varepsilon^2}\right]\,\upsilon_n
\quad{\rm and}\quad\chi_{[\varepsilon/2,\varepsilon]}(x)
=
\begin{cases}
1 & \text{if $x\in[\varepsilon/2,\varepsilon]$,}
\\
0 &\text{if $x\not\in[\varepsilon/2,\varepsilon]$.}
\end{cases}
\end{equation*}
An alternative formula for $h_{R_\varepsilon A}(x)$ is given by 
\begin{equation*}
h_{R_\varepsilon A}(x)
=
\Vert x\Vert\int_{\mathbb R^n} h_A(u)\varphi_\varepsilon\left(\left\Vert u-\frac{x}{\Vert x\Vert}\right\Vert\right)\,\upsilon_n 
\,,
\end{equation*}
from which it follows easily that $h_{R_\varepsilon A}$ is differentiable on $\mathbb R^n\setminus\{0\}$. Observe that, unless $A$ is reduced to a single point, the regularization $h_{R_\varepsilon A}$ is never differentiable at the origin. In order to get a regularized version of $h_A$ which is differentiable on the whole $\mathbb R^n$ one has to drop the requirement of $1$-homogeneity and perform the usual convolution with some regularizing kernel, for example
\begin{equation*}
h_A*\varphi_\varepsilon(x)=\int_{\mathbb R^n}h_A(u)\varphi_\varepsilon(\Vert u-x\Vert )\upsilon_n\,,
\end{equation*}
thus getting, for every $\varepsilon>0$, a smooth function on the whole $\mathbb R^n$ that converges uniformly to $h_A$ on each compact set, as $\varepsilon\to0$. Such a regularization is generally not $1$-homogeneous. The equality
\begin{equation*}
\int_{\mathbb R^n}u_\ell\, \varphi_\varepsilon(\Vert u\Vert)\upsilon_n=0\,
\end{equation*}
for every $1\leqslant\ell\leqslant n$, has two worth-noting consequences: both regularizations converge monotonically to $h_A$ from above; moreover, if $h_A$ is linear on some open subset $\Omega\subset\mathbb R^n$, then 
\begin{equation*}
h_A(x)=h_{R_\varepsilon A}(x)=h_A*\varphi_\varepsilon(x)\,,
\end{equation*}
for any $x\in\Omega\setminus(\partial \Omega)_\varepsilon$.

The $\varepsilon$-neighborhood of $R_\varepsilon A$ belongs to ${\mathcal K}_\infty^2(\mathbb R^n)$ and still converges to $A$, (as $\varepsilon\to0$), so that ${\mathcal K}_\infty^2(\mathbb R^n)$ is a dense subspace of ${\mathcal K}(\mathbb R^n)$, in particular the subspace ${\mathcal K}^2(\mathbb R^n)$ is dense too.

We state the following important result.

\begin{thm}\label{minkowski}
Let~$r\in\mathbb N^*$,~$A_1,\ldots,A_r\in{\mathcal K}(\mathbb R^n)$ and~$\lambda_1,\ldots,\lambda_r\in\mathbb R_{\geq 0}$.  Then $\vol_n\big(\sum_{\ell=1}^r\lambda_\ell A_\ell\big)$ is a homogeneous polynomial of degre~$n$ in $\lambda_1,\ldots,\lambda_r$.
\findim
\end{thm}
This polynomial can be written as
\begin{equation}\label{mink}
\vol_n
\bigg(\sum_{\ell=1}^r\lambda_\ell A_\ell\bigg)
=
\sum_{\rho}
V_n(A_{\rho(1)},\ldots,A_{\rho(n)})
\lambda_{\rho(1)}\cdots\lambda_{\rho(n)}\,,
\end{equation}
where the sum runs on the set of functions~$\rho:\{1,\ldots,n\}\to\{1,\ldots,r\}$ and the coefficients are chosen so as to satisfy, (for any $\rho$), the condition
$$
V_n(A_{\rho(1)},\ldots,A_{\rho(n)})
=
V_n(A_{\varsigma(\rho(1))},\ldots,A_{\varsigma(\rho(n))})\,,
$$
for every permutation~$\varsigma$ of the image of~$\rho$. Notice that, for fixed~$\rho$, the coefficient of the monomial~$\lambda_{\rho(1)}\cdots\lambda_{\rho(n)}$ in the expression~\eqref{mink} depends just on~$A_{\rho(1)},\ldots,A_{\rho(n)}$ as follows by setting~$\lambda_\ell=0$ for every~$\ell$ which does not belong to the image of~$\rho$. If~$n\in\mathbb N^*$ and~$A_1,\ldots,A_n\in{\mathcal K}(\mathbb R^n)$, the $n$-dimensional \textit{(Minkowski) mixed volume\/}~$V_n(A_1,\ldots,A_n)$\label{nvmix} of the convex bodies~$A_1,\ldots,A_n$ is the coefficient in the expression~\eqref{mink} when~$r=n$ and~$\rho$ is the identity. 

For any $k\in\mathbb N$ and any convex body $A\in{\mathcal K}(\mathbb R^n)$, for the sake of notation, we set
\begin{equation*}\label{aktimes}
A[k]
=
\underbrace{A,\ldots,A}_{k\,{\rm times}}\,,
\end{equation*}
meaning that, for $k=0$, the set $A$ is not written at all. For any integer $0\leq k\leq n$, the \textit{$k$-th intrinsic volume} of a convex body $A\in{\mathcal K}(\mathbb R^n)$ is the non-negative real number $\textsl{v}_k(A)$ defined as 
\begin{equation}\label{kth}
\textsl{v}_k(A)=\varkappa_{n-k}^{-1}{n\choose k}V_n(A[k],B_n[n-k])\,.
\end{equation}
For $k=0$ one has $\textsl{v}_0(A)=1$. The number $\textsl{v}_k(A)$ does not depend on the ambient space but just on $A$. For a polytope $\Gamma$, the intrinsic character of $\textsl{v}_k$ is revealed by the equality
\begin{equation*}\label{vintrk}
\textsl{v}_k(\Gamma)=\sum_{\Delta\in{\mathcal B}(\Gamma,k)}\vol_k(\Delta)\psi_\Gamma(\Delta)\,.
\end{equation*}
Notice that the $k$-th intrinsic volume of a $k$-polytope is just its $k$-dimensional Lebesgue measure. Moreover the $k$-th intrinsic volume of a $(k+1)$-polytope equals half the $k$-dimensional volume of its relative boundary, in particular the $1$-st intrinsic volume of a $2$-polytope is the semi-perimeter, whereas the $2$-nd intrinsic volume of a $3$-polytope is half the surface area of its relative boundary. By continuity of intrinsic volumes for the Hausdorff metric, the preceding remarks apply to any convex body.

The notion of valuation on convex bodies will also play an important role in the rest of the paper, so we briefly recall it here.
Let $\mathbb G$ be an abelian group, a \textit{valuation} on ${\mathcal K}(\mathbb R^n)$ with values in $\mathbb G$ is a  mapping $\phi:{\mathcal K}(\mathbb R^n)\cup\{\varnothing\}\to\mathbb G$ such that $\phi(\varnothing)=0$ and $\phi(A_1)+\phi(A_2)=\phi(A_1\cup A_2)+\phi(A_1\cap A_2)$, for every $A_1,A_2\in {\mathcal K}(\mathbb R^n)$ such that $A_1\cup A_2$ is convex. 
If $m\geq 2$, let $A_1,\ldots,A_m\in{\mathcal K}(\mathbb R^n)$ be such that $A_1\cup\ldots\cup A_m$ is convex too. For every non empty subset $I\subseteq\{1,\ldots,m\}$, set 
$$
A_I=\bigcap_{j\in I} A_j\,.
$$
A valuation $\phi$ satisfies the inclusion-exclusion property if 
\begin{equation*}\label{val-I}
\phi(A_1\cup\ldots\cup A_m)
=
\sum_{\varnothing \neq I\subseteq\{1,\ldots,m\}}
(-1)^{\# I -1}
\phi(A_I)\,.
\end{equation*} 

A \textit{weak valuation} on ${\mathcal K}(\mathbb R^n)$ is a  mapping $\phi:{\mathcal K}(\mathbb R^n)\cup\{\varnothing\}\to\mathbb G$ such that $\phi(\varnothing)=0$ and 
\begin{equation*}
\phi(A\cap H^+)+\phi(A\cap H^-)=\phi(A\cap H)+\phi(A)\,,
\end{equation*}
for every $A\in {\mathcal K}(\mathbb R^n)$ and every hyperplane $H\subset\mathbb R^n$, where $H^+$ and $H^-$ are the closed half-spaces bounded by $H$. If $\mathbb G$ has a topology, (weak) valuations can be continuous with respect to Hausdorff metric. A (weak) valuation $\phi$ is $k$-\textit{homogeneous} if $\phi(\lambda A)=\lambda^k\phi(A)$, for every $\lambda\in\mathbb R_{>0}$ and $A\in{\mathcal K}(\mathbb R^n)$. A (weak) valuation $\phi$ is \textit{simple} if $\phi(A)=0$ for every $A\in{\mathcal K}(\mathbb R^n)$ which is not full-dimensional. These notions can be also defined on proper subsets of ${\mathcal K}(\mathbb R^n)$, such as ${\mathcal K}^k(\mathbb R^n)$, ${\mathcal K}_\ell(\mathbb R^n)$, ${\mathcal K}_\ell^k(\mathbb R^n)$ or ${\mathcal P}(\mathbb R^n)$. The following results on valuations are of great interest.
\begin{thm}[\cite{Sal}]\label{Sal}
Every weak valuation on ${\mathcal P}(\mathbb R^n)$ with values in a Hausdorff topological $\mathbb R$-linear space is a valuation.\findim
\end{thm}
\begin{thm}[\cite{Gro}]\label{Gro}
Every continuous weak valuation on ${\mathcal K}(\mathbb R^n)$ with values in a Hausdorff topological $\mathbb R$-linear space is a valuation and satisfies the inclusion-exclusion property.
\findim
\end{thm}

In the following sections we will define the $n$-dimensional Kazarnovski{\v\i} pseudovolume, a notion which generalizes the $n$-dimensional mixed volume to the complex setting. In order to compare these two notions, it is convenient to list here the most important properties of $V_n$. 
\begin{enumerate}
\item 
\textit{Non negativity:} $V_n$ is a function which takes on non-negative real values\,.
\item \textit{Continuity:} The~$n$-dimensional mixed volume is continuous in the Hausdorff metric\,.
\item 
\textit{Symmetry:} $V_n$ is a symmetric function of its arguments\,.
\item 
\textit{Multilinearity:} $V_n$ is multilinear with respect to Minkowski sum and non-negative real multiplication\,.
\item 
\textit{Diagonal property:} $V_n(A[n])=\vol_n(A)\,,$ for every convex body~$A\in{\mathcal K}(\mathbb R^n)$\,.
\item
\textit{Translation invariance:} 
\begin{equation*}
V_n(A_1,\ldots,A_n)=V_n(x_1+A_1,\ldots,x_n+A_n)\,,
\end{equation*}
for any~$A_1,\ldots,A_n\in{\mathcal K}(\mathbb R^n)$ and~$x_1,\ldots,x_n\in\mathbb R^n$.
\item
\textit{Orthogonal invariance:} 
\begin{equation*}
V_n(A_1,\ldots,A_n)=V_n(F(A_1),\ldots,F(A_n))\,,
\end{equation*}
for any~$A_1,\ldots,A_n\in{\mathcal K}(\mathbb R^n)$ and any orthogonal transformation~$F:\mathbb R^n\rightarrow\mathbb R^n$.
\item
\textit{Symmetric formula:} 
\begin{equation*}
V_n(A_1,\ldots,A_n)
=
\frac{1}{n!}\,\frac{\partial^n}{\partial \lambda_1\ldots\partial\lambda_n}\,
\vol_n\left(\sum_{\ell=1}^n\lambda_\ell A_\ell\right)
\,,
\end{equation*}
for any~$A_1,\ldots,A_n\in{\mathcal K}(\mathbb R^n)$.
\item
\textit{Polarization formula:} 
\begin{equation*}
V_n(A_1,\ldots,A_n)
=
\frac{1}{n!}
\sum_{\varnothing\neq I\subseteq\{1,\ldots,n\}}
(-1)^{n-\# I}
\vol_n
\bigg(
\sum_{\ell\in I}
A_\ell
\bigg)
\,,
\end{equation*}
for any~$A_1,\ldots,A_n\in{\mathcal K}(\mathbb R^n)$.
\item 
\textit{Non-degeneracy:} Let~$A_1,\ldots,A_n\in{\mathcal K}(\mathbb R^n)$, then $V_n(A_1,\ldots,A_n)>0$ if and only if there exist segments $\Lambda_1\subseteq A_1,\ldots,\Lambda_n\subseteq A_n$ spanning lines with linearly independent directions.
\item 
\textit{Monotonicity:} The~$n$-dimensional mixed volume is non-decreasing, with respect to inclusion, in each of its arguments.
\item
\textit{Recursive formula:}  
If $\Gamma_1,\ldots,\Gamma_n\in{\mathcal P}(\mathbb R^n)$ are full-dimensional and~$\Gamma^{(2)}=\sum_{\ell=2}^{n}\Gamma_\ell$, then
\begin{equation}\label{inductive}
V_n(\Gamma_1,\ldots,\Gamma_{n})
=
\frac{1}{n}\sum_{\Delta\in{\mathcal B}(\Gamma^{(2)},n-1)}
h_{\Gamma_1}(u_\Delta)V_{n-1}(\Delta_2,\ldots,\Delta_{n})\,,
\end{equation}
where, for any~$\Delta\in{\mathcal B}(\Gamma^{(2)},n-1)$,~$\Delta_2,\ldots,\Delta_{n}$ is the unique sequence of faces summing up to~$\Delta$, $V_{n-1}(\Delta_2,\ldots\Delta_{n})$ is the $(n-1)$-dimensional mixed volume of the projections of $\Delta_2,\ldots,\Delta_{n}$ onto $E_\Delta\simeq\mathbb R^{n-1}$ and $u_\Delta\in E_\Delta^\perp\cap E_{\Gamma^{(2)}}$ is the outer unit normal vector to the facet $\Delta$ of $\Gamma^{(2)}$.
\item 
\textit{Mixed discriminant formula:}
If $A_1,\ldots,A_n\in{\mathcal K}_2(\mathbb R^n)$ then
\begin{equation}\label{mix-D-f}
V_n(A_1,\ldots,A_n)
=
\int_{\partial B_n}
h_{A_1}D(\Hess_\mathbb R h_{A_2},\ldots,\Hess_\mathbb R h_{A_n}) \upsilon_{\partial B_n}\,,
\end{equation}
where, for every $u\in\partial B_n$\label{vformpB} and every $2\leq\ell\leq n$, $\Hess_\mathbb R h_{A_\ell}$\label{real-hess} is the (real) Hessian matrix of $h_{A_\ell}$ at $u$, $\upsilon_{\partial B_n}$ is the standard volume form on $\partial B_n$ and $D(\Hess_\mathbb R h_{A_2},\ldots,\Hess_\mathbb R  h_{A_n})$\label{md} (cf. section \ref{kmd}), equals $1/n!$ times the coefficient of $\lambda_1\lambda_2\ldots\lambda_n$ in the polynomial expression given by $\det(\lambda_1 I_n+\lambda_2\Hess_\mathbb R h_{A_2}+\ldots+\lambda_n\Hess_\mathbb R  h_{A_n})$.
\item
\textit{Intrinsic formula:} 
If~$\Gamma_1,\ldots,\Gamma_k\in{\mathcal P}(\mathbb R^n)$ and~$\Gamma=\sum_{\ell=1}^k\Gamma_\ell$, then
\begin{equation}\label{comb}
\varkappa_{n-k}^{-1}
{n\choose k}
V_n(\Gamma_1,\ldots,\Gamma_k,B_n[n-k])
=
\sum_{\Delta\in{\mathcal B}(\Gamma,k)}
V_k(\Delta_1,\ldots,\Delta_k)
\psi_{\Gamma}(\Delta),
\end{equation}
where, for any~$\Delta\in{\mathcal B}(\Gamma,k)$,~$\Delta_1,\ldots,\Delta_k$ is the unique sequence of faces summing up to~$\Delta$, and~$V_k(\Delta_1,\ldots\Delta_k)$ is the $k$-dimensional mixed volume of the projections of $\Delta_1,\ldots,\Delta_k$ onto $E_\Delta\simeq\mathbb R^k$. 
\item \textit{Rationality on lattice polytopes:} If $\Gamma_1,\ldots,\Gamma_n\in{\mathcal P}(\mathbb R^n)$ are lattice polytopes, then 
\begin{equation*}
n! V_n(\Gamma_1,\ldots,\Gamma_n)\in\mathbb Z\,.
\end{equation*}
\item \textit{Valuation property:} Let $n\in\mathbb N\setminus\{0,1\}$, $m\in\{1,\ldots,n-1\}$ and let $A_{m+1},\ldots,A_n\in{\mathcal K}(\mathbb R^n)$ be fixed. Then the mapping $g_m:{\mathcal K}(\mathbb R^n)\to\mathbb R$, given for any $A\in{\mathcal K}(\mathbb R^n)$ by
\begin{equation*}
g_m(A)=V_n(A[m],A_{m+1},\ldots,A_n)\,,
\end{equation*} 
is a valuation.
\item
\textit{Alexandroff-Fenchel inequality:} If $n\geq2$, for any $A_1,A_2,A_3,\ldots,A_n\in{\mathcal K}(\mathbb R^n)$,
\begin{equation}\label{AF}
V_n(A_1,A_2,A_3,\ldots,A_n)^2\geq V_n(A_1[2],A_3,\ldots,A_n)V_n(A_2[2],A_3,\ldots,A_n)\,.
\end{equation}
\end{enumerate}

The interested reader is referred to \cite{Sch}, \cite{Ewa}, \cite{Sal} and \cite{Gro} for further details on the general theory of convex bodies as well as for the proofs of the results recalled so far.

\section{Mixed {\boldmath $\varphi$}-volume}\label{folume}

Let $\varphi$ be an arbitrary, but fixed, real function defined on the Grassmannian ${\mathcal G}(\mathbb R^n)$. On the basis of \cite{Ka1} we propose the following definition.
\begin{defn}
For any integer $0\leq k\leq n$, the \textit{$k$-th intrinsic $\varphi$-volume} of a polytope $\Gamma\in{\mathcal P}(\mathbb R^n)$ is the real number $\textsl{v}_k^\varphi(A)$ defined as 
\begin{equation}\label{vintrfk}
\textsl{v}_k^\varphi(\Gamma)=\sum_{\Delta\in{\mathcal B}(\Gamma,k)}\varphi(E_\Delta)\vol_k(\Delta)\psi_\Gamma(\Delta)\,.
\end{equation}
\end{defn}
Of course $\textsl{v}_k^\varphi(\Gamma)=0$ as soon as $\dim \Gamma<k$. 
For $k=0$ (resp. $k=n$) one has $\textsl{v}_0^\varphi(\Gamma)=\varphi(\{0\})$ (resp. $\textsl{v}_n^\varphi(\Gamma)=\varphi(E_\Gamma)\vol_n(\Gamma)$). Setting $\textsl{v}_k^\varphi(\varnothing)=0$, the $k$-th intrinsic $\varphi$-volume is $k$-homogeneous, i.e. $\textsl{v}_k^\varphi(\lambda\Gamma)=\lambda^k\textsl{v}_k^\varphi(\Gamma)$, for every $\lambda\in\mathbb R_{\geq0}$. 
The $k$-th intrinsic $\varphi$-volume is just a weighted version of the usual $k$-th intrinsic volume, if $\varphi\equiv 1$, $\textsl{v}_k^\varphi(\Gamma)=\textsl{v}_k(\Gamma)$.
\begin{lem}\label{fi-misto}
Let~$r\in\mathbb N^*$, $\Gamma_1,\ldots,\Gamma_r\in{\mathcal P}(\mathbb R^n)$ and~$\lambda_1,\ldots,\lambda_r\in\mathbb R_{\geq 0}$.  Then $\textsl{v}_k^\varphi\big(\sum_{\ell=1}^r\lambda_\ell \Gamma_\ell\big)$ is a homogeneous polynomial of degree~$k$ in~$\lambda_1,\ldots,\lambda_r$.
\end{lem}
\noindent
\pf Let $\Gamma=\sum_{\ell=1}^r\Gamma_\ell$ and $\Gamma^\prime=\sum_{\ell=1}^r\lambda_\ell\Gamma_\ell$. If $\lambda_1,\ldots,\lambda_r$ are all positive, then there is a bijection from ${\mathcal B}(\Gamma,k)$ to ${\mathcal B}(\Gamma^\prime,k)$ mapping the face $\Delta=\sum_{\ell=1}^r\Delta_\ell$ of $\Gamma$ to the face $\Delta^\prime=\sum_{\ell=1}^r \lambda_\ell\Delta_\ell$ of $\Gamma^\prime$. As a consequence $E_\Delta=E_{\Delta^\prime}$ and $K_{\Delta,\Gamma}=K_{\Delta^\prime,\Gamma^\prime}$, so that
\begin{align*}
\textsl{v}_k^\varphi(\Gamma^\prime)
&=
\sum_{\Delta^\prime\in{\mathcal B}(\Gamma^\prime,k)}\varphi(E_{\Delta^\prime})\vol_k(\Delta^\prime)\psi_{\Gamma^\prime}(\Delta^\prime)
\\
&=
\sum_{\Delta\in{\mathcal B}(\Gamma,k)}\varphi(E_\Delta)\vol_k\left(\sum_{\ell=1}^r \lambda_\ell\Delta_\ell\right)\psi_\Gamma(\Delta)
\,.
\end{align*}
For each $\Delta^\prime\in{\mathcal B}(\Gamma^\prime,k)$, up to a translation carrying $\Delta^\prime$ in $E_{\Delta^\prime}=E_\Delta$, we can use \thmref{minkowski} in $E_\Delta\simeq\mathbb R^k$ thus getting the equality
\begin{equation*}
\vol_k\left(\sum_{\ell=1}^r \lambda_\ell\Delta_\ell\right)
=
\sum_\rho V_n\left(\Delta_{\rho(1)},\ldots,\Delta_{\rho(k)}\right)\lambda_{\rho(1)}\cdots\lambda_{\rho(k)}
\end{equation*}
with $\rho$ running in the set of all the functions $\{1,\ldots,k\}\to\{1,\ldots,r\}$ and the coefficients chosen so as to satisfy, (for any $\rho$), the condition
\begin{equation*}
V_k(\Gamma_{\rho(1)},\ldots,\Gamma_{\rho(n)})
=
V_k(\Gamma_{\varsigma(\rho(1))},\ldots,\Gamma_{\varsigma(\rho(n))})\,,
\end{equation*}
for every permutation~$\varsigma$ of the image of~$\rho$. It follows that
\begin{equation}\label{fi-mix}
\textsl{v}_k^\varphi(\Gamma^\prime)
=
\sum_\rho\left(
\sum_{\Delta\in{\mathcal B}(\Gamma,k)}\varphi(E_\Delta)
V_k\left(\Delta_{\rho(1)},\ldots,\Delta_{\rho(k)}\right)
\psi_\Gamma(\Delta)
\right)
\lambda_{\rho(1)}\cdots\lambda_{\rho(k)}.
\end{equation} 
If there is a $j\in\{1,\ldots,r\}$ such that $\lambda_j=0$, then it is enough to drop $\Gamma_j$ in the sum $\Gamma$.\findim

\begin{defn}
For any integer $1\leq k\leq n$, and polytopes $\Gamma_1,\ldots,\Gamma_k\in{\mathcal P}(\mathbb R^n)$, the \textit{$k$-th mixed $\varphi$-volume} of $\Gamma_1,\ldots,\Gamma_k$ is the coefficient $V_k^\varphi(\Gamma_1,\ldots,\Gamma_k)$ of $\lambda_1\cdots\lambda_k$ in the polynomial expression \eqref{fi-mix}  when~$r=k$ and~$\rho$ is the identity, i.e.
\begin{equation}\label{fi-mixed}
V_k^\varphi(\Gamma_1,\ldots,\Gamma_k)
=
\sum_{\Delta\in{\mathcal B}(\Gamma,k)}\varphi(E_\Delta)
V_k\left(\Delta_1,\ldots,\Delta_k\right)
\psi_\Gamma(\Delta)
\,,
\end{equation} 
with $\Gamma=\sum_{\ell=1}^k\Gamma_\ell$.
\end{defn}

\begin{rmq}
By the non-degeneracy property of $V_k$, the sum in \eqref{fi-mixed} can be performed just on those edge-sum $k$-faces $\Delta\preccurlyeq\sum_{\ell=1}^k\Gamma_k$ such that $\varphi(E_\Delta)\neq0$.
\end{rmq}
\begin{lem}\label{fi-sim}
The $k$-th mixed $\varphi$-volume has the following properties: 
\begin{enumerate}
\item[{$(a)$}] symmetry, i.e. 
$
V_k^\varphi(\Gamma_1,\ldots,\Gamma_k)=V_k^\varphi(\Gamma_{\varsigma(1)},\ldots,\Gamma_{\varsigma(k)})\,,
$
for every sequence $\Gamma_1,\ldots,\Gamma_k\in{\mathcal P}(\mathbb R^n)$ and every permutation $\varsigma$ of $\{1,\ldots,n\}$;
\item[{$(b)$}] diagonal property, i.e. $
\V^\varphi_k(\Gamma[k])=\textsl{v}_k^\varphi(\Gamma)\,,
$
for every $\Gamma\in{\mathcal P}(\mathbb R^n)$;
\item[{$(c)$}] translation invariance.
\end{enumerate}
\end{lem}
\noindent
\pf The claim $(a)$ follows at once from the second of \eqref{utili}, the second of \eqref{utili2} and the definition of $k$-dimensional mixed $\varphi$-volume. 
The statement $(b)$ is a consequence of the diagonal property of $V_k$ and the definition of $V_k^\varphi$ and $\textsl{v}_k^\varphi$. 
The claim $(c)$ follows from the translation invariance of $\varphi$, $V_k$ and the outer angle.\findim

\begin{lem}\label{pol-fi}
The $k$-th mixed $\varphi$-volume satisfies the following formula:
\begin{equation}\label{der}
V_k^\varphi(\Gamma_1,\ldots,\Gamma_k)
=
\frac{1}{k!}\,
\frac{\partial^k}{\partial\lambda_1\ldots\partial\lambda_k}\,
\textsl{v}_k^\varphi\left(\sum_{\ell=1}^k\lambda_\ell \Gamma_\ell\right)
\,,
\end{equation}
for every $\Gamma_1,\ldots,\Gamma_k\in{\mathcal P}(\mathbb R^n)$.
\end{lem}
\noindent
\pf The equality is straightforward if one realizes that the right member of \eqref{der} is the coefficient of $\lambda_1\cdots\lambda_k$ in the polynomial \eqref{fi-mix}.\findim

The idea of the proof of the following theorem comes from \cite{Th}.

\begin{thm}\label{non-sim-fi}
The $k$-th mixed $\varphi$-volume satisfies the following polarization formula:
\begin{equation}\label{pol-fi1}
V_k^\varphi(\Gamma_1,\ldots,\Gamma_k)
=
\frac{1}{k!}
\sum_{\substack{I\subseteq\{1,\ldots,k\}\\ I\neq\varnothing}}(-1)^{k-\#I} \textsl{v}_k^\varphi\left(\sum_{\ell\in I}\Gamma_\ell\right)
\,,
\end{equation}
for every $\Gamma_1,\ldots,\Gamma_k\in{\mathcal P}(\mathbb R^n)$.
\end{thm}
\noindent
\pf Consider the abelian group\footnote{It is in fact an $\mathbb R$-linear space, but we don't need this structure.} $\mathbb R^{{\mathcal P}(\mathbb R^n)}$ of real functionals on ${\mathcal P}(\mathbb R^n)$.
For every $\Gamma\in{\mathcal P}(\mathbb R^n)$, we introduce a shift operator $s_\Gamma\in \text{Hom}_{\mathbb Z}\left(\mathbb R^{{\mathcal P}(\mathbb R^n)},\mathbb R^{{\mathcal P}(\mathbb R^n)}\right)$ by setting, for every $w\in\mathbb R^{{\mathcal P}(\mathbb R^n)}$ and every $\Lambda\in{\mathcal P}(\mathbb R^n)$, $s_\Gamma w(\Lambda)=w(\Lambda+\Gamma)$. Recall that Minkowski addition makes ${\mathcal P}(\mathbb R^n)$ a commutative monoid with $\{0\}$ as zero element, where as $\text{Hom}_{\mathbb Z}\left(\mathbb R^{{\mathcal P}(\mathbb R^n)},\mathbb R^{{\mathcal P}(\mathbb R^n)}\right)$ has a ring structure.\footnote{Recall that for every abelian group $G$, the set $\text{Hom}_{\mathbb Z}(G,G)$ of its endomorphisms has a ring structure: it is an abelian group for the addition of morphisms, a monoid for their composition and these two structures are compatible.} The mapping
$$
s:{\mathcal P}(\mathbb R^n)\rightarrow \text{Hom}_{\mathbb Z}\left(\mathbb R^{{\mathcal P}(\mathbb R^n)},\mathbb R^{{\mathcal P}(\mathbb R^n)}\right)
$$
is a morphism of monoids, i.e. $s_{\{0\}}$ is the identity operator $I$ and $s_{\Gamma_1}\circ s_{\Gamma_2}=s_{\Gamma_1+\Gamma_2}$, for every $\Gamma_1,\Gamma_2\in{\mathcal P}(\mathbb R^n)$. Indeed, for every $w\in\mathbb R^{{\mathcal P}(\mathbb R^n)}$ and every $\Lambda\in{\mathcal P}(\mathbb R^n)$, 
$$
s_{\{0\}}w(\Lambda)=w(\Lambda+\{0\})=w(\Lambda)
$$ 
and 
$$s_{\Gamma_1}\circ s_{\Gamma_2}w(\Lambda)=s_{\Gamma_1}w(\Lambda+\Gamma_2)=w(\Lambda+\Gamma_2+\Gamma_1)=s_{\Gamma_1+\Gamma_2}w(\Lambda)\,,
$$ 
in particular $s_{\Gamma_1}\circ s_{\Gamma_2}=s_{\Gamma_2}\circ s_{\Gamma_1}$, i.e. 
the image of ${\mathcal P}(\mathbb R^n)$ under the morphism $s$ is a commutative submonoid of $\text{Hom}_{\mathbb Z}\left(\mathbb R^{{\mathcal P}(\mathbb R^n)},\mathbb R^{{\mathcal P}(\mathbb R^n)}\right)$. 
Consider the mapping
$$
D:{\mathcal P}(\mathbb R^n)\rightarrow \text{Hom}_{\mathbb Z}\left(\mathbb R^{{\mathcal P}(\mathbb R^n)},\mathbb R^{{\mathcal P}(\mathbb R^n)}\right)\,,
$$
whose value on $\Gamma\in{\mathcal P}(\mathbb R^n)$ is  
$D_\Gamma=s_\Gamma-s_{\{0\}}=s_\Gamma-I$, i.e. $D_\Gamma w(\Lambda)=w(\Lambda+\Gamma)-w(\Lambda)$, for every $w\in\mathbb R^{{\mathcal P}(\mathbb R^n)}$. Since $D_{\{0\}}\neq I$, $D$ is not a morphism of monoids, however the ring structure of $\text{Hom}_{\mathbb Z}\left(\mathbb R^{{\mathcal P}(\mathbb R^n)},\mathbb R^{{\mathcal P}(\mathbb R^n)}\right)$ and the commutativity of $s({\mathcal P}(\mathbb R^n))$ makes $D({\mathcal P}(\mathbb R^n))$ consist of commuting homomorphisms of $\mathbb R^{{\mathcal P}(\mathbb R^n)}$, indeed 
\begin{align*}
D_{\Gamma_1}\circ D_{\Gamma_2}
&=
(s_{\Gamma_1}-I)\circ(s_{\Gamma_2}-I)
\\
&=
s_{\Gamma_1}\circ s_{\Gamma_2}-s_{\Gamma_2}-s_{\Gamma_1}+I
\\
&=
s_{\Gamma_2}\circ s_{\Gamma_1}-s_{\Gamma_1}-s_{\Gamma_2}+I
\\
&=
(s_{\Gamma_2}-I)\circ(s_{\Gamma_1}-I)
\\
&=
D_{\Gamma_2}\circ D_{\Gamma_1}\,.
\end{align*}
An iterated integration of the equality \eqref{der} on $[0,1]^k$ yields
\begin{equation*}
V^\varphi_k(\Gamma_1,\ldots,\Gamma_k)
=
D_{\Gamma_{k}}\circ\cdots\circ D_{\Gamma_1}\textsl{v}_k^\varphi\left(\{0\}\right)\,.
\end{equation*} 
Indeed
\begin{align*}
&V^\varphi_k(\Gamma_1,\ldots,\Gamma_k)
\\&=
\int_{[0,1]^k}
\frac{1}{k!}\,
\frac{\partial^k}{\partial\lambda_1\ldots\partial\lambda_k}\,
\textsl{v}_k^\varphi\left(\sum_{\ell=1}^k\lambda_\ell \Gamma_\ell\right)
\de\lambda_1\ldots\de\lambda_k
\\
&=
\int_{[0,1]^{k-1}}
\frac{1}{k!}\,
\frac{\partial^{k-1}}{\partial\lambda_2\ldots\partial\lambda_k}\,
\left[\textsl{v}_k^\varphi\left(\Gamma_1+\sum_{\ell=2}^k\lambda_\ell \Gamma_\ell\right)
-
\textsl{v}_k^\varphi\left(\sum_{\ell=2}^k\lambda_\ell \Gamma_\ell\right)
\right]
\de\lambda_2\ldots\de\lambda_k
\\
&=
\int_{[0,1]^{k-1}}
\frac{1}{k!}\,
\frac{\partial^{k-1}}{\partial\lambda_2\ldots\partial\lambda_k}\,
\left(D_{\Gamma_1}\textsl{v}_k^\varphi\right)
\left(\sum_{\ell=2}^k\lambda_\ell \Gamma_\ell\right)
\de\lambda_2\ldots\de\lambda_k
\\
&=
\int_{[0,1]^{k-2}}
\frac{1}{k!}\,
\frac{\partial^{k-2}}{\partial\lambda_3\ldots\partial\lambda_k}\,
\left(D_{\Gamma_2}\circ D_{\Gamma_1}\textsl{v}_k^\varphi\right)
\left(\sum_{\ell=3}^k\lambda_\ell \Gamma_\ell\right)
\de\lambda_3\ldots\de\lambda_k\,,
\end{align*}
and $k-3$ more iterations yield
\begin{align*}
&V^\varphi_k(\Gamma_1,\ldots,\Gamma_k)
\\
&=
\int_{[0,1]}
\frac{1}{k!}\,
\frac{\partial}{\partial\lambda_k}\,
\left(D_{\Gamma_{k-1}}\circ\cdots\circ D_{\Gamma_1}\textsl{v}_k^\varphi\right)
\left(\lambda_k \Gamma_k\right)
\de\lambda_k
\\
&=
\frac{1}{k!}
\left[
D_{\Gamma_{k-1}}\circ\cdots\circ D_{\Gamma_1}\textsl{v}_k^\varphi\left(\Gamma_k\right)
-
D_{\Gamma_{k-1}}\circ\cdots\circ D_{\Gamma_1}\textsl{v}_k^\varphi\left(\{0\}\right)
\right]
\\
&=
\frac{1}{k!}
\left[
D_{\Gamma_{k-1}}\circ\cdots\circ D_{\Gamma_1}\textsl{v}_k^\varphi\left(\{0\}+\Gamma_k\right)
-
D_{\Gamma_{k-1}}\circ\cdots\circ D_{\Gamma_1}\textsl{v}_k^\varphi\left(\{0\}\right)
\right]
\\
&=
\frac{1}{k!}
D_{\Gamma_{k}}\circ\cdots\circ D_{\Gamma_1}\textsl{v}_k^\varphi\left(\{0\}\right)\,.
\end{align*} 
It remains to show that $D_{\Gamma_{k}}\circ\cdots\circ D_{\Gamma_1}\textsl{v}_k^\varphi\left(\{0\}\right)$ equals $k!$ times the right hand side of \eqref{pol-fi1}. For every non-empty $J\subseteq\{1,\ldots,k\}$, set $\Gamma_J=\sum_{\ell\in J} \Gamma_\ell$, then
$$
D_{\Gamma_k}\circ\cdots\circ D_{\Gamma_1}
=
\Bigcirc_{j=1}^k(s_{\Gamma_j}-I)
=
\sum_{\substack{J\subseteq\{1,\ldots,k\}\\ J\neq\varnothing}} (-1)^{k-\# J}\Bigcirc_{j\in J}s_{\Gamma_j}
=
\sum_{\substack{J\subseteq\{1,\ldots,k\}\\ J\neq\varnothing}} (-1)^{k-\# J}s_{\Gamma_J}\,.
$$
The conclusion follows because $\sigma_{\Gamma_J}\textsl{v}_k^\varphi(\{0\})=\textsl{v}_k^\varphi(\Gamma_J)$.
\findim

\begin{thm}\label{pol-fi2}
The $k$-th mixed $\varphi$-volume is $k$-linear.
\end{thm}
\noindent
\pf As $V^\varphi_k$ is symmetric, it is enough to prove the linearity of $V_k^\varphi$ in its first argument, i.e. 
\begin{equation*}
V_k^\varphi(\mu_1\Delta_1+\mu_1^\prime\Delta_1^\prime,\Gamma_2,\ldots,\Gamma_k)=\mu_1V_k^\varphi(\Delta_1,\Gamma_2,\ldots,\Gamma_k)+\mu_1^\prime V_k^\varphi(\Delta_1^\prime,\Gamma_2,\ldots,\Gamma_k),
\end{equation*}
for every $\Delta_1,\Delta_1^\prime,\Gamma_2,\ldots,\Gamma_k\in{\mathcal P}(\mathbb R^n)$ and any $\mu_1,\mu_1^\prime\in\mathbb R_{\geq0}$. Set $\Gamma_1=\mu_1\Delta_1+\mu_1^\prime\Delta_1^\prime$ and observe that $\textsl{v}_k^\varphi(\sum_{\ell=1}^k\lambda_\ell\Gamma_\ell)=\textsl{v}_k^\varphi(\lambda_1\mu_1\Delta_1+\lambda_1\mu_1^\prime\Delta_1^\prime+\sum_{\ell=2}^k\lambda_\ell\Gamma_\ell)$, for every $\lambda_1,\lambda_2,\ldots,\lambda_k>0$. Then the coefficients of $\lambda_1\lambda_2\cdots\lambda_k$ in both sides of this equality must agree.
\findim

\begin{coro}\label{uni}
The $k$-th mixed $\varphi$-volume $V_k^\varphi$ is the only $k$-linear symmetric function which agrees with $\textsl{v}_k^\varphi$ on the diagonal.
\end{coro}
\noindent
\pf Let $L$ be a $k$-linear symmetric function on ${\mathcal P}(\mathbb R^n)$ agreeing with $\textsl{v}_k^\varphi$ on the diagonal. Then, on one hand,
\begin{equation*}
\frac{\partial^k}{\partial\lambda_1\ldots\partial\lambda_k}\,
L\left(\sum_{\ell=1}^k\lambda_\ell \Gamma_\ell,\ldots,\sum_{\ell=1}^k\lambda_\ell \Gamma_\ell\right)
=
\frac{\partial^k}{\partial\lambda_1\ldots\partial\lambda_k}\,
\textsl{v}_k^\varphi\left(\sum_{\ell=1}^k\lambda_\ell \Gamma_\ell\right)
\,,
\end{equation*}
for every $\Gamma_1,\ldots,\Gamma_k\in{\mathcal P}(\mathbb R^n)$. On the other hand, as $L$ is multilinear and symmetric, one has
\begin{align*}
&\frac{\partial^k}{\partial\lambda_1\ldots\partial\lambda_k}\,
L\left(\sum_{\ell=1}^k\lambda_\ell \Gamma_\ell,\ldots,\sum_{\ell=1}^k\lambda_\ell \Gamma_\ell\right)
\\
&=
\frac{\partial^k}{\partial\lambda_1\ldots\partial\lambda_k}\,
L\left(
\sum_{\ell_1=1}^k\lambda_{\ell_1} \Gamma_{\ell_1},
\ldots,
\sum_{\ell_k=1}^k\lambda_{\ell_k} \Gamma_{\ell_k},
\right)
\\
&=
\frac{\partial^k}{\partial\lambda_1\ldots\partial\lambda_k}\,
\sum_{\ell_1=1}^k
\ldots
\sum_{\ell_k=1}^k
L(\Gamma_{\ell_1},\ldots,\Gamma_{\ell_k})
\lambda_{\ell_1}\cdots\lambda_{\ell_k}
\\
&=
k!\,L(\Gamma_1,\ldots,\Gamma_k)\,,
\end{align*}
so that the claim follows by \lemref{pol-fi}.\findim

\begin{rmq}\label{qmono}
On the monotonicity of the $k$-th mixed $\varphi$-volume.

\noindent
Let $\Gamma_1,\Gamma_1^\prime,\ldots,\Gamma_k\in{\mathcal P}(\mathbb R^n)$ be such that $\Gamma_1\subset\Gamma_1+\Gamma_1^\prime$. Then, by linearity,
\begin{equation*}
V^\varphi_k(\Gamma_1+\Gamma_1^\prime,\Gamma_2,\ldots,\Gamma_k)=V_k^\varphi(\Gamma_1,\ldots,\Gamma_k)
+V_k^\varphi(\Gamma^\prime_1,\ldots,\Gamma_k)
\geq V_k^\varphi(\Gamma_1,\ldots,\Gamma_k)\,.
\end{equation*}
Without any further hypothesis on the shape of $\varphi$, this is the only monotonicity property $V_k^\varphi$ can satisfy. 
\end{rmq}

\begin{thm}\label{bon}
For every $1\leq k\leq n$, the $k$-th intrinsic $\varphi$-volume $\textsl{v}_k^\varphi $ is a $k$-homogeneous, translation invariant valuation on ${\mathcal P}(\mathbb R^n)$. If $\varphi$ is continuous, $\textsl{v}_k^\varphi $ provides a $k$-homogeneous, translation invariant, continuous valuation on ${\mathcal K}(\mathbb R^n)$.
\end{thm}
\noindent
\pf 
The translation invariance follows at once from the definition of $\textsl{v}^\varphi_k$. 
The $k$-homogeneity is easy, indeed, for every $\Gamma\in {\mathcal P}(\mathbb R^n)$ and $\lambda>0$
\begin{align*}
\textsl{v}^\varphi_k(\lambda\Gamma)
&=
\frac{1}{\varkappa_{n-k}}\sum_{\Delta\in{\mathcal B}(\lambda\Gamma,k)}\varphi(E_\Delta)\vol_k(\Delta)\psi_\Gamma(\Delta)
\\
&=
\frac{1}{\varkappa_{n-k}}\sum_{\Delta\in{\mathcal B}(\Gamma,k)}\varphi(E_{\lambda\Delta})\vol_k(\lambda\Delta)\psi_\Gamma(\lambda\Delta)\\
&=
\frac{1}{\varkappa_{n-k}}\sum_{\Delta\in{\mathcal B}(\Gamma,k)}\varphi(E_\Delta)\lambda^k\vol_k(\Delta)\psi_\Gamma(\Delta)
\\
&=
\lambda^k
\textsl{v}^\varphi_k(\Gamma)\,.
\end{align*}

Since $\mathbb R$ is a Hausdorff space, by virtue of \thmref{Sal}, it is enough to show that $\textsl{v}_k^\varphi $ is a weak $k$-homogeneous valuation. 
Let $H\subset\mathbb R^n$ be a hyperplane intersecting $\Gamma$, then $\Gamma\cap H$, $\Gamma\cap H^+$ and $\Gamma\cap H^-$ are polytopes too. All the $k$-dimensional faces of $\Gamma\cap H$ are also $k$-faces of $\Gamma\cap H^+$ and $\Gamma\cap H^-$. For each such common $k$-face $\Delta$ one has a disjoint union $
K_{\Delta,\Gamma\cap H^+}\cup K_{\Delta,\Gamma\cap H^-}=K_{\Delta,\Gamma\cap H}$. Moreover, $\Gamma\cap H^+$ and $\Gamma\cap H^-$ may respectively have $k$-faces $\Delta_1$ and $\Delta_2$ such that $\Delta_1\cup\Delta_2$ is a $k$-face of $\Gamma$ not included in $H$. In this case, $\varphi(E_{\Delta_1})=\varphi(E_{\Delta_2})=\varphi(E_{\Delta_1\cup\Delta_2})$, $\vol_k(\Delta_1)+\vol_k(\Delta_2)=\vol_k(\Delta_1\cup\Delta_2)$ and $K_{\Delta_1,\Gamma\cap H^+}= K_{\Delta_2,\Gamma\cap H^-}=K_{\Delta_1\cup\Delta_2,\Gamma}$. Any other $k$-face of $\Gamma\cap H^+$ and $\Gamma\cap H^-$ which does not intersect $H$ is also a $k$-face of $\Gamma$. The preceding analysis implies that
\begin{equation*}
\textsl{v}^\varphi_k(\Gamma\cap H^+)+\textsl{v}^\varphi_k(\Gamma\cap H^-)=\textsl{v}^\varphi_k(\Gamma\cap H)+\textsl{v}^\varphi_k(\Gamma)\,.
\end{equation*}
The equality $\textsl{v}^\varphi_{k}(\varnothing)=0$ completes the proof of the first statement. The continuity hypothesis on $\varphi$ implies the continuity of the weak valuation $\textsl{v}^\varphi_k$ on ${\mathcal P}(\mathbb R^n)$. By \thmref{Gro}, this produces a continuous valuation on the whole ${\mathcal K}(\mathbb R^n)$ which is still $k$-homogeneous and translation invariant.
\findim

We have shown that the $V_k^\varphi$ and $V_k$ share some important properties, nevertheless the validity of other properties (like monotonicity, continuity, non-degeneracy, etc) heavily depends on the shape of the function $\varphi$. In the following section we will present a geometrically interesting instance of such a function, playing a fundamental role in the construction of Kazarnovski{\v\i} pseudovolume.


\section{Volume distortion}\label{distorzione}

Let $n\in\mathbb N^*$ and let $E_1,\,E_2\subset\mathbb R^n$ be a couple linear subspaces of positive dimensions  $d_1=\dim E_1$ and $d_2=\dim E_2$. If $\mathbb R^n$ is endowed with the usual scalar product $(\,,)$, consider the $\mathbb R$-linear mapping
\begin{equation*}
\varrho_{E_1,E_2}: E_1\overset{\iota}{\longrightarrow}E_1\oplus E_1^\perp=\mathbb R^n\overset{id}{\longrightarrow}\mathbb R^n=E_2\oplus E_2^\perp\overset{p}{\longrightarrow} E_2\,,
\end{equation*}
where $\iota$ is the inclusion, $id$ is the identity mapping and $p$ is the projection.
To the linear mapping $\varrho_{E_1,E_2}$ one can associate the $\mathbb R$-linear operator 
\begin{equation*}
\tilde\varrho_{E_1,E_2}:E_1\xrightarrow{~\varrho_{E_1,E_2}~} E_2\xrightarrow{~~\psi_2~~} E_2^*\xrightarrow{{}^t(\varrho_{E_1,E_2})} E_1^* \xrightarrow{~~\psi_1^{-1}~~} E_1\,,
\end{equation*}
where $\psi_1$ (resp. $\psi_2$) is the natural isomorphism mapping any element $v$ of $E_1$ (resp. $E_2$) to the linear form $(v,\cdot)$ on $E_1$ (resp. $E_2$).
\textit{The coefficient of volume distortion under the projection of $E_1$ on $E_2$} is the non-negative real number $\varrho(E_1,E_2)$\label{vdist} defined as 
\begin{equation*}
\varrho(E_1,E_2)=\sqrt{\det\tilde\varrho_{E_1,E_2}}\,.
\end{equation*}
The number $\det \tilde\varrho_{E_1,E_2}$ is sometimes called the \emph{Gramian} of the linear mapping $\varrho_{E_1,E_2}$. 
Of course the rank of $\tilde\varrho_{E_1,E_2}$ cannot exceed $\min(d_1,d_2)$, and, in fact, $\varrho(E_1,E_2)=0$ as soon as $d_2<d_1$. 
If $d_1=d_2=d$ and $a_1,\ldots,a_d$ (resp. $b_1,\ldots,b_d$) is an orthonormal basis of $E_1$ (resp. $E_2$), then $\varrho(E_1,E_2)=\vert \det\varrho_{E_1,E_2}\vert=\vert \det( a_\ell,b_j)\vert$, this value equals the $d$-dimensional non-oriented volume of the parallelotope generated in $E_2$ by the row-vectors of the matrix $(( a_\ell,b_j))_{\ell,j}$ and, of course, it does not depend on the choice of the orthonormal basis in $E_1$ and $E_2$. As the norm of each row-vector in this matrix cannot exceed $1$, the same happens to $\varrho(E_1,E_2)$.

By giving $\mathbb C^n$ the euclidean structure induced by the standard scalar product of $\mathbb R^{2n}$, one can define the coefficient of volume distortion for $\mathbb R$-linear subspaces of $\mathbb C^n$. Such a scalar product on $\mathbb C^{n}$ can be realized as the real part of the standard hermitian product of $\mathbb C^n$, i.e. $\re\langle\cdot\, ,\cdot\rangle$.

If $A\subset\mathbb C^n$ is a non-empty subset, 
let $\aff_\mathbb R A$ (resp. $\aff_\mathbb C A$) denote the real (resp. complex) affine subspace of $\mathbb C^n$ spanned by $A$, whereas $E_A$ (still) denotes the $\mathbb R$-linear subspace of $\mathbb C^n$ parallel to $\aff_\mathbb R A$.
 
Let ${\mathcal G}(\mathbb C^n,\mathbb R)$\label{wcgrass} be the Grassmanian of the real linear subspaces of $\mathbb C^n$. For every $E\in {\mathcal G}(\mathbb C^n,\mathbb R)$, the $\mathbb C$-linear subspace spanned by $E$ will be noted $\lin_\mathbb C E$, of course it is equal to $E+i E$. The \textit{dimension} $\dim_\mathbb R E$\label{rdim} of a subspace $E\in{\mathcal G}(\mathbb C^n,\mathbb R)$ is just the real dimension of $E$, whereas the \textit{complex dimension} $\dim_\mathbb C E$\label{cdim} is meant to be the complex dimension of $\lin_\mathbb C E$. In general $\dim_\mathbb R E\leq\dim_\mathbb R(\lin_\mathbb C E)\leq2\dim_\mathbb R E$ and, setting $E^\mathbb C=E\cap iE$\label{cmax},  the equality $\dim_\mathbb R(\lin_\mathbb C E)=2\dim_\mathbb R E$ occurs if and only if and only if $\lin_\mathbb  C E=E\oplus iE$, i.e. $E^\mathbb C=\{0\}$.

By $E^\perp$\label{rperp} (resp. $E^{\perp_\mathbb C}$\label{cperp}) we will denote the orthogonal complement of $E$ respect to the scalar product $\re\langle\,,\rangle$ (resp. the hermitian product $\langle\,,\rangle$). In general one has the relations $(\lin_\mathbb C E)^{\perp_\mathbb C}= E^{\perp_\mathbb C}$, $(E^\perp)^\mathbb C=E^\perp\cap iE^\perp=E^{\perp_\mathbb C}$ and $E^{\perp_\mathbb C}\subseteq E^\perp$, the latter inclusion becoming an equality as soon as $E=\lin_\mathbb C E$. By setting $E^{\,\prime}=E^\perp\cap\lin_\mathbb C E$\label{eprime}, the orthogonal decomposition $\lin_\mathbb C E=E\oplus E^\prime$ yields $\dim_\mathbb R E^\prime\leq \dim_\mathbb R E$, because
\begin{equation*}
\dim_\mathbb R E^\prime=\frac{1}{2}\left(\dim_\mathbb R(\lin_\mathbb C E)-\dim_\mathbb R E^\mathbb C\right)\leq\dim_\mathbb R E\,.
\end{equation*}
In particular, $\dim_\mathbb R E^\prime= \dim_\mathbb R E$ if and only if $E^\mathbb C=\{0\}$, although in general $iE\neq E^\prime$. Anyway $\dim_\mathbb R E^\prime\equiv \dim_\mathbb R E$~mod $2$, since otherwise the complex space $\lin_\mathbb C E$ would have an odd real dimension.

The \textit{Cauchy-Riemann dimension} of $E$ is just the complex dimension of $E^\mathbb C$ . A linear subspace $E\in{\mathcal G}(\mathbb C^n,\mathbb R)$ is \textit{full-dimensional} if its real dimension equals $2n$, it is \textit{purely real} if it is included in $\mathbb R^n$.
A subspace $E\in{\mathcal G}(\mathbb C^n,\mathbb R)$ is \textit{equidimensional} if its real and complex dimensions coincide. Of course $1$-dimensional subspaces are such, as well as purely real ones. If $E_1,E_2\in {\mathcal G}(\mathbb C^n,\mathbb R)$ and $E_1\subset E_2$, then $E_1^\mathbb C\subseteq E_2^\mathbb C$, whence the equidimensionality of $E_2$ is inherited by each of its subspaces. A subspace $E\in {\mathcal G}(\mathbb C^n,\mathbb R)$ with $\dim_{\mathbb R} E>n$ cannot be equidimensional, for in this case $E^\mathbb C\neq\{0\}$. 

The orthogonal decompositions $E\oplus E^{\,\prime}=\lin_\mathbb C E=iE\oplus iE^{\,\prime}$ imply that $iE=(i E^\prime)^\perp\cap\lin_\mathbb C E$, so that one can consider the $\mathbb R$-linear projection $\varrho_{E,iE^\prime}:E\rightarrow iE^{\,\prime}$. If $v\in {\rm ker} \varrho_{E,iE^\prime}$, then $\re \langle v,iE^\prime\rangle=0$ so that $v\in (iE^\prime)^\perp\cap E\subseteq (i E^\prime)^\perp\cap\lin_\mathbb C E= iE$, hence $v\in E^\mathbb C$. On the other hand, $v\in E^\mathbb C$ implies $-iv\in E$, then $\re \langle v,iE^\prime\rangle=\re\langle -iv, E^\prime\rangle=0$, i.e. ${\rm ker} \varrho_{E,iE^\prime}=E^\mathbb C$. It follows that $\varrho_{E,iE^\prime}$ is an isomorphism if and only if $E$ is equidimensional. The \textit{coefficient of volume distortion $\varrho(E)$} is by definition the coefficient of volume distortion under the projection $\varrho_{E,iE^\prime}$ of $E$ onto $iE^\prime$, i.e. $\varrho(E)=\varrho(E,iE^\prime)$\label{roe}.  A subspace  $E\in{\mathcal G}(\mathbb C^n,\mathbb R)$ is \textit{real similar} if $\varrho(E)=1$. The trivial subspace $\{0\}$ is conventionally assumed to be real similar, so that we set $\varrho(\{0\})=1$. Purely real subspaces are real similar as well as 1-dimensional ones, nevertheless there are equidimensional subspaces which are not real similar. Remark that if $E$ is equidimensional, then $E^\prime$ cannot be a complex subspace, since otherwise $E^\prime=iE^\prime$ which yields the contradiction $\varrho({E,iE^\prime})=\varrho({E,E^\prime})=0$. 
If $E$ is equidimensional and $d=\dim_\mathbb R E>0$, let $v_1,\ldots,v_d,w_1,\ldots,w_d$ be a basis of $\lin_\mathbb C E$ over $\mathbb R$ such that $v_1,\ldots,v_d$ is an orthonormal basis of $E$ and $w_1,\ldots,w_d$ an orthonormal basis of $E^{\,\prime}$, then $\varrho(E)=\vert\det(\re\langle v_\ell,i w_j\rangle)_{\ell,j}\vert$. In this case, the equality $\re\langle v_\ell,i w_j\rangle=-\re\langle w_j, iv_\ell\rangle$ implies that $\varrho(E)=\varrho(E^\prime)$. In order to achieve a standard formula for $\varrho(E)$, it is useful to find a standard basis for $E^\prime$.

\begin{lem}\label{lemmat}
Let $E\in{\mathcal G}(\mathbb C^n,\mathbb R)$ with $\dim_\mathbb R E^\mathbb C=2c<d=\dim_\mathbb R E$. If $v_1,\ldots,v_{2c},v_{2c+1},\ldots,v_d$ is an orthonormal basis of $E$ such that $v_1,\ldots,v_{2c}$ span $E^\mathbb C$, then the vectors 
$$
t_\ell=iv_\ell-\sum_{s=1}^d\re\langle iv_\ell,v_s\rangle v_s\,,
$$
for $\ell=1,\ldots,d$ span $E^\prime$ and $t_{2c+1},\ldots,t_d$ provide a basis of $E^\prime$. In particular $E$ is equidimensional if and only if $t_1,\ldots,t_d$ are $\mathbb R$-linearly independent.
\end{lem}
\noindent
\pf
The linear span of $v_{2c+1},\ldots,v_d$ is equidimensional, otherwise $E^\mathbb C$ would not be maximal among the complex subspaces contained in $E$. This means that $\dim_\mathbb R (\lin_\mathbb C E)=2d-2c$, so that $\dim_\mathbb R E^\prime=d-2c$. Observe also that $c>0$ implies $t_1=\ldots=t_{2c}=0$. Suppose, by contradiction, that $\sum_{\ell=2c+1}^d a_\ell t_\ell=0$ is a non trivial linear combination, then 
\begin{equation}\label{ugua}
i\sum_{\ell=2c+1}^d a_\ell v_\ell=\sum_{s=1}^d\left(\sum_{\ell=2c+1}^d a_\ell \re\langle iv_\ell,v_s\rangle \right)v_s\,.
\end{equation}
The vector on the right of \eqref{ugua} cannot be $0$, (otherwise $v_1,\ldots,v_d$ would be linearly dependent), so equality \eqref{ugua} implies that 
\begin{equation}\label{vec}
\sum_{\ell=2c+1}^d a_\ell v_\ell
\in
E^\mathbb C\setminus\{0\}\,.
\end{equation}
If $c=0$ (i.e. $E^\mathbb C=\{0\}$) the relation \eqref{vec} is impossible, if $c>0$ the vector $\sum_{\ell=2c+1}^d a_\ell v_\ell$ must be a non trivial linear combination of $v_1,\ldots,v_{2c}$, which is also impossible since $v_1,\ldots,v_d$ are linearly independent. 

In particular, when $E$ is equidimensional, $t_1,\ldots,t_d\in E^\prime$ provide a basis of $E^\prime$. Conversely, the linear independence of $t_1,\ldots,t_d$  means that $\dim_\mathbb R E=\dim_\mathbb R E^\prime$, i.e. $E$ is equidimensional. \findim

\begin{rmq}{On the natural orientation of complex spaces.}\label{orientazioni}

\noindent
Every complex subspace of $\mathbb C^n$ has a natural orientation coming from the complex structure. This orientation will be referred to as the \textit{positive orientation}. If $e_1,\ldots,e_c$ is a sequence of $\mathbb C$-linearly independent vectors of $\mathbb C^n$, then the corresponding $\mathbb C$-linear span is positively oriented, as a real $2c$-dimensional subspace, by the basis $e_1,ie_1,\ldots,e_c,ie_c$.
\end{rmq}

\begin{coro}\label{pdim}
Let $E\in{\mathcal G}(\mathbb C^n,\mathbb R)$ with $\dim_\mathbb R E^\mathbb C=2c<d=\dim_\mathbb R E$. If $v_1,\ldots,v_{2c},v_{2c+1},\ldots,v_d$ is an orthonormal basis of $E$ such that $v_1,\ldots,v_{2c}$ is an orienting basis of $E^\mathbb C$, then the sequence $v_1,\ldots,v_{2c},v_{2c+1},t_{2c+1},\ldots,v_d,t_d$ is an orienting basis of $E\oplus E^\prime=\lin_{\mathbb C} E$. In particular, the basis $v_1,\ldots,v_d,t_{2c+1},\ldots,t_d$ of $\lin_{\mathbb C} E$ yields the positive orientation if and only if
\begin{equation}\label{pd}
(d-2c)\equiv_4 0\qquad\text{or}\qquad (d-2c)\equiv_4 1\,.
\end{equation}
\end{coro}
\noindent
\pf
Up to a unitary transformation, we may suppose $\lin_\mathbb C E=\mathbb C^{d-c}\times\{0\}^{n-d+c}$ and $E^\mathbb C=\mathbb C^c\times\{0\}^{d-2c}\times \{0\}^{n-d+c}$. As shown by \lemref{unitarion}, from the real point of view, a unitary transformation is an orthogonal transformation of $\mathbb R^{2n}$ with positive determinant, so it preserves the (real) orthonormality of the basis $v_1,\ldots,v_{2c},v_{2c+1},\ldots,v_d$. By virtue of the preceding argument we can carry out the computation in the smaller space $\mathbb C^{d-c}$. As $t_\ell-iv_\ell\in E$, for every $2c+1\leq\ell\leq d$, it follows that
\begin{equation}\label{determ}
\det(v_1,\ldots,v_{2c},v_{2c+1},t_{2c+1},\ldots,v_d,t_d)
=
\det(v_1,\ldots,v_{2c},v_{2c+1},iv_{2c+1},\ldots,v_d,iv_d)\,,
\end{equation}
and of course \eqref{determ} cannot be zero because $v_1,\ldots,v_{2c},v_{2c+1},t_{2c+1},\ldots,v_d,t_d$ is a basis. In fact \eqref{determ} equals
\begin{equation*}
\det 
\left(\begin{array}{ccc|ccccc}
(v_1)_1&\ldots & (v_{2c})_1 & (v_{2c+1})_1 & (iv_{2c+1})_1 & \ldots & (v_d)_1 & (iv_d)_1 \\
 \vdots & \ddots &\vdots & \vdots & \vdots & \ddots & \vdots & \vdots \\
(v_1)_{2c}&\ldots  &(v_{2c})_{2c}  & (v_{2c+1})_{2c} & (iv_{2c+1})_{2c} & \ldots & (v_d)_{2c} & (iv_d)_{2c}  \\\hline
0&\ldots& 0& (v_{2c+1})_{2c+1} & (iv_{2c+1})_{2c+1} & \ldots & (v_d)_{2c+1} & (iv_d)_{2c+1} \\
\vdots& \ddots & \vdots& \vdots & \vdots & \ddots & \vdots & \vdots \\
0&\ldots & 0&(v_{2c+1})_{2d} & (iv_{2c+1})_{2d} & \ldots & (v_{2c+1})_{2d} & (iv_{2c+1})_{2d}\end{array}\right)\,,
\end{equation*}
where both the north-west block and the south-east one cannot have zero determinant. Indeed, the north-west block has positive determinant because $v_1,\ldots,v_{2c}$ is an orienting basis of $E^\mathbb C$, while that of the south-east one equals
\begin{equation*}
\det
\left(\begin{array}{ccccc}
 \re v_{2c+1\,2c+1} &-\im v_{2c+1\,2c+1}  & \ldots  & \re v_{d\,2c+1} & -\im v_{d\,2c+1} \\
 \im v_{2c+1\,2c+1} & \re v_{2c+1\,2c+1} & \ldots & \im v_{d1} & \re v_{d1} \\
 \vdots & \vdots & \ddots & \vdots & \vdots \\
 \re v_{2c+1\,d} &-\im v_{2c+1\,d}  & \ldots  & \re v_{dd} & -\im v_{dd} \\
 \im v_{2c+1\,d} & \re t_{2c+1\,d} & \ldots & \im v_{dd} & \re v_{dd}
\end{array}\right)\,,
\end{equation*}
which, in the notation of \lemref{unitarion}, is just 
\begin{equation*}
\left\vert\det\psi_{d-2c}
\left(\begin{array}{ccc}
v_{2c+1\,2c+1} &\ldots& v_{d\,2c+1} \\
\vdots&\ddots&\vdots\\
v_{2c+1\,d}&\ldots & v_{dd}
\end{array}\right)
\right\vert^2\,.
\end{equation*}
The last statement follows easily by the relation
\begin{equation*}
\det(v_1,\ldots,v_{2c},v_{2c+1},t_{2c+1},\ldots,v_d,t_d)
=
(-1)^{\frac{(d-2c-1)(d-2c)}{2}}
\det(v_1,\ldots,v_d,t_{2c+1},\ldots,t_d)\,.
\end{equation*}
Indeed the number $(d-2c-1)(d-2c)/2$ is even if and only if $(d-2c)\equiv_4 0$ or $(d-2c)\equiv_4 1$.
The proof is thus complete.\findim

\begin{thm}\label{figo}
Let $E\in{\mathcal G}(\mathbb C^n,\mathbb R)$. If $v_1,\ldots,v_d$ is an orthonormal basis of $E$, then 
\begin{equation}\label{formula-figa2}
\varrho(E)=\sqrt{\det(\re\langle t_\ell,t_j\rangle)}\,.
\end{equation}
\end{thm}
\noindent
\pf
Observe, first of all, that $\det(\re\langle t_\ell,t_j\rangle)\geq0$ because it is the square of the volume of the parallelotope spanned by $t_1,\ldots,t_d$.
Let us first suppose $E$ is not equidimensional. In this case the vectors $t_1,\ldots,t_d\in E^\prime$ are $\mathbb R$-linearly dependent and then $\det(\re\langle t_\ell,t_j\rangle)$ vanishes. 
If $E$ is equidimensional, the vectors $t_1,\ldots,t_d$ provide a basis of $E^\prime$. This basis is generally not orthogonal (unless $d=2$) nor orthonormal, so setting $u_1=t_1$ and, for $1<j\leq d$,
\begin{equation*}
u_j=t_j-\sum_{s=1}^{j-1}\frac{\re\langle t_j,u_s\rangle}{\re\langle u_s,u_s\rangle}\,u_s\,,
\end{equation*}
yields an orthogonal basis of $E^\prime$, and the vectors $w_j=u_j/\Vert u_j\Vert$, $1\leq j\leq d$, provide an orthonormal basis of $E^\prime$. As a consequence, $\varrho(E)=\vert\det(\re\langle v_\ell,iw_j\rangle)\vert$. 
Observe that, for any $1\leq j,\ell\leq d$,
\begin{align}
\re\langle t_\ell, t_j\rangle
&=
\delta_{\ell,j}-\sum_{s=1}^d\re\langle iv_j,v_s\rangle\re\langle iv_\ell, v_s\rangle\nonumber
\\
&=
\delta_{\ell,j}-\sum_{s=1}^d\im\langle v_j,v_s\rangle\im\langle v_\ell, v_s\rangle\nonumber
\\
&=
\delta_{\ell,j}+\mathop{\sum_{1\leq s\leq d}^{\phantom{d}}}_{ j\neq s \neq\ell}\langle v_\ell, v_s\rangle\langle v_j,v_s\rangle\label{bua}
\,,
\end{align}
$
\re\langle v_\ell,i t_j\rangle
=
-\re\langle t_\ell, t_j\rangle
$ and, moreover (by induction on $j$)
\begin{equation*}
\re\langle v_\ell,iw_j\rangle
=
-\re\langle t_\ell,w_j\rangle
=
-\Vert u_j\Vert^{-1}\left[\re\langle t_\ell,t_j\rangle-\sum_{s=1}^{j-1}\re\langle v_j,iw_s\rangle\re\langle v_\ell,iw_s\rangle\right].
\end{equation*}
The $j$-th column of the matrix $(\re\langle v_\ell,iw_j\rangle)$ is the sum of $-\Vert u_j\Vert^{-1}\re\langle t_\ell,t_j\rangle$ and a linear combination of the preceding columns, so its determinant equals that of the matrix $(-\Vert u_j\Vert^{-1}\re\langle t_\ell,t_j\rangle)$ whence
\begin{equation*}
\vert\det(\re\langle v_\ell,iw_j\rangle)\vert
=
\det(\Vert u_j\Vert^{-1}\re\langle t_\ell,t_j\rangle)
=
\left(\prod_{s=1}^d \Vert u_s\Vert\right)^{-1}\det(\re\langle t_\ell,t_j\rangle)\,.
\end{equation*}
On the other hand, the real number $\prod_{s=1}^d \Vert u_s\Vert$ equals the $d$-dimensional (non-oriented) volume of the parallelotope of $E^\prime$ generated by the vectors $u_1,\ldots,u_d$ and such a volume can be written as
\begin{equation*}
\prod_{s=1}^d \Vert u_s\Vert=\det(\re\langle t_\ell,w_j\rangle)=\sqrt{\det(\re\langle t_\ell,t_j\rangle)}\,,
\end{equation*}
so that the preceding computation implies that $\varrho(E)=\sqrt{\det(\re\langle t_\ell,t_j\rangle)}$. 
\findim

\begin{rmq}\label{orient}
For every orientation ${\rm or}(E)\in\{-1,1\}$ on $E$ it is possible to choose an orientation ${\rm or} (E^\prime)\in\{-1,1\}$ such that $E\oplus E^\prime=\lin_\mathbb C E$ gets the positive orientation.
This can be achieved, for example, by choosing \emph{any} orthonormal basis $v_1,\ldots,v_d$ of $E$ such that $v_1,\ldots,v_{2c}$ spans $E^\mathbb C$ over $\mathbb R$, followed by the basis $t_{\sigma(2c+1)},\ldots,t_{\sigma(d)}$, where $t_{2c+1},\ldots,t_d$ are the vectors constructed in~\lemref{lemmat} and $\sigma$ is a permutation of the set $\{2c+1,\ldots,d\}$ with $\text{sign}(\sigma)=(-1)^{(d-2c-1)(d-2c)/2}$. In particular, if $(d-2c)\equiv_4 0$ or $(d-2c)\equiv_4 1$, $\sigma$ can simply be the identity, whereas it can be a transposition if $(d-2c)\equiv_4 2$ or $(d-2c)\equiv_4 3$. Even more simply, the basis $v_1,\ldots,v_{2c},(-1)^{(d-2c-1)(d-2c)/2}t_{2c+1},t_{2c+2},\ldots,t_{d}$ will always do. Such an orientation on $E^\prime$ yields an orientation on $E^\perp=E^\prime\oplus E^{\perp_\mathbb C}$, where $E^{\perp_\mathbb C}$ is positively oriented as any complex linear space. It will be referred to as the \emph{coherent} orientation of $E^\prime$ or $E^\perp$.
\end{rmq}

\begin{rmq}\label{2d}
Observe that if $d=2$ and $E$ is equidimensional, the equality \eqref{bua} shows that the vectors $t_1,t_2$ are orthogonal to each-other and that $\Vert t_1\Vert=\Vert t_2\Vert=\sqrt{\varrho(E)}$. By \eqref{pd}, the basis $v_1,v_2,t_2,t_1$ yields the positive orientation of $\lin_\mathbb C E$.
\end{rmq}

\begin{coro}\label{fighissimo}
Let $E\in{\mathcal G}(\mathbb C^n,\mathbb R)$, $d=\dim_\mathbb R E$ and $m=\dim_\mathbb C(\lin_\mathbb C E)$. Let also $v_1,\ldots,v_d$ be an orthonormal basis of $E$ with respect to $\re\langle\,,\rangle$, $e_1,\ldots,e_m$ an orthonormal basis of $ \lin_\mathbb C E$ with respect to $\langle\,,\rangle$. If $A=(\langle v_\ell,v_j\rangle)_{1\leq \ell,j\leq d}$ and $B=(\langle v_\ell,e_k\rangle)_{1\leq \ell\leq d,1\leq k\leq m}$, then 
\begin{equation}\label{formula-figa3}
\varrho(E)=\det A=\det (B\,{}^t\!\bar B).
\end{equation}
In particular, if $d=m$, $\varrho(E)$ equals the real jacobian of the $\mathbb C$-linear operator on $\lin_\mathbb C E$ mapping $e_\ell$ to $v_\ell$, for every $1\leq\ell\leq d$.
\end{coro}
\noindent
\pf
Let $M$ be the $d\times d$ matrix whose $(\ell,j)$-entry is
\begin{equation*}
\mathop{\sum_{1\leq s\leq d}}_{}\langle v_\ell, v_s\rangle\langle v_j,v_s\rangle
=
\delta_{\ell,j}+\mathop{\sum_{1\leq s\leq d}}_{ j\neq s \neq\ell}\langle v_\ell, v_s\rangle\langle v_j,v_s\rangle
\,,
\end{equation*}
then, by \eqref{formula-figa2} and \eqref{bua}, it follows that $\varrho(E)=\sqrt{\det M}$. Observing that $\det A$ is real (since $A$ a hermitian matrix) and $M$ equals the product of $A$ with ${}^t\!A$, it follows that $\varrho(E)=\sqrt{\det M}=\sqrt{(\det A)^2}=\vert\det A\vert$. Nevertheless $\det A\geq 0$. Indeed $A=B\,{}^t\!\bar B$ and $\rk A=\rk B=m$, so  $\det A=\det(B\,{}^t\!\bar B)=0$ if $E$ is not equidimensional (i.e. $d>m$) and $\det A=\det(B\,{}^t\!\bar B)=\vert\det B\vert^2\neq 0$ as soon as $d=m$.
In the latter case, mapping $e_\ell$ to $v_\ell$ yields a $\mathbb C$-linear automorphism, say $L$, on $\lin_\mathbb C E$ whose matrix in the basis $e_1,\ldots,e_d$ is just $(\langle v_\ell,e_j\rangle)$. As happens to every other holomorphic mapping $\mathbb C^d\to\mathbb C^d$, the square of the modulus of the complex jacobian of $L$ (i.e. $\vert\det L\vert^2$) equals its real jacobian, i.e. the jacobian of $L$ regarded as a mapping $\mathbb R^{2d}\to\mathbb R^{2d}$.
\findim
\begin{rmq}
\cororef{fighissimo} shows that $\{E\in{\mathcal G}(\mathbb C^n,\mathbb R)\mid \varrho(E)=0\}$ 
is a Zariski closed subset, thus equidimensionality is a generic property of $\mathbb R$-linear subspaces of $\mathbb C^n$. Moreover it implies that the coefficient of volume distortion is a continuous function.  \findim
\end{rmq}
\begin{rmq}\label{qinv} $\varrho(E)$ is unitarily invariant but not orthogonally invariant.

\noindent
By \cororef{fighissimo}, if $F:\mathbb C^n\to\mathbb C^n$ is a unitary transformation, 
\begin{equation*}
\varrho(F(E))=\det \langle F(v_\ell),F(v_j)\rangle=\det F\det A\det {}^t\overline{F}=\det A=\varrho(E)\,.
\end{equation*}
This is generally false if $F$ is an orthogonal transformation of $\mathbb R^{2n}$. For instance, consider the real subspace $E=\mathbb R^n\subset\mathbb C^n$ and the orthogonal $\mathbb R$-linear mapping $F:\mathbb R^{2n}\to\mathbb R^{2n}$ permuting the elements $ie_1$ and $e_n$ of the canonical basis of $\mathbb R^{2n}$. As $F(E)$ contains the complex line spanned by $e_1$, it follows that $\varrho(F(E))=0$, whereas $\varrho(E)=1$.
\findim
\end{rmq}

\begin{rmq} {Two further formulas for $\varrho(E)$}\label{utile}

\noindent
Up to a renumbering of $v_1,\ldots,v_d$, one can suppose that $v_1,\ldots, v_m$ are $\mathbb C$-linearly independent, so the basis $e_1,\ldots,e_m$ can be chosen as an orthonormalization of $v_1,\ldots, v_m$ with respect to $\langle\,,\rangle$. Setting $e_k=0$, for $m< k\leq d$, a further formula for $\varrho(E)$ would read as follows
\begin{equation*}
\varrho(E)
=
\vert\det(B\mid\bm{0})\vert^2=
\left\vert
\det
\left(
\begin{array}{cccccc}
\langle v_1,e_1\rangle & \cdots & \langle v_1,e_m\rangle&0&\cdots &0 \\
\vdots & \ddots & \vdots &\vdots&\ddots&\vdots\\
\langle v_d,e_1\rangle & \cdots & \langle v_d,e_m\rangle&0&\cdots& 0
 \end{array}
 \right)
 \right\vert^2\,,
\end{equation*}
where the last zero $d\times (d-m)$ block is absent if $d=m$. Using the Gram-Schimdt method one sets $e_1=v_1$, $\omega_\ell=v_\ell-\sum_{k=1}^{\ell-1} \langle v_\ell,e_k\rangle e_k$ and $e_\ell=\omega_\ell/\Vert \omega_\ell\Vert$, for $1<\ell\leq m$. This choice has the advantage of producing some zero entries in the matrix $B$, namely $\langle v_\ell,e_k\rangle=0$, for every $2\leq \ell<k\leq m$, so that, when $d=m$ the matrix $B$ is lower triangular and 
\begin{equation*}
\varrho(E)=\vert\det B\vert^2=\left\vert\prod_{\ell=1}^d\langle v_\ell,e_\ell\rangle\right\vert^2
=
\prod_{\ell=2}^d\vert\langle v_\ell,e_\ell\rangle\vert^2=\prod_{\ell=2}^d\Vert \omega_\ell\Vert^2\,.
\end{equation*}
Another formula for $\varrho(E)$ can be found in \cite{Ale}. Let $B_E=B_{2n}\cap E$ and $iB_E=\{iz\in i E\mid z\in B_E\}$. The subset $B_E+iB_E$ has the structure of a (possibly degenerate) cylinder on the base $B_E$ whose dimension equals $2d$ if and only if $E$ is equidimensional. Since 
\begin{equation*}
0\leq\vol_{2d}(B_E+iB_E)=\varrho(E)\varkappa_d^2\leq\varkappa_d^2\,,
\end{equation*}
one can set $\varrho(E)=\varkappa_d^{-2}\vol_{2d}(B_E+iB_E)$.
\findim
\end{rmq}

\begin{exe} 
The coefficient of length, area and volume distortion.

\noindent
{\rm
When $d=1$ we already know that $\varrho(E)=1$, so there is no length distortion. In fact $\varrho(E)=\langle v_1,v_1\rangle=1$.
If $d=2$, 
$\varrho(E)
=
1+\langle v_1,v_2\rangle^2
=
1-(\im\langle v_1,v_2\rangle)^2
$. 
Since $\im\langle v_1,v_2\rangle=\re\langle iv_1,v_2\rangle$ equals the cosine of the (convex) angle of the vectors $iv_1$ and $v_2$, $\det A$ is the square of the sine of that angle. As a consequence, $E$ is real similar if and only if $v_2$ and $iv_1$ are orthogonal to each other.
For $d=3$, one gets a similar formula
$
\varrho(E)= 1+\langle v_1,v_2\rangle^2+\langle v_1,v_3\rangle^2+\langle v_2,v_3\rangle^2
= 1-(\im\langle v_1,v_2\rangle)^2-(\im \langle v_1,v_3\rangle)^2-(\im\langle v_2,v_3\rangle)^2
$.

So, for $1\leq d\leq 3$, one can write the following simple formula:
\begin{equation*}
\varrho(E)=1+\sum_{1\leq\ell<j\leq d} \langle v_\ell,v_j\rangle^2
=
1-\sum_{1\leq\ell<j\leq d}(\im \langle v_\ell,v_j\rangle)^2
\,,
\end{equation*}
which, unfortunately, does not work anymore as $d$ becomes greater than $3$. 
\findim}
\end{exe}

\begin{exe}
The coefficient of $4$-dimensional volume distortion.

\noindent
{\rm
Unlike the preceding two cases, if $d=4$ the formula for $\varrho(E)$ involves also some mixed products:
\begin{align*}
\varrho(E)
&=
1+\langle v_1,v_2\rangle^2+\langle v_1,v_3\rangle^2+\langle v_1,v_4\rangle^2
+
\langle v_2,v_3\rangle^2+\langle v_2,v_4\rangle^2+\langle v_3,v_4\rangle^2+
\\
&+
\langle v_1,v_2\rangle^2\langle v_3,v_4\rangle^2+\langle v_1,v_3\rangle^2\langle v_2,v_4\rangle^2+\langle v_1,v_4\rangle^2\langle v_2,v_3\rangle^2
\\
&
-2 \langle v_1,v_2\rangle\langle v_1,v_3\rangle\langle v_2,v_4\rangle
\langle v_3,v_4\rangle
\\
&
+2 \langle v_1,v_2\rangle\langle v_1,v_4\rangle\langle v_2,v_3\rangle\langle v_3,v_4\rangle
\\
&
-2\langle v_1,v_3\rangle\langle v_1,v_4\rangle\langle v_2,v_3\rangle\langle v_2,v_4\rangle
\,.
\end{align*}
\findim}
\end{exe}

\begin{exe}\label{realsim}
Real similar $n$-dimensional coordinate subspaces of $\mathbb C^n$.

\noindent
{\rm
The number of real $n$-dimensional coordinate subspaces of $\mathbb C^n$ equals ${2n\choose n}$. If $n=1$ the coordinate axes are both real similar, now suppose $\mathbb C^{n-1}$ admits $2^{n-1}$ real similar $(n-1)$-dimensional real coordinate subspaces. Adding a dimension to each real similar $(n-1)$-dimensional coordinate subspace of $\mathbb C^{n-1}$ yields a real similar $n$-dimensional real coordinate subspace of $\mathbb C^n$ and this can be achieved in only two ways: namely by adding the vector $e_n$ or $ie_n$. So the number of real similar $n$-dimensional real coordinate subspace of $\mathbb C^n$ equals $2^n$.
\findim
}
\end{exe}


\section{Convex bodies in $\mathbb C^n$}\label{corpi-complessi}

The sets of convex bodies, polytopes, piecewise $2$-regular compact subsets and $\varepsilon$-neighborhoods of polytopes of $\mathbb C^n$ are respectively denoted ${\mathcal K}(\mathbb C^n)$, ${\mathcal P}(\mathbb C^n)$, ${\mathcal S}^2(\mathbb C^n)$ and ${\mathcal P}_\circ(\mathbb C^n)$, in the same manner we will use the notation ${\mathcal K}^2_\ell(\mathbb C^n)$ to denote the set of piecewise $2$-regular convex subsets of $\mathbb C^n$ with support function of class ${\mathcal C}^\ell(\mathbb C^n\setminus\{0\})$, $\ell\in\mathbb N^*$\label{konvexc}. For a non-empty subset $A\subseteq \mathbb C^n$, the notions of dimension, complex dimension, full-dimensionality, pure reality, equidimensionality and real similarity as well as the coefficient of volume distortion introduced for $\mathbb R$-linear subspaces of $\mathbb C^n$, can be simply extended by just invoking the real subspace $E_A$. In particular, this can be done for convex bodies. For any $A\in{\mathcal K}(\mathbb C^n)$, set $\dim_\mathbb R A =\dim_\mathbb R E_A$, $\dim_\mathbb C A =\dim_\mathbb C E_A$ and $\varrho(A)=\varrho(E_A)$\label{roA}. In the complex setting ${\mathcal K}(\mathbb R^n)=\{A\in{\mathcal K}(\mathbb C^n)\mid A\subset\mathbb R^n\}$, i.e. the set of purely real convex bodies; in the same way we will use the notation ${\mathcal K}^2(\mathbb R^n)$, ${\mathcal K}_\ell(\mathbb R^n)$, ${\mathcal K}^2_\ell(\mathbb R^n)$, ${\mathcal P}_\circ(\mathbb R^n)$ and ${\mathcal P}(\mathbb R^n)$\label{konvexcr}.  
In the complex setting a \emph{lattice polytope} $\Gamma\in{\mathcal P}(\mathbb C^n)$ is a polytope whose vertices have coordinates in the ring $\mathbb Z+i\mathbb Z$ of Gauss' integers.

If $A\in{\mathcal K}(\mathbb C^n)$, the support function $h_A$ of $A$ is computed with respect to the scalar product $\re\langle\,,\rangle$, i.e.
$
h_A(z)=\sup_{v\in A}\re\langle z,v\rangle
$, so that, given $v\in\mathbb C^n\setminus\{0\}$, the corresponding supporting hyperplane $H_A(v)$ for the convex body $A$ is 
$
H_A(v)=\{z\in\mathbb C^n\mid \re \langle z, v\rangle=h_A(v)\}\,.
$

If $\Gamma\in{\mathcal P}(\mathbb C^n)$, and $\Delta_1\preccurlyeq\Delta_2\preccurlyeq\Gamma$, the inclusions $E_{\Delta_1}\subseteq E_{\Delta_2}\subseteq E_\Gamma$ imply $\lin_\mathbb C E_{\Delta_1}\subseteq \lin_\mathbb C E_{\Delta_1}\subseteq \lin_\mathbb C E_\Gamma$ and $E_{\Delta_1}^{\mathbb C}\subseteq E_{\Delta_2}^{\mathbb C}\subseteq E_\Gamma^{\mathbb C}$. 

For every $0\leq k\leq \dim_\mathbb R\Gamma$, ${\mathcal B}_{\rm ed}(\Gamma,k)$\label{bcked} will denote the subset of ${\mathcal B}(\Gamma)$ consisting of  equidimensional $k$-faces. Notice that, for any (non-empty) polytope, ${\mathcal B}_{\rm ed}(\Gamma,0)={\mathcal B}(\Gamma,0)\neq\varnothing$ and, if $\Gamma$ is not reduced to a single point, ${\mathcal B}_{\rm ed}(\Gamma,1)={\mathcal B}(\Gamma,1)\neq\varnothing$. 

For any proper face $\Delta\prec \Gamma$, one has $\dim_\mathbb C\Delta\leq\dim_\mathbb C\Gamma$, but if $\Gamma$ is purely real the equality cannot occur. The following theorem suggested by Kazarnovski{\v\i} in a private communication clarifies the situation in the general case.

\begin{thm}\label{iok}
Let $\Gamma\in{\mathcal P}(\mathbb C^n)$ be a $d$-polytope and let $2\leq k\leq\dim_\mathbb C\Gamma$ an integer. If ${\mathcal B}_{\rm ed}(\Gamma,k)$ is empty, then the following statements hold true.
\begin{enumerate}
\item For any $\ell\in\{k+1,\ldots,d\}$, the set ${\mathcal B}_{\rm ed}(\Gamma,\ell)$ is empty; in particular $\Gamma$ is not equidimensional. 
\item There exists a  $\Delta\in{\mathcal B}(\Gamma,k)$ such that $\lin_\mathbb C E_\Gamma=\lin_\mathbb C E_\Delta$; in particular $\dim_\mathbb C\Gamma<k$.
\item $\dim_\mathbb C \Gamma=\max\{0\leq j\leq n\mid {\mathcal B}_{\rm ed}(\Gamma,j)\neq\varnothing\}.$
\end{enumerate}
\end{thm}
\noindent
\pf 
If $\Delta^\prime$ is an $\ell$-face of $\Gamma$ and $\Delta$ is a $k$-face of $\Delta^\prime$, then $\{0\}\neq E_\Delta^\mathbb C\subseteq E_{\Delta^\prime}^\mathbb C$, so $\Delta^\prime$ cannot be equidimensional. In particular this occurs to the improper face $\Gamma$, which proves the first statement.

In order to prove the second statement, observe that by first statement one has $d>\dim_\mathbb C \Gamma$ and $E_\Gamma^\mathbb C\neq\{0\}$. Now it may happen that $E_\Gamma=E_\Gamma^\mathbb C$ or $E_\Gamma\supsetneq E_\Gamma^\mathbb C$. In both cases we construct a strictly decreasing sequence $\Delta_{d-1}\succ\ldots\succ\Delta_k$ of faces of $\Gamma$ such that, for every $\ell\in\{1,\ldots,d-k\}$, one has $\dim_\mathbb R\Delta_{d-\ell}=d-\ell$ and $\lin_\mathbb C E_{\Delta_\ell}=\lin_\mathbb C E_\Gamma$. The last term $\Delta_k$ in the sequence will be the desired $k$-face. 

If $E_\Gamma=E_\Gamma^\mathbb C$, then $E_\Gamma=\lin_\mathbb C E_\Gamma=E_\Gamma^\mathbb C$ and $d$ is an even number. For any facet $\Delta_{d-1}\prec\Gamma$ one has $E_{\Delta_{d-1}}\subsetneq\lin_\mathbb C E_{\Delta_{d-1}}\subseteq\lin_\mathbb C E_\Gamma= E_\Gamma$. The first inclusion is strict because $E_{\Delta_{d-1}}$ has odd real dimension whereas $\lin_\mathbb C E_{\Delta_{d-1}}$ (as every complex subspace) has even real dimension. Computing real dimensions yields 
$(d-1)<\dim_\mathbb R\lin_\mathbb C E_{\Delta_{d-1}}\leq \dim_\mathbb R\lin_\mathbb C E_\Gamma=d$, i.e. $\dim_\mathbb R\lin_\mathbb C E_{\Delta_{d-1}}=d$, i.e. $\lin_\mathbb C E_{\Delta_{d-1}}=\lin_\mathbb C E_\Gamma$. 

In the second case, let $v$ be a vertex of $\Gamma$ and consider all the facets of $\Gamma$ admitting $v$ as a vertex. The intersection of the real affine subspaces spanned by such facets is reduced to the single vertex $v$ (otherwise $v$ would not be a vertex), so the intersection of the corresponding real linear subspaces is reduced to $\{0\}$. As $E_\Gamma^\mathbb C\neq\{0\}$, there will be a facet $\Delta_{d-1}$ of $\Gamma$ admitting $v$ as a vertex such that $E_{\Delta_{d-1}}$ does not contain $E_\Gamma^\mathbb C$. This implies that the real linear subspace $E_{\Delta_{d-1}}\cap E_\Gamma^\mathbb C$ has real codimension $1$ in $E_\Gamma^\mathbb C$ and, though it is not a complex subspace, it spans the whole $E_\Gamma^\mathbb C$ over $\mathbb C$. As a consequence $E_\Gamma^\mathbb C$ is a complex subspace of $\lin_\mathbb C E_{\Delta_{d-1}}$ and so $E_{\Delta_{d-1}}$ must be a proper real subspace of $E_\Gamma\cap\lin_\mathbb C E_{\Delta_{d-1}}$, since otherwise it would contain $E_\Gamma^\mathbb C$. Let $v_1,\ldots,v_{d-1}$ be a basis of $E_{\Delta_{d-1}}$ and let $v_d$ be an element of $E_\Gamma\cap\lin_\mathbb C E_{\Delta_{d-1}}$ not belonging to $E_{\Delta_{d-1}}$.
The real subspace spanned by $v_1,\ldots,v_d$ equals $E_\Gamma$ so that $E_\Gamma\subseteq E_\Gamma\cap\lin_\mathbb C E_{\Delta_{d-1}}$, i.e. $E_\Gamma\subseteq\lin_\mathbb C E_{\Delta_{d-1}}$, whence $\lin_\mathbb C E_\Gamma=\lin_\mathbb C E_{\Delta_{d-1}}$. 

Repeating the whole argument on the facet $\Delta_{d-1}$, one finds a ridge $\Delta_{d-2}$ of $\Gamma$, such that $\lin_\mathbb C E_\Gamma=\lin_\mathbb C E_{\Delta_{d-1}}=\lin_\mathbb C E_{\Delta_{d-2}}$. By induction, in $(d-k)$ steps one gets a $k$-face $\Delta_k$ of $\Gamma$ such that $\lin_\mathbb C E_\Gamma=\lin_\mathbb C E_{\Delta_k}$. Like any other $k$-face of $\Gamma$, $\Delta_k$ is not equidimensional, so $\dim_\mathbb C \Gamma=\dim_\mathbb C \Delta_k<k$.

As for the last statement, let $j$ be the minimal integer such that ${\mathcal B}_{\rm ed}(\Gamma,j+1)$ is empty. Then ${\mathcal B}_{\rm ed}(\Gamma,j)$ is non-empty and, by the preceding statement, $j=\dim_\mathbb C\Gamma$.
\findim

Dual cones to the faces of a polytope as well as the corresponding stars are defined like in the real case. 
For any $\lambda\in\mathbb R\setminus\{0\}$ and any polytope $\Gamma$, every face of the polytope $\lambda \Gamma$ has the form $\lambda\Delta$, for a unique face $\Delta$ of $\Gamma$, so that
\begin{equation*}
\qquad
\aff_{\mathbb R}\lambda\Delta=\lambda\aff_{\mathbb R}\Delta\,,
\qquad
\aff_\mathbb C \lambda\Delta=\lambda\aff_\mathbb C \Delta\,,
\qquad
E_{\lambda \Delta}=E_\Delta\,,
\qquad\qquad
\end{equation*}
\begin{equation*}
\lin_\mathbb C E_{\lambda\Delta}=\lin_\mathbb C E_\Delta\,,
\quad
E_{\lambda\Delta}^\perp=E_\Delta^\perp\,,
\quad
E_{\lambda\Delta}^{\perp_\mathbb C}=E_\Delta^{\perp_\mathbb C}\,,
\quad
K_{\lambda\Delta,\lambda \Gamma}=K_{\Delta,\Gamma}\,.
\end{equation*}

If $\Gamma_1,\ldots,\Gamma_s\in{\mathcal P}(\mathbb C^n)$ and  $\Gamma=\sum_{\ell=1}^s\Gamma_s$ there exists a unique sequence of faces $\Delta_\ell\preccurlyeq \Gamma_\ell$, $\ell=1,\ldots,s$, such that $\Delta=\sum_{\ell=1}^s\Delta_\ell$. As a consequence, one gets the equalities
\begin{equation*}
\aff_{\mathbb R}\Delta=\sum_{\ell=1}^s \aff_{\mathbb R}{\Delta_\ell}\,,
\quad
\aff_\mathbb C \Delta=\sum_{\ell=1}^s \aff_\mathbb C \Delta_\ell\,,
\quad
E_\Delta=\sum_{\ell=1}^s E_{\Delta_\ell}\,,
\end{equation*}
\begin{equation*}
\lin_\mathbb C E_\Delta=\sum_{\ell=1}^s \lin_\mathbb C E_{\Delta_\ell}\,,
\quad
E_\Delta^\perp=\bigcap_{\ell=1}^s E_{\Delta_\ell}^\perp\,,
\quad
E_\Delta^{\perp_\mathbb C}=\bigcap_{\ell=1}^s E_{\Delta_\ell}^{\perp_\mathbb C}\,,
\end{equation*}
and
\begin{equation*}
K_{\Delta,\Gamma}=\bigcap_{\ell=1}^s K_{\Delta_\ell,\Gamma_\ell}\,.
\end{equation*}

The mapping $\varrho$ on the real grassmanian ${\mathcal G}(\mathbb C^n,\mathbb R)$ yields a corresponding $k$-th intrinsic $\varrho$-volume $\textsl{v}_k^\varrho$, for $0\leq k\leq 2n$, nevertheless $\textsl{v}_k^\varrho\equiv 0$ if $n< k\leq 2n$ since there is no equidimensional real subspace in $\mathbb C^n$ of (real) dimension greater than $n$. Moreover, as $\varrho\leq 1$, one has \begin{equation}\label{ineq}
\textsl{v}_k^\varrho(\Gamma)
\leqslant
\left[\max_{\Delta\in{\mathcal B}(\Gamma,k)}
\varrho(\Delta)\right]
\textsl{v}_k(\Gamma)
\leqslant
\textsl{v}_k(\Gamma), \end{equation}
for every $\Gamma\in{\mathcal P}(\mathbb C^n)$ and, by a continuity and multilinearity argument,
\begin{equation}\label{diseg}
\textsl{V}_k^\varrho(A_1,\ldots,A_k)
\leqslant
\varkappa_{2n-k}^{-1}{2n\choose k} V_{2n}(A_1,\ldots,A_k,B_{2n}[2n-k]),
\end{equation} 
for every $A_1,\ldots,A_k\in{\mathcal K}(\mathbb C^n)$. Since $\varrho=1$ on $1$-dimensional real subspaces, then $\textsl{v}^\varrho_1
=\textsl{v}_1$ and, by virtue of (\ref{kth}), a continuity argument yields
\begin{equation}\label{meta}
\textsl{v}_1^\varrho(A)=\textsl{v}_1(A)=\dfrac{2n}{\varkappa_{2n-1}}V_{2n}(A,B_{2n}[2n-1]),
\end{equation}
for every $A\in{\mathcal K}(\mathbb C^n)$.

The following result provides a non vanishing condition for the $k$-th mixed $\varrho$-volume.

\begin{coro}\label{ndeg}
Let $\Gamma_1,\ldots,\Gamma_k\in{\mathcal P}(\mathbb C^n)$ and $\Gamma=\sum_{j=1}^k\Gamma_j$. Then $\dim_\mathbb C\Gamma<k$ if and only if $\textsl{V}_k^\varrho(\Gamma_1,\ldots,\Gamma_k)=0$.
\end{coro}
\noindent
\pf
Suppose $\dim_\mathbb C\Gamma<k$ and, for every non empty subset  $I\subseteq\{1,\ldots,k\}$, set $\Gamma_I=\sum_{j\in I}\Gamma_j$. Then, for every $I$, $\dim_\mathbb C\Gamma_I\leq \dim_\mathbb C\Gamma<k
$
and by \thmref{iok} ${\mathcal B}_{\rm ed}(\Gamma_I,k)=\varnothing$ whence $\textsl{v}_k^\varrho(\Gamma_I)=0$ and, by \eqref{pol-fi1}, $\textsl{V}_k^\varrho(\Gamma_1,\ldots,\Gamma_k)=0$. On the other hand, if 
$\textsl{V}_k^\varrho(\Gamma_1,\ldots,\Gamma_k)=0$, then ${\mathcal B}_{\rm ed}(\Gamma,k)=\varnothing$ and again by \thmref{iok} we obtain $\dim_\mathbb C\Gamma<k$.
\findim

\begin{coro}\label{ndeg2}
Let $\Gamma_1,\ldots,\Gamma_k\in{\mathcal P}(\mathbb C^n)$. The following assertions are equivalent.
\begin{enumerate}
\item[{$(a)$}] $V_k^\varrho(\Gamma_1,\ldots,\Gamma_k)>0$;
\item[{$(b)$}] there are line segments $\Lambda_1\subset \Gamma_1,\ldots,\Lambda_k\subset \Gamma_k$ with linearly independent complex directions;
\item[{$(c)$}] For every non empty subset $I\subseteq\{1,\ldots,k\}$, $\dim_\mathbb C \left(\sum_{j\in I} \Gamma_j\right)\geq \# I$.
\end{enumerate}
\end{coro}
\noindent\pf 
$(a)\Rightarrow (b)$. By~\cororef{ndeg}, the Minkowski sum $\Gamma=\Gamma_1+\ldots+\Gamma_k$ admits at least an equidimensional $k$-face $\Delta=\Delta_1+\ldots+\Delta_k$ such that $V_k(\Delta_1,\ldots,\Delta_k)>0$. By the non-degeneracy property of mixed volume, this implies that $\Delta$ is an edge-sum $k$-face, i.e. there exist line segments $\Lambda_1\subseteq \Delta_1,\ldots,\Lambda_k\subseteq \Delta_k$ spanning over $\mathbb R$ lines with $\mathbb R$-linearly independent directions. Among all such sequences, there must be at least one sequence of segments spanning over $\mathbb C$ lines with $\mathbb C$-linearly independent directions, otherwise $\Delta$ would not be equidimensional.

\noindent
$(b)\Rightarrow (a)$. As $k=\dim_\mathbb C (\Lambda_1+\ldots+\Lambda_k)\leq \dim_\mathbb C (\Gamma)$, by~\cororef{ndeg}, $V_k^\varrho(\Gamma_1,\ldots,\Gamma_k)>0$.

\noindent
$(b)\Rightarrow (c)$. For every $I\subseteq\{1,\ldots,k\}$, $\dim_\mathbb C \left(\sum_{j\in I} \Gamma_j\right)\geq \dim_\mathbb C \left(\sum_{j\in I} \Lambda_j\right)\geq \# I$.

\noindent
$(c)\Rightarrow (b)$ is a consequence of Lemma~5.1.8 in~\cite{Sch}. 
\findim

For any $\varepsilon\in\mathbb R_{>0}$, the $\varepsilon$-regularization of a convex body $A\in{\mathcal K}(\mathbb C^n)$ is defined as in the real case, its support function is
\begin{equation*}
h_{R_\varepsilon A}(z)=\int_{\mathbb C^n} h_A(z+u \Vert z\Vert )\varphi_\varepsilon(\Vert u\Vert)\,\upsilon_{2n} \,,
\end{equation*}
where 
\begin{eqnarray}\label{vform2n}
\upsilon_{2n}=\bigwedge_{\ell=1}^n(i/2)\de u_\ell\wedge\de {\bar u}_\ell=\bigwedge_{\ell=1}^n\de(\re u_\ell)\wedge\de(\im u_\ell)
\end{eqnarray}
is the usual standard volume form on $\mathbb C^n=\mathbb R^{2n}$ and 
$\varphi_\varepsilon\in{\mathcal C}^\infty(\mathbb R,\mathbb R)$ is a non-negative function with (compact) support included in the interval $[\varepsilon/2,\varepsilon]$ such that 
$
\int_{\mathbb C^n} \varphi_\varepsilon(\Vert u\Vert)\,\upsilon_{2n}=1\,.
$
The regularization $h_{R_\varepsilon A}$ is convex but not differentiable at the origin, whereas the usual one, $h_A*\varphi_\varepsilon$, provides a plurisubharmonic smooth function on the whole $\mathbb C^n$.

\section{Kazarnovski{\v\i} pseudovolume on the class {\boldmath${\mathcal S}^2(\mathbb C^n)$}}\label{pS}

We will use the identification $\mathbb C^n= (\mathbb R+i\mathbb R)^n\simeq\mathbb R^{2n}$, so that $z_\ell=x_\ell+iy_\ell=\xi_{2\ell-1}+i\xi_{2\ell}$ and $i:\mathbb R^{2n}\to \mathbb R^{2n}$ acts on the vector $\xi$ in the usual way, i.e. $i\xi=(-\xi_2,\xi_1,\ldots,-\xi_{2n},\xi_{2n-1})$. Moreover 
\begin{align*}
\frac{\partial}{\partial z_\ell}
&=
\frac{1}{2}\left(\frac{\partial}{\partial x_\ell}-i\frac{\partial}{\partial y_\ell}\right)
\,,
&
\de z_\ell
&=
\de x_\ell+i\de y_\ell
\,,
\\
\frac{\partial}{\partial \bar z_\ell}
&=
\frac{1}{2}\left(\frac{\partial}{\partial x_\ell}+i\frac{\partial}{\partial y_\ell}\right)
\,,
&
\de \bar z_\ell
&=
\de x_\ell-i\de y_\ell
\,,
\end{align*}
for every $\ell=1,\ldots,n$. In consequence, the operators $\partial$, $\bar\partial$, $\de=\partial+\bar\partial$ and $\dc=i(\bar\partial-\partial)$\label{dc} are given by
\begin{equation*}
\partial
=
\frac{1}{2}\sum_{\ell=1}^n\left(\frac{\partial}{\partial x_\ell}\de x_\ell+\frac{\partial}{\partial y_\ell}\de y_\ell\right)+i\left(\frac{\partial}{\partial x_\ell}\de y_\ell-\frac{\partial}{\partial y_\ell}\de x_\ell\right)
=
\sum_{\ell=1}^n\frac{\partial}{\partial z_\ell}\de z_\ell\,,
\end{equation*}
\begin{equation*}
\bar\partial
=
\frac{1}{2}\sum_{\ell=1}^n\left(\frac{\partial}{\partial x_\ell}\de x_\ell+\frac{\partial}{\partial y_\ell}\de y_\ell\right)-i\left(\frac{\partial}{\partial x_\ell}\de y_\ell-\frac{\partial}{\partial y_\ell}\de x_\ell\right)
=
\sum_{\ell=1}^n\frac{\partial}{\partial \bar z_\ell}\de \bar z_\ell\,.
\end{equation*}
\begin{equation*}
\de
=
\sum_{\ell=1}^n\frac{\partial}{\partial x_\ell}\de x_\ell+\frac{\partial}{\partial y_\ell}\de y_\ell
=
\sum_{\ell=1}^n\frac{\partial}{\partial z_\ell}\de z_\ell+\frac{\partial}{\partial \bar z_\ell}\de \bar z_\ell\,,
\end{equation*}
\begin{equation*}
\dc
=
\sum_{\ell=1}^n \frac{\partial}{\partial x_\ell}\de y_\ell-\frac{\partial}{\partial y_\ell}\de x_\ell
=
i\sum_{\ell=1}^n\frac{\partial}{\partial \bar z_\ell}\de \bar z_\ell-\frac{\partial}{\partial z_\ell}\de z_\ell\,.
\end{equation*}
Since $\dc$ equals its conjugate operator, it follows that $\dc$ and $\ddc$\label{ddc} (often referred to as \emph{the Bott-Chern operator}) are real operators, moreover $\ddc=2i\partial\bar\partial$ is given by
\begin{align}\label{ddcreale}
\ddc
&=
\sum_{\ell,j=1}^n
\frac{\partial^2}{\partial x_\ell\partial x_j}\,\de x_\ell\wedge\de y_j
+
\frac{\partial^2}{\partial y_\ell\partial x_j}\,\de y_\ell\wedge\de y_j
-
\frac{\partial^2}{\partial x_\ell\partial y_j}\,\de x_\ell\wedge\de x_j
-
\frac{\partial^2}{\partial y_\ell\partial y_j}\,\de y_\ell\wedge\de x_j\,
\\
&=
\sum_{\ell,j=1}^n
\frac{\partial^2}{\partial x_\ell\partial x_j}\,\de x_\ell\wedge\de y_j
+
\frac{\partial^2}{\partial y_\ell\partial y_j}\,\de x_j\wedge\de y_\ell
+
\frac{\partial^2}{\partial y_\ell\partial x_j}\,\de y_\ell\wedge\de y_j
-
\frac{\partial^2}{\partial y_j\partial x_\ell}\,\de x_\ell\wedge\de x_j\,
\\
&=
\sum_{\ell=1}^n
\left(\frac{\partial^2}{\partial x_\ell^2}
+
\frac{\partial^2}{\partial y_\ell^2}\right)
\de x_\ell\wedge\de y_\ell
\\
&+
\sum_{1\leqslant\ell\neq j\leqslant n}
\frac{\partial^2}{\partial x_\ell\partial x_j}\,\de x_\ell\wedge\de y_j
+
\frac{\partial^2}{\partial y_\ell\partial y_j}\,\de x_j\wedge\de y_\ell
+
\frac{\partial^2}{\partial y_\ell\partial x_j}\,\de y_\ell\wedge\de y_j
-
\frac{\partial^2}{\partial y_j\partial x_\ell}\,\de x_\ell\wedge\de x_j
\,
\end{align}
or
\begin{equation}\label{ddcrealexi}
\begin{split}
\ddc
&=
\sum_{\ell,j=1}^n
\frac{\partial^2}{\partial \xi_{2\ell-1}\partial \xi_{2j-1}}\,\de \xi_{2\ell-1}\wedge\de \xi_{2j}
+
\frac{\partial^2}{\partial \xi_{2\ell}\partial \xi_{2j-1}}\,\de \xi_{2\ell}\wedge\de \xi_{2j}
\\
&-
\frac{\partial^2}{\partial \xi_{2\ell-1}\partial \xi_{2j}}\,\de \xi_{2\ell-1}\wedge\de \xi_{2j-1}
-
\frac{\partial^2}{\partial \xi_{2\ell}\partial \xi_{2j}}\,\de \xi_{2\ell}\wedge\de \xi_{2j-1}\,
\end{split}
\end{equation}
or, in complex coordinates, by
\begin{equation*}
\ddc
=
2i\sum_{\ell,j=1}^n \frac{\partial^2}{\partial z_\ell\partial\bar z_j} \de z_\ell\wedge\de\bar z_j\,.
\end{equation*}

It should be mentioned that very often in the literature the operator $\dc$ is defined in different ways, namely $(i/2)(\bar\partial-\partial)$ or even $(i/2\pi)(\bar\partial-\partial)$, nevertheless we will conform to the choice of \cite{Ka1}, where $\dc=i(\bar\partial-\partial)$. 

Remark that for every $h\in{\mathcal C}^2(\mathbb C^n,\mathbb C)$, $\ddc h$ is a $2$-form so that the $2n$-form $(\ddc h)^{\wedge n}$ is identically zero as soon as $h$ depends on less than $n$-complex variables.

The integration with respect to Lebesque measure in $\mathbb C^n$ will be usually performed by integrating the standard volume form
\begin{equation*}
\upsilon_{\mathbb C^n}=\upsilon_{2n}=\frac{i^n}{2^n}\bigwedge_{\ell=1}^n\de z_\ell\wedge\de {\bar z}_\ell
=
\bigwedge_{\ell=1}^n\de x_\ell\wedge\de y_\ell
=
(-1)^{\frac{n(n-1)}{2}}\de x_1\wedge\ldots\wedge\de x_n\wedge\de y_1\wedge\ldots\wedge\de y_n\,,
\end{equation*}
then, for $h\in{\mathcal C}^2(\mathbb C^n,\mathbb C)$, 
\begin{equation}\label{MA}
(\ddc h)^{\wedge n}=(2 i\partial\bar\partial h)^{\wedge n}=4^n n!\det\left(\frac{\partial^2 h}{\partial z_\ell\partial\bar z_k}\right)\upsilon_{2n}\,.
\end{equation}
If $f:\mathbb C^n\to\mathbb C^n$ is a holomorphic mapping, then 
\begin{equation*}
(\ddc (h\circ f))^{\wedge n}=\vert\det J_f\vert^{2n} (\ddc h)^{\wedge n}\,,
\end{equation*}
where $J_f$ is the (complex) jacobian determinant of $f$.

For every differential form $\omega$ on $\mathbb C^n$, by $*\omega$ we will denote the unique form such that $\omega\wedge *\omega=\upsilon_{2n}$. The operator $*$\label{*-op} on forms is well defined and known as \textit{Hodge $*$-operator}. It can be proved that it defines a linear operator on $(p,q)$-forms with values on $(n-p,n-q)$-forms.

For any $j=1,\ldots,n$, we will also use the $(2n-1)$-forms
\begin{align}
\upsilon_{2n}[x_j]&=*\de x_j=\de y_j\wedge\bigwedge_{\substack{\ell=1 \\ \ell\neq j}}^{n}\de x_\ell\wedge\de y_\ell\,,\label{basex}
\\
\upsilon_{2n}[y_j]&=-*\de y_j=\de x_j\wedge\bigwedge_{\substack{\ell=1 \\ \ell\neq j}}^{n}\de x_\ell\wedge\de y_\ell\,,\label{basey}
\\
\upsilon_{2n}[z_j]&=*\de z_j= \frac{i^n}{2^n}\,\de \bar z_j\wedge\bigwedge_{\substack{\ell=1 \\ \ell\neq j}}^{n}\de z_\ell\wedge\de \bar z_\ell\,,\nonumber
\\
\upsilon_{2n}[\bar z_j]&=-*\de \bar z_j=\frac{i^n}{2^n}\,\de  z_j\wedge\bigwedge_{\substack{\ell=1 \\ \ell\neq j}}^{n}\de z_\ell\wedge\de \bar z_\ell\,,\nonumber
\end{align}
with the obvious relations
\begin{equation*}
\upsilon_{2n}[z_j]-\upsilon_{2n}[\bar z_j]=\upsilon_{2n}[x_j]
\quad{\rm and}\quad
\upsilon_{2n}[z_j]+\upsilon_{2n}[\bar z_j]=i\upsilon_{2n}[y_j]\,.
\end{equation*}

For every $\mathbb R$-linear oriented subspace $E\subset\mathbb C^n$ let $\upsilon_E$ denote the volume form on $E$.\label{vformE} 
Observe that, for every $v\in\mathbb C^n$ the $\dc$ of the linear form $\re\langle v, z\rangle$ is given by 
\begin{align*}
\dc \re\langle v,z\rangle
&=
i(\bar\partial-\partial)\re\langle v, z\rangle
\\
&=
\frac{i}{2}(\bar\partial-\partial)(\langle v, z\rangle+\langle z,v\rangle)
\\
&=
\frac{i}{2}\left(\sum_{\ell=1}^n v_\ell\de\bar z_\ell-\bar v_\ell\de z_\ell\right)
\end{align*}
so that, for any $w\in\mathbb C^n$,
\begin{align*}
(\dc \re\langle v,z\rangle)(w)
&=
\frac{i}{2}\sum_{\ell=1}^n v_\ell\bar w_\ell-\bar v_\ell w_\ell
\\
&=
\frac{i}{2}\left(\langle v,w\rangle-\langle w,v\rangle\right)
\\
&=
\frac{i}{2}\left(\langle v,w\rangle-\overline{\langle v,w\rangle}\right)
\\
&=
-\im\langle v,w\rangle
\\
&=
\re\langle iv,w\rangle\,.
\end{align*}
This shows that, if $\Vert v\Vert=1$ and $E$ is the line spanned by $v$, the $1$-form $\dc  \re\langle v,z\rangle$ is a volume form $\upsilon_{E^\prime}$ on the line $E^\prime=iE$ spanned by $iv$.

Unlike the one dimensional case, if $E\subset\mathbb C^n$ is an equidimensional real subspace with $d>1$ and $v_1,\ldots,v_d$ is an orthonormal basis of $E$ one gets the following equalities valid just on $E^\perp$:
\begin{equation}\label{forme}
\varrho(E)\upsilon_{E^\prime}=(-1)^{d(d-1)/2}\upsilon_{i E}=(-1)^{d(d-1)/2}\dc\re\langle v_1,z\rangle\wedge\ldots\wedge\dc\re\langle v_d,z\rangle\,.
\end{equation} 
Indeed, by using the notation of \thmref{figo}, 
\begin{align*}
\upsilon_{E^\prime}
&=
(-1)^{d(d-1)/2}\re\langle w_1,-\rangle\wedge\ldots\wedge\re\langle w_d,-\rangle
\\
&=
\frac{(-1)^{d(d-1)/2}}{\varrho(E)}\re\langle u_1,-\rangle\wedge\ldots\wedge\re\langle u_d,-\rangle
\\
&=
\frac{(-1)^{d(d-1)/2}}{\varrho(E)}\re\langle t_1,-\rangle\wedge\ldots\wedge\re\langle t_d,-\rangle\,.
\end{align*}
When $a_1,\ldots,a_d\in E^\perp$, for any $1\leq j,\ell\leq d$, 
\begin{align*}
\re\langle t_j,a_\ell\rangle
&=
\re\langle iv_j-\sum_{s=1}^d\re\langle iv_\ell,v_s\rangle v_s,a_\ell\rangle
\\
&=
\re\langle iv_j,a_\ell\rangle-\sum_{s=1}^d\re\langle iv_\ell,v_s\rangle\re\langle v_s,a_\ell\rangle
\\
&=
\re\langle iv_j,a_\ell\rangle\,,
\end{align*}
because $\re\langle v_s,a_\ell\rangle=0$, then
\begin{align*}
\upsilon_{E^\prime}(a_1,\ldots,a_d)
&=
\frac{(-1)^{d(d-1)/2}}{\varrho(E)}\left(\re\langle iv_1,-\rangle\wedge\ldots\wedge\re\langle iv_d,-\rangle\right)
(a_1,\ldots,a_d)
\\
&=
\frac{(-1)^{d(d-1)/2}}{\varrho(E)}\upsilon_{i E}(a_1,\ldots,a_d)
\\
&=
\frac{(-1)^{d(d-1)/2}}{\varrho(E)}\left(\dc\re\langle v_1,z\rangle\wedge\ldots\wedge\dc\re\langle v_d,z\rangle\right) (a_1,\ldots,a_d)\,.
\end{align*}
Moreover, as $E^\perp=E^\prime\oplus E^{\perp_\mathbb C}$ and $iE\perp E^{\perp_\mathbb C}$, it follows that $iE\cap E^\perp=iE\cap(E^\prime\oplus E^{\perp_\mathbb C})=iE\cap E^\prime$, whence
the ``restriction'' of $\upsilon_{iE}$ to $E^\perp$ coincides with the ``restriction'' to $E^\prime$.

Let $A\in{\mathcal S}^2(\mathbb C^n)$, then $\partial A$ admits a \textit{Gauss map} $\nu_{\partial A}:\partial A\rightarrow\partial B_{2n}$\label{Gaussmap}, given by $\nu_{\partial A}=\nabla\rho_A/\Vert\nabla\rho_A\Vert$, for every fixed defining function $\rho_A$. The Gauss map $\nu_{\partial A}$ is of class ${\mathcal C}^1$ and it does not depend on the choice of the defining function used to represent it. Infact, if $\tilde\rho_A$ is another defining function for $A$, then there exists a positive function $g\in{\mathcal C}^1(\mathbb C^n,\mathbb R)$ such that $\tilde\rho_A=g\rho_A$ on $\mathbb C^n$, $g\neq0$ on $\partial A$, and $\nabla \tilde\rho_A=g\nabla\rho_A$ on $\partial A$. 
As a consequence, we may choose a defining function $\rho_A$ such that $\Vert\nabla\rho_A\Vert\equiv 1$ on $\partial A$, but we will not assume this normalization.

It follows that
\begin{equation*}
\nabla\rho_A(z)=\left(\frac{\partial\rho_A}{\partial x_1},\frac{\partial\rho_A}{\partial y_1},\ldots,\frac{\partial\rho_A}{\partial x_n},\frac{\partial\rho_A}{\partial y_n}\right)=2\left(\frac{\partial\rho_A}{\partial \bar z_1},\ldots,\frac{\partial\rho_A}{\partial \bar z_n}\right)\,,
\end{equation*}
\begin{equation*}
\nu_{\partial A}(z)=\frac{1}{\Vert\nabla\rho_A\Vert}\left(\frac{\partial\rho_A}{\partial x_1},\frac{\partial\rho_A}{\partial y_1},\ldots,\frac{\partial\rho_A}{\partial x_n},\frac{\partial\rho_A}{\partial y_n}\right)
=
\frac{2}{\Vert\nabla\rho_A\Vert}\left(\frac{\partial\rho_A}{\partial \bar z_1},\ldots,\frac{\partial\rho_A}{\partial \bar z_n}\right)\,.
\end{equation*}
so that , if $\iota_{\partial A}:\partial A\to\mathbb C^n$ is the inclusion mapping, the volume form $\upsilon_{\partial A}$\label{vformpA} on $\partial A$ is given by
\begin{align*}
\upsilon_{\partial A}
&=
\iota_{\partial A}^*\left(\frac{*(\de \rho_A)}{\Vert\nabla\rho_A\Vert}\right)
\\
&=
\iota_{\partial A}^*\left(\frac{1}{\Vert\nabla\rho_A\Vert}\sum_{\ell=1}^n \frac{\partial\rho_A}{\partial x_\ell}\upsilon_{2n}[x_\ell]-\frac{\partial\rho_A}{\partial y_\ell}\upsilon_{2n}[y_\ell]\right)
\\
&=
\iota_{\partial A}^*\left(\frac{2}{\Vert\nabla\rho_A\Vert}\sum_{\ell=1}^n \frac{\partial\rho_A}{\partial\bar z_\ell}\upsilon_{2n}[z_\ell]-\frac{\partial\rho_A}{\partial z_\ell}\upsilon_{2n}[\bar z_\ell]\right)\,.
\end{align*}
The unit vector $\nu_{\partial A}(z)$ is, by its definition, orthogonal to the tangent hyperplane $T_z\partial A$ and the latter subspace admits a maximal complex subspace $T_z^\mathbb C\partial A$, which has real codimension $1$ in $T_z\partial A$. It follows that there is a unit vector in $T_z\partial A$ spanning the orthocomplement of $T_z^\mathbb C\partial A$ in $T_z\partial A$, i.e. spanning the so called \textit{characteristic direction} of $T_z\partial A$. The real plane spanned by $\nu_{\partial A}(z)$ and the characteristic direction of $T_z\partial A$ is just the complex line $(T_z\partial A^\mathbb C)^{\perp_\mathbb C}$, so that the characteristic direction of $T_z\partial A$ is indeed spanned by $i\nu_{\partial A}(z)$. 

The central object in the definition of Kazarnovski{\v\i} pseudovolume on the class ${\mathcal S}^2(\mathbb C^n)$ is the differential form $\alpha_{\partial A}$\label{alphaf} defined, for every $z\in\partial A$ and every $\zeta\in T_z\partial A$, as
\begin{equation*}
(\alpha_{\partial A})_z(\zeta)=\re\langle\zeta,i\nu_{\partial A}(z)\rangle=\frac{i}{2}\,\left(\langle\nu_{\partial A}(z),\zeta\rangle-\langle\zeta,\nu_{\partial A}(z)\rangle\right)\,.
\end{equation*}
The real number $(\alpha_{\partial A})_z(\zeta)$ measures the projection of $-i\zeta$ along the outward unit normal vector $\nu_{\partial A}(z)$ and, as such, it vanishes if and only if $\zeta\in T^\mathbb C_z\partial A$. Equivalently, $(\alpha_{\partial A})_z(\zeta)$ measures the projection of $\zeta$ on the characteristic direction of $T_z\partial A$.

In global coordinates, this form is given by
\begin{align*}
\alpha_{\partial A}
&=
\iota_{\partial A}^*\left(
\frac{1}{\Vert\nabla\rho_A\Vert}
\sum_{\ell=1}^n\frac{\partial\rho_A}{\partial x_\ell}\de y_\ell-\frac{\partial\rho_A}{\partial y_\ell}\de x_\ell\right)
\\
&=
\iota_{\partial A}^*\left(
\frac{\dc\rho_A}{\Vert\nabla\rho_A\Vert}\right)
\\
&=
\iota_{\partial A}^*\left(
\frac{i(\bar\partial\rho_A -\partial\rho_A)}{\Vert\nabla\rho_A\Vert}\right)
\\
&=
\iota_{\partial A}^*\left(
\frac{2i\bar\partial\rho_A }{\Vert\nabla\rho_A\Vert}\right)
\,,
\end{align*}
where the last equality follows from the relation $\iota_{\partial A}^*\de\rho_A=\iota_{\partial A}^*(\partial\rho_A+\partial\rho_A)=0$, on $T_z\partial A$.
The preceding expressions of $\alpha_{\partial A}$ have some remarkable consequences, namely
\begin{equation}\label{bforma}
\de\alpha_{\partial A}
=
\iota_{\partial A}^*\left(
\Vert\nabla\rho_A\Vert\de\left(\frac{1}{\Vert\nabla\rho_A\Vert}\right)\wedge\alpha_{\partial A}+\frac{\de\dc\rho_A}{\Vert\nabla\rho_A\Vert}\right)\,,
\end{equation}
\begin{align}
\alpha_{\partial A}\wedge \left(\de\alpha_{\partial A}\right)^{\wedge(n-1)}
&=
\iota_{\partial A}^*\left[\alpha_{\partial A}\wedge\left(\frac{\de\dc\rho_A}{\Vert\nabla\rho_A\Vert}\right)^{\wedge(n-1)}
\right]
\nonumber
\\
&=
\iota_{\partial A}^*\left[
\alpha_{\partial A}\wedge
\left(\frac{2i}{\Vert\nabla\rho_A\Vert}\right)^{n-1}
(\partial\bar\partial\rho_A)^{\wedge(n-1)}\right]\label{de-alfa}
\end{align}
and
\begin{equation}\label{ottima}
\upsilon_{\partial A}
=
\frac{1}{(n-1)!}\,
\iota_{\partial A}^*
\left[
\alpha_{\partial A}\wedge
\left(\frac{i}{2}\sum_{\ell=1}^n\de z_\ell\wedge\de\bar z_\ell\right)^{\wedge(n-1)}
\right].
\end{equation}
Observe that $\de\alpha_{\partial A}$ is defined only almost everywhere on $\partial A$ because, although $\partial A$ is a ${\mathcal C}^1$-hyper\-surface, it admits only a piecewise ${\mathcal C}^2$ structure.
\begin{rmq}\label{autovalori}
Observe that, as $(\alpha_{\partial A})_z$ vanishes on $T_z^\mathbb C\partial A$, the only part of $(\de\alpha_{\partial A})_z$ giving a non zero contribution to the product $(\alpha_{\partial A}\wedge(\de\alpha_{\partial A})^{\wedge(n-1)})_z$ is just its restriction to $T_z^\mathbb C\partial A$. It follows that the value of the $(2n-1)$-form $\alpha_{\partial A}\wedge(\de\alpha_{\partial A})^{\wedge(n-1)}$ at the point $z\in\partial A$ is proportional to the value of the volume form $\upsilon_{\partial A}$ at the same point and the proportionality constant equals $(n-1)!(2i/\Vert\nabla\rho_A(z)\Vert)^{n-1}$ times the product $k_1(z)\cdots k_{n-1}(z)$ of the eighenvalues of the complex hessian matrix $\partial\bar\partial \rho_A$ restricted to $T_z^\mathbb C\partial A$.
\end{rmq}

\begin{defn}\label{prima}
Let $A\in {\mathcal S}^2(\mathbb C^n)$. The \textbf{$\bm n$-dimensional Kazarnovski{\v\i} pseudovolume}\label{nkaza} (or simply the $n$-pseudovolume) of $A$ is the real number $P_n(A)$ defined as
\begin{equation}\label{alfa}
P_n(A)=\frac{1}{n!\varkappa_n}\int_{\partial A} \alpha_{\partial A}\wedge (\de\alpha_{\partial A})^{\wedge(n-1)}\,.
\end{equation}
\end{defn}

\begin{lem}\label{semip}
Let $A\in {\mathcal S}^2(\mathbb C)$. Then the $1$-dimensional pseudovolume is the semi-perimeter of $A$.
\end{lem}
\noindent
\pf
In this case we have
\begin{equation*}
\alpha_{\partial A}
=
\iota_{\partial A}^*\left[\frac{i}{\Vert\nabla\rho_A\Vert}\left(\frac{\partial\rho_A}{\partial\bar z}\de\bar z-\frac{\partial \rho_A}{\partial z}\de z\right)\right]
=
\upsilon_{\partial A}\,.
\end{equation*}
As $\varkappa_1=2$, $P_1(A)$ equals the half of the integral of $\upsilon_{\partial A}$ over $\partial A$, i.e the semi-perimeter.
\findim

\begin{rmq} The coefficient in the definition of $P_n$.

\noindent
In the english version of the original article \cite{Ka1}, the coefficient in the formula \eqref{alfa} is denoted $c_n$ and about this coefficient one reads: \textit{``..... $c_n=1/(\sigma_n n!)$, where $\sigma_n$ is the volume of the unit $n$-sphere." } 

\noindent
In this point there might be a problem with the translation from Russian into English. Infact $\sigma_n$ must be the volume of the unit $n$-ball (i.e. our $\varkappa_n$), otherwise \lemref{semip}, (whose statement is present in the first page of~\cite{Ka1}), should have claimed that $P_1(A)$ equals $(1/\pi)$ times the semi-perimeter of $A$.
\end{rmq}

\begin{lem}\label{transl}
The $1$-form $\alpha_{\partial A}$ is translation invariant. In particular, the pseudovolume is such.
\end{lem}
\noindent
\pf 
If $F:\mathbb C^n\rightarrow\mathbb C^n$ is a translation, for every $z\in\partial A$ and $\zeta\in T_z\partial A$, one has $\nu_{\partial F(A)}(F(z))=\nu_{\partial A}(z)$ and
\begin{align*}
(F^*\alpha_{\partial F(A)})_{z}(\zeta)
&=
(\alpha_{\partial F(A)})_{F(z)}(\de F_z(\zeta))
\\
&=
(\alpha_{\partial F(A)})_{F(z)}(\zeta)
\\
&=
\re\langle \zeta,i\nu_{\partial F(A)}(F(z))\rangle
\\
&=
\re\langle \zeta,i\nu_{\partial A}(z)\rangle
\\
&=
(\alpha_{\partial A})_z(\zeta)\,,
\end{align*}
so that
\begin{align*}
P_n(T(A))
&=
\frac{1}{n!\varkappa_n}\int_{\partial F(A)}
\alpha_{\partial F(A)}\wedge(\de\alpha_{\partial F(A)})^{\wedge(n-1)}
\\
&=
\frac{1}{n!\varkappa_n}\int_{\partial A}
F^*\left(\alpha_{\partial F(A)}\wedge(\de\alpha_{\partial F(A)})^{\wedge(n-1)}\right)
\\
&=
P_n(A)\,.
\end{align*}
The proof is complete.
\findim

\begin{lem}\label{buo}
The $1$-form $\alpha_{\partial A}$ is unitarily invariant. In particular, the pseudovolume is such.
\end{lem}
\noindent\pf 
If $F:\mathbb C^n\rightarrow\mathbb C^n$ is a unitary transformation, one has $\rho_{F(A)}=\rho_A\circ F^{-1}=\rho_A\circ {}^t\overline{F}$. If $w_1,\ldots,w_n$ are the coordinates in the range of $F$, by the chain rule, one gets 
\begin{equation*}
\nabla\rho_{F(A)}(F(z))
=
\overline{\nabla\rho_A}(z)\cdot\left(\frac{\partial(F^{-1})_j}{\partial\bar w_\ell}\left(F(z)\right)\right)_{1\leq j,\ell\leq n}
+
\nabla\rho_A(z)\cdot\left(\frac{\partial(\overline{F^{-1}})_j}{\partial \bar w_\ell}\left(F(z)\right)\right)_{1\leq j,\ell\leq n}\,.
\end{equation*}
As $F^{-1}$ is $\mathbb C$-linear, (in particular holomorphic), by confusing $F$ with its jacobian matrix one gets
\begin{equation*}
\frac{\partial(F^{-1})_j}{\partial\bar w_\ell}=0
\qquad
{\rm and}
\qquad
\frac{\partial\big(\overline{F^{-1}}\big)_j}{\partial \bar w_\ell}
=
\overline{\left(\frac{\partial(F^{-1})_j}{\partial w_\ell}\right)}=(\overline{F^{-1}})_{j,\ell}=({}^t\!F)_{j,\ell}
\,,
\end{equation*}
for every $1\leq j,\ell\leq n$, then
\begin{equation*}
\nabla\rho_{F(A)}(F(z))
=
\nabla\rho_A(z)\cdot{}^t\!F
=
F\cdot{}^t\nabla\rho_A(z)
=
F(\nabla\rho_A(z))\,,
\end{equation*}
and so 
\begin{align}
\nu_{\partial F(A)}(F(z))
&=
\frac{\nabla\rho_{F(A)}(F(z))}{\Vert \nabla\rho_{F(A)}(F(z))\Vert}\nonumber
\\
&=
\frac{F(\nabla\rho_A(z))}{\Vert F(\nabla\rho_A(z))\Vert}\nonumber
\\
&=
\frac{F(\nabla\rho_A(z))}{\Vert \nabla\rho_A(z)\Vert}\nonumber
\\
&=
F\left(\frac{\nabla\rho_A(z)}{\Vert \nabla\rho_A(z)\Vert}\right)\nonumber
\\
&=
F(\nu_{\partial A}(z))\label{bu}\,,
\end{align}
for every $z\in\partial A$.
It follows that, for any $\zeta\in T_z\partial A$,
\begin{align*}
(F^*\alpha_{\partial F(A)})_{z}(\zeta)
&=
(\alpha_{\partial F(A)})_{F(z)}(\de F_z(\zeta))
\\
&=
(\alpha_{\partial F(A)})_{F(z)}(F(\zeta))
\\
&=
\re\langle F(\zeta),i\nu_{\partial F(A)}(F(z))\rangle
\\
&=
\re\langle F(\zeta),iF(\nu_{\partial A}(z))\rangle
\\
&=
\re\langle F(\zeta),F(i\nu_{\partial A}(z))\rangle
\\
&=
\re\langle \zeta,i\nu_{\partial A}(z)\rangle
\\
&=
(\alpha_{\partial A})_z(\zeta)\,,
\end{align*}
so that
\begin{align*}
P_n(F(A))
&=
\frac{1}{n!\varkappa_n}\int_{\partial F(A)}
\alpha_{\partial F(A)}\wedge(\de\alpha_{\partial F(A)})^{\wedge(n-1)}
\\
&=
\frac{1}{n!\varkappa_n}\int_{\partial A}
F^*\left(\alpha_{\partial F(A)}\wedge(\de\alpha_{\partial F(A)})^{\wedge(n-1)}\right)
\\
&=
P_n(A)\,.
\end{align*}
The proof is complete.
\findim

\begin{lem}
The $1$-form $\alpha_{\partial A}$ is (positively) $1$-homogeneous and so the pseudovolume is (positively) $n$-homogeneous.
\end{lem}
\noindent\pf
Let $\lambda\in\mathbb R_{>0}$ and let $F:\mathbb C^n\rightarrow\mathbb C^n$ be the homothety corresponding to the multiplication by $\lambda$. For every $z\in\partial A$ and $\zeta\in T_z\partial A$, one has $\nu_{\partial F(A)}(F(z))=\nu_{\partial \lambda A}(\lambda z)=\nu_{\partial A}(z)$ and
\begin{align*}
(F^*\alpha_{\partial F(A)})_{z}(\zeta)
&=
(\alpha_{\partial F(A)})_{F(z)}(\de F_z(\zeta))
\\
&=
(\alpha_{\partial \lambda A})_{\lambda z}(\lambda\zeta)
\\
&=
\re\langle \lambda\zeta,i\nu_{\partial \lambda A}(\lambda z)\rangle
\\
&=
\re\langle \lambda\zeta,i\nu_{\partial A}(z)\rangle
\\
&=
\lambda\re\langle \zeta,i\nu_{\partial A}(z)\rangle
\\
&=
\lambda(\alpha_{\partial A})_z(\zeta)\,,
\end{align*}
so that
\begin{align*}
P_n(\lambda A)
&=
\frac{1}{n!\varkappa_n}\int_{\partial \lambda(A)}
\alpha_{\partial F(A)}\wedge(\de\alpha_{\partial F(A)})^{\wedge(n-1)}
\\
&=
\frac{1}{n!\varkappa_n}\int_{\partial A}
F^*\left(\alpha_{\partial F(A)}\wedge(\de\alpha_{\partial F(A)})^{\wedge(n-1)}\right)
\\
&=
\lambda^n P_n(A)\,.
\end{align*}
The proof is complete.
\findim

When $n>1$, an alternative formula for the pseudovolume involves the second quadratic form of $\partial A$. Recall that, the \textit{second quadratic form} of $\partial A$ is the twice covariant tensor $\sec_{\partial A}$\label{second} on $T\partial A$ whose value at a point $z\in\partial A$ is given, for any $\zeta,\eta\in T_z\partial A$, by
\begin{equation*}
(\sec_{\partial A})_z(\zeta,\eta)=\re\langle (\de\nu_{\partial A})_z(\zeta),\eta\rangle\,,
\end{equation*}
with
\begin{align*}
(\de\nu_{\partial A})_z(\zeta)
&=
\de\left(
\frac{\nabla\rho_A}{\Vert\nabla\rho_A\Vert}
\right)(\zeta)
\\
&=
\left(
\de\left(
\frac{2}{\Vert\nabla\rho_A\Vert}\frac{\partial\rho_A}{\partial \bar z_1}
\right)(\zeta),
\ldots,
\de\left(
\frac{2}{\Vert\nabla\rho_A\Vert}\frac{\partial\rho_A}{\partial \bar z_n}
\right)(\zeta)
\right).\qquad
\end{align*} 
Notice that the second quadratic form of $\partial A$ behaves well with respect to unitary transformations. Infact, if $F:\mathbb C^n\rightarrow \mathbb C^n$ is such a mapping, then, by virtue of \eqref{bu},
\begin{align*}
F^*(\sec_{\partial F(A)})_z(\zeta,\eta)
&=
(\sec_{\partial F(A)})_{F(z)}(F(\zeta),F(\eta))
\\
&=
\re\langle(\de\nu_{\partial F(A)})_{F(z)}(F(\zeta)),F(\eta)\rangle
\\
&=
\re\langle F\circ(\de\nu_{\partial A})_z(\zeta),F(\eta)\rangle
\\
&=
\re\langle(\de\nu_{\partial A})_z(\zeta),\eta\rangle\
\\
&=
(\sec_{\partial A})_z(\zeta,\eta)\,.
\end{align*}
By restricting $(\sec_{\partial A})_z$ to complex tangent vectors $\zeta,\eta\in T_z^{\mathbb C}\partial A$, one can deduce an $\mathbb R$-bilinear symmetric form $\sec_{\partial A}^{\mathbb C}$\label{secondc} by setting
\begin{equation*}
(\sec_{\partial A}^{\mathbb C})_z(\zeta,\eta)=\frac{1}{2}\big[(\sec_{\partial A})_z(\zeta,\eta)+(\sec_{\partial A})_z(i\zeta,i\eta)\big]\,.
\end{equation*}
Since $(\sec_{\partial A}^{\mathbb C})_z(\zeta,\eta)=(\sec_{\partial A}^{\mathbb C})_z(i\zeta,i\eta)$, it follows that
\begin{equation*}
{\mathcal L}_{\partial A,z}(\zeta,\eta)=(\sec_{\partial A}^{\mathbb C})_z(\zeta,\eta)+i(\sec_{\partial A}^{\mathbb C})_z(\zeta,i\eta)\,,
\end{equation*}
defines a hermitian form on $T_z^{\mathbb C}\partial A$ known as the \textit{Levi form} of $\partial A$ at the point $z$. Notice that, for $n=1$, $T_z^\mathbb C\partial A=\{0\}$, so ${\mathcal L}_{\partial A,z}\equiv 0$ in this case.

In order to find a simpler expression of ${\mathcal L}_{\partial A,z}$, let us carry on some computation:
\begin{align*}
2{\mathcal L}_{\partial A,z}(\zeta,\eta)
&=
\re\langle(\de\nu_{\partial A})_z(\zeta),\eta\rangle+\re\langle(\de\nu_{\partial A})_z(i\zeta),i\eta\rangle
\\
&+
i\re\langle(\de\nu_{\partial A})_z(\zeta),i\eta\rangle+i\re\langle(\de\nu_{\partial A})_z(i\zeta),i^2\eta\rangle
\qquad\;
\\
&=
\re\langle(\de\nu_{\partial A})_z(\zeta),\eta\rangle+\re\langle(\de\nu_{\partial A})_z(i\zeta),i\eta\rangle
\\
&+
i\im\langle(\de\nu_{\partial A})_z(\zeta),\eta\rangle+i\im\langle(\de\nu_{\partial A})_z(i\zeta),i\eta\rangle
\\
&=
\langle(\de\nu_{\partial A})_z(\zeta),\eta\rangle+\langle(\de\nu_{\partial A})_z(i\zeta),i\eta\rangle
\\
&=
\langle(\de\nu_{\partial A})_z(\zeta)-i(\de\nu_{\partial A})_z(i\zeta),\eta\rangle
\\
&=
\langle(\partial\nu_{\partial A})_z(\zeta)+(\bar\partial\nu_{\partial A})_z(\zeta)-i(\partial\nu_{\partial A})_z(i\zeta)-i(\bar\partial\nu_{\partial A})_z(i\zeta),\eta\rangle
\\
&=
\langle(\partial\nu_{\partial A})_z(\zeta)+(\bar\partial\nu_{\partial A})_z(\zeta)-i^2(\partial\nu_{\partial A})_z(\zeta)+i^2(\bar\partial\nu_{\partial A})_z(\zeta),\eta\rangle
\\
&=
2\langle(\partial\nu_{\partial A})_z(\zeta),\eta\rangle\,.
\end{align*}
For any $k=1,\ldots,n$, the $k$-th component of $\partial\nu_{\partial A}$ equals
\begin{align*}
\partial
\left(
\frac{2}{\Vert\nabla\rho_A\Vert}\frac{\partial\rho_A}{\partial \bar z_k}
\right)
&=
2
\sum_{\ell=1}^n
\left(
\frac{\partial\Vert\nabla\rho_A\Vert^{-1}}{\partial z_\ell}
\frac{\partial\rho_A}{\partial\bar z_k}
+
\frac{1}{\Vert\nabla\rho_A\Vert}
\frac{\partial^2\rho_A}{\partial z_\ell\partial\bar z_k}
\right)\de z_\ell\,,
\end{align*}
so that
\begin{align*}
{\mathcal L}_{\partial A,z}(\zeta,\eta)
&=
2
\sum_{k=1}^n
\sum_{\ell=1}^n
\left(
\frac{\partial\Vert\nabla\rho_A\Vert^{-1}}{\partial z_\ell}
\frac{\partial\rho_A}{\partial\bar z_k}
+
\frac{1}{\Vert\nabla\rho_A\Vert}
\frac{\partial^2\rho_A}{\partial z_\ell\partial\bar z_k}
\right)\zeta_\ell\bar\eta_k
\\
&=
\langle\nabla\rho_A,\eta\rangle
\sum_{\ell=1}^n
\frac{\partial\Vert\nabla\rho_A\Vert^{-1}}{\partial z_\ell}
\zeta_\ell
+
\frac{2}{\Vert\nabla\rho_A\Vert}
\sum_{\ell,k=1}^n
\frac{\partial^2\rho_A}{\partial z_\ell\partial\bar z_k}
\zeta_\ell\bar\eta_k
\\
&=
\frac{2}{\Vert\nabla\rho_A\Vert}
\sum_{\ell,k=1}^n
\frac{\partial^2\rho_A}{\partial z_\ell\partial\bar z_k}
\zeta_\ell\bar\eta_k\,,
\end{align*}
as $\langle\nabla\rho_A,\eta\rangle=0$ for $\eta\in T_z^\mathbb C\partial A$. It follows that the Levi form of $\partial A$ has the following expression:
\begin{equation*}
{\mathcal L}_{\partial A}
=
\frac{2}{\Vert\nabla\rho_A\Vert}
\sum_{\ell,k=1}^n
\frac{\partial^2\rho_A}{\partial z_\ell\partial\bar z_k}
\de z_\ell\otimes\de\bar z_k\,.
\end{equation*}
It is important to remark that, though involving global coordinates, the preceding expression is always restricted to $T_z^\mathbb C \partial A$, on which $\partial \rho_A=\bar\partial \rho_A=0$. It follows that ${\mathcal L}_{\partial A}$ actually depends on at most $(n-1)$ complex coordinates. 

Thanks to \eqref{bu} and the unitary invariance of $\sec_{\partial A}$, it follows that
$
F^*{\mathcal L}_{\partial F(A)}
=
{\mathcal L}_{\partial A}\,,
$ for any unitary transformation $F:\mathbb C^n\rightarrow\mathbb C^n$.

For any $z\in\partial A$, let ${\mathcal E}_{\partial A}(z)$\label{kappalevi} denote\footnote{In the paper \cite{Ka1} the author emploies the notation $K_{\partial A}$ but we prefer ${\mathcal E}_{\partial A}$ to avoid confusion with dual cones.} the product of the eigenvalues of ${\mathcal L}_{\partial A,z}$. Remark that ${\mathcal E}_{\partial A}(z)$ equals $(2/\Vert\nabla\rho_A(z)\Vert)^{n-1}$ times the product $k_1(z)\cdots k_{n-1}(z)$ of the eigenvalues of the complex Hessian of $\rho_A$ restricted to $T_z^\mathbb C\partial A$. We recall that a subset $A\in {\mathcal S}^2(\mathbb C^n)$ is \textit{(Levi) pseudoconvex} if,  for every $z\in\partial A$, the Levi form ${\mathcal L}_{\partial A,z}$ is semi-positive definite, i.e. ${\mathcal L}_{\partial A,z}(\zeta,\eta)\geq 0$, for every $z\in\partial A$ and any $\zeta,\eta\in T_z^\mathbb C\partial A$. It is also useful to recall that a subset $A\in {\mathcal S}^2(\mathbb C^n)$ is \textit{strictly (Levi) pseudoconvex} if,  for every $z\in\partial A$, the Levi form ${\mathcal L}_{\partial A,z}$ is positive definite, i.e. ${\mathcal L}_{\partial A,z}(\zeta,\eta)> 0$, for every $z\in\partial A$ and any $\zeta,\eta\in T_z^\mathbb C\partial A\setminus\{0\}$. Any (strictly) convex subset $A\in{\mathcal S}^2(\mathbb C^n)$ is (strictly) pseudoconvex, but the converse is generally false already for $n=1$, in which case the class of $k$-regular pseudoconvex subsets coincide with the space ${\mathcal S}^k(\mathbb C)$. Of course, if $A\in {\mathcal S}^2(\mathbb C^n)$ is pseudoconvex, then ${\mathcal E}_{\partial A}(z)\geq 0$ for every $z\in\partial A$, with a strict inequality when $A$ is strictly pseudoconvex.

\begin{thm}\label{thmlevi}
Let $n\in\mathbb N\setminus\{0,1\}$ and $A\in {\mathcal S}^2(\mathbb C^n)$. Then 
\begin{equation}
P_n(A)=\frac{2^{n-1}(n-1)!}{n!\varkappa_n}\int_{\partial A}{\mathcal E}_{\partial A}\upsilon_{\partial A}\,.\label{levi}
\end{equation}
In particular, if $A$ is pseudoconvex then $P_n(A)\geq 0$. 
\end{thm}
\noindent\pf
Observe that both the forms $\alpha_{\partial A}\wedge(\de\alpha_{\partial A})^{\wedge(n-1)}$ and ${\mathcal E}_{\partial A}(z)\upsilon_{\partial A}$ have maximal degree, so they have to be proportional. 
By \eqref{ottima}, \rmqref{autovalori} and \eqref{de-alfa}, one has
\begin{align*}
&
2^{n-1}(n-1)! {\mathcal E}_{\partial A}\upsilon_{\partial A}
\\
&=
i^{n-1} {\mathcal E}_{\partial A}\,\iota^*_{\partial A}\left[\alpha_{\partial A}\wedge\left(\sum_{\ell=1}^n \de z_\ell\wedge\de \bar z_\ell\right)^{\wedge(n-1)}\right]
\\
&=
\iota^*_{\partial A}\left[
\alpha_{\partial A}\wedge
\left(\frac{2i}{\Vert\nabla\rho_A\Vert}\right)^{n-1}
k_1(z)\cdots k_{n-1}(z)
\left(\sum_{\ell=1}^n \de z_\ell\wedge\de \bar z_\ell\right)^{\wedge(n-1)}\right]
\\
&=
\iota^*_{\partial A}\left[
\alpha_{\partial A}\wedge
\left(\frac{2i}{\Vert\nabla\rho_A\Vert}\right)^{n-1}
(\partial\bar\partial\rho_A)^{\wedge(n-1)}
\right]
\\
&=
\alpha_{\partial A}\wedge(\de\alpha_{\partial A})^{\wedge(n-1)}
\,.
\end{align*}
When $A$ is pseudoconvex, the integrand in \eqref{levi} is non negative, so that $P_n(A)\geq 0$.
\findim

\begin{rmq}
The dimensional coefficient $(n-1)!$ in formula \eqref{levi}.

\noindent
Kazarnovski{\v\i} in his paper \cite{Ka1} states formula \eqref{levi} without the dimensional coefficient $(n-1)!$. The reason should be in the definition of the volume form $\upsilon_{\partial A}$ or in that of the Levi form ${\mathcal L}_{\partial A}$, on which the reference \cite{Ka1} lacks a fully detailed description. Any way, as shown by \exeref{pseudoball}, our presentation is correct.
\findim
\end{rmq}

\begin{thm}\label{continuo}
The Kazarnovski{\v\i} $n$-pseudovolume is continuous on ${\mathcal K}^2(\mathbb C^n)$ for the Hausdorff metric.
\end{thm}
\noindent\pf Let $A\in{\mathcal K}^2(\mathbb C^n)$ and let $\{A_m\}_{m\in \mathbb N}$ be a sequence of convex bodies from ${\mathcal K}^2(\mathbb C^n)$ converging to $A$ in the Hausdorff metric. We want to show that $\lim_{m\to+\infty} P_n( A_m)=P_n(A)$. 
Let $z_o$ an interior point of $A$ and $A_m$, for sufficiently big $m\in\mathbb N$. For any such $m$, let $\gamma_m:\partial A\to\partial A_m$ map the point $z\in\partial A$ to the point $\gamma_m(z)\in\partial A_m$ at which the ray issuing from $z_o$ and pointing to $z$ intersects $\partial A_m$. The $2$-regularity of $A$ and $A_m$ imply that $\gamma_m$ is an orientation preserving diffeomorphism, so that
\begin{align*}
\int_{\partial A_m} \alpha_{\partial A_m}\wedge (\de \alpha_{\partial A_m})^{\wedge(n-1)} 
&=
\int_{\partial A} \gamma_m^*
\left(
\alpha_{\partial A_m}\wedge (\de \alpha_{\partial A_m})^{\wedge(n-1)} 
\right)
\\
&=
\int_{\partial A} \gamma_m^*\alpha_{\partial A_m}\wedge (\de \gamma_m^* \alpha_{\partial A_m})^{\wedge(n-1)} \,.
\end{align*}
For every $z\in\partial A$ and any $\zeta\in T_z\partial A$, 
\begin{equation*}
(\gamma_m^*\alpha_{\partial A_m})_z(\zeta)
=
\alpha_{\partial A_m,\gamma_m(z)}((\de\gamma_m)_z(\zeta))
=
\re\langle (\de \gamma_m)_z(\zeta),i\nu_{\partial A_m}(\gamma_m(z))\rangle
\end{equation*}
As $m\to+\infty$, the mapping $\gamma_m$ approaches the identity mapping on $\partial A$, hence $\gamma_m^*\alpha_{\partial A_m}$ approaches $\alpha_{\partial A}$ and the theorem follows.
\findim

\begin{rmq} The Kazarnovski{\v\i} $n$-pseudovolume is continuous on ${\mathcal S}^2(\mathbb C^n)$ for the Hausdorff metric.

\noindent
\thmref{continuo} can be extended to the whole ${\mathcal S}^2(\mathbb C^n)$. The proof of this more general result would differ from the one presented above just in the choice of the diffeomorphism $\gamma_m$. 
\findim
\end{rmq}

The following theorem provides an easy formula for computing the $n$-pseudovolume on ${\mathcal P}_\circ(\mathbb C^n)$.

\begin{thm}\label{azzok0}
Let $n\in\mathbb N^*$, $\Gamma\in {\mathcal P}(\mathbb C^n)$ and $\varepsilon>0$. For each $0\leq k<2n$ and any $\Delta\in{\mathcal B}(\Gamma,k)$,  
\begin{equation}\label{f-azzok0}
\frac{1}{n!2^{n-k}\varkappa_{2n-k}}
\int_{\Delta+(K_\Delta\cap\partial \varepsilon B_{2n})}
\alpha_{\partial(\Gamma)_\varepsilon}\wedge(\de \alpha_{\partial(\Gamma)_\varepsilon})^{\wedge(n-1)}
=
\varepsilon^{n-k}
\varrho(\Delta)\vol_k(\Delta)\psi_\Gamma(\Delta)\,,
\end{equation}
then, by neglecting non equidimensional faces,
\begin{equation}\label{f-azzok1}
P_n((\Gamma)_\varepsilon)
=
\sum_{k=0}^n\frac{2^{n-k}\varkappa_{2n-k}}{\varkappa_n}\,{\textsl v}^\varrho_k(\Gamma)\varepsilon^{n-k}\,.
\end{equation}
\end{thm}
\noindent\pf
The hypersurface $\Delta+(K_\Delta\cap\partial \varepsilon B_{2n})$ has the structure of a cylinder on the $(2n-k-1)$-dimensional base $K_\Delta\cap\partial \varepsilon B_{2n}$. As the outward unit normal vector to $\Delta+(K_\Delta\cap\partial \varepsilon B_{2n})$ belongs to $E_\Delta^\perp$, it follows that $\alpha_{\partial(\Gamma)_\varepsilon}$ depends on $(2n-k)$ variables.
If $k>n$, one has $(2n-k)=n+(n-k)<n$ so that the $(2n-1)$-form $\alpha_{\partial(\Gamma)_\varepsilon}\wedge(\de \alpha_{\partial(\Gamma)_\varepsilon})^{\wedge(n-1)}$ depends on $(2n-k)<n\leq(2n-1)$ variables, i.e. 
$\alpha_{\partial(\Gamma)_\varepsilon}\wedge(\de \alpha_{\partial(\Gamma)_\varepsilon})^{\wedge(n-1)}\equiv 0$.
On the other hand, $\varrho(\Delta)=0$ if $k>n$, then the statement holds true for $n<k<2n$.

Let us now suppose $0\leq k\leq n $. In order to compute the integral on the left of \eqref{f-azzok0}, we first consider the case of a non equidimensional face $\Delta$. In this case $k=\dim_\mathbb R E_\Delta\geq 2$ and the Cauchy-Riemann dimension $k^\prime=\dim_\mathbb C E_\Delta^\mathbb C$ satisfies the relations $k\geq 2 k^\prime>k^\prime>0$. Up to a non singular $\mathbb R$-linear transformation, (affecting our computation by a non zero multiplicative constant), we may suppose that 
\begin{align*}
E_\Delta^\perp
&=
\{z\in\mathbb C^n\mid z_1=\ldots=z_{k^\prime}=\re z_{k^\prime+1}=\ldots=\re z_{k-k^\prime}=0\}
\\
&=
\{
z\in\mathbb C^n\mid x_1=\ldots=x_{k-k^\prime}=y_1=\ldots=y_{k^\prime}=0\}\,.
\end{align*}
For the sake of notation, set $\nu=\nu_{\partial(\Gamma)_\varepsilon}$. As $\nu(z)\in E_\Delta^\perp$, for every $z\in \Delta+ (K_\Delta+\partial\varepsilon B_{2n})$, we deduce that
\begin{equation*}
\alpha_{\partial(\Gamma)_\varepsilon}
=
\sum_{\ell=k-k^\prime+1}^n(\re\nu_\ell)\,\de y_\ell-\sum_{\ell=k^\prime+1}^n (\im\nu_\ell)\, \de x_\ell\,,
\end{equation*}
on $\Delta+ (K_\Delta\cap\partial\varepsilon B_{2n})$. Moreover $\nu(z)=\nu(z+u)$, for every $u\in E_\Delta$, so the partial derivatives of the components of $\nu$ with respect to $x_1,\ldots
x_{k-k^\prime}$ and $y_1,\ldots,y_{k^\prime}$ are all zero on $\Delta+ (K_\Delta\cap\partial\varepsilon B_{2n})$. This implies that $\de\alpha_{\partial(\Gamma)_\varepsilon}$ has the following expression
\begin{align*}
&
\sum_{\ell=k-k^\prime+1}^n
\left(\sum_{j=k-k^\prime+1}^n
\frac{\partial(\re\nu_\ell)}{\partial x_j}\de x_j\wedge\de y_\ell
+
\sum_{j=k^\prime+1}^n
\frac{\partial(\re\nu_\ell)}{\partial y_j}\de y_j\wedge\de y_\ell
\right)
\\
&-
\sum_{\ell=k^\prime+1}^n
\left(\sum_{j=k-k^\prime+1}^n
\frac{\partial(\im\nu_\ell)}{\partial x_j}\de x_j\wedge\de x_\ell
+
\sum_{j=k^\prime+1}^n
\frac{\partial(\im\nu_\ell)}{\partial y_j}\de y_j\wedge\de x_\ell
\right)\,,
\end{align*}
with $(k-k^\prime)>0$ and $k^\prime>0$. As a consequence the $(2n-1)$-form $\alpha_{\partial(\Gamma)_\varepsilon}\wedge(\de\alpha_{\partial(\Gamma)_\varepsilon})^{\wedge(n-1)}$ vanishes identically on $\Delta+ (K_\Delta\cap\partial\varepsilon B_{2n})$, because it depends on no more than $(2n-2k^\prime)$ variables and, of course, $2k^\prime>1$. We deduce that the statement of the theorem holds true for every face which is not equidimensional. 

Finally, if $\Delta\in{\mathcal B}_{\rm ed}(\Gamma,k)$, up to a non singular $\mathbb R$-linear transformation, say $L_\Delta$, we may suppose 
\begin{align*}
E_\Delta^\perp
&=
\{z\in\mathbb C^n\mid \re z_1=\ldots=\re z_k=0\}
\\
&=
\{z\in\mathbb C^n\mid x_1=\ldots=x_k=0\}\,.
\end{align*}
Indeed, by the unitary invariance of $\alpha_{\partial(\Gamma)_\varepsilon}$, we can make $\lin_\mathbb C E_\Delta=\{z\in\mathbb C^n\mid z_{k+1}=\ldots=z_n=0\}$, so we need to define  $L_\Delta$ just on the smaller subspace $\lin_\mathbb C E_\Delta$. Let $e_1,\ldots,e_k$ be the canonical basis of $\lin_\mathbb C E_\Delta=\mathbb C^k\times\{0\}^{n-k}$, let also $v_1,\ldots,v_k$ be an orthonormal basis of $E_\Delta$ (with respect to the scalar product $\re\langle\,,\rangle$) and consider the operator $L_\Delta:\lin_\mathbb C E_\Delta\to\lin_\mathbb C E_\Delta$ mapping $e_\ell$ to $v_\ell$. As $\Delta$ is equidimensional, $L_\Delta$ is an isomorphism of complex linear spaces, with $\det L_\Delta=\det(\langle v_\ell,e_j\rangle)$. If we look at $L_\Delta$ as a mapping on $\mathbb R^{2k}$, its jacobian equals $\vert \det(\langle v_\ell,e_j\rangle)\vert^2$ which, by virtue of \cororef{fighissimo}, is nothing but $\varrho(\Delta)$. 

From now on we should systematically consider the pullback of $\alpha_{\partial(\Gamma)_\varepsilon}$ via $L_\Delta$ and  change the names of variables, however, for the sake of notation, we will content ourselves to append the jacobian $\varrho(\Delta)$ to the computation of the integral on the left of \eqref{f-azzok0}, leaving the names of variables unchanged.

Arguing as before, in the simplified coordinates we now get
\begin{equation*}
\alpha_{\partial(\Gamma)_\varepsilon}
=
\sum_{\ell=k+1}^n(\re\nu_\ell)\,\de y_\ell-\sum_{\ell=1}^n (\im\nu_\ell)\, \de x_\ell\,
\end{equation*}
so that, on $\Delta+ (K_\Delta\cap\partial\varepsilon B_{2n})$, the form $\de\alpha_{\partial(\Gamma)_\varepsilon}$ is given by
\begin{align*}
&
\sum_{\ell=k+1}^n
\left(\sum_{j=k+1}^n
\frac{\partial(\re\nu_\ell)}{\partial x_j}\de x_j\wedge\de y_\ell
+
\sum_{j=1}^n
\frac{\partial(\re\nu_\ell)}{\partial y_j}\de y_j\wedge\de y_\ell
\right)
\\
&-
\sum_{\ell=1}^n
\left(\sum_{j=k+1}^n
\frac{\partial(\im\nu_\ell)}{\partial x_j}\de x_j\wedge\de x_\ell
+
\sum_{j=1}^n
\frac{\partial(\im\nu_\ell)}{\partial y_j}\de y_j\wedge\de x_\ell
\right)\,.
\end{align*}
Notice that in the present situation the $(2n-1)$-form $\alpha_{\partial(\Gamma)_\varepsilon}\wedge(\de\alpha_{\partial(\Gamma)_\varepsilon})^{\wedge(n-1)}$ depends on the full array of $2n$ variables. In order to complete the computation we need to specify a defining function of $(\Gamma)_\varepsilon$ or at least its restriction to the relatively open subset $\Delta+ K_\Delta$. Up to a translation of $\Gamma$ (which cannot affect the overall computation) we may suppose that $p_\Delta=0$ so, according to (\ref{funzione-def}), the restriction  $\rho_\Delta$ of $\rho_{(\Gamma)_\varepsilon}$ to $\Delta+ K_\Delta$, can be chosen as the function  
\begin{equation*}
\rho_\Delta(z)=-\varepsilon^2+\sum_{\ell=k+1}^n x_\ell^2+\sum_{\ell=1}^n y_\ell^2\,,
\end{equation*}
so that $\Vert\nabla\rho_\Delta(z)\Vert=2\varepsilon$ and 
$\nu(z)=
\varepsilon^{-1}
(0,y_1,\ldots,0,y_k,x_{k+1},y_{k+1},\ldots,x_n,y_n)
$, for every $z\in\Delta+(K_\Delta\cap\partial\varepsilon B_{2n})$. This yields
\begin{equation*}
\alpha_{\partial(\Gamma)_\varepsilon}
=
\varepsilon^{-1}
\left(
\sum_{\ell=k+1}^nx_\ell\de y_\ell-\sum_{\ell=1}^n y_\ell \de x_\ell
\right)
\,
\end{equation*}
and
\begin{equation*}
\de\alpha_{\partial(\Gamma)_\varepsilon}
=
\varepsilon^{-1}
\left(
\sum_{\ell=1}^k
\de x_\ell\wedge\de y_\ell
+
2\sum_{j=k+1}^n
\de x_\ell\wedge\de y_\ell
\right)
\,,
\end{equation*}
for every $z\in\Delta+(K_\Delta\cap\partial\varepsilon B_{2n})$. Observe that
\begin{equation*}
(\de\alpha_{\partial(\Gamma)_\varepsilon})^{\wedge n}
=
\varepsilon^{-n}2^{n-k}n!
\bigwedge_{\ell=1}^n\de x_\ell\wedge \de y_\ell
=
\varepsilon^{-n}2^{n-k}n!\upsilon_{2n}
\,,
\end{equation*}
so, by Stokes' theorem,
\begin{align*}
&
\frac{1}{n!2^{n-k}\varkappa_{2n-k}}
\int_{\Delta+(K_\Delta\cap\partial \varepsilon B_{2n})}
\alpha_{\partial(\Gamma)_\varepsilon}\wedge(\de \alpha_{\partial(\Gamma)_\varepsilon})^{\wedge(n-1)}
\\
&=
\frac{1}{n!2^{n-k}\varkappa_{2n-k}}
\int_{\Delta+(K_\Delta\cap\varepsilon B_{2n})}
(\de \alpha_{\partial(\Gamma)_\varepsilon})^{\wedge n}
\\
&=
\varrho(\Delta)
\frac{\varepsilon^{-n}}{\varkappa_{2n-k}}
\vol_{2n}(\Delta+(K_\Delta\cap\varepsilon B_{2n}))
\\
&=
\varrho(\Delta)
\frac{\varepsilon^{-n}}{\varkappa_{2n-k}}
\vol_k(\Delta)\vol_{2n-k}(K_\Delta\cap\varepsilon B_{2n})
\\
&=
\varrho(\Delta)
\frac{\varepsilon^{-n}}{\varkappa_{2n-k}}
\vol_k(\Delta)\varepsilon^{2n-k}\vol_{2n-k}(K_\Delta\cap B_{2n})
\\
&=
\varepsilon^{n-k}
\varrho(\Delta)\,
\vol_k(\Delta)\psi_\Gamma(\Delta)\,,
\end{align*}  
hence \eqref{f-azzok0} is proved. Notice that \eqref{f-azzok0} is a non trivial equality if and only if $\Delta$ is equidimensional. Moreover
\begin{equation*}
\partial(\Gamma)_\varepsilon=\bigcup_{k=0}^{2n-1}\bigcup_{\Delta\in{\mathcal B}(\Gamma,k)}\Delta+(K_\Delta\cap\partial \varepsilon B_{2n})\,,
\end{equation*}
so, by summing, as $k$ runs in the set $\{0,\ldots,2n-1\}$, the relations \eqref{f-azzok0} and neglecting non equidimensional faces, one gets \eqref{f-azzok1}. The proof is thus complete.
\findim

The notion of pseudovolume applies in particular to convex bodies from ${\mathcal K}^2_2(\mathbb C^n)\subset{\mathcal S}^2(\mathbb C^n)$. The following results of \cite{Ka1} imply a further formula for $P_n$ on ${\mathcal K}^2_1(\mathbb C^n)$.

\begin{lem}
Let $A\in {\mathcal K}^2_1(\mathbb C^n)$. Then the Gauss map $\nu_{\partial A}$ and the restriction  to $\partial B_{2n}$ of the gradient mapping $\nabla h_A$ are inverse diffeomorphisms. 
\end{lem}
\noindent\pf
Let us start with $\nu_{\partial A}:\partial A\rightarrow\partial B_{2n}$. We already know that the Gauss map is ${\mathcal C}^1$. For any $u\in\partial B_{2n}$, the fiber $\nu_{\partial A}^{-1}(u)$ is non-empty since it equals the exposed face $A\cap H_A(u)$. By the strict convexity of $A$, $A\cap H_A(u)$ is reduced to a single point, so $\nu_{\partial A}$ is bijective. The exposed face $A\cap H_A(u)$ is nothing but the sub-gradient of $h_A$ at the point $u$, so that
\begin{equation*}
\nabla h_A(u)\in\Subd h_A(u)=A\cap H_A(u)=\{\nu_{\partial A}^{-1}(u)\}\,,
\end{equation*} 
i.e. $\nu_{\partial A}^{-1}=\nabla h_A$. It follows that $\nu_{\partial A}^{-1}$ is also differentiable. \findim

\begin{lem}\label{koush}
For every $A\in {\mathcal K}^2_2(\mathbb C^n)$, the forms $\dc h_A$ and  $(\nu_{\partial A}^{-1})^*\alpha_{\partial A}$ are cohomologous.
\end{lem}
\noindent\pf 
Let $u\in\partial B_{2n}$ and $\zeta\in T_u\partial B_{2n}$ be fixed. Then
\begin{align*}
\left(\dc h_A\right)_u(\zeta)
&=
i\sum_{\ell=1}^n
\frac{\partial h_A}{\partial \bar u_\ell}(u)\bar\zeta_\ell
-
\frac{\partial h_A}{\partial u_\ell}(u)\zeta_\ell
\\
&=
(i/2)\left(\langle\nabla h_A(u),\zeta\rangle-\langle\zeta,\nabla h_A(u)\rangle\right)\,,
\end{align*}
whereas
\begin{align}
\left((\nu_{\partial A}^{-1})^*\alpha_{\partial A}\right)_u
&=
\left(\alpha_{\partial A}\right)_{\nu_{\partial A}^{-1}(u)}\left((\de\nu_{\partial A}^{-1})_u(\zeta)\right)
\\
&=
\re\langle \de(\nabla h_A)_u(\zeta),iu\rangle
\\
&=
\im \langle \de(\nabla h_A)_u(\zeta),u\rangle
\\
&=
(i/2)\left(\langle u,\de(\nabla h_A)_u(\zeta)\rangle-\langle \de(\nabla h_A)_u(\zeta),u\rangle\right)\label{koushnirenko}
\,. 
\end{align}
It follows that
\begin{align*}
&
\left(\dc h_A-(\nu_{\partial A}^{-1})^*\alpha_{\partial A}\right)_u(\zeta)
\\
&=
(i/2)\left[
\langle\nabla h_A(u),\zeta\rangle-\langle\zeta,\nabla h_A(u)\rangle-\langle u,\de(\nabla h_A)_u(\zeta)\rangle+\langle \de(\nabla h_A)_u(\zeta),u\rangle\right]\qquad\qquad
\\
&=
(i/2)\left[
\langle\nabla h_A(u),\zeta\rangle
+
\langle \de(\nabla h_A)_u(\zeta),u\rangle
-
\langle u,\de(\nabla h_A)_u(\zeta)\rangle
-
\langle\zeta,\nabla h_A(u)\rangle\right]
\\
&=
(i/2)\left[
\de\left(\langle\nabla h_A(u),u\rangle\right)-\de\left(\langle u,\nabla h_A(u)\rangle\right)
\right](\zeta)
\\
&=
\de\left[
\im\langle u,\nabla h_A(u)\rangle\right](\zeta)\,.
\end{align*}
The lemma is thus proved.\findim

\begin{rmq}
\lemref{koush} can be found in~\cite{Sha} pag.~199 but, although I cannot provide a precise reference, Kazarnovsk{\v\i} told me it's originally due to Koushnirenko.
\end{rmq}

\begin{lem}\label{hinv}
For every $A\in {\mathcal K}_1(\mathbb C^n)$, the form $\dc h_A$ is unitarily invariant.
\end{lem}
\noindent\pf
Let $F:\mathbb C^n\rightarrow\mathbb C^n$ a unitary transformation. Then $h_{F(A)}=h_A\circ F^{-1}$, $F(\partial B_{2n})=\partial B_{2n}$ and (arguing as in the proof of \lemref{buo}) $\nabla (h_{A}\circ F^{-1})(F(z))=F(\nabla h_{A}(z))$, so that, for every $z\in \partial B_{2n}$ and any $\zeta\in T_z\partial B_{2n}$,
\begin{align*}
F^*(\dc h_{F(A)})_z(\zeta)
&=
(\dc h_{F(A)})_{F(z)}(F(\zeta))
\\
&=
\frac{i}{2}
\left[
\langle\nabla h_{F(A)}(F(z)),F(\zeta)\rangle-\langle F(\zeta),\nabla h_{F(A)}(F(z))\rangle
\right]
\\
&=
\frac{i}{2}
\left[
\langle\nabla (h_{A}\circ F^{-1})(F(z)),F(\zeta)\rangle-\langle F(\zeta),\nabla (h_{A}\circ F^{-1})(F(z))\rangle
\right]
\\
&=
\frac{i}{2}
\left[
\langle F(\nabla h_{A}(z)),F(\zeta)\rangle-\langle F(\zeta),F(\nabla h_{A}(z))\rangle
\right]
\\
&=
\frac{i}{2}
\left[
\langle \nabla h_{A}(z),\zeta\rangle-\langle \zeta,\nabla h_{A}(z)\rangle
\right]
\\
&=
(\dc h_A)_z(\zeta)\,.
\end{align*} 
The proof is complete.
\findim

\begin{coro}\label{kaza}
For any $A\in{\mathcal K}_2^2(\mathbb C^n)$, one has
\begin{equation}\label{f-koush}
P_n(A)=\frac{1}{n!\varkappa_n}\int_{\partial B_{2n}}\dc h_A\wedge(\ddc h_A)^{\wedge(n-1)}=\frac{1}{n!\varkappa_n}\int_{B_{2n}}(\ddc h_A)^{\wedge n}\,.
\end{equation}
In particular $P_n(A)> 0$.
\end{coro}
\noindent\pf
By \lemref{koush}, 
\begin{equation*}
\frac{1}{n!\varkappa_n}\int_{\partial B_{2n}}\dc h_A\wedge(\ddc h_A)^{\wedge(n-1)}
=
\frac{1}{n!\varkappa_n}\int_{\partial B_{2n}}
\left(\nu_{\partial A}^{-1}\right)^*\left(\alpha_{\partial A}\wedge\left(\de\alpha_{\partial A}\right)^{\wedge(n-1)}\right)
=
P_n(A)\,.
\end{equation*}
Moreover, the only singularity of $h_A$ is at the origin, where the second derivatives of $h_A$ are $O(1/\Vert z\Vert)$, (because $h_A$ is positively homogeneous of degree $1$). It follows that the form $(\ddc h_A)^{\wedge n}$ has a locally integrable density and so, by Stokes' theorem, one obtains the right-most term in the equality \eqref{f-koush}. Since a convex body is pseudoconvex, \thmref{thmlevi} implies that $P_n(A)\geq 0$. However, the strict convexity of $A$ implies its strict pseudoconvexity, hence $P_n(A)$ is strictly positive.
\findim

\section{Kazarnovski{\v\i} pseudovolume on the class {\boldmath ${\mathcal K}(\mathbb C^n)$}}\label{ksez}

\cororef{kaza} suggests a way to extend the definition of $P_n$ to all convex bodies, even those with the worst boundary structure. The passage from ${\mathcal K}^\infty_1(\mathbb C^n)$ to ${\mathcal K}(\mathbb C^n)$ requires the regularization of the complex Monge-Amp\`ere operator as described in \cite{BT} or \cite{Dem}.

\begin{defn}\label{mkaza}
Let $A_1,\ldots,A_n\in{\mathcal K}(\mathbb C^n)$. The \textbf{Kazarnovski{\v\i} $\bm n$-dimensio\-nal mixed pseudovolume} of $A_1,\ldots,A_n$ is the (non-negative) real number, noted $Q_n(A_1,\ldots,A_n)$, given by
\begin{equation}\label{def-pseudo}
Q_n(A_1,\ldots,A_n)=\frac{1}{n!\varkappa_n}\int_{B_{2n}}\ddc h_{A_1}\wedge\ldots\wedge \ddc h_{A_n}\,.
\end{equation}
If $A=A_1,\ldots,A_n$, set 
\begin{equation}
P_n(A)=Q_n(A[n])
=
\frac{1}{n!\varkappa_n}\int_{B_{2n}}(\ddc h_{A})^{\wedge n}
\,,
\end{equation}
and call it the \textbf{Kazarnovski{\v\i} $\bm n$-dimensional pseudovolume} of $A$.
\end{defn}

\begin{rmq}
Kazarnovski{\v\i} $n$-dimensional mixed pseudovolume is well defined.

\noindent
All the currents involved in \defnref{mkaza} are positive, so that their wedge product is legitimate. Indeed, it is defined inductively by setting 
\begin{equation}\label{defMA}
\ddc h_{A_1}\wedge\ddc h_{A_2}\wedge\ldots\wedge \ddc h_{A_n}=\ddc(h_{A_1}\ddc h_{A_2}\wedge\ldots\wedge \ddc h_{A_n})\,.
\end{equation}
Although this definition is not symmetric in $A_1,\ldots,A_n$, it can be shown that the resulting current is such. Indeed, let $\varepsilon>0$ be fixed and, for any $1\leq\ell\leq n$, let us replace in the definition \eqref{defMA} the function $h_{A_\ell}$ with its regularization $h_{A_\ell}*\varphi_\varepsilon$. As shown in \cite{BT} or \cite{Dem}, the resulting relation yields an equality of smooth differential forms the left side of which is clearly symmetric. Since the form $\ddc (h_{A_1}*\varphi_\varepsilon)\wedge\ldots\wedge \ddc (h_{A_n}*\varphi_\varepsilon)$ is symmetric and weakly converges to the current $\ddc h_{A_1}\wedge\ldots\wedge \ddc h_{A_n}$, as $\varepsilon\to0$, it follows that the limiting current is symmetric too.

Positive currents have measure coefficients, so the current $\kappa=\ddc h_{A_1}\wedge\ldots\wedge \ddc h_{A_n}$ is a positive measure and the integral on the right side of \eqref{def-pseudo} stands for the $\kappa$-measure of the unit ball $B_{2n}$. In the sequel, the pairing between a current $T$\label{current} and a test form $\varsigma$\label{tform} will be denoted $\langle\!\langle T,\varsigma\rangle\!\rangle$\label{pairing}. The mixed pseudovolume of $A_1,\ldots,A_n$, can be computed as
\begin{equation*}
\frac{1}{n!\varkappa_n}
\lim_{m\to\infty}
\langle\!\langle 
\ddc h_{A_1}\wedge\ldots\wedge \ddc h_{A_n},\varsigma_m
\rangle\!\rangle\,,
\end{equation*}
where $(\varsigma_m)$ is a regularization of $\chi_{B_{2n}}$. It follows that $Q_n$ is non-negative, symmetric, multilinear and translation invariant.
\end{rmq}

\begin{rmq}\label{continuo2}
Kazarnovski{\v\i} $n$-dimensional mixed pseudovolume is continuous for the Hausdorff metric.

\noindent
If, for every $1\leq\ell\leq n$, $(A_{m,\ell})_{m\in\mathbb N}$ is a sequence of convex bodies belonging to some dense subset of ${\mathcal K}(\mathbb C^n)$ and converging to $A_\ell$ in the Hausdorff metric, then the corresponding sequence $h_{A_{m,\ell}}$ of support functions converges uniformly to $h_{A_\ell}$ on each compact set of $\mathbb C^n$ and so the current $\ddc h_{A_{m,\ell}}$ converges weakly to the current $\ddc h_{A_\ell}$, as $m\to+\infty$. This fact implies the continuity of $Q_n$ for the Hausdorff metric, so that $Q_n(A_1,\ldots,A_n)$ can be computed as the limit
\begin{equation*}
Q_n(A_1,\ldots,A_n)=\lim_{m\to\infty}Q_n(A_{m,1},\ldots,A_{m,n})\,.
\end{equation*}
As a consequence, any property of $Q_n$ on some dense subspace of ${\mathcal K}(\mathbb C^n)$ can be extended to the whole ${\mathcal K}(\mathbb C^n)$. 
\end{rmq}

\begin{rmq} Kazarnovski{\v\i} $n$-dimensional pseudovolume on ${\mathcal K}^2(\mathbb C^n)$.\label{continuo3}

\noindent
On the dense subspace ${\mathcal K}^2_2(\mathbb C^n)\subset {\mathcal K}^2(\mathbb C^n)$ both \defnref{mkaza} and \defnref{prima} may apply and, in fact, by virtue of \cororef{kaza}, the two definitions agree, in particular
\begin{equation}\label{magnifica0}
\int_{B_{2n}}(\ddc h_{(R_\varepsilon A)_\varepsilon})^{\wedge n}
=
\int_{\partial (R_\varepsilon A)_\varepsilon} \alpha_{\partial(R_\varepsilon A)_\varepsilon}\wedge (\de \alpha_{\partial(R_\varepsilon A)_\varepsilon})^{\wedge(n-1)}
\,,
\end{equation}
for every $A\in{\mathcal K}^2(\mathbb C^n)$ and any $\varepsilon>0$. \thmref{continuo} and \rmqref{continuo2} imply that the two definitions agree on the whole ${\mathcal K}^2(\mathbb C^n)$, indeed passing to the limit in \eqref{magnifica0}, as $\varepsilon\to0$, yields
\begin{equation}\label{magnifica1}
\int_{B_{2n}}(\ddc h_{A})^{\wedge n}
=
\int_{\partial A} \alpha_{\partial A}\wedge (\de \alpha_{\partial A})^{\wedge(n-1)}
\,.
\end{equation}
In particular,  \eqref{magnifica1} holds true on ${\mathcal P}_\circ(\mathbb C^n)$, so for every $\Gamma\in{\mathcal P}(\mathbb C^n)$ and any $\varepsilon>0$,  
\begin{equation}\label{magnifica1,5}
\int_{B_{2n}}(\ddc h_{(\Gamma)_\varepsilon})^{\wedge n}
=
\int_{\partial (\Gamma)_\varepsilon} \alpha_{\partial (\Gamma)_\varepsilon}\wedge (\de \alpha_{\partial (\Gamma)_\varepsilon})^{\wedge(n-1)}
\,.
\end{equation}
Then, as $\varepsilon\to0$, \rmqref{continuo2}, \eqref{magnifica1,5} and  \thmref{azzok0}  yield the equality
\begin{equation}\label{magnifica2}
\frac{1}{n!\varkappa_n}\int_{B_{2n}}(\ddc h_{\Gamma})^{\wedge n}
=
\textsl{v}_n^\varrho(\Gamma)\,,
\end{equation}
which reveals a remarkable combinatorial feature of Kazarnovski{\v\i} pseudovolume.
\findim
\end{rmq}

Inspired by a privately communicated idea of Kazarnovski{\v\i}'s, the following \thmref{azzok1} and \cororef{azzok2} generalise \eqref{magnifica2}. 

\begin{thm}\label{azzok1}
Let $n\in\mathbb N^*$, $\Gamma\in {\mathcal P}(\mathbb C^n)$ and let $0\leqslant k\leqslant n$ be an integer. Then
\begin{equation}\label{f-azzok}
\frac{\varkappa_n}{2^{n-k}\varkappa_{2n-k}}{n\choose k}Q_n(\Gamma[k],B_{2n}[n-k])
=
\textsl{v}_k^\varrho(\Gamma)
\,.
\end{equation}

\end{thm}
\noindent\pf
Let $\varepsilon>0$. By virtue of \ref{magnifica1} we can compute $P_n((\Gamma)_\varepsilon)$ in two ways. On one hand
\begin{align}
P_n((\Gamma)_\varepsilon)
&=
\frac{1}{n!\varkappa_n}\int_{B_{2n}}(\ddc h_{(\Gamma)_\varepsilon})^{\wedge n}\nonumber
\\
&=
\frac{1}{n!\varkappa_n}\int_{B_{2n}}\left(\ddc (h_\Gamma+\varepsilon h_{B_{2n}})\right)^{\wedge n}\nonumber
\\
&=
\frac{1}{n!\varkappa_n}\int_{B_{2n}}(\ddc h_\Gamma+\varepsilon\ddc h_{B_{2n}})^{\wedge n}\nonumber
\\
&=
\sum_{k=0}^n {n\choose k}\varepsilon^{n-k}
\frac{1}{n!\varkappa_n}
\int_{B_{2n}}(\ddc h_\Gamma)^k\wedge(\ddc h_{B_{2n}})^{\wedge(n-k)}\nonumber
\\
&=
\sum_{k=0}^n {n\choose k}\varepsilon^{n-k}
Q_n(\Gamma[k],B_{2n}[n-k])\,,\label{k-uno}
\end{align}
on the other hand, (thanks to \thmref{azzok0}),
\begin{align}
P_n((\Gamma)_\varepsilon)
&=
\frac{1}{n!\varkappa_n}
\int_{\partial (\Gamma)_\varepsilon }
\alpha_{\partial (\Gamma)_\varepsilon }\wedge (\de\alpha_{\partial (\Gamma)_\varepsilon} )^{\wedge(n-1)}\nonumber
\\
&=
\sum_{k=0}^n\frac{2^{n-k}\varkappa_{2n-k}}{\varkappa_n}\,\textsl{v}^\varrho_k(\Gamma)\varepsilon^{n-k}\,.\label{k-due}
\end{align}
The theorem follows by comparing coefficients in \eqref{k-uno} and \eqref{k-due}. \findim

\begin{coro}\label{azzok2}
For any integer $0\leq k\leq n$ and $A_1,\ldots,A_k\in {\mathcal K}(\mathbb C^n)$, 
\begin{equation}\label{f-azzokk}
\frac{\varkappa_n}{2^{n-k}\varkappa_{2n-k}}{n\choose k}Q_n(A_1,\ldots,A_k,B_{2n}[n-k])
=
\textsl{V}_k^\varrho(A_1,\ldots,A_k)
\,.
\end{equation}
\end{coro}
\noindent\pf 
If $k=0$, by definition of $\textsl{v}_0^\varrho$ and by \exeref{pseudoball}, both size of (\ref{f-azzokk}) equal $1$. If $k>0$, 
both sides of \eqref{f-azzokk} are continuous in each of the $k$ arguments $A_1,\ldots,A_k$, so it's enough to prove the equality in the case of polytopes. In this case we know that both sides of \eqref{f-azzokk} are symmetric and $k$-linear, so, by \cororef{uni}, they agree as soon as they do it on the diagonal, i.e. for $A_1=\ldots=A_k=A\in{\mathcal P}(\mathbb C^n)$. The claim follows from \thmref{azzok1}.\findim

\begin{coro}
For any $A_1,\ldots,A_n\in {\mathcal K}(\mathbb R^n)$, 
\begin{equation*}
Q_n(A_1,\ldots,A_n)
=
V_n(A_1,\ldots,A_k)\,.
\end{equation*}
\end{coro}
\noindent\pf By continuity, it's enough to prove the equality for polytopes. In this case the statement is an easy consequence of \cororef{azzok2} with $k=n$ and purely real polytopes.\findim

\begin{rmq} $Q_n$ is unitarily invariant and satisfies both the polarization formula both the symmetric one.

\noindent
By the continuity of the mixed $n$-pseudovolume, it follows that $Q_n$ is unitarily invariant and satisfies the usual polarization and symmetric formulas,
\begin{equation}
Q_n(A_1,\ldots,A_n)
=
\frac{1}{n!}\sum_{\varnothing\neq I\subseteq\{1,\ldots,n\}} (-1)^{n-\#I}P_n\left(\sum_{\ell\in I} A_\ell\right)\,,\label{polf}
\end{equation}
\begin{equation}
Q_n(A_1,\ldots,A_n)
=
\frac{1}{n!}\frac{\partial^n}{\partial\lambda_1\ldots\partial\lambda_n} P_n\left(\sum_{\ell=1}^n\lambda_\ell A_\ell\right)\,.\label{sym2}
\end{equation}
Indeed, by \lemref{hinv}, $Q_n$ is unitarily invariant on ${\mathcal K}^\infty_1(\mathbb C^n)$, whereas by \thmref{non-sim-fi},  \lemref{pol-fi} and \cororef{azzok2}, it satisfies \eqref{polf} and \eqref{sym2} on ${\mathcal P}(\mathbb C^n)$.\findim
\end{rmq}

\section{The current {\boldmath $(\ddc h_\Gamma)^{\wedge k}$}}\label{chk}

The present section is devoted to a closer study of the current $(\ddc h_\Gamma)^{\wedge k}$, for any $\Gamma\in{\mathcal P}(\mathbb C^n)$. The case in which the polytope $\Gamma$ is replaced by a general convex body is much harder, but it can be viewed as a weak limit of the polytopal case. Let us first consider the case $k=1$ and observe that $\dc h_\Gamma$ is simply a $1$-form the coefficients of which are locally constant functions, so that the current $\ddc h_\Gamma$ is supported on the corner locus of $h_\Gamma$. Such a corner locus equals the set of points where the coefficients of $\dc h_\Gamma$ are not continuous, i.e. the closure of the $1$-star of ${\Gamma}$, hence
\begin{equation}\label{supp1}
\Supp \ddc h_\Gamma=\overline{\Sigma}_{1,\Gamma}\,.
\end{equation}
Moreover, if $A_1,A_2\in{\mathcal K}(\mathbb C^n)$ are such that $A_1\cup A_2$ is convex, the well-known equality 
\begin{equation*}
A_1+A_2=(A_1\cup A_2)+(A_1\cap A_2)
\end{equation*}
implies that $h_{A_1}+h_{A_2}=h_{A_1\cup A_2}+h_{A_1\cap A_2}$, whence 
\begin{equation}\label{val-1}
\ddc h_{A_1}+\ddc h_{A_2}=\ddc h_{A_1\cup A_2}+\ddc h_{A_1\cap A_2}\,.
\end{equation}
Setting $h_\varnothing\equiv 0$, the relation \eqref{val-1} shows that taking the $\ddc$ of support functions yields a $1$-homogeneous valuation on ${\mathcal K}(\mathbb C^n)$ with values in the space of $(1,1)$-currents. In this section we will prove the following generalizations, for $k>1$, of \eqref{supp1} and \eqref{val-1}:
\begin{equation}\label{suppk}
\Supp (\ddc h_\Gamma)^{\wedge k}=\overline{\Sigma}_{k,\Gamma}\,,
\end{equation}
\begin{equation}\label{vale-k}
(\ddc h_{A_1})^{\wedge k}+(\ddc h_{A_2})^{\wedge k}=(\ddc h_{A_1\cup A_2})^{\wedge k}+(\ddc h_{A_1\cap A_2})^{\wedge k}\,.
\end{equation}

If $\Gamma\in{\mathcal P}(\mathbb C^n)$ is oriented, for every $1\leq k\leq n$ and any $\Delta\in{\mathcal B}_{\rm ed}(\Gamma,k)$, let $\lambda_\Delta$ be the $(k,k)$-current of measure type on $\mathbb C^n$ acting on any continuous and compactly supported $(n-k,n-k)$-form $\varsigma$ as
\begin{equation}\label{lambdak}
\langle\!\langle \lambda_\Delta,\varsigma\rangle\!\rangle   
=
\int_{K_\Delta}\iota_\Delta^*(\upsilon_{E_\Delta^\prime}\wedge\varsigma)\,,
\end{equation}
where $\iota_\Delta:K_\Delta \rightarrow \mathbb C^n$ is the inclusion, $K_\Delta$ is given the coherent orientation as an open submanifold of $E^\perp_\Delta$ and $\upsilon_{E^\prime_\Delta}$ is the corresponding volume form on $E^\prime$.
Let us point out some remarks.
\begin{itemize}
\item The definition of $\lambda_\Delta$ makes sense if and only if $\Delta$ is equidimensional. Indeed, the degree of the form $\iota_\Delta^*(\upsilon_{E_\Delta^\prime}\wedge\varsigma)$ equals the dimension of $K_\Delta$ if and only if $\Delta$ is equidimensional. 
\item For any equidimensional face $\Delta$ of $\Gamma$, one has $E_\Delta\oplus E_\Delta^\prime\oplus E_\Delta^{\perp_\mathbb C}=\mathbb C^n$ and $K_\Delta\subset E_\Delta^\perp= E_\Delta^\prime\oplus E_\Delta^{\perp_\mathbb C}$, so the orientation of $K_\Delta$ comes from that of $E_\Delta^\prime\oplus E_\Delta^{\perp_\mathbb C}$. Since the mutually orthogonal $\mathbb C$-linear subspaces $E_\Delta\oplus E_\Delta^\prime$ and $E_\Delta^{\perp_\mathbb C}$ are naturally positively oriented, it follows that the orientation of $E_\Delta^\perp$ ultimately depends on the orientation of $E_\Delta^\prime$. As shown in \rmqref{orient}, choosing an orientation on $E_\Delta$ (resp. $E_\Delta^\prime$) implies a corresponding choice of orientation on $E_\Delta^\prime$ (resp. $E_\Delta$) yielding the natural positive orientation on $E_\Delta\oplus E_\Delta^\prime$. The volume form $\upsilon_{E_\Delta^\prime}$ involved in the definition of the current $\lambda_\Delta$ conforms to the chosen orientation on $E_\Delta^\prime$. This way the integral in \eqref{lambdak} is not affected by the choice of such an orientation and the current $\lambda_\Delta$ is thus well defined. According to de Rham~\cite{dR}, the form $\upsilon_{E^\prime_\Delta}$ is an example of \emph{odd form}, i.e. a form changing its sign as $E_\Delta^\prime$ (an thus $E_\Delta$) changes its orientation. A convenient way to make a coherent choice is to fix an orthonormal basis $v_1,\ldots,v_k$ of $E_\Delta$ and set 
\begin{equation}\label{convenient}
\upsilon_{E_\Delta^\prime}
=
(-1)^{k(k-1)/2}\re\langle w_1,-\rangle\wedge\ldots\wedge \re\langle w_k,-\rangle
\,,
\end{equation}
where $w_1,\ldots,w_k$ is the basis of $E_\Delta^\prime$ defined in~\thmref{figo}. Accordingly, the space $E_\Delta^\perp$ (and hence $K_\Delta$) will get the orientation associated to the basis $$(-1)^{k(k-1)/2}w_1,\ldots,w_k,b_{2k+1},\ldots,b_{2n}\,,$$ where $b_{2k+1},\ldots,b_{2n}$ is a positive basis of $E_\Delta^{\perp_\mathbb C}$.

If $x_1,\ldots,x_k,y_1,\ldots,y_k,x_{k+1},y_{k+1},\ldots,x_n,y_n$ are the coordinates corresponding to the basis $v_1,\ldots,v_k,(-1)^{k(k-1)/2}w_1,\ldots,w_k,b_{2k+1},\ldots,b_{2n}$ of $\mathbb C^n$ and $[{K}_\Delta]$ denotes the integration current on the open oriented submanifold  ${K}_\Delta\subset\mathbb C^n$, the 
expressions of $[{K}_\Delta]$, $\upsilon_{E^\prime_\Delta}$ and $\lambda_\Delta$ in the these coordinates are as follows
\begin{align}
[{K}_\Delta]&=({\mathcal H}^{2n-k}\,\llcorner\, K_\Delta)\,\de x_1\wedge\ldots\wedge\de x_k\,,
\\
\upsilon_{E^\prime_\Delta}&=(-1)^{k(k-1)/2}\,\de y_1\wedge\ldots\wedge\de y_k\,,
\\
\lambda_\Delta&=[{K}_\Delta]\wedge\upsilon_{E^\prime_\Delta}
\\
&=
({\mathcal H}^{2n-k}\,\llcorner\, K_\Delta)\,\de x_1\wedge\de y_1\ldots\wedge\de x_k\wedge\de y_k\,,
\end{align}
where ${\mathcal H}^{2n-k}$ is the $(2n-k)$-dimensional Hausdorff measure on $\mathbb C^n$ and ${\mathcal H}^{2n-k}\,\llcorner\, K_\Delta$ denotes its restriction to $K_\Delta$.

\item The only terms of $\varsigma$ giving a non zero contribution to $\iota_\Delta^*(\upsilon_{E_\Delta^\prime}\wedge\varsigma)=\iota_\Delta^*\upsilon_{E_\Delta^\prime}\wedge\iota_\Delta^*\varsigma$ are those involving the differentials of the coordinate functions of $E_\Delta^{\perp_\mathbb C}$, nevertheless the coefficients of such terms do depend on the full array of coordinates of $K_\Delta\subset E_\Delta^\perp$. 
\item The support of the current $\lambda_\Delta$ is the closure of $K_\Delta$ but $\lambda_\Delta$ is  concentrated just on $K_\Delta$.
\end{itemize}

\subsection[The case $k=1$]{The case {\boldmath $k=1$}}

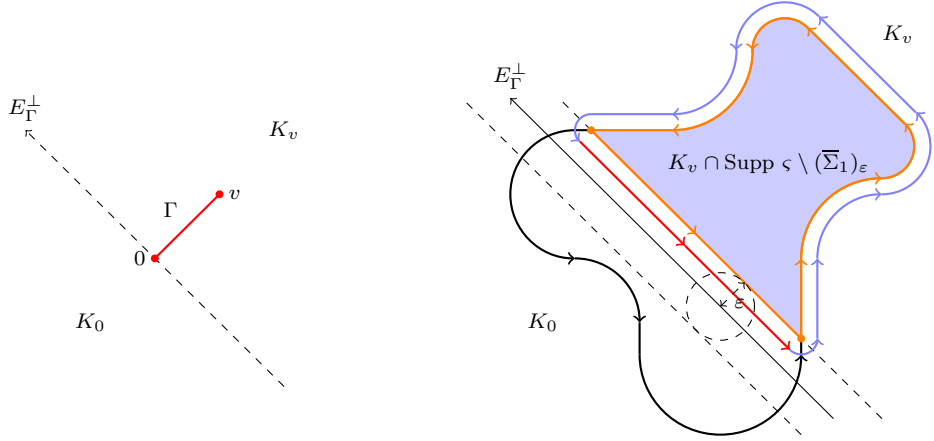
\begin{figure}[ht]

\begin{center}
\begin{tikzpicture}
\draw (3.5,3.5) node{\footnotesize$K_v$};
\draw (-2,-1) node{\footnotesize$K_0$};
\draw (-9,-1) node{\footnotesize$K_0$};
\draw (-6,2) node{\footnotesize$K_v$};

\draw[dashed,->] (-6,-2) -- (-10,2) node[above]{\footnotesize$E_\Gamma^\perp$}; 
\draw[red,thick] (-8,0) -- (-7,1);
\draw (-7.5,0.5) node[above left]{\footnotesize$\Gamma$};
\filldraw[red] (-7,1) circle(1.5pt) node[black, right]{\footnotesize$v$};
\filldraw[red] (-8,0) circle(1.5pt) node[black,left]{\footnotesize$0$};
\draw[->] (2.5,-2.5) -- (-2.5,2.5) node[above]{\footnotesize$E_\Gamma^\perp$}; 
\draw[dashed] (3.25,-2.5)--(2,-1.25);
\draw[dashed] (-1.25,2)--(-1.75,2.5);
\draw[dashed] (2,-2.75)--(-3.25,2.5);
\draw[thick] (-1.25,2) -- (-1.5,2);
\draw[->, thick] (-1.5,2) arc (90:270:1);
\draw[->, thick] (-1.5,0) arc (90:0:1);
\draw[thick] (-0.5,-1)--(-0.5,-1.5);
\draw[->, thick] (-0.5,-1.5) arc (180:360:1.25);
\draw[thick] (2,-1.5) -- (2,-1.25);
\draw[thick,orange] (0,2) -- (-1.25,2);
\filldraw[orange] (-1.25,2) circle(1.5pt);
\draw[thick,blue!50!white] (0.05,2.25) -- (-1.25,2.25);
\draw[->,thick,blue!50!white] (-1.25,2.25) arc(90:225:0.25);

\draw[orange,thick] (0.375,0.375)--(2,-1.25);
\filldraw[orange] (2,-1.25) circle(1.5pt);

\filldraw[thick,draw=orange,fill=blue!20!white] (2,0) arc (180:90:1.25)--(3.25,1.25) arc (270:405:0.5)-- (2.1,3.6) arc (45:180:0.5)-- (1.25,3.25) arc (360:270:1.25)--(-1.25,2)--(2,-1.25)--cycle;

\draw[orange,thick,->] (2,0) arc (180:90:1.25); 
\draw[blue!50!white,thick,->] (2.25,0) arc (180:90:1); 

\draw[->, thick,orange] (3.25,1.25) arc (270:405:0.5);
\draw[->, thick,blue!50!white] (3.25,1) arc (270:405:0.75);

\draw[thick,orange,->] (3.6,2.1) --(2.1,3.6);
\draw[thick,blue!50!white,->] (3.775,2.275) --(2.275,3.775);

\draw[->, thick,orange] (2.1,3.6) arc (45:180:0.5);
\draw[->, thick,blue!50!white] (2.275,3.775) arc (45:180:0.75);

\draw[->, thick,orange] (1.25,3.25) arc (360:270:1.25);
\draw[->, thick,blue!50!white] (1,3.25) arc (360:270:1);

\draw[->,thick,orange] (2,-1.25) -- (2,0);
\draw[->,thick,blue!50!white] (2.25,-1.25) -- (2.25,0);

\draw[->,orange,thick] (-1.25,2)--(0.375,0.375);
\draw[->,red,thick] (-1.425,1.825)--(0.2,0.2);
\draw[->,red,thick] (0.2,0.2)--(1.825,-1.425);

\draw[->,thick,blue!50!white] (1.825,-1.425) arc(225:360:0.25);

\draw (1.5,1.5) node{\footnotesize$K_v\cap \Supp\varsigma\setminus(\overline{\Sigma}_1)_\varepsilon$};
%
\draw[dashed] (0.75,-0.75) circle(15pt);
\draw[dashed,<->] (0.75,-0.75)--(1.125,-0.375);
\draw (1.05,-0.7) node{\footnotesize$\varepsilon$};
\end{tikzpicture}
\caption{On the left the 1-polytope $\Gamma$. The half-spaces bounded by the hyperplane $E_\Delta^\perp$ are $K_0$ and $K_v$. On the right, $\Supp\varsigma$ is depicted as the compact set with smooth rounded boundary, $(\overline{\Sigma}_1)_\varepsilon$ is the strip with dashed boundary. The blue part is $K_v\cap \Supp\varsigma\setminus(\overline{\Sigma}_1)_\varepsilon$, the components of its boundary are drawn in orange. The blue/red smooth path around $K_v\cap \Supp\varsigma\setminus(\overline{\Sigma}_1)_\varepsilon$ is the boundary of $S_{v,\varepsilon}$. The red segment is the only piece of $\partial S_{v,\varepsilon}$ that counts.}
\label{1-Gamma}
\end{center}
\end{figure}

\begin{lem}\label{lem1}
Let $\Gamma\in{\mathcal P}(\mathbb C)$, then $\ddc h_\Gamma=\sum_{\Lambda\in{\mathcal B}(\Gamma,1)}\vol_1(\Lambda)\lambda_\Lambda$. In particular $\Supp \ddc h_\Gamma=\overline{\Sigma}_1$.
\end{lem}
\noindent\pf 
If $\Gamma$ is a point, the lemma is trivial. If $\Gamma$ is not a point, $h_\Gamma$ is piece-wise linear, the current $\ddc h_\Gamma$ is zero on $\mathbb C\setminus\overline{\Sigma}_{1}$ and it acts on test functions, so let $\varsigma$ be a smooth function with compact support.
If $d=1$, $\Gamma$ is a line segment and, up to a translation (which do not affect $\dc h_\Gamma$), we may suppose that $\Gamma=[0,v]$ with $v\in \mathbb C^*$. Then $E_\Gamma$ is spanned by $v$ and, according to our convention, $iv$ is a basis of $E_\Gamma^\perp$ such that $v,iv$ gives $E_\Gamma\oplus E_\Gamma^\perp=\mathbb C$ the natural orientation. The closure of the $1$-star of $\Gamma$ is the closure of the dual cone to $\Gamma$ and, in this case, it is nothing but $E_\Gamma^\prime=E_\Gamma^\perp=iE_\Gamma=\{ie^{i\Arg v}\rho\in \mathbb C\mid \rho\in\mathbb R\}$.
Then 
\begin{align*}
\langle\!\langle \ddc h_\Gamma,\varsigma\rangle\!\rangle
&=
\langle\!\langle \dc h_\Gamma,\de\varsigma\rangle\!\rangle
\\
&=
\int_{\Supp\varsigma}\dc h_\Gamma\wedge\de \varsigma
\\
&=
\lim_{\varepsilon\to0}
\int_{\Supp\varsigma\setminus(\overline{\Sigma}_{1})_\varepsilon}\dc h_\Gamma\wedge\de \varsigma
\\
&=
\lim_{\varepsilon\to0}
\int_{K_{0}\cap\Supp\varsigma\setminus(\overline{\Sigma}_{1})_\varepsilon}\dc h_\Gamma\wedge\de \varsigma
\\
&+
\lim_{\varepsilon\to0}
\int_{K_{v}\cap\Supp\varsigma\setminus(\overline{\Sigma}_{1})_\varepsilon}\dc h_\Gamma\wedge\de \varsigma
\,.
\end{align*}
Since $h_\Gamma$ vanishes on $K_0$,
\begin{align*}
\langle\!\langle \ddc h_\Gamma,\varsigma\rangle\!\rangle
&=
\lim_{\varepsilon\to0}
\int_{K_{v}\cap\Supp\varsigma\setminus(\overline{\Sigma}_{1})_\varepsilon}
(i/2)(v\de\bar z-\bar v\de z)\wedge\de \varsigma
\\
&=
\lim_{\varepsilon\to0}
\int_{K_{v}\cap\Supp\varsigma\setminus(\overline{\Sigma}_{1})_\varepsilon}
(i/2)\de[(\bar v\de z-v\de\bar z)\varsigma]
\\
&=
\lim_{\varepsilon\to0}
\int_{S_{v,\varepsilon}}
(i/2)\de[(\bar v\de z-v\de\bar z)\varsigma]
\\
&=
\lim_{\varepsilon\to0}
\int_{\partial S_{v,\varepsilon}}
(i/2)(\bar v\de z-v\de\bar z)\varsigma\,.
\end{align*}
where $S_{v,\varepsilon}$ is the $(\varepsilon/2)$-neighbourhood of $K_{v}\cap\Supp\varsigma\setminus(\overline{\Sigma}_{1})_\varepsilon$. Of course $K_{v}\cap\Supp\varsigma\setminus(\overline{\Sigma}_{1})_\varepsilon\subsetneq S_{v,\varepsilon}\subsetneq K_v$ and both $K_{v}\cap\Supp\varsigma\setminus(\overline{\Sigma}_{1})_\varepsilon$ and $S_{v,\varepsilon}$ approach the set $K_v\cap\Supp\varsigma$, as $\varepsilon\to 0^+$. However, unlike $K_{v}\cap\Supp\varsigma\setminus(\overline{\Sigma}_{1})_\varepsilon$, the set $S_{v,\varepsilon}$ has smooth boundary. This boundary can be decomposed into three parts: the pieces which belong to $\partial \Supp\varsigma$ (on which $\varsigma$ is zero), the spherical pieces over the singular points of $\partial[K_{v}\cap\Supp\varsigma\setminus(\overline{\Sigma}_{1})_\varepsilon]$, (which shrink to the corresponding singular points as $\varepsilon\to0$) and the flat pieces parallel and near to $K_{\Gamma}$ (in Figure \ref{1-Gamma} the set $\partial S_{v,\varepsilon}$ admits just one of such parts). Each of the latter components is a relatively open line segment inheriting from $\partial S_{v,\varepsilon}$ an orientation which is opposite to that of $E_{\Gamma}^\perp$. As $\varepsilon\to 0$, each of these segments  approaches the subset
$\overline{K}_{\Gamma}\cap\Supp\varsigma=K_\Gamma\cap\Supp\varsigma$, on which $\iota^*_{\Gamma}\de z=ie^{i\Arg v}\de \rho$, $\iota^*_{\Gamma}\de\bar z=-ie^{-i\Arg v}\de\rho$ and $\iota^*_{\Gamma}v_{E_\Gamma^\prime}=\de\rho$, thus 
\begin{align*}
\langle\!\langle \ddc h_\Gamma,\varsigma\rangle\!\rangle
&=
\int_{\overline{K}_{\Gamma}\cap\Supp\varsigma}
-(i/2)(\bar vie^{i\Arg v}+v ie^{-i\Arg v})\varsigma\de\rho
\\
&=
\int_{\overline{K}_{\Gamma}\cap\Supp\varsigma}
\vert v\vert \varsigma\de\rho
\\
&=
\int_{K_{\Gamma}\cap\Supp\varsigma}
\vert v\vert \varsigma\de\rho
\\
&=
\vol_1(\Gamma) \langle\!\langle \lambda_\Gamma,\varsigma\rangle\!\rangle\,,
\end{align*}
i.e. $\ddc h_\Gamma=\vol_1(\Gamma)\lambda_\Gamma$. 

Observe that giving $E_\Gamma$ the opposite orientation means choosing $-v$ as a basis of $E_\Gamma$ and, consequently, $-v,-iv$ as a basis of $E_\Gamma\oplus E_\Gamma^\perp=\mathbb C$. In both cases the two basis of the latter space are positively oriented and the proof proceeds in the same way.

\begin{figure}[ht]
\begin{center}
\begin{tikzpicture}[scale=1]
\filldraw[blue!20!white] (-9,0)--(-7,0)--(-10,3)--cycle;
\filldraw[orange] (-9,0) circle(1.5pt);
\draw (-9,-0.1) node[black,below left]{$v_0$};
\filldraw[orange] (-7,0) circle(1.5pt) node[black,below right]{$v_1$};
\filldraw[orange] (-10,3) circle(1.5pt) node[black,above left]{$v_2$};
\draw[thick,orange,line cap=round, line join=round] (-9,0)--(-7,0)--(-10,3)--cycle;
\draw[thick,->,red,line cap=round, line join=round] (-9,0) -- (-8,0);
\draw[thick,red,line cap=round, line join=round] (-8,0) -- (-7,0);

\draw[thick,->,green,line cap=round, line join=round] (-7,0) -- (-8.5,1.5);
\draw[thick,green,line cap=round, line join=round] (-7,0) -- (-10,3);

\draw[thick,->,brown,line cap=round, line join=round] (-10,3) -- (-9.5,1.5);
\draw[thick,brown,line cap=round, line join=round] (-9.5,1.5)-- (-9,0);

\draw (-8.75,1) node{$\Gamma$};
\draw[dashed,line cap=round, line join=round] (-9,0) -- (-9,-2);
\draw[dashed,line cap=round, line join=round] (-9,0) -- (-11,-0.666);
\draw[dashed,line cap=round, line join=round] (-10,3) -- (-12,2.333);
\draw[dashed,line cap=round, line join=round] (-10,3) -- (-8.5,4.5);
\draw[dashed,line cap=round, line join=round] (-7,0) -- (-5.5,1.5);
\draw[dashed,line cap=round, line join=round] (-7,0) -- (-7,-2);
\draw (-10,-1) node{$v_0+K_{v_0}$};
\draw (-6,-1) node{$v_1+K_{v_1}$};
\draw (-11,4) node{$v_2+K_{v_2}$};
\draw (-8,0) node[below]{$\Lambda_0$};
\draw (-8.5,1.5) node[above right]{$\Lambda_1$};
\draw (-9.5,1.5) node[left]{$\Lambda_2$};
\end{tikzpicture}
\caption{The 2-polytope $\Gamma$ is a triangle.}
\label{2-Gamma}
\end{center}
\end{figure}
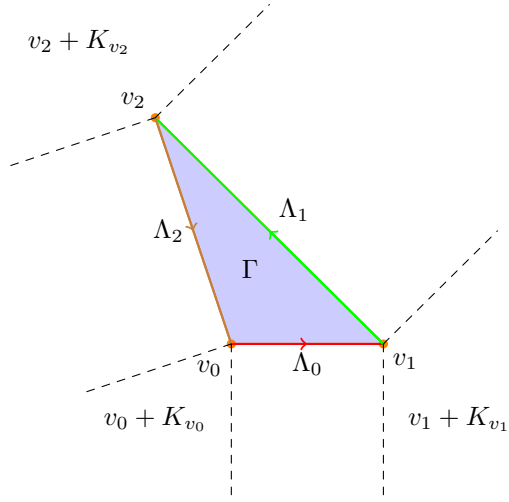

Let us now suppose $d=\dim\Gamma=2$, (cf. Figure~\ref{2-Gamma}). If $\Gamma$ is oriented and  $v_0,v_1,\ldots,v_s$ is a numbering of the vertices of $\Gamma$ which agrees with the orientation of $\partial\Gamma$ induced by unit outer normal vectors, the sides of $\Gamma$ are the oriented line segments $\Lambda_\ell=[v_\ell,v_{\ell+1}]$, $-1\leq\ell\leq s$, with $v_{s+1}=v_0$ and $v_{-1}=v_s$.

\begin{figure}[ht]
\begin{center}
\begin{tikzpicture}[scale=1.5]
\draw[thick,->,line cap=round, line join=round] (0,0) -- (0,-4) node[below]{$K_{\Lambda_0}$};
\draw[thick,->,line cap=round, line join=round] (0,0) -- (3.5,3.5)  node[above right]{$K_{\Lambda_1}$};
\draw[thick,->,line cap=round, line join=round] (0,0) -- (-4.5,-1.5257)  node[below left]{$K_{\Lambda_2}$};

\draw (-2,-2) node{$K_{v_0}\setminus(\overline\Sigma_1)_\varepsilon$};
\draw (3.5,0) node{$K_{v_1}\setminus(\overline\Sigma_1)_\varepsilon$};
\draw (-2,3) node{$K_{v_2}\setminus(\overline\Sigma_1)_\varepsilon$};
\draw[dashed] (0,0) circle(0.5cm);

\draw[dashed,line cap=round, line join=round] (-0.5,0) -- (-0.5,-4);
\draw[dashed,line cap=round, line join=round] (0.5,0)--(0.5,-4);
\draw[dashed,line cap=round, line join=round] (2.25,3)--(3,3.707);
\draw[dashed,line cap=round, line join=round] (3,2.25)--(4,3.292);
\draw[dashed,line cap=round, line join=round] (-1.5,0.0257)--(-4.52,-1.032);
\draw[dashed,line cap=round, line join=round] (0.1581,-0.4743) -- (-4.26,-1.95);
\filldraw[thick,draw=orange,fill=blue!20!white,line cap=round, line join=round] (0,2) -- (-1.5,2.0257) arc (90:270:1) -- (-0.2701,0.447)--(2.25,3)--(2.25,3.25) arc (0:180:0.5) --(1.25,3.25) arc (360:270:1.25)--cycle;
\filldraw[thick,draw=orange,fill=blue!20!white,line cap=round, line join=round] (0.5,-0.21)--(0.5,-2.7)--(0.5,-2.725) arc (258.5:360:1.25) -- (2,0) arc (180:90:1.25)--(3.25,1.25) arc (270:450:0.5)--(3,2.25)--cycle;
\draw[->, thick,orange,line cap=round, line join=round] (0,2) -- (-1.5,2.0257) ; 
\draw[->, thick,blue!50!white,line cap=round, line join=round] (0,2.25) -- (-1.5,2.25); 

\draw[->, thick,orange,line cap=round, line join=round] (-1.5,2.0257)  arc (90:270:1); 
\draw[->, thick,blue!50!white,line cap=round, line join=round] (-1.5,2.25) arc (90:270:1.25); 

\draw[thick,line cap=round, line join=round] (-1.5,0.0257)-- (-1.25,0.0257); 

\draw[thick,blue!50!white,line cap=round, line join=round] (-1.5,-0.25) arc(270:288:0.25); 

\draw[thick,brown,->,line cap=round, line join=round] (-1.4226,-0.2377)-- (-0.19,0.2); 

\draw[thick,blue!50!white,line cap=round, line join=round]  (-0.19,0.2) arc(288:315:0.25); 
\draw[thick,green,->,line cap=round, line join=round] (-0.09,0.26)-- (1.1225,1.536); 
\draw[thick,green,->,line cap=round, line join=round] (1.1225,1.536)-- (2.425,2.825); 

\draw[->, thick,line cap=round, line join=round] (-1.25,0.0257) arc (90:0:0.75);
\draw[thick,->,line cap=round, line join=round] (-0.5,-0.7243)--(-0.5,-1.5);
\draw[thick,->,line cap=round, line join=round] (-0.5,-1.5) arc (180:258:1.25);
\draw[thick,orange,->,line cap=round, line join=round] (0.5,-2.725) arc (258.5:360:1.25);
\draw[thick,blue!50!white,line cap=round, line join=round] (2.25,-1.5) arc (0:-101.5:1.5);

\draw[->,thick,orange,line cap=round, line join=round] (2,-1.5) -- (2,0);
\draw[->,thick,blue!50!white] (2.25,-1.5) -- (2.25,0);

\draw[->, thick,orange,line cap=round, line join=round] (2,0) arc (180:90:1.25);
\draw[->, thick,blue!50!white,line cap=round, line join=round] (2.25,0) arc (180:90:1);

\draw[->, thick,blue!50!white,line cap=round, line join=round] (1,3.25) arc (360:270:1);
\draw[->,thick,blue!50!white,line cap=round, line join=round] (2.5,3.25) arc (0:180:0.75);

\draw[->, thick,orange,line cap=round, line join=round] (3.25,2.25)--(3,2.25);
\draw[->, thick,blue!50!white,line cap=round, line join=round] (3.25,2.5)--(3,2.5);
\draw[thick,blue!50!white,line cap=round, line join=round] (3,2.5) arc(90:135:0.25);

\draw[->,thick,line cap=round, line join=round] (3,2.25) arc (270:180:0.75); 
\draw[thick,orange,line cap=round, line join=round] (2.25,3)--(2.25,3.25);
\draw[thick,blue!50!white,line cap=round, line join=round] (2.5,3)--(2.5,3.25);

\draw[thick,blue!50!white,line cap=round, line join=round] (2.5,3) arc(0:-45:0.25);

\draw[thick,blue!50!white,line cap=round, line join=round] (3.25,1) arc (270:450:0.75);
\draw[thick,orange,->,line cap=round, line join=round] (-1.5,0.0257)--(-0.2701,0.447);
\draw[thick,orange,->,line cap=round, line join=round] (-0.2701,0.447)--(0.99,1.72);
\draw[thick,orange,line cap=round, line join=round] (0.99,1.72)--(2.25,3);
\draw[thick,orange,->,line cap=round, line join=round] (3,2.25)--(1.75,1.02);
\draw[thick,green,->,line cap=round, line join=round] (2.825,2.425)--(1.575,1.195);

\draw[thick,orange,->,line cap=round, line join=round] (1.75,1.02)--(0.5,-0.21);
\draw[thick,green,->,line cap=round, line join=round] (1.575,1.195)--(0.325,-0.035);

\draw[thick,orange,->,line cap=round, line join=round] (0.5,-0.21)--(0.5,-1.455);
\draw[thick,blue!50!white,line cap=round, line join=round] (0.25,-0.21) arc(180:135:0.25);
\draw[thick,red,->,line cap=round, line join=round] (0.25,-0.21)--(0.25,-1.455);

\draw[thick,orange,line cap=round, line join=round] (0.5,-1.455)--(0.5,-2.7);
\draw[thick,red,->,line cap=round, line join=round] (0.25,-1.455)--(0.25,-2.725);
\draw[thick,blue!50!white,line cap=round, line join=round] (0.25,-2.725) arc (180:258.65:0.25);
\draw[dashed,<->,line cap=round, line join=round] (0,-3.25)--(0.5,-3.25);
\draw (0.25,-3.25) node[below]{\footnotesize$\varepsilon$};
\draw[dotted] (0.25,-2.725)--(0.25,-4);
\draw[dashed,<->,line cap=round, line join=round] (0.25,-3.6)--(0.5,-3.6);
\draw (0.375,-3.6) node[below]{\footnotesize$\frac{\varepsilon}{2}$};
\end{tikzpicture}
\caption{The set $\Supp\varsigma$ is depicted as the compact set with smooth rounded boundary, $(\overline{\Sigma}_1)_\varepsilon$ is the set with black, dashed boundary. The blue $2$-dimensional part of $\Supp\varsigma$ is the union of $K_{v_1}\cap \Supp\varsigma\setminus(\overline{\Sigma}_1)_\varepsilon$ and $K_{v_2}\cap \Supp\varsigma\setminus(\overline{\Sigma}_1)_\varepsilon$, the components of $\partial [K_{v_1}\cap \Supp\varsigma\setminus(\overline{\Sigma}_1)_\varepsilon]$ and $\partial [K_{v_2}\cap \Supp\varsigma\setminus(\overline{\Sigma}_1)_\varepsilon]$ are coloured according to  the colours of the edges of $\Gamma$. The coloured smooth path all around $K_{v_1}\cap \Supp\varsigma\setminus(\overline{\Sigma}_1)_\varepsilon$ is the boundary of $S_{v_1,\varepsilon}$. The red, brown and green segments are the only pieces of $\partial S_{v_1,\varepsilon}$ which count. The same occurs to $S_{v_2,\varepsilon}$.}
\label{3-Gamma}
\end{center}
\end{figure}
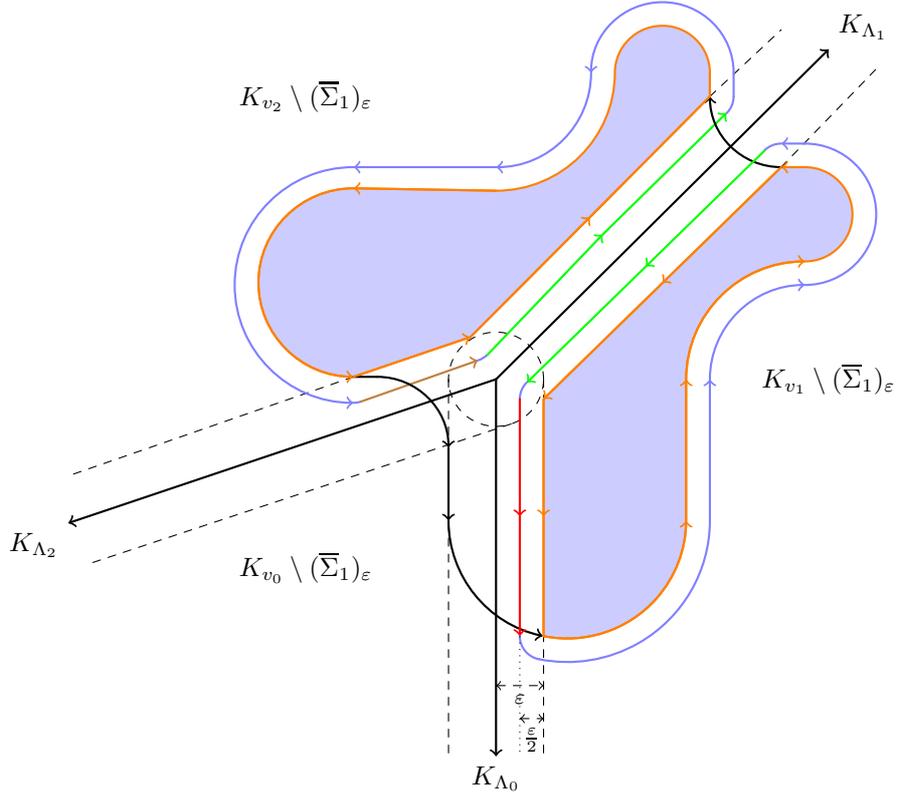

The dual cone to the side $\Lambda_\ell$ is spanned by the unit outer normal vector $u_{\Lambda_\ell,\Gamma}=-ie^{i\Arg(v_{\ell+1}-v_\ell)}$ hence $K_{\Lambda_\ell}=\{-ie^{i\Arg(v_{\ell+1}-v_\ell)}\rho\in\mathbb C\mid \rho>0\}$ so that $\iota^*_{\Lambda_\ell}\de z=-ie^{i\Arg(v_{\ell+1}-v_\ell)}\de\rho$, $\iota^*_{\Lambda_\ell}\de\bar z=ie^{-i\Arg(v_{\ell+1}-v_\ell)}\de \rho$ and $v_{E_{\Lambda_\ell}^\prime}=\de\rho$. 
\begin{align*}
\langle\!\langle \ddc h_\Gamma,\varsigma\rangle\!\rangle
&=
\langle\!\langle \dc h_\Gamma,\de\varsigma\rangle\!\rangle
\\
&=
\int_{\Supp\varsigma}\dc h_\Gamma\wedge\de \varsigma
\\
&=
\lim_{\varepsilon\to0}
\int_{\Supp\varsigma\setminus(\overline{\Sigma}_{1})_\varepsilon}\dc h_\Gamma\wedge\de \varsigma
\\
&=
\lim_{\varepsilon\to0}
\sum_{\ell=0}^s
\int_{K_{v_\ell}\cap\Supp\varsigma\setminus(\overline{\Sigma}_{1})_\varepsilon}\dc h_\Gamma\wedge\de \varsigma\,.
\end{align*}
The domains of linearity of $h_\Gamma$ are the dual cones to the vertices of $\Gamma$, so
\begin{align}
\langle\!\langle \ddc h_\Gamma,\varsigma\rangle\!\rangle
&=
\sum_{\ell=0}^s
\lim_{\varepsilon\to0}
\int_{K_{v_\ell}\cap\Supp\varsigma\setminus(\overline{\Sigma}_{1})_\varepsilon}
(i/2)(v_\ell\de\bar z-\bar v_\ell\de z)\wedge\de \varsigma\nonumber
\\
&=
-\sum_{\ell=0}^s
\lim_{\varepsilon\to0}
\int_{K_{v_\ell}\cap\Supp\varsigma\setminus(\overline{\Sigma}_{1})_\varepsilon}
(i/2)\de[(v_\ell\de\bar z-\bar v_\ell\de z)\varsigma]\nonumber
\\
&=
-\sum_{\ell=0}^s
\lim_{\varepsilon\to0}
\int_{S_{v_\ell,\varepsilon}}
(i/2)\de[(v_\ell\de\bar z-\bar v_\ell\de z)\varsigma]\nonumber
\\
&=
\sum_{\ell=0}^s
\lim_{\varepsilon\to0}
\int_{\partial S_{v_\ell,\varepsilon}}
(i/2)(\bar v_\ell\de z-v_\ell\de\bar z)\varsigma\,,\label{sum}
\end{align}
where $S_{v_\ell,\varepsilon}$ is the $(\varepsilon/2)$-neighbourhood of $K_{v_\ell}\cap\Supp\varsigma\setminus(\overline{\Sigma}_{1})_\varepsilon$. 
As before, the only parts of the $\partial S_{v_\ell,\varepsilon}$ which really count are those which are parallel and near to $K_{\Lambda_\ell}$ and $K_{\Lambda_{\ell-1}}$ respectively, with the first one inheriting from $\partial S_{v_\ell,\varepsilon}$ an orientation which is opposite to that of $K_{\Lambda_\ell}$ and the second one inheriting from $\partial S_{v_\ell,\varepsilon}$ an orientation which is equal to that of $K_{\Lambda_{\ell-1}}$, (cf. Figure \ref{3-Gamma}). Observe that giving $\Gamma$ the opposite orientation simply switches the orientations of the latter boundary pieces.

As $\varepsilon\to 0$, such boundary pieces are line segments approaching respectively 
$\overline{K}_{\Lambda_\ell}\cap\Supp\varsigma$ and $\overline{K}_{\Lambda_{\ell-1}}\cap\Supp\varsigma$, thus \eqref{sum} becomes
\begin{align*}
&
\sum_{\ell=0}^s
\int_{\overline{K}_{\Lambda_\ell}\cap\Supp\varsigma}
(i/2)(\bar v_\ell ie^{i\Arg(v_{\ell+1}-v_\ell)}+v_\ell ie^{-i\Arg(v_{\ell+1}-v_\ell)})\varsigma\de\rho
\\
&+
\sum_{\ell=0}^s
\int_{\overline{K}_{\Lambda_{\ell-1}}\cap\Supp\varsigma}
(i/2)(-\bar v_\ell ie^{i\Arg(v_{\ell}-v_{\ell-1})}- v_\ell ie^{-i\Arg(v_\ell-v_{\ell-1})})\varsigma\de\rho
\\
&=
(1/2)
\sum_{\ell=0}^s
\int_{K_{\Lambda_\ell}\cap\Supp\varsigma}
(-\bar v_\ell e^{i\Arg(v_{\ell+1}-v_\ell)}- v_\ell e^{-i\Arg(v_{\ell+1}-v_\ell)})\varsigma\de\rho
\\
&+
(1/2)
\sum_{\ell=0}^s
\int_{K_{\Lambda_{\ell-1}}\cap\Supp\varsigma}
(\bar v_\ell e^{i\Arg(v_{\ell}-v_{\ell-1})}+ v_\ell e^{-i\Arg(v_\ell-v_{\ell-1})})\varsigma\de\rho
\\
&=
(1/2)
\sum_{\ell=0}^s
\int_{K_{\Lambda_\ell}\cap\Supp\varsigma}
(-\bar v_\ell e^{i\Arg(v_{\ell+1}-v_\ell)}- v_\ell e^{-i\Arg(v_{\ell+1}-v_\ell)})\varsigma\de\rho
\\
&+
(1/2)
\sum_{j=-1}^{s-1}
\int_{K_{\Lambda_{j}}\cap\Supp\varsigma}
(\bar v_{j+1} e^{i\Arg(v_{j+1}-v_{j})}+ v_{j+1} e^{-i\Arg(v_{j+1}-v_{j})})\varsigma\de\rho\,.
\end{align*}
Since $\Lambda_{-1}=\Lambda_s$, we simply get
\begin{equation*}
(1/2)
\sum_{\ell=0}^s
\int_{K_{\Lambda_\ell}\cap\Supp\varsigma}
[(\bar v_{\ell+1} -\bar v_{\ell})e^{i\Arg(v_{\ell+1}-v_\ell)}+( v_{\ell+1} - v_{\ell}) e^{-i\Arg(v_{\ell+1}-v_\ell)}]\varsigma\de\rho\,,
\end{equation*}
where the term in brackets is nothing but $2\vert v_{\ell+1}-v_\ell\vert=2\vol_1(\Lambda_\ell)$. It then follows that 
\begin{equation}
\langle\!\langle \ddc h_\Gamma,\varsigma\rangle\!\rangle
=
\sum_{\ell=0}^s \vert v_{\ell+1}-v_\ell\vert\int_{K_{\Lambda_\ell}}\varsigma\de\rho
=
\sum_{\Lambda\in{\mathcal B}(\Gamma,1)} 
\vol_1(\Lambda)
\langle\!\langle \lambda_\Lambda,\varsigma\rangle\!\rangle\,,
\end{equation}
i.e. $\ddc h_\Gamma=\sum_{\Lambda\in{\mathcal B}(\Gamma,1)}\vol_1(\Lambda)\lambda_\Lambda$. The support of $\ddc h_\Gamma$ is the union of the supports of the currents $\lambda_\Lambda$, as $\Lambda\in{\mathcal B}(\Gamma,1)$, i.e. $\overline{\Sigma}_1$.
\findim

\begin{thm}\label{supp-1}
Let $\Gamma\in {\mathcal P}(\mathbb C^n)$, then 
\begin{equation}\label{current1}
\ddc h_\Gamma
=
\sum_{\Lambda\in{\mathcal B}(\Gamma,1)}
\vol_1(\Lambda)
\lambda_\Lambda
\,.
\end{equation}
In particular $\Supp \ddc h_\Gamma=\overline{\Sigma}_1$.
\end{thm}
\noindent\pf 
If $\Gamma$ is reduced to a single point, $h_\Gamma$ is linear and both sides of \eqref{current1} are zero. If $\Gamma$ is not reduced to a single point, the current $\dc h_\Gamma$ is just a smooth form almost everywhere, i.e. on $\mathbb C^n\setminus\overline\Sigma_1$, so, for any test $(n-1,n-1)$-form $\varsigma$,
\begin{align*}
\langle\!\langle \ddc h_\Gamma,\varsigma\rangle\!\rangle
&=
\langle\!\langle \dc h_\Gamma,\de\varsigma\rangle\!\rangle
\\
&=
\lim_{\varepsilon\to0}
\int_{\Supp\varsigma\setminus(\overline\Sigma_1)_\varepsilon}
\dc h_\Gamma\wedge\de\varsigma
\\
&=
\lim_{\varepsilon\to0}
\sum_{v\in{\mathcal B}(\Gamma,0)}
\int_{K_v\cap\Supp\varsigma\setminus(\overline\Sigma_1)_\varepsilon}
\dc h_\Gamma\wedge\de \varsigma
\\
&=
\lim_{\varepsilon\to0}
\sum_{v\in{\mathcal B}(\Gamma,0)}
\int_{S_{v,\varepsilon}}
\dc h_\Gamma\wedge\de \varsigma\,,
\end{align*}
where $S_{v,\varepsilon}$ denotes the $(\varepsilon/2)$-neighbourhood of $K_v\cap\Supp\varsigma\setminus(\overline\Sigma_1)_\varepsilon$. 
Of course $K_{v}\cap\Supp\varsigma\setminus(\overline{\Sigma}_{1})_\varepsilon\subsetneq S_{v,\varepsilon}\subsetneq K_v$ and both $K_{v}\cap\Supp\varsigma\setminus(\overline{\Sigma}_{1})_\varepsilon$ and $S_{v,\varepsilon}$ approach the set $K_v\cap\Supp\varsigma$, as $\varepsilon\to 0^+$. However, unlike $K_{v}\cap\Supp\varsigma\setminus(\overline{\Sigma}_{1})_\varepsilon$, the set $S_{v,\varepsilon}$ has smooth boundary. This boundary can be decomposed into three parts: the pieces which belong to $\partial \Supp\varsigma$ (on which $\varsigma$ is zero), the spherical pieces over the singular points of $\partial[K_{v}\cap\Supp\varsigma\setminus(\overline{\Sigma}_{1})_\varepsilon]$, (which shrink to the corresponding singular points of $\partial[K_{v}\cap\Supp\varsigma\setminus(\overline{\Sigma}_{1})_\varepsilon]$ as $\varepsilon\to0$) and the flat pieces parallel and near to $\partial K_{v}$

For every $v\in{\mathcal B}(\Gamma,0)$, we have
\begin{equation*}
\dc h_\Gamma
=
(i/2)\sum_{\ell=1}^n v_\ell\de\bar z_\ell-\bar v_\ell\de z_\ell
\end{equation*}
on $K_v$, so that $\dc h_\Gamma\wedge\de\varsigma=-\de(\dc h_\Gamma\wedge\varsigma)$ on such a cone, then by Stokes' theorem $\langle\!\langle \ddc h_\Gamma,\varsigma\rangle\!\rangle$ equals the limit, as $\varepsilon\to0$, of
\begin{equation}\label{somma}
-\sum_{v\in{\mathcal B}(\Gamma,0)}
\int_{\partial S_{v,\varepsilon}}
(i/2)\left(\sum_{\ell=1}^n v_\ell\de\bar z_\ell-\bar v_\ell\de z_\ell\right)\wedge\varsigma\,.
\end{equation}
For any fixed $v\in{\mathcal B}(\Gamma,0)$ the only parts of the $\partial S_{v_\ell,\varepsilon}$ which really count are those which are parallel and near to the $1$-codimensional parts of $\partial K_v$, if any. Since, for any vertex $v$,
\begin{equation*}
\partial K_v=\bigcup_{\substack{\Lambda\in{\mathcal B}(\Gamma,1)\\ v\prec\Lambda}}
\overline K_\Lambda\,,
\end{equation*}
it follows that, for any fixed side $\Lambda$ of $\Gamma$, the set 
\begin{equation}\label{insieme}
\bigcup_{v\in{\mathcal B}(\Gamma,0)}
\partial S_{v,\varepsilon}
\end{equation}
may admit components parallel to $K_\Lambda$ and belonging to different half-spaces. If $\Lambda$ is the oriented line segment $[w_\Lambda,v_\Lambda]$, then $a_\Lambda=(v_\Lambda-w_\Lambda)\vol_1(\Lambda)^{-1}$ is an orienting orthonormal basis of $E_\Lambda$. Let $X^{\varepsilon}_{v(\Lambda)}$ and $X^{\varepsilon}_{w(\Lambda)}$ be, respectively, the subsets of \eqref{insieme} parallel and near to $K_\Lambda$ and belonging to different half-spaces, then $X^{\varepsilon}_{v(\Lambda)}$ inherits from \eqref{insieme} an orientation which is opposite to that of $E_\Lambda^\perp$, whereas $X^{\varepsilon}_{w(\Lambda)}$ inherits from \eqref{insieme} an orientation which equals that of $E_\Lambda^\perp$, so that \eqref{somma} becomes 
\begin{align*}
&
\lim_{\varepsilon\to0}
\sum_{\Lambda\in{\mathcal B}(\Gamma,1)}
\int_{X^{\varepsilon}_{v(\Lambda)}}
(i/2)\left(\sum_{\ell=1}^n v_{\Lambda,\ell}\de\bar z_\ell-\bar v_{\Lambda,\ell}\de z_\ell\right)\wedge\varsigma
\\
&-
\lim_{\varepsilon\to0}
\sum_{\Lambda\in{\mathcal B}(\Gamma,1)}
\int_{X^{\varepsilon}_{w(\Lambda)}}
(i/2)\left(\sum_{\ell=1}^n w_{\Lambda,\ell}\de\bar z_\ell-\bar w_{\Lambda,\ell}\de z_\ell\right)\wedge\varsigma\,.
\end{align*}
As $\varepsilon\to0$, both $X^{\varepsilon}_{v(\Lambda)}$ and $X^{\varepsilon}_{w(\Lambda)}$ approach $\overline{K}_\Lambda\cap\Supp\varsigma$, then the continuity of the integrands and the compactness of their supports imply that the sum of the preceding terms becomes 
\begin{equation*}
\sum_{\Lambda\in{\mathcal B}(\Gamma,1)}
\int_{\overline{K}_\Lambda\cap\Supp\varsigma}
\vol_1(\Lambda)\left((i/2)\sum_{\ell=1}^n a_{\Lambda,\ell}\de\bar z_\ell-\bar a_{\Lambda,\ell}\de z_\ell\right)\wedge\varsigma\,,
\end{equation*}
whence
\begin{equation*}
\langle\!\langle \ddc h_\Gamma,\varsigma\rangle\!\rangle
=
\sum_{\Lambda\in{\mathcal B}(\Gamma,1)}
\vol_1(\Lambda) 
\int_{K_\Lambda}
\iota^*_\Lambda(\upsilon_{E_\Lambda^\prime}\wedge\varsigma)
=
\sum_{\Lambda\in{\mathcal B}(\Gamma,1)}
\vol_1(\Lambda) 
\langle\!\langle \lambda_\Lambda,\varsigma\rangle\!\rangle
\,.
\end{equation*}
Observe that giving $\Lambda$ the opposite orientation does not affect the computation because both $\vol_1(\Lambda)$ and $\lambda_\Lambda$ do not depend on the chosen orientation.

The support of $\ddc h_\Gamma$ is the union of the supports of the currents $\lambda_\Lambda$, as $\Lambda\in{\mathcal B}(\Gamma,1)$, i.e. $\overline{\Sigma}_1$. 
The proof is thus complete.
\findim

\begin{rmq}
\lemref{lem1} and \thmref{current1} can be viewed as a generalisation of the following formula known as  \textit{la formule des sautes} in the french literature.
Let $S\subset \mathbb R$ be a discrete subset and $g:\mathbb R\to\mathbb R$ be a continuous piecewise smooth function admitting a corner point at each $s\in S$, (i.e. $g^\prime_-(s)=\lim_{x\to s^-}g^\prime(x)$ and $g^\prime_+(s)=\lim_{x\to s^+} g^\prime(x)$ are both finite but distinct). Then, it is well-known that the second derivative of $g$ is a measure $\mu$ with the following decomposition
\begin{equation}\label{jump}
\mu=g^{\prime\prime}\chi_{\mathbb R\setminus S}+\sum_{s\in S}\left[g^\prime_+(s)-g^\prime_-(s)\right]\delta_s\,,
\end{equation}
where $\chi_{\mathbb R\setminus S}$ is the indicator function of the set $\mathbb R\setminus S$, (i.e. $\chi_{\mathbb R\setminus S}(x)=1$ if $x\in \mathbb R\setminus S$ and $\chi_{\mathbb R\setminus S}(x)=0$ if $x\in S$) and $\delta_s$ is the Dirac delta on the point $s\in S$. 
The first term in \eqref{jump} is the ordinary second derivative of $g$ on the points where $f$ is smooth, whereas the second one gives the contribution of the singular points, each singular point producing a Dirac delta multiplied by the corresponding jump. 
If $g$ is convex then $\mu$ is a positive measure, indeed both the second derivative of $g$ (on the points where $g^\prime$ is smooth) and the jumps of $g^\prime$ (at singular points) are non-negative. When $g$ is a convex piecewise linear function, the absolutely continuous part of $\mu$ vanishes, and $\mu$ becomes a positive Borel measure which is singular with respect to the Lebesgue measure.
\end{rmq}

\subsection[The case $k=2$]{The case {\boldmath $k=2$}}

\begin{lem}\label{lemma-2}
Let $\Delta\in {\mathcal P}(\mathbb C^n)$ be an oriented polygon, then 
\begin{equation*}
\sum_{\Lambda\in{\mathcal B}(\Delta,1)}\vol_1(\Lambda)
\upsilon_{E_\Lambda^\prime}=0\,.
\end{equation*}
\end{lem}
\noindent
\pf Let us first suppose that $\Delta$ is a triangle with vertices $a,b,c$. If we denote $\Lambda_1$, $\Lambda_2$ and $\Lambda_3$ the oriented sides $[a,b]$, $[b,c]$ and $[c,a]$ respectively, we get
\begin{align*}
\vol_1(\Lambda_1)\upsilon_{E_{\Lambda_1}^\prime}
&=
\frac{i}{2}\sum_{\ell=1}^n (b_\ell-a_\ell)\de \bar z_\ell-(\bar b_\ell-\bar a_\ell)\de z_\ell\,,
\\
\vol_1(\Lambda_2)\upsilon_{E_{\Lambda_2}^\prime}
&=
\frac{i}{2}\sum_{\ell=1}^n(c_\ell-b_\ell)\de \bar z_\ell-(\bar c_\ell-\bar b_\ell)\de z_\ell\,,
\\
\vol_1(\Lambda_3)\upsilon_{E_{\Lambda_3}^\prime}
&=
\frac{i}{2}\sum_{\ell=1}^n(a_\ell-c_\ell)\de \bar z_\ell-(\bar a_\ell-\bar c_\ell)\de z_\ell\,,
\end{align*}
and consequently 
\begin{align*}
&
\vol_1(\Lambda_1)\upsilon_{E_{\Lambda_1}^\prime}
+
\vol_1(\Lambda_2)\upsilon_{E_{\Lambda_2}^\prime}
+
\vol_1(\Lambda_3)\upsilon_{E_{\Lambda_3}^\prime}
\\
&=
\frac{i}{2}\sum_{\ell=1}^n (b_\ell-a_\ell+c_\ell-b_\ell+a_\ell-c_\ell)\de \bar z_\ell-(\bar b_\ell-\bar a_\ell+\bar c_\ell-\bar b_\ell+\bar a_\ell-\bar c_\ell)\de z_\ell=0\,.
\end{align*}
If $\Delta$ is given the opposite orientation, a minus sign will affect each of the three terms involved in the preceding sum without changing the result.
If $\Delta$ is not a triangle, it can be decomposed into two or more triangles and by applying the preceding argument to each triangle one realises that the contributions of common sides cancel out for orientation reasons. \findim 

\begin{lem}\label{lemma-s2-1}
Let $\Delta\in {\mathcal P}(\mathbb C^n)$ an equidimensional oriented polygon then
\begin{equation*}
\sum_{\Lambda\in{\mathcal B}(\Delta,1)}\vol_1(\Lambda)
\upsilon_{E_\Lambda^\prime}\wedge\iota_\Lambda^*\dc h_\Gamma
=
\sum_{\Lambda\in{\mathcal B}(\Delta,1)}\vol_1(\Lambda)
h_\Delta(u_{\Lambda,\Delta})\upsilon_{E_\Lambda^\prime}\wedge \re\langle i u_{\Lambda,\Delta},-\rangle\,,
\end{equation*} 
where, for every $\Lambda\in{\mathcal B}(\Delta,1)$, $u_{\Lambda,\Delta}$ denotes the unit outer normal vector to the side $\Lambda$ of $\Delta$ and $\iota_\Lambda:K_\Lambda\to \mathbb C^n$ is the inclusion.
\end{lem}
\noindent\pf
Let $\Lambda\in{\mathcal B}(\Delta,1)$ be fixed. Recall that $E_\Lambda^\perp\cap\aff_\mathbb R \Lambda=\{p_\Lambda\}$ and $E_\Delta^\perp\cap\aff_\mathbb R \Delta=\{p_\Delta\}$. Then, by virtue of (\ref{decomposition}), $p_\Lambda=p_\Delta+h_\Delta(u_{\Lambda,\Delta}) u_{\Lambda,\Delta}$ and 
\begin{align*}
\upsilon_{E_\Lambda^\prime}\wedge\iota_\Lambda^* \dc h_\Gamma
&=
\upsilon_{E_\Lambda^\prime}\wedge\re\langle i p_\Lambda,-\rangle
\\
&=
\upsilon_{E_\Lambda^\prime}\wedge\re\langle i p_\Delta,-\rangle+h_\Delta(u_{\Lambda,\Delta})\upsilon_{E_\Lambda^\prime}\wedge \re\langle i u_{\Lambda,\Delta},-\rangle\,,
\end{align*}
By \lemref{lemma-2}, 
\begin{equation*}
\sum_{\Lambda\in{\mathcal B}(\Delta,1)}\vol_1(\Lambda)
\upsilon_{E_\Lambda^\prime}\wedge\re\langle i p_\Delta,-\rangle
=
-\re\langle i p_\Delta,-\rangle
\wedge
\sum_{\Lambda\in{\mathcal B}(\Delta,1)}\vol_1(\Lambda)\upsilon_{E_\Lambda^\prime}
=
0
\end{equation*}
whence
\begin{equation*}
\sum_{\Lambda\in{\mathcal B}(\Delta,1)}\vol_1(\Lambda)
\upsilon_{E_\Lambda^\prime}\wedge\iota_\Lambda^* \dc h_\Gamma
=
\sum_{\Lambda\in{\mathcal B}(\Delta,1)}\vol_1(\Lambda)
h_\Delta(u_{\Lambda,\Delta})\upsilon_{E_\Lambda^\prime}\wedge \re\langle i u_{\Lambda,\Delta},-\rangle
\end{equation*}
which proves the lemma.\findim

\begin{lem}\label{lemma-s2}
Let $\Delta\in {\mathcal P}(\mathbb C^n)$ an equidimensional oriented polygon then, for every $\Lambda\in{\mathcal B}(\Delta,1)$,
\begin{equation}\label{s2-2}
\vol_1(\Lambda)
h_\Delta(u_{\Lambda,\Delta})\upsilon_{E_\Lambda^\prime}\wedge \re\langle i u_{\Lambda,\Delta},-\rangle
=
\varrho(\Delta)
\vol_1(\Lambda)h_\Delta(u_{\Lambda,\Delta})
\upsilon_{E^\prime_\Delta}
\,.
\end{equation} 
In particular 
\begin{equation}\label{s2}
\sum_{\Lambda\in{\mathcal B}(\Delta,1)}
\vol_1(\Lambda)
h_\Delta(u_{\Lambda,\Delta})\upsilon_{E_\Lambda^\prime}\wedge \re\langle i u_{\Lambda,\Delta},-\rangle
=
2\varrho(\Delta)
\vol_2(\Delta)
\upsilon_{E^\prime_\Delta}\,.
\end{equation}
\end{lem}
\noindent\pf
As $\Delta$ is oriented, each of its facets gets the orientation induced by the outer unit normal vector.  
If $a_\Lambda$ is an orienting basis of $E_\Lambda$, then $(u_{\Lambda,\Delta},a_\Lambda)$ is an orienting orthonormal basis of $E_\Delta$, $\upsilon_{E_\Lambda^\prime}=\re\langle i a_\Lambda,-\rangle$ and $\vol_1(\Lambda)\re\langle i a_\Lambda,-\rangle\wedge \re\langle i u_{\Lambda,\Delta},-\rangle$ is a non-zero $2$-form.
Consider the vectors 
\begin{align*}
t_{\Lambda,1}
&=
iu_{\Lambda,\Delta}-\re\langle iu_{\Lambda,\Delta},a_\Lambda\rangle a_\Lambda\,, 
\\
t_{\Lambda,2}
&=
ia_\Lambda-\re\langle ia_\Lambda,u_{\Lambda,\Delta}\rangle u_{\Lambda,\Delta}
\end{align*}
$w_{\Lambda,1}=t_{\Lambda,1}/\sqrt{\varrho(\Delta)}$ and $w_{\Lambda,2}=t_{\Lambda,2}/\sqrt{\varrho(\Delta)}$. According to \thmref{figo}, since $\Delta$ is equidimensional, the sequence  $u_{\Lambda,\Delta},a_\Lambda,w_{\Lambda,1},w_{\Lambda,2}$ is an orthonormal basis of $\lin_\mathbb C E_\Delta$ over $\mathbb R$ such that $w_{\Lambda,1},w_{\Lambda,2}$ span $E_\Delta^\prime$. For dimensional reasons, ($k=2$), this basis does not belong to the positive orientation of $\lin_\mathbb C E_\Delta$ so, by~\eqref{convenient}, we have $\upsilon_{E_\Delta^\prime}=-\re\langle w_{\Lambda,1},-\rangle\wedge\re\langle w_{\Lambda,2},-\rangle$. Let us complete the preceding basis by adding an orienting orthonormal basis $b_5,\ldots,b_{2n}$ of $E_\Delta^{\perp_\mathbb C}$ over $\mathbb R$. 
Since $ia_\Lambda$ and $iu_{\Lambda,\Delta}$ belong to $\lin_\mathbb C E_\Delta$, both the linear forms $\re\langle ia_\Lambda,-\rangle$ and $\re\langle iu_{\Lambda,\Delta},-\rangle$ vanish on $E_\Delta^{\perp_\mathbb C}$, i.e. they both vanish on each of the vectors $b_5,\ldots,b_{2n}$.
A direct computation shows that $\re\langle ia_\Lambda,w_{\Lambda,2}\rangle=\re\langle i u_{\Lambda,\Delta},w_{\Lambda,1}\rangle=\sqrt{\varrho(\Delta)}$ and $\re\langle ia_\Lambda,w_{\Lambda,1}\rangle=0$.
Indeed,
\begin{align*}
\re\langle ia_\Lambda,w_{\Lambda,2}\rangle
&=
\frac{\re\langle ia_\Lambda, 
ia_\Lambda-\re\langle i a_\Lambda,u_{\Lambda,\Delta}\rangle u_{\Lambda,\Delta}\rangle}{\sqrt{\varrho(\Delta)}}
\\
&=
\frac{1-(\re\langle i a_\Lambda,u_{\Lambda,\Delta}\rangle)^2}{\sqrt{\varrho(\Delta)}}
=
\sqrt{\varrho(\Delta)}\,,
\\
\re\langle i u_{\Lambda,\Delta},w_{\Lambda,1}\rangle
&=
\frac{\re\langle iu_{\Lambda,\Delta},iu_{\Lambda,\Delta}-\re\langle i u_{\Lambda,\Delta},a_\Lambda\rangle a_\Lambda\rangle}{\sqrt{\varrho(\Delta)}}
\\
&=
\frac{1-(\re\langle i u_{\Lambda,\Delta},a_\Lambda\rangle)^2}{\sqrt{\varrho(\Delta)}}
\\
&=
\sqrt{\varrho(\Delta)}\,,
\\
\re\langle ia_\Lambda,w_{\Lambda,1}\rangle
&=
\frac{\re\langle ia_\Lambda,iu_{\Lambda,\Delta}-\re\langle i u_{\Lambda,\Delta},a_\Lambda\rangle a_\Lambda\rangle}{\sqrt{\varrho(\Delta)}}
\\
&=
\frac{\re\langle a_\Lambda,u_{\Lambda,\Delta}\rangle-\re\langle i u_{\Lambda,\Delta},a_\Lambda\rangle \re\langle ia_\Lambda,a_\Lambda\rangle}{\sqrt{\varrho(\Delta)}}
\\
&=
0\,.
\end{align*}
It follows that
\begin{align*}
&
\left[
\vol_1(\Lambda)h_\Delta(u_{\Lambda,\Delta})
\re\langle ia_\Lambda,-\rangle\wedge \re\langle i u_{\Lambda,\Delta},-\rangle\right]
(w_{\Lambda,1},w_{\Lambda,2})
\\
&=
\vol_1(\Lambda)h_\Delta(u_{\Lambda,\Delta})
\left[
-\re\langle i u_{\Lambda,\Delta},-\rangle\wedge \re\langle ia_\Lambda,-\rangle\right]
(w_{\Lambda,1},w_{\Lambda,2})
\\
&=
2\varrho(\Delta)
\left(\frac{\vol_1(\Lambda)h_\Delta(u_{\Lambda,\Delta})}{2}\right)
(-1)
\\
&=
2\varrho(\Delta)
\left(\frac{\vol_1(\Lambda)h_\Delta(u_{\Lambda,\Delta})}{2}\right)
\left[
-\re\langle w_{\Lambda,1},-\rangle\wedge \re\langle w_{\Lambda,2},-\rangle\right]
(w_{\Lambda,1},w_{\Lambda,2})
\\
&=
2\varrho(\Delta)
\left(\frac{\vol_1(\Lambda)h_\Delta(u_{\Lambda,\Delta})}{2}\right)
\upsilon_{E^\prime_\Delta}(w_{\Lambda,1},w_{\Lambda,2})
\,,
\end{align*} 
which implies (\ref{s2-2}). Observe that changing the orientation of $\Delta$, and hence of $\Lambda$, affects by a minus sign the vectors $a_\Lambda,w_{\Lambda,1}$ and the forms $\upsilon_{E_\Lambda^\prime},\upsilon_{E_\Delta^\prime}$, however it changes neither the vectors $u_{\Lambda,\Delta},w_{\Lambda,2}$ nor the preceding computations. By summing the relations (\ref{s2-2}) as $\Lambda$ runs in ${\mathcal B}(\Delta,1)$, the equality (\ref{piramidi}) yields the desired relation (\ref{s2}). The lemma is thus completely proved.
\findim

\begin{rmq}\label{onesta}
{Remark that the basis 
$$
u_{\Lambda,\Delta},a_\Lambda,-w_{\Lambda,1},w_{\Lambda,2},b_5,\ldots,b_{2n}\qquad \text{and} \qquad a_\Lambda,-u_{\Lambda,\Delta},-w_{\Lambda,1},w_{\Lambda,2},b_5,\ldots,b_{2n}
$$ 
are both positively oriented, so $\mathbb C^n$ gives $E_\Lambda^\perp$ the orientation induced by $a_\Lambda$ and $E_\Delta^\perp$ gets from $E_\Lambda^\perp$ the orientation induced by $-u_{\Lambda,\Delta}$. This means that $-u_{\Lambda,\Delta}$ orients $E_\Delta^\perp$ in the right way.}
\end{rmq}

\begin{thm}\label{supp-2}
Let $\Gamma\in {\mathcal P}(\mathbb C^n)$, then 
\begin{equation}\label{current-2}
(\ddc h_\Gamma)^{\wedge 2}
=
2\sum_{\Delta\in{\mathcal B}(\Gamma,2)}
\varrho(\Delta)
\vol_2(\Delta)
\lambda_\Delta
\,.
\end{equation}
In particular $\Supp (\ddc h_\Gamma)^{\wedge 2}=\overline{\Sigma}_2$.
\end{thm}
\noindent\pf By \thmref{supp-1} we already know that
\begin{equation}
\ddc h_\Gamma
=
\sum_{\Lambda\in{\mathcal B}(\Gamma,1)}
\vol_1(\Lambda)
\lambda_\Lambda
\,,
\end{equation}
so, for any $(n-2,n-2)$-test form $\varsigma$,
\begin{align*} 
\langle\!\langle (\ddc h_\Gamma)^{\wedge 2},\varsigma\rangle\!\rangle
&=
\langle\!\langle \ddc(h_\Gamma\ddc h_\Gamma),\varsigma\rangle\!\rangle\\
&=
\langle\!\langle h_\Gamma\ddc h_\Gamma,\ddc\varsigma\rangle\!\rangle\\
&=
\langle\!\langle \ddc h_\Gamma,h_\Gamma\ddc\varsigma\rangle\!\rangle\\
&=
\sum_{\Lambda\in{\mathcal B}(\Gamma,1)}\vol_1(\Lambda)\langle\!\langle \lambda_{\Lambda},h_\Gamma\ddc\varsigma\rangle\!\rangle\\
&=
\sum_{\Lambda\in{\mathcal B}(\Gamma,1)}\vol_1(\Lambda)\int_{K_\Lambda} \upsilon_{E_\Lambda^\prime}\wedge h_\Gamma\ddc \varsigma\\
&=
\lim_{\varepsilon\to0}
\sum_{\Lambda\in{\mathcal B}(\Gamma,1)}\vol_1(\Lambda)\int_{K_\Lambda\cap\Supp\varsigma\setminus(\overline{\Sigma}_2)_\varepsilon} \upsilon_{E_\Lambda^\prime}\wedge h_\Gamma\ddc \varsigma
\\
&=
\lim_{\varepsilon\to0}
\sum_{\Lambda\in{\mathcal B}(\Gamma,1)}\vol_1(\Lambda)\int_{S_{\Lambda,\varepsilon}}\upsilon_{E_\Lambda^\prime}\wedge h_\Gamma\ddc \varsigma
\,,
\end{align*}
where $S_{\Lambda,\varepsilon}=K_\Lambda\cap (K_\Lambda\cap\Supp\varsigma\setminus(\overline{\Sigma}_2)_\varepsilon)_{\varepsilon/2}$, i.e. the intersection of $K_\Lambda$ with the $(\varepsilon/2)$-neighbourhood of $K_\Lambda\cap\Supp\varsigma\setminus(\overline{\Sigma}_2)_\varepsilon$. Unlike $K_\Lambda\cap\Supp\varsigma\setminus(\overline{\Sigma}_2)_\varepsilon$, the set $S_{\Lambda,\varepsilon}$ has a smooth relative boundary.

The 1-form $\upsilon_{E_\Lambda^\prime}$ is a closed on $E_\Lambda^\prime\subset E_\Lambda^\perp$, moreover the $(2n-2)$-forms $-\de h_\Gamma\wedge\dc\varsigma$ and $\dc h_\Gamma\wedge \de \varsigma$ have the same $(n-1,n-1)$-parts on $E^{\perp_\mathbb C}_\Lambda$ so that on $K_\Lambda$ one has
\begin{align*}
\de(v_{E^\prime_\Lambda}\wedge h_\Gamma\dc\varsigma)
&=
\de \upsilon_{E_\Lambda^\prime}\wedge h_\Gamma\dc\varsigma-\upsilon_{E_\Lambda^\prime}\wedge \de(h_\Gamma\dc\varsigma)\\
&=
0-\upsilon_{E_\Lambda^\prime}\wedge \de h_\Gamma\wedge\dc\varsigma-\upsilon_{E_\Lambda^\prime}\wedge h_\Gamma\ddc\varsigma
\\
&=
\upsilon_{E_\Lambda^\prime}\wedge \dc h_\Gamma\wedge\de\varsigma-
\upsilon_{E_\Lambda^\prime}\wedge h_\Gamma\ddc\varsigma\,.
\end{align*}
Since on $K_\Lambda$ the $1$-form $\dc h_\Gamma$ has constant coefficients, it follows that
\begin{equation*}
\upsilon_{E_\Lambda^\prime}\wedge h_\Gamma\ddc \varsigma
=
\de (\upsilon_{E_\Lambda^\prime}\wedge\dc h_\Gamma\wedge\varsigma)-\de (\upsilon_{E_\Lambda^\prime}\wedge h_\Gamma\dc\varsigma)\,,
\end{equation*}
whence
\begin{align}
\langle\!\langle (\ddc h_\Gamma)^{\wedge 2},\varsigma\rangle\!\rangle
&=
\lim_{\varepsilon\to0}
\sum_{\Lambda\in{\mathcal B}(\Gamma,1)}\vol_1(\Lambda)\int_{\relbd S_{\Lambda,\varepsilon}}
\upsilon_{E_\Lambda^\prime}\wedge\dc h_\Gamma\wedge\varsigma\label{t1}
\\
&-
\lim_{\varepsilon\to0}
\sum_{\Lambda\in{\mathcal B}(\Gamma,1)}\vol_1(\Lambda)\int_{\relbd S_{\Lambda,\varepsilon}}
\upsilon_{E_\Lambda^\prime}\wedge h_\Gamma\dc\varsigma\,.\label{t2}
\end{align}
For any fixed side $\Lambda$, if $\Supp\varsigma\cap \relbd K_\Lambda\neq\varnothing$, one has 
$K_{\Lambda}\cap\Supp\varsigma\setminus(\overline{\Sigma}_{2})_\varepsilon\subsetneq S_{\Lambda,\varepsilon}\subsetneq K_\Lambda$ and both $K_{\Lambda}\cap\Supp\varsigma\setminus(\overline{\Sigma}_{1})_\varepsilon$ and $S_{\Lambda,\varepsilon}$ approach the set $K_\Lambda\cap\Supp\varsigma$, as $\varepsilon\to 0^+$. The relative boundary of $S_{\Lambda,\varepsilon}$ can be decomposed into three parts: the pieces belonging to $\partial \Supp\varsigma$ (on which $\varsigma$ is zero), the spherical pieces over the singular points of $\partial[K_{\Lambda}\cap\Supp\varsigma\setminus(\overline{\Sigma}_{2})_\varepsilon]$, (which shrink to the corresponding singular points of $\partial[K_{\Lambda}\cap\Supp\varsigma\setminus(\overline{\Sigma}_{2})_\varepsilon]$ as $\varepsilon\to0$) and the flat pieces parallel and near to the $(2n-2)$-dimensional parts of $\relbd K_{\Lambda}$, if any.
The latter are the only parts of $\relbd S_{\Lambda,\varepsilon}$ which really count. The $(2n-2)$-dimensional parts of $\relbd K_\Lambda$ are nothing but the cones which are dual to the $2$-faces $\Delta$ of $\Gamma$ admitting $\Lambda$ as a side. 

\begin{figure}[ht]
\begin{center}
\begin{tikzpicture}[scale=2.85,path fading=fade right]
\draw[dashed,thick,line cap=round, line join=round] (1,1,0)--(2,1,0);
\draw[dashed,thick,line cap=round, line join=round] (1,1,1)--(1,1,0)--(1,2,0);
\fill[fill=orange,nearly opaque,path fading,fading transform={rotate=45},line cap=round, line join=round] (2,2,0)--(4.2,2,0)--(4.2,3.75,0)--(2,3.75,0)--cycle;
\draw (4.2,3.75,0) node[below left]{$\mathbb R^3\cap K_{\Lambda_1}$};

\fill[fill=green,nearly opaque,path fading,fading transform={rotate=135},line cap=round, line join=round] (2,2,0)--(0.5,2,0)--(0.5,3.75,0)--(2,3.75,0)--cycle;
\draw (0.8,3.75,0) node[below]{$\mathbb R^3\cap K_{\Lambda_2}$};

\draw[draw=black,dotted,line cap=round, line join=round] (2,3.75,0)--(2,4.1,0);

\draw[draw=black,dotted,line cap=round, line join=round] (2.2,2.2,0)--(2.2,4.2,0);
\draw[draw=black,dotted,line cap=round, line join=round] (2.1,2.2,0)--(2.1,4.2,0);

\draw[<->,dashed] (1.8,4.1,0)--(2,4.1,0);
\draw (1.9,4.1,0) node[above]{$\varepsilon$};
\draw[<->,dashed] (2.1,4.2,0)--(2.2,4.2,0);
\draw (2.15,4.2,0) node[above]{$\varepsilon/2$};
\draw[red,very thick,->,line cap=round, line join=round] (2,2,0)--(2,2.75,0);
\draw[red,very thick,line cap=round, line join=round] (2,2.75,0)--(2,3.75,0) node[black,above]{$\mathbb R^3\cap K_{\Delta}$};

\draw[violet!20!white,thick,->,line cap=round, line join=round] (2,2,0)--(2.9,2,0);
\draw[violet!20!white, thick,line cap=round, line join=round] (2.9,2,0)--(4.2,2,0);
\draw[<->,dashed] (4.4,2,0)--(4.4,2.2,0);
\draw (4.4,2.1,0) node[left]{$\varepsilon$};
\draw[<->,dashed] (4.6,2.1,0)--(4.6,2.2,0);
\draw (4.6,2.15,0) node[right]{$\varepsilon/2$};
\draw[dotted] (4.2,2,0)--(4.4,2,0);
\draw[dotted] (3.6,2.2,0)--(4.6,2.2,0);
\draw[dotted] (3.6,2.1,0)--(4.6,2.1,0);

\draw[thick,brown] (2.2,2.2,0.4)--(2.2,2.2,0);
\draw[thick,brown] (1,2.2,0.4)--(1,2.2,0);
\draw[thick,brown,dashed] (2.2,2.2,0)--(2.2,2.2,-0.4);
\draw[thick, brown] (2.2,3.25,0.4)--(2.2,3.25,0);
\draw[thick, brown,dashed] (2.2,3.25,0)--(2.2,3.25,-0.4);
\draw[thick,brown] (3.6,2.2,0.4)--(3.6,2.2,0);
\draw[thick,brown,dashed] (3.6,2.2,0)--(3.6,2.2,-0.4);
\draw[thick,brown] (1.8,2.2,0.4)--(1.8,2.2,0);
\draw[thick,brown,dashed] (1.8,2.2,0)--(1.8,2.2,-0.4);
\draw[thick,brown,dashed] (1.8,3.25,-0.4)--(1.8,2.2,-0.4);
\draw[thick,brown,dashed] (1.8,3.25,0)--(1.8,3.25,-0.4);
\draw[thick,brown] (1.8,3.25,0)--(1.8,3.25,0.4);
\draw[thick,brown] (1.8,2.2,0.4)--(1.8,3.25,0.4);
\draw[thick,brown] (1,2.2,0.4)--(1.8,2.2,0.4);
\draw[thick,brown,dashed] (1,2.2,-0.4)--(1.8,2.2,-0.4);

\draw[brown] (3.25,1,0) node{$S$};

\draw[thick, brown,dashed,line cap=round, line join=round] (2.2,3.25,-0.4) -- (2.2,2.2,-0.4)--(3.6,2.2,-0.4)--(3.6,2.5,-0.4)--(3.6,2.5,-0.4) arc(0:90:0.75)--(1,3.25,-0.4)-- (1,1.5,-0.4) arc(180:270:0.5)-- (2.6,1,-0.4) arc (270:360:1) -- (3.6,2.2,-0.4);

\draw[thick, brown]  (3.314,1.3,-0.4) arc (315.5:355:1);

\draw[thick, brown] (3.314,1.3,0.4)--(3.314,1.3,-0.4);
\draw[thick, brown] (1,3.25,0.4)--(1,3.25,0);
\draw[thick, brown,dashed] (1,3.25,0)--(1,3.25,-0.4);


\fill[fill=blue!30!white, very thick,line cap=round, line join=round] (3.6,2.2,0)--(3.6,2.5,0)--(3.6,2.5,0) arc(0:90:0.75)--(2.2,3.25,0)-- (2.2,2.2,0) --(3.6,2.2,0);

\draw[->,very thick, brown,line cap=round, line join=round] (3.6,2.2,0)--(3.6,2.5,0);
\draw[->,very thick,draw=black,line cap=round, line join=round] (3.7,2.2,0)--(3.7,2.5,0);

\draw[->,very thick, brown,line cap=round, line join=round] (3.6,2.5,0) arc(0:90:0.75);
\draw[->,very thick,draw=black,line cap=round, line join=round] (3.7,2.5,0) arc(0:90:0.85);

\draw[very thick, brown,line cap=round, line join=round] (2.9,3.25,0) --(2.2,3.25,0);
\draw[very thick,draw=black,line cap=round, line join=round] (2.9,3.35,0) --(2.2,3.35,0);

\draw[blue!50!white,very thick,line cap=round, line join=round] (2.2,3.35,0) arc (90:180:0.1);
\fill[blue!50!white] (2.2,3.25,0) circle(0.5pt);

\draw[violet!20!white,very thick,->,line cap=round, line join=round] (2.2,2.2,0)--(2.9,2.2,0);
\draw[violet!20!white,very thick,->,line cap=round, line join=round] (2.2,2.1,0)--(2.9,2.1,0);

\draw[violet!20!white,very thick,line cap=round, line join=round] (2.9,2.2,0)--(3.6,2.2,0);
\draw[violet!20!white,very thick,line cap=round, line join=round] (2.9,2.1,0)--(3.6,2.1,0);

\draw[blue!50!white,very thick,line cap=round, line join=round] (3.6,2.1,0) arc(270:360:0.1);
\fill[blue!50!white] (3.6,2.2,0) circle(0.5pt);

\draw (2.2,2.725,0) node[black, right]{\footnotesize$X_{\Lambda_1(\Delta)}^\varepsilon$};

\draw[thick, brown,line cap=round, line join=round] (2.2,3.25,0.4) -- (2.2,2.2,0.4)--(3.6,2.2,0.4)--(3.6,2.5,0.4)--(3.6,2.5,0.4) arc(0:90:0.75)--(1,3.25,0.4) -- (1,2,0.4) ;
\draw[dashed,thick, brown,line cap=round, line join=round]  (1,2,0.4)-- (1,1.5,0.4) arc(180:270:0.5) -- (2,1,0.4) ;
\draw[thick, brown,line cap=round, line join=round] (2,1,0.4)-- (2.6,1,0.4) arc (270:360:1) -- (3.6,2.2,0.4);

\draw (2.9,2.5,0) node{\footnotesize$K_{\Lambda_1}\cap S\setminus(\overline\Sigma_2)_\varepsilon$};

\draw[red,->,very thick,line cap=round, line join=round] (2.2,3.25,0)-- (2.2,2.725,0) ;
\draw[red,->,very thick,line cap=round, line join=round] (2.1,3.25,0)-- (2.1,2.725,0) ;

\draw[red,very thick,line cap=round, line join=round] (2.2,2.725,0)-- (2.2,2.2,0) ;
\draw[red,very thick,line cap=round, line join=round] (2.1,2.725,0)-- (2.1,2.2,0) ;

\draw[blue!50!white,very thick,line cap=round, line join=round] (2.1,2.2,0) arc(180:270:0.1) ;
\fill[blue!50!white] (2.2,2.2,0) circle(0.5pt);

\draw[thick,brown] (1,2,0.4)--(1,2.25,0.4);
\draw[thick,brown] (2,1,0.4)--(2.25,1,0.4);

\fill[blue!50!white] (2.2,3.25,0) circle(0.5pt);

\fill[blue!30!white,line cap=round] (1,3.25,0)--(1,2.2,0)--(1.8,2.2,0)--(1.8,3.25,0)--cycle;
\draw[red,very thick,line cap=round] (1.8,2.75,0) --(1.8,2,0);
\draw[red,->,very thick,line cap=round] (1.8,3.25,0)--(1.8,2.725,0) ;
\draw[brown,very thick,line cap=round] (1.8,3.25,0) --(1.4,3.25,0);
\draw[brown,very thick,line cap=round,->] (1,3.25,0) --(1.4,3.25,0);
\draw[brown,very thick,line cap=round,->] (1,2.2,0) --(1,2.725,0);
\draw[brown,very thick,line cap=round] (1,3.25,0)--(1,2.725,0) ;
\draw[gray!20!white,very thick,line cap=round,->] (1.8,2.2,0)--(1.4,2.2,0) ;
\draw[gray!20!white,very thick,line cap=round] (1,2.2,0)--(1.4,2.2,0) ;
\fill[blue!50!white] (1.8,3.25,0) circle(0.5pt);
\fill[blue!50!white] (1.8,2.2,0) circle(0.5pt);
\fill[blue!50!white] (1,2.2,0) circle(0.5pt);
\fill[blue!50!white] (1,3.25,0) circle(0.5pt);
%
\draw[red,very thick,line cap=round] (1.9,2.2,0) --(1.9,2.75,0);
\draw[red,->,very thick,line cap=round] (1.9,3.25,0) --(1.9,2.725,0);
\draw[black,very thick,line cap=round] (1.8,3.35,0) --(1.4,3.35,0);
\draw[black,very thick,line cap=round,->] (1,3.35,0) --(1.4,3.35,0);
\draw[black,very thick,line cap=round,->] (0.9,2.2,0)-- (0.9,2.725,0);
\draw[black,very thick,line cap=round] (0.9,3.25,0)--(0.9,2.725,0) ;
\draw[gray!20!white,very thick,line cap=round,->] (1.8,2.1,0)--(1.4,2.1,0) ;
\draw[gray!20!white,very thick,line cap=round] (1,2.1,0)--(1.4,2.1,0) ;
\draw[very thick,blue!30!white] (1.9,3.25,0) arc(0:90:0.1);
\draw[very thick,blue!30!white] (1,3.35,0) arc(90:180:0.1);
\draw[very thick,blue!30!white] (0.9,2.2,0) arc(180:270:0.1);
\draw[very thick,blue!30!white] (1.8,2.1,0) arc(270:360:0.1);

\draw[gray!20!white,thick,line cap=round] (0.5,2,0)--(1,2,0);


\draw[very thick] (1,1,1)--(2,1,1);
\filldraw[fill=gray!20!white, very thick,nearly opaque,line cap=round, line join=round] (1,1,1)--(1,1,0)--(1,2,0)--(1,2,1)--cycle;
\filldraw[fill=violet!20!white, very thick,nearly opaque,line cap=round, line join=round] (2,2,1)--(2,1,1)--(2,1,0)--(2,2,0)--cycle;
\filldraw[fill=red!40!white, very thick,nearly opaque,line cap=round, line join=round] (1,2,0)--(1,2,1)--(2,2,1)--(2,2,0)--cycle;
\draw (1.5,2,0.5) node{$\Delta$};
\draw (2,2,0.5) node[black, right]{$\Lambda_1$};
\draw (0.95,2,0.5) node[black, left]{$\Lambda_2$};
\draw[thick,brown] (1.8,2.2,0.4)--(1.8,3.25,0.4);
\draw[thick,brown] (1.8,3.25,0.4)--(1.8,3.25,0);
\draw (1.4,2.5,0) node{\footnotesize$K_{\Lambda_2}\cap S\setminus(\overline\Sigma_2)_\varepsilon$};
\draw (1.8,2.725,0) node[black, left]{\footnotesize$X_{\Lambda_2(\Delta)}^\varepsilon$};
\draw[orange,very thick,line cap=round, line join=round] (2,2,0)--(2,2,0.5);
\draw[orange,very thick,line cap=round, line join=round] (2,2,0.5)--(2,2,1);

\fill (2,2,0) circle(0.5pt) node[below right]{$o$};
\draw[very thick,green,line cap=round] (1,2,1)--(1,2,0.5);
\draw[very thick,green,line cap=round] (1,2,0.5)--(1,2,0);

\draw[dotted,line cap=round] plot[smooth,variable=\t,domain=0:360,samples=360] ({2+0.2*cos(\t)},{2.2},{0.2*sin(\t)});
\draw[dotted,line cap=round] plot[smooth,variable=\t,domain=0:360,samples=360] ({2+0.2*cos(\t)},{3.25},{0.2*sin(\t)});
\draw[dotted,line cap=round] plot[smooth,variable=\t,domain=0:360,samples=360] ({1.8+0.1*cos(\t)},{3.25},{0.1*sin(\t)});
\draw[dotted,line cap=round] plot[smooth,variable=\t,domain=0:360,samples=360] ({1.8+0.1*cos(\t)},{2.2},{0.1*sin(\t)});
\draw[dotted,line cap=round] plot[smooth,variable=\t,domain=0:360,samples=360] ({1+0.1*cos(\t)},{3.25},{0.1*sin(\t)});
\draw[draw=black,dotted,line cap=round, line join=round] (1.8,3.25,0)--(1.8,4.1,0);
\draw[dotted,line cap=round] plot[smooth,variable=\t,domain=0:360,samples=360] ({2.2+0.1*cos(\t)},{3.25},{0.1*sin(\t)},);
\draw[dotted,line cap=round] plot[smooth,variable=\t,domain=0:360,samples=360] ({1+0.1*cos(\t)},2.2,{0.1*sin(\t)},);
\draw[dotted,line cap=round] plot[smooth,variable=\t,domain=0:360,samples=360] ({2.2+0.1*cos(\t)},2.2,{0.1*sin(\t)},);

\fill (1,2,1) circle(0.5pt) node[left]{$q$};

\draw[very thick,->,line cap=round] (2.5,4.5,0)--(2.5,4.5,-0.5) node[above right]{$u_{o,\Lambda_1}$};
\draw[->,very thick,red,line cap=round] (2.5,4.5,0)--(2.5,5,0) node[above]{$u_{\Delta,\Gamma}$};
\draw[->,very thick,orange,line cap=round] (2.5,4.5,0)--(3,4.5,0) node[right]{$u_{\Lambda_1,\Delta}$};

\draw[->,very thick,red,line cap=round] (1.5,4.5,0)--(1.5,5,0) node[above]{$u_{\Delta,\Gamma}$};
\draw[->,very thick,green,line cap=round] (1.5,4.5,0)--(1,4.5,0) node[left]{$u_{\Lambda_2,\Delta}$};
\draw[very thick,->,line cap=round] (1.5,4.5,0)--(1.5,4.5,0.5) node[below left]{$u_{q,\Lambda_2}$};

\draw[thick,brown] (1,3.25,0.4)--(2.2,3.25,0.4);

\end{tikzpicture}
\caption{$S=\Supp\varsigma$ is the transparent cylinder with brown contour, $\Gamma$ is a cube with a vertex in the origin, $\Delta$ is the red face, $\Lambda_1$ is the orange edge of $\Delta$ with a vertex in the origin $o$, whereas $\Lambda_2$ is the green edge parallel to $\Lambda_1$. The space $\mathbb C^3$ is spanned by the orthonormal basis $(u_{\Delta,\Gamma},u_{\Lambda_1,\Delta},u_{o,\Lambda_1})$ or $(u_{\Delta,\Gamma},u_{\Lambda_2,\Delta},u_{q,\Lambda_2})$ and three more unit vectors $w_1,w_2,w_3$ spanning $E_\Gamma^\perp$ thus getting in both cases positively oriented basis of $\mathbb C^3$. This choice implies that $E_\Delta$ is oriented by $(u_{\Lambda_1,\Delta},u_{o,\Lambda_1})$ or $(u_{\Lambda_2,\Delta},u_{q,\Lambda_1})$. The set $K_{\Lambda_1}\cap\Supp\varsigma\setminus(\overline\Sigma_2)_\varepsilon$ is the blue slice (on the right) of $\Supp\varsigma$ along $K_{\Lambda_1}\setminus(\overline\Sigma_2)_\varepsilon$ while  $S_{\Lambda_1,\varepsilon}$ is the set bounded by the smooth path all around $K_{\Lambda_1}\cap\Supp\varsigma\setminus(\overline\Sigma_2)_\varepsilon$. $X_{\Lambda_1(\Delta)}^\varepsilon$ is the red segment on the boundary of $S_{\Lambda_1,\varepsilon}$.  
Along the black part of the boundary of $S_{\Lambda_1,\varepsilon}$ the form $\varsigma$ is zero because this part is included in $\partial\Supp\varsigma$. On the blue circular (in fact cylindrical!) pieces of $\partial S_{\Lambda_1,\varepsilon}$, over the singular points of $K_{\Lambda_1}\cap\Supp\varsigma\setminus(\overline\Sigma_2)_\varepsilon$, the integral vanishes in the limit. The same considerations are valid for $\Lambda_2$.}
\label{corrente}
\end{center}
\end{figure}
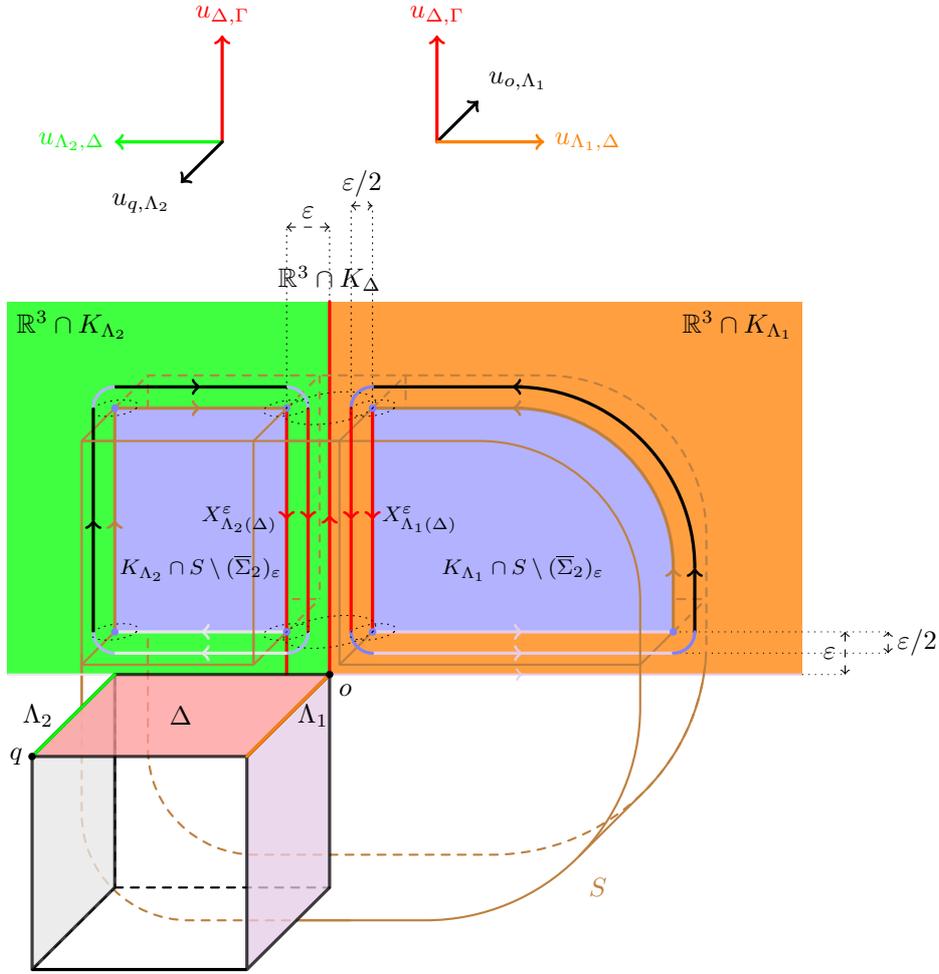

Let us denote $X^\varepsilon_{\Lambda(\Delta)}$ the subset of $\relbd S_{\Lambda,\varepsilon}$ parallel and near to $K_\Delta$ and observe, by the way, that $X^\varepsilon_{\Lambda(\Delta)}$ approaches $\overline{K}_\Delta\cap \Supp\varsigma$ as $\varepsilon\to0$. 
As a part of the manifold $\relbd S_{\Lambda,\varepsilon}$, the set $X^\varepsilon_{\Lambda(\Delta)}$ has the orientation of $\relbd S_{\Lambda,\varepsilon}$, namely the orientation induced by the outer unit normal vector to $S_{\Lambda,\varepsilon}$. On $X^\varepsilon_{\Lambda(\Delta)}$ the outer unit normal vector to $S_{\Lambda,\varepsilon}$ equals the vector pointing outside $K_\Lambda\cap\Supp\varsigma\setminus(\overline{\Sigma}_2)_\varepsilon$ i.e. $-u_{\Lambda,\Delta}$, which, by~\rmqref{onesta}, is coherent with the choice of the basis made in the proof of~\lemref{lemma-s2}. 
The situation is depicted in Figure~\ref{corrente} for a real $3$-dimensional cube.\footnote{Of course the picture refers to the intersection of $K_\Lambda\cap\Supp\varsigma\setminus(\overline{\Sigma}_2)_\varepsilon$ with the $3$-dimensional subspace $E_\Gamma=\mathbb R^3\subset\mathbb C^3$. In fact $\Supp\varsigma$ can be $6$-dimensional while $K_{\Lambda_1}$ and $K_{\Lambda_2}$ are $4$-dimensional, so there are three more dimensions that cannot be represented in the picture.}

As a last remark, observe that the form $\upsilon_{E_\Lambda^\prime}$ may vanish on $K_\Delta$, indeed this happens precisely when $E^\prime_\Lambda$ is orthogonal to $K_\Delta\subset E_\Delta^\perp $ and since the ortho-complement of $E_\Delta^\perp$ is $E_\Delta$, it follows that the pullback of $\upsilon_{E_\Lambda^\prime}$ to $K_\Delta$ is zero if and only if $E^\prime_\Lambda\cap E_\Delta\neq\{0\}$, which means that $ia_\Lambda=\pm u_{\Lambda,\Delta}$. We deduce that such a pullback is non-zero if and only if $\Delta$ is equidimensional.
Whence
\begin{align*}
\langle\!\langle (\ddc h_\Gamma)^{\wedge 2},\varsigma\rangle\!\rangle
&=
\lim_{\varepsilon\to0}
\sum_{\Lambda\in{\mathcal B}(\Gamma,1)}
\vol_1(\Lambda)
\sum_{\substack{\Delta\in{\mathcal B}(\Gamma,2)\\\Lambda\prec\Delta}}
\int_{X^\varepsilon_{\Lambda(\Delta)}}
\upsilon_{E_\Lambda^\prime}\wedge\dc h_\Gamma\wedge\varsigma
\\
&-
\lim_{\varepsilon\to0}
\sum_{\Lambda\in{\mathcal B}(\Gamma,1)}
\vol_1(\Lambda)
\sum_{\substack{\Delta\in{\mathcal B}(\Gamma,2)\\\Lambda\prec\Delta}}
\int_{X^\varepsilon_{\Lambda(\Delta)}}
\upsilon_{E_\Lambda^\prime}\wedge h_\Gamma\dc\varsigma
\\
&=
\sum_{\Delta\in{\mathcal B}_{\rm ed}(\Gamma,2)}
\sum_{\Lambda\in{\mathcal B}(\Delta,1)}
\vol_1(\Lambda)
\int_{\overline{K}_\Delta\cap\Supp\varsigma}
\upsilon_{E_\Lambda^\prime}\wedge\dc h_\Gamma\wedge\varsigma
\\
&-
\sum_{\Delta\in{\mathcal B}_{\rm ed}(\Gamma,2)}
\sum_{\Lambda\in{\mathcal B}(\Delta,1)}\vol_1(\Lambda)
\int_{\overline{K}_\Delta\cap\Supp\varsigma}
\upsilon_{E_\Lambda^\prime}\wedge h_\Gamma\dc\varsigma\,.
\end{align*}
Now, for every $\Delta\in{\mathcal B}_{\rm ed}(\Gamma,2)$, by \lemref{lemma-s2-1} and \lemref{lemma-s2}
\begin{align*}
&
\sum_{\Lambda\in{\mathcal B}(\Delta,1)}\vol_1(\Lambda)
\int_{K_\Delta\cap\Supp\varsigma}
\upsilon_{E_\Lambda^\prime}\wedge\dc h_\Gamma\wedge\varsigma
\\
&=
\int_{K_\Delta\cap\Supp\varsigma}
\left(
\sum_{\Lambda\in{\mathcal B}(\Delta,1)}\vol_1(\Lambda)
\upsilon_{E_\Lambda^\prime}\wedge\dc h_\Gamma
\right)
\wedge\varsigma
\\
&=
2\varrho(\Delta)\vol_2(\Delta)\int_{K_\Delta\cap\Supp\varsigma}
\upsilon_{E^\prime_\Delta}\wedge\varsigma
\\
&=
2\varrho(\Delta)
\vol_2(\Delta)
\langle\!\langle \lambda_\Delta,\varsigma\rangle\!\rangle\,,
\end{align*} 
whereas by \lemref{lemma-2}
\begin{align*}
&
\sum_{\Lambda\in{\mathcal B}(\Delta,1)}\vol_1(\Lambda)
\int_{K_\Delta\cap\Supp\varsigma}
\upsilon_{E_\Lambda^\prime}\wedge h_\Gamma\dc\varsigma\,
\\
&=
\int_{K_\Delta\cap\Supp\varsigma}
\left(
\sum_{\Lambda\in{\mathcal B}(\Delta,1)}\vol_1(\Lambda)
\upsilon_{E_\Lambda^\prime}\right)\wedge h_\Gamma\dc\varsigma
\\
&=
0\,.
\end{align*}
It follows that
\begin{equation}
(\ddc h_\Gamma)^{\wedge 2}
=
2\sum_{\Delta\in{\mathcal B}_{\rm ed}(\Gamma,2)}
\varrho(\Delta)
\vol_2(\Delta)
\lambda_\Delta
=
2\sum_{\Delta\in{\mathcal B}(\Gamma,2)}
\varrho(\Delta)
\vol_2(\Delta)
\lambda_\Delta\,
\end{equation}
and, consequently, $\Supp (\ddc h_\Gamma)^{\wedge 2}=\overline{\Sigma}_2$, 
which completes the proof.
\findim

\begin{coro}\label{val2}
Let $A_1,A_2\in{\mathcal K}(\mathbb C^n)$ be such that $A_1\cup A_2$ is convex. Then
\begin{equation}\label{val-2}
(\ddc h_{A_1})^{\wedge 2}+(\ddc h_{A_2})^{\wedge 2}=(\ddc h_{A_1\cup A_2})^{\wedge 2}+(\ddc h_{A_1\cap A_2})^{\wedge 2}\,.
\end{equation}
In particular
\begin{equation}\label{pseudo2}
\ddc h_{A_1}\wedge \ddc h_{A_2}=\ddc h_{A_1\cup A_2}\wedge\ddc h_{A_1\cap A_2}\,.
\end{equation}
\end{coro}
\noindent\pf 
Since the space of $4$-currents is Hausdorff, by virtue of \thmref{Sal}, it is enough to show that, for every hyperplane $H\subset\mathbb C^n$ and any polytope $\Gamma$, 
\begin{equation}\label{wval-2}
(\ddc h_{\Gamma\cap H})^{\wedge 2}+(\ddc h_{\Gamma})^{\wedge 2}=(\ddc h_{\Gamma\cap H^+})^{\wedge 2}+(\ddc h_{\Gamma\cap H^-})^{\wedge 2}\,.
\end{equation}
Remark that all the $2$-dimensional faces of $\Gamma\cap H$ are also $2$-faces of $\Gamma\cap H^+$ and $\Gamma\cap H^-$. For each such common $2$-face $\Delta$ one has a disjoint union $
K_{\Delta,\Gamma\cap H^+}\cup K_{\Delta,\Gamma\cap H^-}=K_{\Delta,\Gamma\cap H}$. Moreover, $\Gamma\cap H^+$ and $\Gamma\cap H^-$ may respectively have $2$-faces $\Delta_1$ and $\Delta_2$ such that $\Delta_1\cup\Delta_2$ is a $2$-face of $\Gamma$ not included in $H$. In this case, $\varrho({\Delta_1})=\varrho({\Delta_2})=\varrho({\Delta_1\cup\Delta_2})$, $\vol_2(\Delta_1)+\vol_2(\Delta_2)=\vol_2(\Delta_1\cup\Delta_2)$ and $K_{\Delta_1,\Gamma\cap H^+}= K_{\Delta_2,\Gamma\cap H^-}=K_{\Delta_1\cup\Delta_2,\Gamma}$. Any other $2$-face of $\Gamma\cap H^+$ and $\Gamma\cap H^-$ which does not intersect $H$ is also a $2$-face of $\Gamma$. The preceding analysis and \thmref{supp-2} imply \eqref{wval-2}. Moreover, by \eqref{val-1}
\begin{equation*}
(\ddc h_{A_1}+\ddc h_{A_2})^{\wedge 2}=(\ddc h_{A_1\cup A_2}+\ddc h_{A_1\cap A_2})^{\wedge 2}\,,
\end{equation*}
i.e.
\begin{align*}
&
(\ddc h_{A_1})^{\wedge 2}+2\,\ddc h_{A_1}\wedge \ddc h_{A_2}+(\ddc h_{A_2})^{\wedge 2}
\\
&=
(\ddc h_{A_1\cup A_2})^{\wedge 2}+2\,\ddc h_{A_1\cup A_2}\wedge\ddc h_{A_1\cap A_2}+(\ddc h_{A_1\cap A_2})^{\wedge 2}\,,
\end{align*}
which, thanks to \eqref{val-2}, yields \eqref{pseudo2}.
\findim

\subsection[The case $k>2$]{The case {\boldmath $k>2$}}

For $k>2$, the valuation property of the current $(\ddc h_A)^{\wedge k}$ can be proved by adapting the argument used in the case $k=2$, however the following corollary shows that the general case is indeed a consequence of the case $k=2$.

\begin{coro}\label{valk}
Let $A_1,A_2\in{\mathcal K}(\mathbb C^n)$ be such that $A_1\cup A_2$ is convex. Then
\begin{equation}\label{val-k}
(\ddc h_{A_1})^{\wedge k}+(\ddc h_{A_2})^{\wedge k}=(\ddc h_{A_1\cup A_2})^{\wedge k}+(\ddc h_{A_1\cap A_2})^{\wedge k}\,.
\end{equation}
\end{coro}
\noindent\pf 
We prove the statement by induction on $k$. For the sake of notation let us set $a=\ddc h_{A_1}$, $b=\ddc h_{A_2}$, $c=\ddc h_{A_1\cup A_2}$ and $d=\ddc h_{A_1\cap A_2}$. By \eqref{val-1} and \cororef{val2},  we know that $a+b=c+d$, $a^{\wedge 2}+b^{\wedge 2}=c^{\wedge 2}+d^{\wedge 2}$ and $a\wedge b=c\wedge d$ (whence $(a\wedge b)^{\wedge \ell}=(c\wedge d)^{\wedge \ell}$, for every $\ell\in\mathbb N$). By induction, suppose that, for any $1\leq s< k$, one has $a^{\wedge s}+b^{\wedge s}=c^{\wedge s}+d^{\wedge s}$. 
It then follows, for $1\leq \ell\leq k-1$, that
\begin{align*}
a^{\wedge(k-\ell)}\wedge b^{\wedge \ell}+a^{\wedge \ell}\wedge b^{\wedge (k-\ell)}
&=
\begin{cases}
(a\wedge b)^{\wedge \ell}\wedge \left(a^{\wedge (k-2\ell)}+b^{\wedge (k-2\ell)}\right) &\text{if } \ell< k-\ell \\
(a\wedge b)^{\wedge (k-\ell)}\wedge \left(a^{\wedge (2\ell-k)}+b^{\wedge (2\ell-k)}\right) &\text{if } \ell> k-\ell
\end{cases}\nonumber
\\
&=
\begin{cases}
(c\wedge d)^{\wedge \ell}\wedge \left(c^{\wedge (k-2\ell)}+d^{\wedge (k-2\ell)}\right) &\text{if } \ell< k-\ell \\
(c\wedge d)^{\wedge (k-\ell)}\wedge \left(c^{\wedge (2\ell-k)}+d^{\wedge (2\ell-k)}\right) &\text{if } \ell> k-\ell
\end{cases}\nonumber
\\
&=
c^{\wedge(k-\ell)}\wedge d^{\wedge \ell}+c^{\wedge \ell}\wedge d^{\wedge (k-\ell)}\,.\nonumber
\end{align*}
The equality $(a+b)^{\wedge k}=(c+d)^{\wedge k}$ and the preceding computation show that
\begin{align*}
&
a^{\wedge k}+b^{\wedge k}
\\
&=
(a+b)^{\wedge k}
-
\sum_{\ell=1}^{k-1}
{k\choose \ell} a^{\wedge (k-\ell)}\wedge b^{\wedge \ell}
\\
&=
(a+b)^{\wedge k}
-
\begin{cases}
\displaystyle{k\choose k/2} (a\wedge b)^{\wedge k/2}+
\sum_{\ell=1}^{(k/2)-1}
{k\choose \ell} \left(a^{\wedge (k-\ell)}\wedge b^{\wedge \ell} + a^{\wedge \ell}\wedge b^{\wedge (k-\ell)}\right),
&\text{if $k$ is even} \\
\displaystyle\sum_{\ell=1}^{\lfloor k/2\rfloor}
{k\choose \ell} \left(a^{\wedge (k-\ell)}\wedge b^{\wedge \ell} + a^{\wedge \ell}\wedge b^{\wedge (k-\ell)}\right),
&\text{if $k$ is odd} 
\end{cases}\nonumber
\\
&=
(c+d)^{\wedge k}
-
\begin{cases}
\displaystyle{k\choose k/2} (c\wedge d)^{\wedge k/2}+
\sum_{\ell=1}^{(k/2)-1}
{k\choose \ell} \left(c^{\wedge (k-\ell)}\wedge d^{\wedge \ell} + c^{\wedge \ell}\wedge d^{\wedge (k-\ell)}\right),
&\text{if $k$ is even} \\
\displaystyle\sum_{\ell=1}^{\lfloor k/2\rfloor}
{k\choose \ell} \left(c^{\wedge (k-\ell)}\wedge d^{\wedge \ell} + c^{\wedge \ell}\wedge d^{\wedge (k-\ell)}\right),
&\text{if $k$ is odd} 
\end{cases}\nonumber
\\
&=
c^{\wedge k}+d^{\wedge k}
\,.
\end{align*}
The proof is thus complete.\findim

\begin{coro}\label{valutazione-completa}
Let $m\in\{1,\ldots,n\}$ and, if $m<n$ let $A_{m+1},\ldots,A_n\in{\mathcal K}(\mathbb C^n)$ be fixed. Then the mapping $g_m:{\mathcal K}(\mathbb C^n)\to\mathbb R$ given, for any $A\in{\mathcal K}(\mathbb C^n)$, by
\begin{equation*}
g_m(A)=Q_n(A[m],A_{m+1},\ldots,A_n)
\,,
\end{equation*} 
is a continuous, translation invariant and unitarily invariant valuation. In particular $P_n$ is a continuous, translation invariant and unitarily invariant valuation.
\end{coro}
\noindent
\pf By \cororef{valk}, the mapping $A\mapsto (\ddc h_A)^{\wedge m}$ is a valuation, then so is the mapping $A\mapsto(\ddc h_A)^{\wedge m}\wedge\ddc h_{A_{m+1}}\wedge\ldots\wedge \ddc h_{A_n}$. It follows that $g_m$ is a valuation too. Continuity as well as translation and unitary invariance follow from the analogous properties of $Q_n$. In particular, by choosing $m=n$ one gets the statement for $P_n$.
\findim

\begin{lem}\label{utilem}
Let $\Gamma\in{\mathcal P}(\mathbb C^n)$, $\Delta\in{\mathcal B}(\Gamma,k)$ and $\Lambda\in{\mathcal B}_{\rm ed}(\Delta, k-1)$. Then $\Delta$ is equidimensional if and only if $E_\Lambda^\prime\cap E_\Delta=\{0\}$.
\end{lem}
\noindent\pf 
If $E_\Lambda^\prime\cap E_\Delta=\{0\}$, the outer normal vector $u_{\Lambda,\Delta}$, spanning $E_\Delta\cap E_\Lambda^\perp$, does not belong to $E_\Lambda^\prime$. Observe that $\dim_\mathbb R E_\Delta^\prime\neq k-1$ since otherwise $\lin_\mathbb C E_\Delta$ would have odd real dimension. Since $E_\Lambda^\prime\subset\lin_\mathbb C E_\Lambda\subseteq\lin_\mathbb C E_\Delta=E_\Delta\oplus E_\Delta^\prime=E_\Lambda\oplus (E_\Delta\cap E_\Lambda^\perp)\oplus E_\Delta^\prime$, it follows that $E_\Lambda^\prime\subseteq(E_\Delta\cap E_\Lambda^\perp)\oplus E_\Delta^\prime$, whence $\dim_\mathbb R E_\Delta^\prime= k-2$ or $\dim_\mathbb R E_\Delta^\prime= k$. If $\dim_\mathbb R E_\Delta^\prime=k-2$, then $E_\Lambda^\prime=(E_\Delta\cap E_\Lambda^\perp)\oplus E_\Delta^\prime$, but this means that $u_{\Lambda,\Delta}$ belongs to $E_\Lambda^\prime$, which is contrary to our assumption. It follows that $\dim_\mathbb R E_\Delta^\prime= k$, i.e. $\Delta$ is equidimensional. 
On the other hand, suppose $E_\Lambda^\prime\cap E_\Delta\neq\{0\}$. As $E_\Lambda^\prime\subset E_\Lambda^\perp$ and since $u_{\Lambda,\Delta}$ spans $E_\Lambda^\perp\cap E_\Delta$, it follows that $u_{\Lambda,\Delta}\in E_\Lambda^\prime\cap E_\Delta$. Let $v_1,\ldots,v_{k-1}$ an orthonormal basis of $E_\Lambda$, by \lemref{lemmat} the vectors $t_1,\ldots,t_{k-1}$ defined by $t_\ell=iv_\ell-\sum_{s=1}^{k-1}\re\langle i v_\ell,v_s\rangle v_s$, for $1\leq \ell\leq k-1$, provide a basis of $E_\Lambda^\prime$. Since $u_{\Lambda,\Delta}\in E_\Lambda^\prime$, we have $u_{\Lambda,\Delta}=\sum_{\ell=1}^{k-1} c_\ell t_\ell$, for some non trivial choice of coefficients $c_1,\ldots,c_{k-1}\in\mathbb R$. This representation yields the equality
 \begin{equation*}
 u_{\Lambda,\Delta}+\sum_{\ell=1}^{k-1} c_\ell \sum_{s=1}^{k-1}\re\langle i v_\ell,v_s\rangle v_s=\sum_{\ell=1}^{k-1}c_\ell i v_\ell\,,
 \end{equation*}
 where the vector on the left hand side belongs to $E_\Delta$ and that on the right hand side belongs to $iE_\Lambda\setminus\{0\}$. If $a$ denotes such a vector, then $ia\in E_\Lambda\cap iE_\Delta\subset E_\Delta\cap iE_\Delta$, whence $\varrho(\Delta)=0$.
\findim

\begin{lem}\label{lemma-k}
Let $\Delta\in{\mathcal P}(\mathbb C^n)$ be an equidimensional oriented $k$-polytope. Then
\begin{equation}\label{somma-lemma-k}
\iota_\Delta^*\left(\sum_{\Lambda\in{\mathcal B}(\Delta,k-1)}
\varrho(\Lambda)\vol_{k-1}(\Lambda)\upsilon_{E^\prime_\Lambda}\right)=0\,,
\end{equation}
where $\iota_\Delta:K_\Delta\to \mathbb C^n$ is the inclusion mapping and every $\Lambda\in{\mathcal B}(\Delta,k-1)$ is oriented by the unit outer normal vector $u_{\Lambda,\Delta}$.
\end{lem}
\noindent\pf 
Since $\Delta$ is equidimensional, any of its facets is such. As $\Delta$ is oriented each facet gets the  orientation corresponding to the outer unit normal vector. For every $\Lambda\in{\mathcal B}(\Delta,k-1)$, we have $E_\Delta^\perp=u_{\Lambda,\Delta}^\perp\cap E_\Lambda^\perp$ and by \eqref{forme}, on $E_\Lambda^\perp$, $\varrho(\Lambda)\upsilon_{E_\Lambda^\prime}=(-1)^{(k-1)(k-2)/2}\upsilon_{iE_\Lambda}$. 
Observe that
\begin{equation}\label{hodg}
\de\re\langle i u_{\Lambda,\Delta},z\rangle\wedge \upsilon_{iE_\Lambda}\wedge\upsilon_{iE_\Delta^\perp}=\upsilon_{2n}
\end{equation}
so we get the equality 
$
*\de\re\langle i u_{\Lambda,\Delta},z\rangle
=
\upsilon_{iE_\Lambda}\wedge\upsilon_{iE_\Delta^\perp}
$
of $(2n-1)$-forms on $\mathbb C^n$. 
In particular, if $b_{k+1},\ldots,b_{2n}$ is a fixed  positive basis of $E_\Delta^\perp$, consider the form 
$$
(-1)^{(k-1)(2n-k)}
\left[ib_{k+1}\,\ctr\ldots\ctr\, ib_{2n}\ctr *\de\re\langle i u_{\Lambda,\Delta},z\rangle\right]
$$
obtained, up to the sign, as a sequence of iterated interior products of the vectors $-ib_{2n},\ldots,-ib_{k+1}$ with $*\de\re\langle i u_{\Lambda,\Delta},z\rangle$. 
Then the relation (\ref{hodg}) yields the equalities of $(k-1)$-forms on $\mathbb C^n$ 
\begin{align} 
&
(-1)^{(k-1)(2n-k)}
\left[ib_{k+1}\ctr\ldots\ctr\, ib_{2n}\ctr *\de\re\langle i u_{\Lambda,\Delta},z\rangle\right](-,\ldots,-)
\\
&=
(-1)^{(k-1)(2n-k)}
*\de\re\langle i u_{\Lambda,\Delta},z\rangle(ib_{k+1},\ldots,ib_{2n},-,\ldots,-)
\\
&=
*\de\re\langle i u_{\Lambda,\Delta},z\rangle(-,\ldots,-,ib_{k+1},\ldots,ib_{2n})
\\
&=
\upsilon_{iE_\Lambda}\wedge\upsilon_{iE_\Delta^\perp}(-,\ldots,-,ib_{k+1},\ldots,ib_{2n})
\\
&=
\upsilon_{iE_\Lambda}(-,\ldots,-)
\end{align}
whence
\begin{equation}\label{contrazione}
\iota_\Delta^*\upsilon_{iE_\Lambda}
=
\iota_\Delta^*
\left[
(-1)^{(k-1)(2n-k)}
\left(
ib_{2k-1}\ctr\ldots\ctr\, ib_{2n}\ctr *\de\re\langle i u_{\Lambda,\Delta},z\rangle\right)\right].
\end{equation}
Remark, by the way, that the form $\iota^*_\Delta\upsilon_{iE_\Lambda}$ vanishes on $E^\perp_\Delta\cap iE_\Delta^\perp=E_\Delta^{\perp_\mathbb C}$ but it doesn't on $E_\Delta^\prime$.
By linearity, 
\begin{align*}
&\iota_\Delta^*\left(\sum_{\Lambda\in{\mathcal B}(\Delta,k-1)}
\varrho(\Lambda)\vol_{k-1}(\Lambda)\upsilon_{E^\prime_\Lambda}\right)
\\
&=
(-1)^{\frac{(k-1)(k-2)}{2}}
\sum_{\Lambda\in{\mathcal B}(\Delta,k-1)}
\vol_{k-1}(\Lambda)\iota_\Delta^*\upsilon_{iE_\Lambda}
\\
&=
(-1)^{\frac{(k-1)(k-2)}{2}}
\sum_{\Lambda\in{\mathcal B}(\Delta,k-1)}
\vol_{k-1}(\Lambda)\iota_\Delta^*
(-1)^{(k-1)(2n-k)}
\left[ib_{2k-1}\ctr\ldots\ctr\, ib_{2n}\ctr *\de\re\langle i u_{\Lambda,\Delta},z\rangle\right]
\\
&=
(-1)^{\frac{(k-1)(4n-k-2)}{2}}
\iota_\Delta^*
\left[ib_{2k-1}\ctr\ldots\ctr\, ib_{2n}\ctr *\left(\de\re\left\langle i\sum_{\Lambda\in{\mathcal B}(\Delta,k-1)}
\vol_{k-1}(\Lambda)u_{\Lambda,\Delta},z\right\rangle\right)
\right],
\end{align*}
now thanks to \eqref{piramidi} one obtains the desired result. 
Observe that changing the orientation of $\Delta$ changes the orientation of each of its facets, so that (\ref{contrazione}) becomes
$$\iota_\Delta^*\upsilon_{iE_\Lambda}
=
-\iota_\Delta^*
\left[
(-1)^{(k-1)(2k-1)}
ib_{2k-1}\ctr\ldots\ctr\, ib_{2n}\ctr *\de\re\langle i u_{\Lambda,\Delta},z\rangle\right],$$
however this change of sign does not affect the equality \eqref{somma-lemma-k}.\findim

\begin{lem}\label{lemma-sk-1}
Let $\Gamma\in {\mathcal P}(\mathbb C^n)$ and $\Delta\in{\mathcal B}_{\rm ed}(\Gamma,k)$. If $\Delta$ is oriented, then
\begin{align*}
&
\iota_\Delta^*\left(\sum_{\Lambda\in{\mathcal B}(\Delta,k-1)}\varrho(\Lambda)\vol_{k-1}(\Lambda)
\upsilon_{E_\Lambda^\prime}\wedge\dc h_\Gamma\right)
\\
&=
\iota_\Delta^*\left(\sum_{\Lambda\in{\mathcal B}(\Delta,k-1)}\varrho(\Lambda)\vol_{k-1}(\Lambda)
h_\Delta(u_{\Lambda,\Delta})\upsilon_{E_\Lambda^\prime}\wedge \re\langle i u_{\Lambda,\Delta},-\rangle\right)\,,
\end{align*} 
where $\iota_\Delta:K_\Delta\to \mathbb C^n$ is the inclusion and every $\Lambda\in{\mathcal B}(\Delta,k-1)$ is oriented by the unit outer normal vector $u_{\Lambda,\Delta}$.
\end{lem}
\noindent\pf
Let $\Lambda\in{\mathcal B}(\Delta,k-1)$ be fixed. Recall that $E_\Lambda^\perp\cap\aff_\mathbb R \Lambda=\{p_\Lambda\}$ and $E_\Delta^\perp\cap\aff_\mathbb R \Delta=\{p_\Delta\}$. Then, by virtue of (\ref{decomposition}), $p_\Lambda=p_\Delta+h_\Delta(u_{\Lambda,\Delta}) u_{\Lambda,\Delta}$ and 
\begin{align*}
\upsilon_{E_\Lambda^\prime}\wedge \dc h_\Gamma
&=
\upsilon_{E_\Lambda^\prime}\wedge\re\langle i p_\Lambda,-\rangle
\\
&=
\upsilon_{E_\Lambda^\prime}\wedge\re\langle i p_\Delta,-\rangle+h_\Delta(u_{\Lambda,\Delta})\upsilon_{E_\Lambda^\prime}\wedge \re\langle i u_{\Lambda,\Delta},-\rangle\,,
\end{align*}
By \lemref{lemma-k}, 
\begin{align*}
&
\iota_\Delta^*\left(\sum_{\Lambda\in{\mathcal B}(\Delta,k-1)}\varrho(\Lambda)\vol_{k-1}(\Lambda)
\upsilon_{E_\Lambda^\prime}\wedge\re\langle i p_\Delta,-\rangle\right)
\\
&=
\iota_\Delta^*\left(
-\re\langle i p_\Delta,-\rangle
\right)
\wedge
\iota_\Delta^*\left(
\sum_{\Lambda\in{\mathcal B}(\Delta,k-1)}\varrho(\Lambda)\vol_{k-1}(\Lambda)\upsilon_{E_\Lambda^\prime}\right)
=
0
\end{align*}
whence
\begin{align*}
&
\iota_\Delta^*\left(\sum_{\Lambda\in{\mathcal B}(\Delta,k-1)}\varrho(\Lambda)\vol_{k-1}(\Lambda)
\upsilon_{E_\Lambda^\prime}\wedge\iota_\Lambda^* \dc h_\Gamma\right)
\\
&=
\iota_\Delta^*\left(\sum_{\Lambda\in{\mathcal B}(\Delta,k-1)}\varrho(\Lambda)\vol_{k-1}(\Lambda)
h_\Delta(u_{\Lambda,\Delta})\upsilon_{E_\Lambda^\prime}\wedge \re\langle i u_{\Lambda,\Delta},-\rangle\right)
\end{align*}
which proves the lemma.\findim

\begin{lem}\label{lemma-sk}
Let $\Gamma\in {\mathcal P}(\mathbb C^n)$ and $\Delta\in{\mathcal B}_{\rm ed}(\Gamma,k)$. If $\Delta$ is oriented and $\Lambda\in{\mathcal B}(\Delta,k-1)$ is oriented by the unit outer normal vector $u_{\Lambda,\Delta}$, then
\begin{align*}
&
\iota_\Delta^*\left(
\varrho(\Lambda)\vol_{k-1}(\Lambda)
h_\Delta(u_{\Lambda,\Delta})\upsilon_{E_\Lambda^\prime}\wedge \re\langle i u_{\Lambda,\Delta},-\rangle\right)
\\
&=
\varrho(\Delta)
\vol_{k-1}(\Lambda)h_\Delta(u_{\Lambda,\Delta})
\upsilon_{E^\prime_\Delta}
\,.
\end{align*} 
In particular 
\begin{align*}
\iota_\Delta^*\left(
\sum_{\Lambda\in{\mathcal B}(\Delta,k-1)}
\varrho(\Lambda)\vol_{k-1}(\Lambda)
h_\Delta(u_{\Lambda,\Delta})\upsilon_{E_\Lambda^\prime}\wedge \re\langle i u_{\Lambda,\Delta},-\rangle
\right)
=
k\varrho(\Delta)
\vol_k(\Delta)
\upsilon_{E^\prime_\Delta}\,.
\end{align*}
\end{lem}
\noindent\pf 
By \eqref{forme}, on $E_\Delta^\prime\subset E_\Delta^\perp\subset E_\Lambda^\perp$ we have 
$
\iota_\Delta^*\left(
\varrho(\Lambda)
\upsilon_{E_\Lambda^\prime}
\right)
=
(-1)^{(k-1)(k-2)/2}
\iota_\Delta^*\left(
\upsilon_{iE_\Lambda}
\right)
$, then 
\begin{align*}
\iota_\Delta^*\left(
\varrho(\Lambda)
\upsilon_{E_\Lambda^\prime}\wedge \re\langle i u_{\Lambda,\Delta},-\rangle
\right)
&=
(-1)^{(k-1)(k-2)/2}
\iota_\Delta^*\left(
\upsilon_{iE_\Lambda}\wedge \re\langle i u_{\Lambda,\Delta},-\rangle
\right)
\\
&=
(-1)^{(k-1)(k-2)/2}(-1)^{(k-1)}
\iota_\Delta^*\left(
\re\langle i u_{\Lambda,\Delta},-\rangle\wedge \upsilon_{iE_\Lambda}
\right)
\\
&=
(-1)^{k(k-1)/2}
\iota_\Delta^*\left(
v_{i E_\Delta}
\right)
\\
&=
\varrho(\Delta)v_{E_\Delta^\prime}
\,.
\end{align*}

As $\Delta$ is oriented, each of its facets gets the orientation induced by the outer unit normal vector and, by \eqref{vol-piramidi}, it follows that 
\begin{align*}
&
\iota_\Delta^*\left(
\sum_{\Lambda\in{\mathcal B}(\Delta,k-1)}
\varrho(\Lambda)\vol_{k-1}(\Lambda)
h_\Delta(u_{\Lambda,\Delta})\upsilon_{E_\Lambda^\prime}\wedge \re\langle i u_{\Lambda,\Delta},-\rangle
\right)
\\
&=
k\varrho(\Delta)
\left(
\sum_{\Lambda\in{\mathcal B}(\Delta,k-1)}
\frac{\vol_{k-1}(\Lambda)
h_\Delta(u_{\Lambda,\Delta})}{k}
\right)
v_{E_\Delta^\prime}
\\
&=
k\varrho(\Delta)\vol_k(\Delta)
v_{E_\Delta^\prime}\,.
\end{align*}
The proof is thus complete.\findim

\begin{thm}\label{supp-k}
Let $\Gamma\in {\mathcal P}(\mathbb C^n)$ and $1< k\leq n$, then 
\begin{equation}\label{current-k}
(\ddc h_\Gamma)^{\wedge k}
=
k!\sum_{\Delta\in{\mathcal B}(\Gamma,k)}
\varrho(\Delta)
\vol_k(\Delta)
\lambda_\Delta
\,.
\end{equation}
In particular $\Supp (\ddc h_\Gamma)^{\wedge k}=\overline{\Sigma}_k$.
\end{thm}
\noindent\pf 
We prove the theorem by induction on $k$. The case $k=1$ and $k=2$ have been already proved by \thmref{supp-1} and \thmref{supp-2}, respectively. By induction, suppose that 
\begin{equation*}
(\ddc h_\Gamma)^{\wedge (k-1)}
=
(k-1)!\sum_{\Lambda\in{\mathcal B}(\Gamma,k-1)}
\varrho(\Lambda)
\vol_{k-1}(\Lambda)
\lambda_\Lambda
\,.
\end{equation*}
Then, for any $(n-k,n-k)$-test form $\varsigma$,
\begin{align*} 
\langle\!\langle (\ddc h_\Gamma)^{\wedge k},\varsigma\rangle\!\rangle
&=
\langle\!\langle \ddc\big(h_\Gamma(\ddc h_\Gamma)^{\wedge(k-1)}\big),\varsigma\rangle\!\rangle\\
&=
\langle\!\langle h_\Gamma(\ddc h_\Gamma)^{\wedge(k-1)},\ddc\varsigma\rangle\!\rangle\\
&=
(k-1)!\sum_{\Lambda\in{\mathcal B}_{\rm ed}(\Gamma,k-1)}\varrho(\Lambda)\vol_{k-1}(\Lambda)\langle\!\langle \lambda_{\Lambda},h_\Gamma\ddc\varsigma\rangle\!\rangle\\
&=
(k-1)!\sum_{\Lambda\in{\mathcal B}_{\rm ed}(\Gamma,k-1)}\varrho(\Lambda)\vol_{k-1}(\Lambda)\int_{K_\Lambda} \upsilon_{E_\Lambda^\prime}\wedge h_\Gamma\ddc \varsigma\\
&=
(k-1)!
\lim_{\varepsilon\to0}
\sum_{\Delta\in{\mathcal B}_{\rm ed}(\Gamma,k-1)}\varrho(\Lambda)\vol_{k-1}(\Lambda)\int_{K_\Lambda\cap\Supp\varsigma\setminus(\overline{\Sigma}_k)_\varepsilon} \upsilon_{E_\Lambda^\prime}\wedge h_\Gamma\ddc \varsigma
\\
&=
(k-1)!
\lim_{\varepsilon\to0}
\sum_{\Delta\in{\mathcal B}_{\rm ed}(\Gamma,k-1)}\varrho(\Lambda)\vol_{k-1}(\Lambda)\int_{S_{\Lambda,\varepsilon}} \upsilon_{E_\Lambda^\prime}\wedge h_\Gamma\ddc \varsigma\,,
\end{align*}
where $S_{\Lambda,\varepsilon}=(K_\Lambda\cap(K_\Lambda\cap\Supp\varsigma\setminus(\overline{\Sigma}_k)_\varepsilon)_{\varepsilon/2}$, i.e. the $(\varepsilon/2)$-neighbourhood of $K_\Lambda\cap\Supp\varsigma\setminus(\overline{\Sigma}_k)_\varepsilon$.
Since $\upsilon_{E_\Lambda^\prime}$ is a closed $(k-1)$-form on $E_\Lambda^\prime\subset E_\Lambda^\perp$ and as the $(2n-2k+2)$-forms $-\de h_\Gamma\wedge\dc\varsigma$ and $\dc h_\Gamma\wedge \de \varsigma$ have the same $(n-k+1,n-k+1)$-parts, it follows that on the manifold $S_{\Lambda,\varepsilon}$ 
\begin{align*}
\de(v_{E^\prime_\Lambda}\wedge h_\Gamma\dc\varsigma)
&=
\de \upsilon_{E_\Lambda^\prime}\wedge h_\Gamma\dc\varsigma  + (-1)^{k-1}\upsilon_{E_\Lambda^\prime}\wedge \de(h_\Gamma\dc\varsigma)\\
&=
0+(-1)^{k-1}\upsilon_{E_\Lambda^\prime}\wedge \de h_\Gamma\wedge\dc\varsigma+(-1)^{k-1}\upsilon_{E_\Lambda^\prime}\wedge h_\Gamma\ddc\varsigma
\\
&=
(-1)^k \upsilon_{E_\Lambda^\prime}\wedge \dc h_\Gamma\wedge\de\varsigma+(-1)^{k-1}
\upsilon_{E_\Lambda^\prime}\wedge h_\Gamma\ddc\varsigma\,.
\end{align*}
Observe that, on $K_\Lambda$, the $1$-form $\dc h_\Gamma$ has constant coefficients, so the $k$-form $\upsilon_{E_\Lambda^\prime}\wedge \dc h_\Gamma$ is closed on  $K_\Lambda$ and $(-1)^k \upsilon_{E_\Lambda^\prime}\wedge \dc h_\Gamma\wedge\de\varsigma=\de (\upsilon_{E_\Lambda^\prime}\wedge\dc h_\Gamma\wedge\varsigma)$. It follows that
\begin{equation*}
(-1)^{k}\upsilon_{E_\Lambda^\prime}\wedge h_\Gamma\ddc \varsigma
=
\de (\upsilon_{E_\Lambda^\prime}\wedge\dc h_\Gamma\wedge\varsigma)-\de (\upsilon_{E_\Lambda^\prime}\wedge h_\Gamma\dc\varsigma)\,,
\end{equation*}
whence 
\begin{align}
&
\langle\!\langle (\ddc h_\Gamma)^{\wedge k},\varsigma\rangle\!\rangle\nonumber
\\
&=
(k-1)!\lim_{\varepsilon\to0}
\sum_{\Lambda\in{\mathcal B}_{\rm ed}(\Gamma,k-1)}
(-1)^{k}\varrho(\Lambda)\vol_{k-1}(\Lambda)\int_{\relbd S_{\Lambda,\varepsilon}}
\upsilon_{E_\Lambda^\prime}\wedge\dc h_\Gamma\wedge\varsigma\label{ppt}
\\
&-
(k-1)!\lim_{\varepsilon\to0}
\sum_{\Lambda\in{\mathcal B}_{\rm ed}(\Gamma,k-1)}(-1)^{k}\varrho(\Lambda)\vol_{k-1}(\Lambda)\int_{\relbd S_{\Lambda,\varepsilon}}
\upsilon_{E_\Lambda^\prime}\wedge h_\Gamma\dc\varsigma\,.\label{sst}
\end{align}
Let us first consider the term \eqref{sst}. For any fixed equidimensional $(k-1)$-face $\Lambda$ such that $\Supp\varsigma\cap \relbd K_\Lambda\neq\varnothing$, the only parts of $\relbd S_{\Lambda,\varepsilon}$
which really count are the flat parts parallel and near to the $(2n-k)$-dimensional pieces of $\relbd K_\Lambda$ on which $h_\Gamma\dc\varsigma$ and $\upsilon_{E_\Lambda^\prime}$ are linearly independent. Such $(2n-k)$-dimensional pieces of $\relbd K_\Lambda$ are nothing but the cones which are dual to the $k$-faces $\Delta$ of $\Gamma$ admitting $\Lambda$ as a facet. Let us denote $X^\varepsilon_{\Lambda(\Delta)}$ the subset of $\relbd S_{\Lambda,\varepsilon}$ parallel and near to $K_\Delta$ and observe, by the way, that $X^\varepsilon_{\Lambda(\Delta)}$ approaches $K_\Delta\cap \Supp\varsigma$ as $\varepsilon\to0$. The form $\upsilon_{E_\Lambda^\prime}$ may vanish on $K_\Delta$, indeed this happens precisely when $E^\prime_\Lambda$ admits a non-zero vector that is orthogonal to $K_\Delta\subset E_\Delta^\perp $. Since the ortho-complement of $E_\Delta^\perp$ is $E_\Delta$, it follows that the pullback of $\upsilon_{E_\Lambda^\prime}$ to $K_\Delta$ is zero if and only if $E^\prime_\Lambda\cap  E_\Delta\neq\{0\}$, which, by \lemref{utilem}, means that $\Delta$ is not equidimensional.

It follows that \eqref{sst} becomes
\begin{equation}
(-1)^{k-1}(k-1)!\lim_{\varepsilon\to0}
\sum_{\Lambda\in{\mathcal B}_{\rm ed}(\Gamma,k-1)}\varrho(\Lambda)\vol_{k-1}(\Lambda)\sum_{\substack{\Delta\in{\mathcal B}_{\rm ed}(\Gamma,k)\\ \Lambda\prec\Delta}}\int_{X^\varepsilon_{\Lambda(\Delta)}}
\upsilon_{E_\Lambda^\prime}\wedge h_\Gamma\dc\varsigma\,.\label{ttt}
\end{equation}

Let us now fix $\Delta\in{\mathcal B}_{\rm ed}(\Gamma,k)$ and give it an orientation. This provides each of its facets $\Lambda\in{\mathcal B}(\Delta,k-1)$ the orientation induced by the outer unit normal vector  $u_{\Lambda,\Delta}\in E_\Lambda^\perp\cap E_\Delta$. In particular, if $v_{\Lambda,1},\ldots,v_{\Lambda,k-1}$ spans $E_\Lambda$ and
\begin{equation*}
u_{\Lambda,\Delta},v_{\Lambda,1},\ldots,v_{\Lambda,k-1},(-1)^{(k-1)k/2} w_{\Delta,1},\ldots,w_{\Delta,k}
\end{equation*}
is the positive basis of $\lin_\mathbb C E_\Delta$ defined in~\rmqref{orient}, up to completing it to a positive basis 
\begin{equation*}
u_{\Lambda,\Delta},v_{\Lambda,1},\ldots,v_{\Lambda,k-1},(-1)^{(k-1)k/2} w_{\Delta,1},\ldots,w_{\Delta,k},b_{2k+1},\ldots,b_{2n}
\end{equation*}
of $\mathbb C^n$, it is clear that the basis
\begin{equation*}
v_{\Lambda,1},\ldots,v_{\Lambda,k-1},(-1)^{k-1}u_{\Lambda,\Delta},(-1)^{(k-1)k/2} w_{\Delta,1},\ldots,w_{\Delta,k},,b_{2k+1},\ldots,b_{2n}
\end{equation*}
is positive too and the subspace $E_\Lambda^\perp$ is oriented by the vector $(-1)^{k-1}u_{\Lambda,\Delta}$. On the other hand $X^\varepsilon_{\Lambda(\Delta)}$ has the orientation coming from $\relbd S_{\Lambda,\varepsilon}$, this orientation corresponds to the unit normal vector to $X^\varepsilon_{\Lambda(\Delta)}$ pointing outside $K_\Lambda\cap\Supp\varsigma\setminus(\overline{\Sigma}_k)_\varepsilon$, i.e. $-u_{\Lambda,\Delta}$, 
so that the orientation of $X^\varepsilon_{\Lambda(\Delta)}$ is $(-1)^k$ times that of $K_\Delta$ as a subset of $K_\Lambda$. Thus the expression~\eqref{ttt} becomes
\begin{align*}
&
(-1)(k-1)!\sum_{\Delta\in{\mathcal B}_{\rm ed}(\Gamma,k)}
\sum_{\Lambda\in{\mathcal B}_{\rm ed}(\Delta,k-1)}\varrho(\Lambda)\vol_{k-1}(\Lambda)
\int_{K_\Delta\cap\Supp\varsigma}
\upsilon_{E_\Lambda^\prime}\wedge h_\Gamma\dc\varsigma
\\
&=
(-1)(k-1)!
\sum_{\Delta\in{\mathcal B}_{\rm ed}(\Gamma,k)}
\int_{K_\Delta\cap\Supp\varsigma}
\left(\sum_{\Lambda\in{\mathcal B}_{\rm ed}(\Delta,k-1)}\varrho(\Lambda)\vol_{k-1}(\Lambda)
\upsilon_{E_\Lambda^\prime}\right)\wedge h_\Gamma\dc\varsigma\,.
\end{align*}
By \lemref{lemma-k}, the latter term is zero. 

We now turn our attention to the term \eqref{ppt}. Again, for any fixed equidimensional $(k-1)$-face $\Lambda$, if $\Supp\varsigma\cap \relbd K_\Lambda\neq\varnothing$, the only parts of $\relbd S_{\Lambda,\varepsilon}$ on which $\iota_\Lambda^*(\upsilon_{E_\Lambda^\prime}\wedge \dc h_\Gamma\wedge\varsigma)$ gives a non-zero contribution are the flat pieces on which $\iota_\Lambda^*\varsigma\neq0$ and $\iota_\Lambda^*(\upsilon_{E_\Lambda^\prime}\wedge\dc h_\Gamma)$ is a non-zero $k$-form. The first condition produces the same conclusions as in the case of the term \eqref{sst}, i.e. the integral in \eqref{ppt} has to be performed on the subsets of the form $X^\varepsilon_{\Lambda(\Delta)}$. The second condition requires that $\Delta$ has to be an equidimensional $k$-face and $\iota_\Lambda^*(\upsilon_{E_\Lambda^\prime}\wedge\dc h_\Gamma)$ a volume form on $E_\Delta^\prime$. So 
\begin{align*}
&
(k-1)!\lim_{\varepsilon\to0}
\sum_{\Lambda\in{\mathcal B}_{\rm ed}(\Gamma,k-1)}(-1)^{k}\varrho(\Lambda)\vol_{k-1}(\Lambda)\sum_{\substack{\Delta\in{\mathcal B}_{\rm ed}(\Gamma,k)\\ \Lambda\prec\Delta}}\int_{X^\varepsilon_{\Lambda(\Delta)}}
\upsilon_{E_\Lambda^\prime}\wedge \dc h_\Gamma\wedge\varsigma\\
&=
(k-1)!\sum_{\Delta\in{\mathcal B}_{\rm ed}(\Gamma,k)}
\sum_{\Lambda\in{\mathcal B}_{\rm ed}(\Delta,k-1)}\varrho(\Lambda)\vol_{k-1}(\Lambda)
\int_{K_\Delta\cap\Supp\varsigma}
\upsilon_{E_\Lambda^\prime}\wedge \dc h_\Gamma\wedge\varsigma\\
&=
(k-1)!\sum_{\Delta\in{\mathcal B}_{\rm ed}(\Gamma,k)}
\int_{K_\Delta\cap\Supp\varsigma}
\left(\sum_{\Lambda\in{\mathcal B}(\Delta,k-1)}\varrho(\Lambda)\vol_{k-1}(\Lambda)
\upsilon_{E_\Lambda^\prime}\wedge \dc h_\Gamma\right)\wedge\varsigma\,.\label{qt}
\end{align*}
By \lemref{lemma-sk}, $\iota_\Delta^*\sum_{\Lambda\in{\mathcal B}(\Delta,k-1)}\varrho(\Lambda)\vol_{k-1}(\Lambda)
\upsilon_{E_\Lambda^\prime}\wedge \dc h_\Gamma=k\varrho(\Delta)\vol_k(\Delta)v_{E_\Delta^\prime}$,
hence
\begin{equation*}
(\ddc h_\Gamma)^{\wedge k}
=
k!\sum_{\Delta\in{\mathcal B}_{\rm ed}(\Gamma,k)}
\varrho(\Delta)
\vol_k(\Delta)
\lambda_\Delta
=
k!\sum_{\Delta\in{\mathcal B}(\Gamma,k)}
\varrho(\Delta)
\vol_k(\Delta)
\lambda_\Delta\,.
\end{equation*}
Of course, $\Supp (\ddc h_\Gamma)^{\wedge k}=\overline{\Sigma}_k$ and the proof is thus complete.
\findim

\begin{coro}\label{combinatorio}
Let $1\leq k\leq n$, $\Gamma_1,\ldots,\Gamma_k\in{\mathcal P}(\mathbb C^n)$ and $\Gamma=\sum_{\ell=1}^k\Gamma_\ell$. Then
\begin{equation}
\ddc h_{\Gamma_1}\wedge\ldots\wedge\ddc h_{\Gamma_k}=k!\sum_{\Delta\in{\mathcal B}(\Gamma,k)}
\varrho(\Delta)V_k(\Delta_1,\ldots,\Delta_k)\lambda_\Delta\,,
\end{equation}
where, for every $\Delta\in{\mathcal B}(\Gamma,k)$ and any $1\leq\ell\leq k$, $\Delta_\ell\preccurlyeq\Gamma_\ell$ is such that $\Delta=\sum_{\ell=1}^k\Delta_\ell$. In particular $\Supp \ddc h_{\Gamma_1}\wedge\ldots\wedge\ddc h_{\Gamma_k}=\overline\Sigma_{k,\Gamma}$.
\end{coro}
\noindent\pf
Let $t_1,\ldots,t_k>0$. By applying \thmref{supp-k} to the polytope $\widetilde\Gamma=\sum_{\ell=1}^k t_\ell\Gamma_\ell$ one gets
\begin{equation*}
\left(\sum_{\ell=1}^k t_\ell \ddc h_{\Gamma_\ell}\right)^{\wedge k}
=
k!\sum_{\widetilde\Delta\in{\mathcal B}(\widetilde\Gamma,k)}
\varrho(\widetilde\Delta)\vol_k\left(\sum_{\ell=1}^kt_\ell\Delta_\ell\right)\lambda_{\widetilde\Delta}\,,
\end{equation*}
where $\widetilde \Delta=\sum_{\ell=1}^kt_\ell\Delta_\ell$ is the representation of $\widetilde\Delta$ as a sum of faces $\Delta_\ell\preccurlyeq\Gamma_\ell$, for $1\leq\ell\leq k$. Set $\Delta=\sum_{\ell=1}^k\Delta_\ell$, with the same face summands as $\widetilde \Delta$. The mapping $\widetilde\Delta\mapsto\Delta$ is a bijection from ${\mathcal B}(\widetilde\Gamma,k)$ to ${\mathcal B}(\Gamma,k)$ such that $E_{\widetilde\Delta}=E_\Delta$, whence $E_{\widetilde\Delta}^\prime=E_\Delta^\prime$, $\varrho(\widetilde\Delta)=\varrho(\Delta)$ and $\upsilon_{E_{\widetilde\Delta}^\prime}=\upsilon_{E_\Delta^\prime}$. By comparing the coefficients of the monomial $t_1\cdots t_k$ on both sides one obtains
\begin{equation*}
k!\bigwedge_{\ell=1}^k \ddc h_{\Gamma_\ell}
=
(k!)^2
\sum_{\Delta\in{\mathcal B}(\Gamma,k)}
\varrho(\Delta)V_k\left(\Delta_1,\ldots\Delta_k\right)\lambda_{\Delta}\,,
\end{equation*}
whence
\begin{equation*}
\ddc h_{\Gamma_1}\wedge\ldots\wedge\ddc h_{\Gamma_k}
=
k!
\sum_{\Delta\in{\mathcal B}(\Gamma,k)}
\varrho(\Delta)V_k\left(\Delta_1,\ldots\Delta_k\right)\lambda_{\Delta}\,.
\end{equation*}
The statement about the support of $\ddc h_{\Gamma_1}\wedge\ldots\wedge\ddc h_{\Gamma_k}$ is a consequence of \thmref{supp-k}.\findim

For a fixed polytope $\Gamma\in{\mathcal P}(\mathbb C^n)$ and a fixed integer $1\leqslant k\leqslant n$, the \emph{trace measure}  
of the current $(\ddc h_\Gamma)^{\wedge k}$ is the current (of  measure type) given by 
\begin{align}\label{trace-measure}
\dfrac{1}{4^{n-k}(n-k)!}(\ddc h_\Gamma)^{\wedge k}\wedge (\ddc h_{B_{2n}}^2)^{\wedge (n-k)}.
\end{align}
More generally, if $\Gamma_1,\ldots,\Gamma_k\in{\mathcal P}(\mathbb C^n)$ and $\Gamma=\sum_{\ell=1}^k \Gamma_\ell$, the trace measure of $\ddc h_{\Gamma_1}\wedge\ldots\wedge\ddc h_{\Gamma_k}$ is given by
$$
\dfrac{1}{4^{n-k}(n-k)!}
\left(
\bigwedge_{\ell=1}^k
\ddc h_{\Gamma_\ell}
\right)
\wedge
(\ddc h^2_{B_{2n}})^{\wedge (n-k)}.
$$
By~\cororef{combinatorio}, the trace measure of $\ddc h_{\Gamma_1}\wedge\ldots\wedge\ddc h_{\Gamma_k}$ has the following simple expression:
\begin{align*}
&
\dfrac{1}{4^{n-k}(n-k)!}
\left(
\bigwedge_{\ell=1}^k
\ddc h_{\Gamma_\ell}
\right)
\wedge
(\ddc h^2_{B_{2n}})^{\wedge (n-k)}
\\
&=
\dfrac{k!}{4^{n-k}(n-k)!}
\sum_{\Delta\in{\mathcal B}_{\rm ed}(\Gamma,k)}
\varrho(\Delta)V_k(\Delta_1,\ldots,\Delta_k)
\lambda_\Delta\wedge(\ddc h_{B_{2n}}^2)^{\wedge (n-k)}
\\
&=
k!
\sum_{\Delta\in{\mathcal B}_{\rm ed}(\Gamma,k)}
\varrho(\Delta)V_k(\Delta_1,\ldots,\Delta_k)
[K_\Delta]
\wedge
\upsilon_{E_\Delta^\prime}
\wedge
\upsilon_{E_\Delta^{\perp_\mathbb C}}
\\
&=
k!
\sum_{\Delta\in{\mathcal B}_{\rm ed}(\Gamma,k)}
\varrho(\Delta)V_k(\Delta_1,\ldots,\Delta_k)
[K_\Delta]
\wedge
\upsilon_{E_\Delta^\perp}.
\end{align*}

\begin{lem}\label{stranissimo}
Let $0\leqslant k\leqslant n$ and $0\leqslant \ell\leqslant n-k$ be integers. For every $A\in{\mathcal K}(\mathbb C^n)$
\begin{align}\label{stranissima}
\int_{B_{2n}}
(\ddc h_A)^{\wedge k}
\wedge
(\ddc h^2_{B_{2n}})^{\wedge (n-k)}
=
2^\ell
\int_{B_{2n}}
(\ddc h_A)^{\wedge k}
\wedge
(\ddc h_{B_{2n}})^{\wedge \ell}
\wedge
(\ddc h^2_{B_{2n}})^{\wedge (n-k-\ell)}
\end{align}
\end{lem}

\noindent
\pf
Let $0\leqslant k\leqslant n$ be fixed and let us show the claim by induction on $\ell$. If $\ell=0$ the statement is trivial. Now suppose that (\ref{stranissima}) holds true for some $0< \ell<n-k$. Let $(\varsigma_m)$ be a sequence 
of radial smooth cut-off function equal to $1$ on $B_{2n}$ and decresing to $\chi_{B_{2n}}$, let also $\varphi_\varepsilon$ be a regularizing kernel such that $h_A*\varphi_\varepsilon$ decreses uniformly to $h_A$ on compacta. Then
\begin{align*}
&\int_{B_{2n}}
(\ddc h_A)^{\wedge k}
\wedge
(\ddc h^2_{B_{2n}})^{\wedge (n-k)}
\\
&=
2^\ell
\int_{B_{2n}}
(\ddc h_A)^{\wedge k}
\wedge
(\ddc h_{B_{2n}})^{\wedge \ell}
\wedge
(\ddc h^2_{B_{2n}})^{\wedge (n-k-\ell)}
\\
&=
2^\ell
\int_{B_{2n}}
(\ddc h_A)^{\wedge k}\wedge(\ddc h_{B_{2n}})^{\wedge \ell}\wedge(\ddc h_{B_{2n}}^2)\wedge(\ddc h_{B_{2n}}^2)^{\wedge(n-k-\ell-1)}
\\
&=
2^\ell
\lim_{m\to\infty}
\langle\!\langle 
(\ddc h_A)^{\wedge k}\wedge(\ddc h_{B_{2n}})^{\wedge \ell}\wedge(\ddc h_{B_{2n}}^2)\wedge(\ddc h_{B_{2n}}^2)^{\wedge(n-k-\ell-1)},
\varsigma_m
\rangle\!\rangle
\\
&=
2^\ell
\lim_{m\to\infty}
\lim_{\varepsilon\to0}
\int_{\Supp\varsigma_m}
[\ddc (h_A*\varphi_\varepsilon)]^{\wedge k}\wedge(\ddc h_{B_{2n}})^{\wedge \ell}\wedge(\ddc h_{B_{2n}}^2)\wedge(\ddc h_{B_{2n}}^2)^{\wedge(n-k-\ell-1)}
\varsigma_m
\\
&=
2^\ell
\lim_{\varepsilon\to0}
\lim_{m\to\infty}
\int_{\Supp\varsigma_m}
[\ddc (h_A*\varphi_\varepsilon)]^{\wedge k}\wedge(\ddc h_{B_{2n}})^{\wedge \ell}\wedge(\ddc h_{B_{2n}}^2)\wedge(\ddc h_{B_{2n}}^2)^{\wedge(n-k-\ell-1)}\varsigma_m
\\
&=
2^\ell
\lim_{\varepsilon\to0}
\int_{B_{2n}}
[\ddc (h_A*\varphi_\varepsilon)]^{\wedge k}\wedge(\ddc h_{B_{2n}})^{\wedge \ell}\wedge(\ddc h_{B_{2n}}^2)\wedge(\ddc h_{B_{2n}}^2)^{\wedge(n-k-\ell-1)}\\
&=
2^\ell
\lim_{\varepsilon\to0}
\lim_{\delta\to 0}
\int_{B_{2n}\setminus\delta B_{2n}}
[\ddc (h_A*\varphi_\varepsilon)]^{\wedge k}\wedge(\ddc h_{B_{2n}})^{\wedge \ell}\wedge(\ddc h_{B_{2n}}^2)\wedge(\ddc h_{B_{2n}}^2)^{\wedge(n-k-\ell-1)}\\
&=
2^\ell
\lim_{\varepsilon\to0}
\int_{\partial B_{2n}}
[\ddc (h_A*\varphi_\varepsilon)]^{\wedge k}\wedge(\ddc h_{B_{2n}})^{\wedge \ell}\wedge(\dc h_{B_{2n}}^2)\wedge(\ddc h_{B_{2n}}^2)^{\wedge(n-k-\ell-1)}
\\
&-
2^\ell
\lim_{\varepsilon\to0}
\lim_{\delta\to 0}
\int_{\partial (\delta B_{2n})}
[\ddc (h_A*\varphi_\varepsilon)]^{\wedge k}\wedge(\ddc h_{B_{2n}})^{\wedge \ell}\wedge(\dc h_{B_{2n}}^2)\wedge(\ddc h_{B_{2n}}^2)^{\wedge(n-k-\ell-1)}
\\
&=
2^\ell
\lim_{\varepsilon\to0}
\int_{\partial B_{2n}}
[\ddc (h_A*\varphi_\varepsilon)]^{\wedge k}\wedge(\ddc h_{B_{2n}})^{\wedge \ell}\wedge(2h_{B_{2n}}\dc h_{B_{2n}})\wedge(\ddc h_{B_{2n}}^2)^{\wedge(n-k-\ell-1)}
\\
&-
2^\ell
\lim_{\varepsilon\to0}
\lim_{\delta\to 0}
\int_{\partial (\delta B_{2n})}
[\ddc (h_A*\varphi_\varepsilon)]^{\wedge k}\wedge(\ddc h_{B_{2n}})^{\wedge \ell}\wedge(2h_{B_{2n}}\dc h_{B_{2n}})\wedge(\ddc h_{B_{2n}}^2)^{\wedge(n-k-\ell-1)}
\end{align*}
Now $h_{B_{2n}}(z)= 1$ for every $z\in\partial B_{2n}$, $h_{B_{2n}}(z)= \delta$ for every $z\in\partial \delta B_{2n}$ and 
\begin{equation}\label{integrand}
\lim_{\delta\to 0}
\int_{\partial (\delta B_{2n})}
[\ddc (h_A*\varphi_\varepsilon)]^{\wedge k}\wedge(\ddc h_{B_{2n}})^{\wedge \ell}\wedge\dc h_{B_{2n}}\wedge(\ddc h_{B_{2n}}^2)^{\wedge(n-k-\ell-1)}
=0
\end{equation}
because the integrand is locally summable, 
so we get that the left side of (\ref{stranissima}) equals
\begin{align}\label{fine1}
2^{\ell+1}
\lim_{\varepsilon\to0}
\int_{\partial B_{2n}}
[\ddc (h_A*\varphi_\varepsilon)]^{\wedge k}\wedge(\ddc h_{B_{2n}})^{\wedge \ell}\wedge\dc h_{B_{2n}}\wedge(\ddc h_{B_{2n}}^2)^{\wedge(n-k-\ell-1)}.
\end{align}
On the other hand
\begin{align}
&2^{\ell+1}
\int_{B_{2n}}
(\ddc h_A)^{\wedge k}\wedge(\ddc h_{B_{2n}})^{\wedge (\ell+1)}\wedge(\ddc h_{B_{2n}}^2)^{\wedge(n-k-\ell-1)}\nonumber
\\
&=
2^{\ell+1}
\lim_{m\to\infty}
\langle\!\langle 
(\ddc h_A)^{\wedge k}\wedge(\ddc h_{B_{2n}})^{\wedge (\ell+1)}\wedge(\ddc h_{B_{2n}}^2)^{\wedge(n-k-\ell-1)},
\varsigma_m
\rangle\!\rangle\nonumber
\\
&=
2^{\ell+1}
\lim_{m\to\infty}
\lim_{\varepsilon\to0}
\int_{\Supp\varsigma_m}
[\ddc (h_A*\varphi_\varepsilon)]^{\wedge k}\wedge(\ddc h_{B_{2n}})^{\wedge (\ell+1)}\wedge(\ddc h_{B_{2n}}^2)^{\wedge(n-k-\ell-1)}\varsigma_m\nonumber
\\
&=
2^{\ell+1}
\lim_{\varepsilon\to0}
\lim_{m\to\infty}
\int_{\Supp\varsigma_m}
[\ddc (h_A*\varphi_\varepsilon)]^{\wedge k}\wedge(\ddc h_{B_{2n}})^{\wedge (\ell+1)}\wedge(\ddc h_{B_{2n}}^2)^{\wedge(n-k-\ell-1)}\varsigma_m\nonumber
\\
&=
2^{\ell+1}
\lim_{\varepsilon\to0}
\int_{B_{2n}}
[\ddc (h_A*\varphi_\varepsilon)]^{\wedge k}\wedge(\ddc h_{B_{2n}})^{\wedge (\ell+1)}\wedge(\ddc h_{B_{2n}}^2)^{\wedge(n-k-\ell-1)}\nonumber
\\
&=
2^{\ell+1}
\lim_{\varepsilon\to0}
\lim_{\delta\to 0}
\int_{B_{2n}\setminus\delta B_{2n}}
[\ddc (h_A*\varphi_\varepsilon)]^{\wedge k}\wedge(\ddc h_{B_{2n}})^{\wedge (\ell+1)}\wedge(\ddc h_{B_{2n}}^2)^{\wedge(n-k-\ell-1)}\nonumber
\\
&=
2^{\ell+1}
\lim_{\varepsilon\to0}
\int_{\partial B_{2n}}
[\ddc (h_A*\varphi_\varepsilon)]^{\wedge k}\wedge(\ddc h_{B_{2n}})^{\wedge \ell}\wedge\dc h_{B_{2n}}\wedge(\ddc h_{B_{2n}}^2)^{\wedge(n-k-\ell-1)}\nonumber
\\
&-
2^{\ell+1}
\lim_{\varepsilon\to0}
\lim_{\delta\to 0}
\int_{\partial(\delta B_{2n})}
[\ddc (h_A*\varphi_\varepsilon)]^{\wedge k}\wedge(\ddc h_{B_{2n}})^{\wedge \ell}\wedge\dc h_{B_{2n}}\wedge(\ddc h_{B_{2n}}^2)^{\wedge(n-k-\ell-1)}\nonumber
\\
&=
2^{\ell+1}
\lim_{\varepsilon\to0}
\int_{\partial B_{2n}}
[\ddc (h_A*\varphi_\varepsilon)]^{\wedge k}\wedge(\ddc h_{B_{2n}})^{\wedge \ell}\wedge\dc h_{B_{2n}}\wedge(\ddc h_{B_{2n}}^2)^{\wedge(n-k-\ell-1)}\label{fine2}
\end{align}
because of (\ref{integrand}). Since (\ref{fine1}) and (\ref{fine2}) coincide, the induction is complete. 
\findim

\begin{rmq}\label{strano}
The case $\ell=n-k$ in \lemref{stranissimo} can be deduced by \thmref{supp-k} and \thmref{azzok1}. Indeed if $\Gamma\in{\mathcal P}(\mathbb C^n)$, by \thmref{supp-k}
\begin{align}
&\int_{B_{2n}}
(\ddc h_\Gamma)^{\wedge k}\wedge (\ddc h_{B_{2n}}^2)^{\wedge (n-k)}\nonumber
\\
&=
4^{n-k}(n-k)!\,k!
\lim_{m\to+\infty}
\sum_{\Delta\in{\mathcal B}_{\rm ed}(\Gamma,k)}
\varrho(\Delta)\vol_k(\Delta)
\langle\!\langle
[K_\Delta]\wedge \upsilon_{E^\perp_\Delta},\varsigma_m
\rangle\!\rangle\nonumber
\\
&=
4^{n-k}(n-k)!\,k!
\sum_{\Delta\in{\mathcal B}_{\rm ed}(\Gamma,k)}
\varrho(\Delta)\vol_k(\Delta)
\vol_{2n-k}(K_\Delta\cap B_{2n})\nonumber
\\
&=4^{n-k}(n-k)!\,k!\varkappa_{2n-k}
{\textsl v}^\varrho_k(\Gamma).\label{great}
\end{align}
On the other hand, by~\thmref{azzok1},
\begin{align*}
{\textsl v}^\varrho_k(\Gamma)
&=
\dfrac{\varkappa_n}{2^{n-k}\varkappa_{2n-k}}{n\choose k}
Q_n(\Gamma[k],B_{2n}[n-k])
\\
&=
\dfrac{1}{2^{n-k}\varkappa_{2n-k}k!(n-k)!}
\int_{B_{2n}}
(\ddc h_\Gamma)^{\wedge k}\wedge(\ddc h_{B_{2n}})^{\wedge(n-k)},
\end{align*}
so
\begin{align*}
\int_{B_{2n}}
(\ddc h_\Gamma)^{\wedge k}\wedge(\ddc h^2_{B_{2n}})^{\wedge(n-k)}
=
2^{n-k}
\int_{B_{2n}}
(\ddc h_\Gamma)^{\wedge k}\wedge(\ddc h_{B_{2n}})^{\wedge(n-k)}.
\end{align*}
For a general $A\in{\mathcal K}(\mathbb C^n)$, the conclusion follows by continuity.
\end{rmq}

\begin{coro}\label{strano1}
Let $0\leqslant k\leqslant n$ and $0\leqslant \ell\leqslant n-k$ be integers. For every $A_1,\ldots,A_k\in{\mathcal K}(\mathbb C^n)$,  
\begin{align*}
\int_{B_{2n}}
\left(
\bigwedge_{\ell=1}^k\ddc h_{A_\ell}
\right)
\wedge(\ddc h_{B_{2n}}^2)^{\wedge(n-k)}
=
2^\ell
\int_{B_{2n}}
\left(
\bigwedge_{\ell=1}^k\ddc h_{A_\ell}
\right)
\wedge(\ddc h_{B_{2n}})^{\wedge \ell}
\wedge(\ddc h^2_{B_{2n}})^{\wedge(n-k-\ell)}.
\end{align*}
In particular,
\begin{align}\label{strana2}
Q_n(A_1,\ldots,A_k,B_{2n}[n-k])
=
\dfrac{1}{2^{n-k}n!\varkappa_n}
\int_{B_{2n}}
\left(
\bigwedge_{\ell=1}^k\ddc h_{A_\ell}
\right)
\wedge(\ddc h^2_{B_{2n}})^{\wedge (n-k)}\,.
\end{align}
\end{coro}

\noindent
\pf Both sides are multilinear with respect to $A_1,\ldots,A_k$, so they agree as soon as they do it when $A_1=\ldots,A_k=A$, which boils down to \lemref{stranissimo}. Formula (\ref{strana2}) is just the first statement for $\ell=n-k$.\findim

\begin{rmq}\label{vcorrk}

\cororef{combinatorio} implies a proof of \cororef{azzok2} not depending on \thmref{azzok1}. Indeed, if $\Gamma_1,\ldots,\Gamma_k\in{\mathcal P}(\mathbb C^n)$ and $\Gamma=\sum_{\ell=1}^k \Gamma_\ell$, \cororef{combinatorio} and \cororef{strano1} with $\ell=n-k$ yield
\begin{align*}
&\dfrac{\varkappa_n}{2^{n-k}\varkappa_{2n-k}}{n\choose k}Q_n(\Gamma_1,\ldots,\Gamma_k,B_{2n}[n-k])
\\
&=
\dfrac{1}{2^{n-k}\varkappa_{2n-k}k!(n-k)!}
\int_{B_{2n}}
\left(
\bigwedge_{\ell=1}^k
\ddc h_{\Gamma_\ell}
\right)\wedge(\ddc h_{B_{2n}})^{\wedge(n-k)}
\\
&=
\dfrac{1}{4^{n-k}\varkappa_{2n-k}k!(n-k)!}
\int_{B_{2n}}
\left(
\bigwedge_{\ell=1}^k
\ddc h_{\Gamma_\ell}
\right)\wedge(\ddc h^2_{B_{2n}})^{\wedge(n-k)}
\\
&=
\lim_{m\to 0}
\dfrac{1}{\varkappa_{2n-k}}
\sum_{\Delta\in{\mathcal B}_{\rm ed}(\Gamma,k)}
\varrho(\Delta)
V_k(\Delta_1,\ldots,\Delta_k)
\langle\!\langle 
[K_\Delta]\wedge\upsilon_{E_\Delta^\perp},\varsigma_m
\rangle\!\rangle
\\
&=
\dfrac{1}{\varkappa_{2n-k}}
\sum_{\Delta\in{\mathcal B}_{\rm ed}(\Gamma,k)}
\varrho(\Delta)
V_k(\Delta_1,\ldots,\Delta_k) 
\vol_{2n-k}(K_\Delta\cap B_{2n})
\\
&=
\sum_{\Delta\in{\mathcal B}_{\rm ed}(\Gamma,k)}
\varrho(\Delta)
V_k(\Delta_1,\ldots,\Delta_k) 
\psi_\Gamma(\Delta)
\\
&=
V^\varrho_k(\Gamma_1,\ldots,\Gamma_k)\,.
\end{align*}
By continuity, the conclusion of \cororef{azzok2} follows for arbitrary convex bodies. 
\end{rmq}

\section{Kazarnovsk{\v\i}i pseudovolume and mixed discriminants}\label{kmd}

In this section we investigate the relations between mixed pseudovolume and mixed discriminants\footnote{It should be noted that the terminology \emph{mixed discriminant} refers at least to two different notions in mathematics: one is related to the elimination theory of Laurent polynomials, the other belongs to linear algebra and is the one we are going to deal with here.} in order to get an analogue of \eqref{mix-D-f} for the mixed pseudovolume. As a result of our investigation we will obtain a few new formulas for $Q_n$ and $V_n$.

\subsection{Mixed discriminants}
Let ${\mathfrak S}_n$\label{gruppo-simm} be the symmetric group over the set $\{1,\ldots,n\}$. For $\ell=1,\ldots,n$, let $M_\ell=(m^{(\ell)}_{j,k})$ be an $n$-by-$n$ complex matrix and, for every $\sigma\in{\mathfrak S}_n$, let $\text{mix}_\sigma(M_1,\ldots,M_n)$ and $\text{mix}^\sigma(M_1,\ldots,M_n)$\label{mix} denote the $n$-by-$n$ matrices respectively defined by
\begin{equation*}
\text{mix}_\sigma(M_1,\ldots,M_n)
=
\left(m^{(\ell)}_{j,\sigma(\ell)}\right)=
\left(
\begin{array}{ccc}
m^{(1)}_{1,\sigma(1)} & \cdots & m^{(n)}_{1,\sigma(n)} \\
\vdots & \ddots & \vdots \\
m^{(1)}_{n,\sigma(1)} & \cdots & m^{(n)}_{n,\sigma(n)}
\end{array}
\right)\,,
\end{equation*}
\begin{equation*}
\text{mix}^\sigma(M_1,\ldots,M_n)
=
\left(m^{(\sigma(\ell))}_{j,\ell}\right)
=
\left(
\begin{array}{ccc}
m^{(\sigma(1))}_{1,1} & \cdots & m^{(\sigma(n))}_{1,n} \\
\vdots & \ddots & \vdots \\
m^{(\sigma(1))}_{n,1} & \cdots & m^{(\sigma(n))}_{n,n}
\end{array}
\right)\,.
\end{equation*}
The matrices $\text{mix}_\sigma(M_1,\ldots,M_n)$ and $\text{mix}^\sigma(M_1,\ldots,M_n)$ are mixed versions of $M_1,\ldots,M_n$, indeed  the $\ell$-th column of $\text{mix}_\sigma(M_1,\ldots,M_n)$ is the $\sigma(\ell)$-th column of $M_\ell$ and the $\ell$-th column of $\text{mix}^\sigma(M_1,\ldots,M_n)$ is $\ell$-th column of $M_{\sigma(\ell)}$. In general $\text{mix}_\sigma(M_1,\ldots,M_n)\neq\text{mix}^\sigma(M_1,\ldots,M_n)$ unless $\sigma$ is the identity. When $M_1=\ldots=M_n=M$ one has $\text{mix}^\sigma(M,\ldots,M)=M$ too. Of course $\text{mix}_\sigma(\overline{M}_1,\ldots,\overline{M}_n)=\overline{\text{mix}_\sigma(M_1,\ldots,M_n)}$. 

If $\text{sgn}(\sigma)$ denotes the sign of $\sigma\in{\mathfrak S}_n$, then it is easy to see that 
\begin{equation*}
\text{sgn}(\sigma)
\det \left(\text{mix}_\sigma(M_1,\ldots,M_n)\right)
=
\det \left(\text{mix}^{\sigma^{-1}}(M_1,\ldots,M_n)\right)
\,.
\end{equation*}
The \emph{mixed discriminant} of $M_1,\ldots,M_n$ is the complex number $D_n(M_1,\ldots,M_n)$\label{dm} given by 
\begin{align}
D_n(M_1,\ldots,M_n)
&=
\dfrac{1}{n!}
\sum_{\sigma\in{\mathfrak S}_n}
\sum_{\tau\in{\mathfrak S}_n}
\text{sgn}(\tau)
\prod_{\ell=1}^n
m^{(\sigma\cdot\tau(\ell))}_{\ell,\tau(\ell)}
\\
&=\dfrac{1}{n!}
\sum_{\sigma\in{\mathfrak S}_n}
\det\left(
\text{mix}^\sigma(M_1,\ldots,M_n)
\right)\label{disc}
\\
&=
\dfrac{1}{n!}\sum_{\sigma\in{\mathfrak S}_n}\text{sgn}(\sigma) \det \left(\text{mix}_\sigma(M_1,\ldots,M_n)\right)\,.\label{disc2}
\end{align}
From this definition it follows that, for every $n$-by-$n$ complex matrices $M,N,M_1,M_2,\ldots,M_n$,
\begin{enumerate}
\item\label{sigmaetau} $D_n(M_1,\ldots,M_n)
=
\dfrac{1}{n!}
\displaystyle\sum_{\sigma,\tau\in{\mathfrak S}_n}
\text{sgn}(\sigma\cdot\tau)
\prod_{\ell=1}^n
m^{(\ell)}_{\sigma(\ell),\tau(\ell)}
$;
\item $n!D_n(M_1,\ldots,M_n)$ is the coefficient of $t_1\ldots t_n$ in the polynomial $\det(t_1M_1+\ldots+t_nM_n)$;
\item $n!D_n(M_1,\ldots,M_n)=\dfrac{\partial}{\partial t_1\ldots\partial t_n}\det(t_1M_1+\ldots+t_nM_n)$;
\item $
n!D_n(M_1,\ldots,M_n)
=
\displaystyle\sum_{\varnothing\neq I\subseteq\{1,\ldots,n\}}
(-1)^{n-\# I}
\det
\bigg(
\sum_{\ell\in I}
M_\ell
\bigg)
$;
\item $D_n(M,\ldots,M)=\det M$;
\item $D_n(I_n,\ldots,I_n)=1$;
\item $D_n(M_{\sigma(1)},\ldots,M_{\sigma(n)})=D_n(M_1,\ldots,M_n)$, for every $\sigma\in{\mathfrak S}_n$;
\item $D_n(aM+bN,M_2,\ldots,M_n)=aD_n(M,M_2,\ldots,M_n)+bD_n(N,M_2,\ldots,M_n)$ for every $a,b\in\mathbb C$;
\item $D_n({}^t\!M_1,\ldots,{}^t\!M_n)=D_n(M_1,\ldots,M_n)$ and $D_n(\overline{M}_1,\ldots,\overline{M}_n)=\overline{D_n(M_1,\ldots,M_n)}$, in particular\footnote{As usual $M^*={}^t\!\overline{M}$.} $D_n(M^*_1,\ldots,M^*_n)=\overline{D_n(M_1,\ldots,M_n)}$;
\item If $M_1,\ldots,M_n$ are hermitian positive semidefinite (resp. definite), then $D_n(M_1,\ldots,M_n)\geq 0$ (resp. $D_n(M_1,\ldots,M_n)> 0$);
\item $D_n(NM_1,\ldots,NM_n)=(\det N)\cdot D_n(M_1,\ldots,M_n)=D_n(M_1N,\ldots,M_nN)$, in particular one has $D_n(N^*M_1N,\ldots,N^*M_nN)=\vert\det N\vert^2 D_n(M_1,\ldots,M_n)$ and so $D_n(N^*M_1N,\ldots,N^*M_nN)=D_n(M_1,\ldots,M_n)$ in case $\vert\det N\vert=1$.
\end{enumerate}
Il $M$ is an $n$-by-$n$ matrix and $1\leq j,k\leq n$, let $M^{[j,-]}$, $M^{[-,k]}$ and $M^{[j,k]}$\label{sottomatrice} denote the submatrices respectively obtained by deleting the $j$-th row, the $k$-th column and both the $j$-th row and the $k$-th column of $M$. For every $1\leq \ell\leq n$, by Laplace expansion along the $\ell$-th column of the matrices in (\ref{disc2}) it follows that $D_n(M_1,\ldots,M_n)$ is a linear form in the entries of $M_\ell$. The formula (\ref{comoda}) proved by the following lemma will be referred to as the \textit{Laplace expansion of $D_n(M_1,\ldots,M_n)$ along the $\ell$-th matrix}.

\begin{lem}\label{laplace1}
If $M_1=(m^{(1)}_{j,k}),\ldots,M_n=(m^{(n)}_{j,k})$ are $n$-by-$n$ matrices, then, for every $1\leqslant\ell\leqslant n$,
\begin{equation}\label{comoda}
D_n(M_1,\ldots,M_n)
=
\dfrac{1}{n}
\sum_{j,k=1}^n
(-1)^{j+k}
m_{j,k}^{(\ell)}\,
 D_{n-1}(M_1^{[j,k]},\ldots,M_{\ell-1}^{[j,k]},M_{\ell+1}^{[j,k]},\ldots,M_n^{[j,k]})
\end{equation}
\end{lem}
\noindent
\pf
By (\ref{disc2}), 
\begin{align}
D_n(M_1,\ldots,M_n)
&=
\dfrac{1}{n!}
\sum_{\sigma\in{\mathfrak S}_n}
\text{sgn}(\sigma)
\sum_{j=1}^n
(-1)^{j+\ell}
m^{(\ell)}_{j,\sigma(\ell)}
\det
\text{mix}_\sigma(M_1,\ldots,M_n)^{[j,\ell]}
\\
&=
\dfrac{1}{n!}
\sum_{j=1}^n
\sum_{k=1}^n
\sum_{\sigma\in{\mathfrak S}_n (\ell,k)}
(-1)^{j+\ell}
\text{sgn}(\sigma)
m^{(\ell)}_{j,k}
\det
\text{mix}_\sigma(M_1,\ldots,M_n)^{[j,\ell]}\,,\label{quasiquasi}
\end{align}
where, for any fixed $1\leqslant \ell,k\leqslant n$, we set ${\mathfrak S}_n(\ell,k)=\{\sigma\in{\mathfrak S}_n\mid\sigma(\ell)=k\}$. The subset ${\mathfrak S}_n(\ell,k)\subset{\mathfrak S}_n$ and the subgroup ${\mathfrak S}_{n-1}=\{\tau\in{\mathfrak S}_n\mid \tau(n)=n\}\simeq{\mathfrak S}_{n-1}$ have the same number of elements, an explicit bijection $\gamma_{\ell,k}:{\mathfrak S}_n(\ell,k)\to {\mathfrak S}_{n-1}$ is given by
\begin{align}
\gamma_{\ell,k}(\sigma)(r)
&=
\begin{cases}
\sigma(r)&\text{if }\sigma(r)<k
\\
\sigma(r)-1&\text{if }\sigma(r)>k
\end{cases}
&\text{for } 1\leqslant r\leqslant\ell-1
\end{align}
and
\begin{align}
\gamma_{\ell,k}(\sigma)(r)
&=
\begin{cases}
\sigma(r+1)&\text{if }\sigma(r+1)< k
\\
\sigma(r+1)-1&\text{if }\sigma(r+1)>k\end{cases}
&\text{for }\ell\leqslant r\leqslant n-1\,.
\end{align}
The permutation $\gamma_{\ell,k}(\sigma)$ is such that
\begin{equation}
\text{mix}_\sigma(M_1,\ldots,M_n)^{[j,\ell]}
=
\text{mix}_{\gamma_{\ell,k}(\sigma)}
(M_1^{[j,k]},\ldots,
M_{\ell-1}^{[j,k]},M_{\ell+1}^{[j,k]},\ldots,M_n^{[j,k]})
\end{equation}
and $\text{sgn}(\sigma)=(-1)^{\ell+k}\text{sgn}(\gamma_{\ell,k}(\sigma))$, it follows that (\ref{quasiquasi}) becomes 
\begin{align*}
&\dfrac{1}{n!}
\sum_{j,k=1}^n
\sum_{\sigma\in{\mathfrak S}_{n}(\ell,k)}
(-1)^{j+2\ell+k}
\text{sgn}(\gamma_{\ell,k}(\sigma))
m^{(\ell)}_{j,k}
\det
\text{mix}_{\gamma_{\ell,k}(\sigma)}
(M_1^{[j,k]},\ldots,
M_{\ell-1}^{[j,k]},M_{\ell+1}^{[j,k]},\ldots,M_n^{[j,k]})
\\
&\dfrac{1}{n!}
\sum_{j,k=1}^n
\sum_{\tau\in{\mathfrak S}_{n-1}}
(-1)^{j+k}
\text{sgn}(\tau)
m^{(\ell)}_{j,k}
\det
\text{mix}_{\tau}
(M_1^{[j,k]},\ldots,
M_{\ell-1}^{[j,k]},M_{\ell+1}^{[j,k]},\ldots,M_n^{[j,k]})
\\
&=
\dfrac{1}{n!}
\sum_{j,k=1}^n
(-1)^{j+k}
m^{(\ell)}_{j,k}
\sum_{\tau\in{\mathfrak S}_{n-1}}
\text{sgn}(\tau)
\det\text{mix}_\tau
(M_1^{[j,k]},\ldots,M_{\ell-1}^{[j,k]},M_{\ell+1}^{[j,k]},\ldots,M_n^{[j,k]})
\\
&=
\dfrac{1}{n}
\sum_{j,k=1}^n
(-1)^{j+k}
m^{(\ell)}_{j,k}\,
D_{n-1}(M_1^{[j,k]},\ldots,M_{\ell-1}^{[j,k]},M_{\ell+1}^{[j,k]},\ldots,M_n^{[j,k]}).
\end{align*}
The proof is complete.
\findim
The coefficient $(-1)^{j+k}(1/n)D_{n-1}(M_1^{[j,k]},\ldots,M_{\ell-1}^{[j,k]},M_{\ell+1}^{[j,k]},\ldots,M_n^{[j,k]})$ will be referred to as the \emph{mixed minor} of $m_{j,k}^{(\ell)}$ in the Laplace expansion of $D_n(M_1,\ldots,M_n)$ along $M_\ell$; a \emph{principal mixed minor} will be the mixed minor of a diagonal entry of $M_\ell$.

In the special case of a hermitian matrix $M$, although the submatrices $M^{[j,k]}$ are not hermitian (unless $j=k$), they still satisfy the identities
\begin{equation}
M^{[j,k]}
=
M^{[k,j]*}
\,,
\end{equation}
for $1\leqslant j,k\leqslant n$. 
If $M_1,\ldots,M_{\ell-1},M_{\ell},\ldots,M_n$ are hermitian matrices, the preceding relations yield 
\begin{align}
&D_{n-1}\left(M_1^{[j,k]},\ldots,M_{\ell-1}^{[j,k]},M_{\ell+1}^{[j,k]},\ldots,M_n^{[j,k]}\right)
\\
&=
D_{n-1}\left(M_1^{[k,j]*},\ldots,M_{\ell-1}^{[k,j]*},M_{\ell+1}^{[k,j]*},\ldots,M_n^{[k,j]*}\right)
\\
&=
\overline{
D_{n-1}
\left(
M_1^{[k,j]},\ldots,
M_{\ell-1}^{[k,j]},
M_{\ell+1}^{[k,j]},\ldots
M_n^{[k,j]}
\right)\,,
}
\end{align}
where the last equality follows by property 9 of mixed discriminants. This means that the $n$-by-$n$ matrix $Z^{\langle\ell\rangle}$ consisting of the mixed minors 
\begin{equation}\label{mZ}
Z^{\langle\ell\rangle}_{j,k}=\frac{(-1)^{j+k}}{n}D_{n-1}\left(M_1^{[j,k]},\ldots,M_{\ell-1}^{[j,k]},M_{\ell+1}^{[j,k]},\ldots,M_n^{[j,k]}\right)
\end{equation}
is a hermitian too. When $M_\ell$ is also hermitian, formula (\ref{comoda}) gets the following simpler form
\begin{equation}\label{comoda2}
\begin{split}
D_n(M_1,\ldots,M_n)
&=
\sum_{1\leqslant j\leqslant n}
m^{(\ell)}_{j,j}\,
Z^{\langle\ell\rangle}_{j,j}
+
\sum_{1\leqslant j<k\leqslant n}
2
\re\left(
m^{(\ell)}_{j,k}
Z^{\langle\ell\rangle}_{j,k}\right)\,,
\end{split}
\end{equation}
and, if $M_\ell$ or $Z^{\langle\ell\rangle}$ is in diagonal form, the preceding formula becomes even simpler:
\begin{align}\label{comoda3}
D_n(M_1,\ldots,M_n)
&=
\sum_{j=1}^n
m^{(\ell)}_{j,j}\,
Z^{\langle\ell\rangle}_{j,j}\,.
\end{align}

If, moreover, $M_1,\ldots,M_{\ell-1},M_{\ell+1},\ldots,M_n$ are positive definite (resp. semidefined), then the matrix $Z^{\langle\ell\rangle}$ is positive definite (resp. semidefined) too. Indeed, up to a unitary transformation of coordinates, $Z^{\langle\ell\rangle}$ can be diagonalized and 
\begin{equation}
Z^{\langle\ell\rangle}_{j,j}=\frac{1}{n}D_{n-1}\left(M_1^{[j,j]},\ldots,M_{\ell-1}^{[j,j]},M_{\ell+1}^{[j,j]},\ldots,M_n^{[j,j]}\right)
\end{equation}
is positive (resp. non-negative) by property 10 of the mixed discriminant applied to the hermitian matrices $M_1^{[j,j]},\ldots,M_{\ell-1}^{[j,j]},M_{\ell+1}^{[j,j]},\ldots,M_n^{[j,j]}$, which are positive definite (resp. semidefined) by Sylvester's criterion. 

For the sake of notation, we set $Z=Z^{\langle1\rangle}$ and we state the next lemma for $\ell=1$, but it holds true for every $1\leqslant\ell\leqslant n$.

\begin{lem}\label{Alex0}
Let $M_2,\ldots,M_n$ be $n$-by-$n$ positive definite hermitian matrices and let $U$ be a non singular $n$-by-$n$ matrix. If 
\begin{equation}
Z^\prime_{j,k}=\dfrac{(-1)^{j+k}}{n}D_{n-1}\left((U^*M_2U)^{[j,k]},\ldots,(U^*M_nU)^{[j,k]}\right)
\end{equation}
then 
$Z^\prime=
\vert\det U\vert^2\;
\overline{U^{-1}}\;Z\;{}^tU^{-1}
$. In particular $\det Z^\prime=\det Z$.
\end{lem}
\noindent
\pf
Since $M_2,\ldots,M_n$ are hermitian and positive definite, the same happens to the matrix $Z$, in particular $Z$ is non singular. Suppose that $M_2=\ldots=M_n=M$, then 
\begin{equation}
\overline{Z}_{j,k}=Z_{k,j}=\dfrac{(-1)^{k+j}}{n}\det M^{[k,j]}=\frac{\det M}{n}(M^{-1})_{j,k},
\end{equation}
i.e. 
\begin{equation}\label{binet}
 Z=\dfrac{\det M}{n}\overline{M^{-1}}.
 \end{equation} 
Formula~\ref{binet} for $U^*MU$ reads
\begin{align}
Z^\prime
&=
\dfrac{\det(U^*MU)}{n}\;\overline{(U^*MU)^{-1}}
=
\vert \det U\vert^2\;\dfrac{\det M}{n}\;\overline{U^{-1}}\;\overline{M^{-1}}\;{}^tU^{-1}
=
\vert \det U\vert^2\;
\overline{U^{-1}}\;Z\;{}^tU^{-1},
\end{align}
or equivalently
\begin{equation}
Z^\prime_{j,k}
=
\vert\det U\vert^2
\sum_{r,s=1}^n
\left(\overline{U^{-1}}\right)_{j,r}Z_{r,s}\left(U^{-1}\right)_{k,s},
\end{equation}
for every $1\leqslant j,k\leqslant n$. 
We now apply the preceding argument to the matrix $N=\sum_{\ell=2}^n t_\ell M_\ell$, so 
\begin{align}
&\dfrac{(-1)^{j+k}}{n}
D_{n-1}\left((U^*NU)^{[j,k]},\ldots,(U^*NU)^{[j,k]}\right)
\\
&=
\vert\det U\vert^2
\dfrac{(-1)^{j+k}}{n}
\sum_{r,s=1}^n
\left(\overline{U^{-1}}\right)_{j,r}
D_{n-1}\left(N^{[r,s]},\ldots,N^{[r,s]}\right)
\left(U^{-1}\right)_{k,s}
\end{align}
and by the multilinearity of $D_{n-1}$ we get
\begin{align}
&\dfrac{(-1)^{j+k}}{n}
\sum_{\ell_2,\ldots,\ell_{n}=2}^{n}
t_{\ell_2}\cdots t_{\ell_n}
D_{n-1}\left((U^*M_{\ell_2}U)^{[j,k]},\ldots,(U^*M_{\ell_n}U)^{[j,k]}\right)
\\
&=
\vert\det U\vert^2
\dfrac{(-1)^{j+k}}{n}
\sum_{\ell_2,\ldots,\ell_{n}=2}^{n}
t_{\ell_2}\cdots t_{\ell_n}
\sum_{r,s=1}^n
\left(\overline{U^{-1}}\right)_{j,r}
D_{n-1}\left(M_{\ell_2}^{[r,s]},\ldots,M_{\ell_n}^{[r,s]}\right)
\left(U^{-1}\right)_{k,s}.
\end{align}
Now by selecting the coefficients of $t_2\cdots t_n$ in the preceding equality we obtain 
\begin{equation}
Z^\prime_{j,k}
=
\vert\det U\vert^2\sum_{r,s=1}^n
\left(\overline{U^{-1}}\right)_{j,r}
Z_{r,s}
\left(U^{-1}\right)_{k,s},
\end{equation}
i.e. $Z^\prime=\vert\det U\vert^2\; \overline{U^{-1}}\;Z\;{}^tU^{-1}$. As a consequence $\det Z^\prime=\vert\det U\vert^2\overline{\det U^{-1}}\det Z\det U^{-1}=\det Z$.
\findim

\begin{coro}\label{Alex3}
Let $M_1,M_2,\ldots,M_n$ be complex hermitian $n$-by-$n$ matrices such that $M_2,\ldots,M_n$ are positive definite. Then  $a\tr M_1\leqslant D_n(M_1,M_2,\ldots,M_n)\leqslant b \tr M_1$, for some positive constants $a$ and $b$.
\end{coro}
\noindent
\pf
Since $M_2,\ldots,M_n$ are hermitian and positive definite, the matrix $Z=Z^{\langle1\rangle}$ is hermitian and positive definite too. Let $U$ be a unitary matrix such that $U^*M_1U$ is diagonal, then by~\lemref{Alex0} 
\begin{align}
D_n(M_1,M_2,\ldots,M_n)
&=
\vert \det U\vert^2
D_n(M_1,M_2,\ldots,M_n)
\\
&=
D_n\left(U^*M_1U,U^*M_2U,\ldots,U^*M_nU\right)
\\
&=
\sum_{j=1}^n
\left(U^*M_1U\right)_{j,j}
(\overline{U^{-1}}\;Z\;{}^tU^{-1})_{j,j}.
\end{align}
Now set $a=\min_{1\leqslant j\leqslant n}
(\overline{U^{-1}}\;Z\;{}^tU^{-1})_{j,j}$, $b=\max_{1\leqslant j\leqslant n}
(\overline{U^{-1}}\;Z\;{}^tU^{-1})_{j,j}$ and observe that both $a$ and $b$ are positive since $\overline{U^{-1}}\;Z\;{}^tU^{-1}$ is positive definite. The conclusion follows from the equality $\tr U^*M_1U=\tr M_1$. \findim

The mixed discriminant has several other properties, among which \thmref{Ada} in Section~\ref{open} will recall one of the most important. The interested reader is referred to~\cite{A2}, \cite{Bap} or~\cite{FMS} for further information. 

\subsection{Kazarnovski{\v\i} mixed pseudovolume through mixed discriminant}

We can now realize that (\ref{mix-D-f}) is an instance of (\ref{comoda3}) with $\ell=1$, $M_1=I_n$ and, for $2\leqslant k\leqslant n$, $M_k=\Hess_\mathbb R h_{h_{A_k}}$, $A_k\in{\mathcal K}_2(\mathbb R^n)$. Indeed, by (\ref{mix-D-f}), one obtains 
\begin{align}
V_n(A_1,A_2\ldots,A_n)
&=
\int_{\partial B_n}
h_{A_1}
D(\Hess_\mathbb R h_{A_2},\ldots,\Hess_\mathbb R h_{A_n})\upsilon_{\partial B_n}
\\
&=
\int_{\partial B_n}
h_{A_1}
D_n(I_n,\Hess_\mathbb R h_{A_2},\ldots,\Hess_\mathbb R h_{A_n})\upsilon_{\partial B_n}\label{mix-D-f2}
\\
&=
\dfrac{1}{n}
\int_{\partial B_n}
h_{A_1}
\sum_{j=1}^n
D_{n-1}\left(\Hess_\mathbb R h_{A_2}^{[j,j]},\ldots,\Hess_\mathbb R h_{A_n}^{[j,j]}\right)\upsilon_{\partial B_n}\label{mix-D-f3}\,,
\end{align}
so that the term $D(\Hess_\mathbb R h_{A_2},\ldots,\Hess_\mathbb R h_{A_n})$ in (\ref{mix-D-f}) is nothing but the sum of the principal mixed minors corresponding to the nonzero entries of $I_n$ in the Laplace expansion of the mixed discriminant $D_n(I_n,\Hess_\mathbb R h_{A_2},\ldots,\Hess_\mathbb R h_{A_n})$ along $I_n$.

Moreover, if $A_1,\ldots,A_n\in{\mathcal K}_2(\mathbb C^n)$ and $\Hess_\mathbb C h_{A_1},\ldots,\Hess_\mathbb C h_{A_n}$\label{c-hess} are the respective complex hessian matrices of the support functions $h_{A_1},\ldots,h_{A_n}$, then 
\begin{equation*}
\ddc h_{A_1}\wedge\ldots\wedge \ddc h_{A_n}= 4^n n! D_n(\Hess_\mathbb C h_{A_1},\ldots,\Hess_\mathbb C h_{A_n})\,\upsilon_{2n}\,.
\end{equation*}
Indeed, both sides depend multilinearly on $A_1,\ldots,A_n$ and, for $A_1=\ldots=A_n=A$, by~\eqref{MA} and the equality $D_n(\Hess_\mathbb C h_{A},\ldots,\Hess_\mathbb C h_{A})=\det(\Hess_\mathbb C h_A)$, one gets $4^nn!\det (\Hess_\mathbb C h_A) \upsilon_{2n}$ on both sides. It follows that
\begin{equation}
Q_n(A_1,\ldots,A_n)=\dfrac{4^n}{\varkappa_n}\int_{B_{2n}} D_n(\Hess_\mathbb C h_{A_1},\ldots,\Hess_\mathbb C h_{A_n})\,\upsilon_{2n}\,,\label{kmix-D-f}
\end{equation}
which provides a further formula for $Q_n$ on ${\mathcal K}_2(\mathbb C^n)$.

As shown by the corollary to the following lemma, the preceding formula~(\ref{kmix-D-f}) is false if the bodies $A_1,\ldots,A_n$ belong to ${\mathcal K}_2(\mathbb R^n)$.

\begin{lem}\label{Eulero}
Let $f\in{\mathcal C}^2(\mathbb C^n\setminus\{0\})$ be a positively homogeneous function of degree $1$ depending on the real variables $x_1,\ldots,x_n$ only. If $\Hess_\mathbb R f$ denotes its hessian matrix with respect to $x_1,\ldots,x_n$ at the point $z$, then
\begin{equation}
\Hess_\mathbb R f(z) \re z=0
\end{equation} 
for every $z\in\mathbb C^n\setminus\{0\}$, in particular $\det \Hess_\mathbb R f(z)=0$.
\end{lem}

\noindent
\pf
By Euler's relation and thanks to the independence of $f$ on $y_1,\ldots,y_n$
\begin{equation}
f(z)
=
\sum_{s=1}^n
\dfrac{\partial f}{\partial x_s}(z) x_s
+
\dfrac{\partial f}{\partial y_s}(z) y_s
=
\sum_{s=1}^n
\dfrac{\partial f}{\partial x_s}(z) x_s,
\end{equation}
for every $z=x+iy$. A differentiation with respect to $x_r$, $1\leqslant r\leqslant n$, yields
\begin{equation}
\sum_{s=1}^n
\dfrac{\partial^2 f}{\partial x_s\partial x_r}(z) x_s=0\,,
\end{equation} 
in particular $\det\Hess_\mathbb R f(z)=0$.
\findim

\begin{coro}
Let $A_1,\ldots,A_n\in{\mathcal K}_2(\mathbb R^n)$, then 
\begin{equation}
D_n(\Hess_\mathbb R h_{A_1},\ldots,\Hess_\mathbb R h_{A_n})=0\,.
\end{equation}
In particular, if $V_n(A_1,\ldots,A_n)\neq0$, the equality (\ref{kmix-D-f}) is no longer correct.
\end{coro}

\noindent
\pf 
For any $\varnothing\neq I\subseteq \{1,\ldots,n\}$, set $A_I=\sum_{k\in I}A_k$, then 
\begin{equation}
\det\left(\sum_{k\in I}\Hess_\mathbb R h_{A_k}\right)=\det\left(\Hess_\mathbb R h_{A_I}\right)=0,
\end{equation} where the second equality is due to \lemref{Eulero}. Then $D_n(\Hess_\mathbb R h_{A_1},\ldots,\Hess_\mathbb R h_{A_n})=0$ by property 4 of mixed discriminant.
Now suppose that $V_n(A_1,\ldots,A_n)\neq0$. 
As a subset of $\mathbb C^n$, $A_k$ is not strictly convex and
\begin{equation}\label{rchess}
\Hess_\mathbb C h_{A_k}=\dfrac{1}{4}\Hess_\mathbb R h_{A_k},
\end{equation} 
for every $1\leqslant k\leqslant n$. 
The expression~(\ref{ddcreale}) of $\ddc$ in real coordinates reads
\begin{equation}
\ddc h_{A_k}(z)=\sum_{\ell,j=1}^n\frac{\partial^2 h_{A_k}}{\partial x_\ell\partial x_j}\,\de x_\ell\wedge\de y_j\,,
\end{equation} 
whence
\begin{align}
\dfrac{4^n}{\varkappa_n}\int_{B_{2n}} 
\ddc h_{A_1}\wedge\ldots\wedge\ddc h_{A_n}\nonumber
=
\dfrac{1}{\varkappa_n}\int_{B_{2n}} 
D_n(\Hess_\mathbb R h_{A_1},\ldots,\Hess_\mathbb R h_{A_n})\upsilon_{2n}
=0
\,,
\end{align}
however $Q_n(A_1,\ldots,A_n)=V_n(A_1,\ldots,A_n)\neq0$.\findim

Now we are going to manipulate (\ref{kmix-D-f}) in order to get an equivalent representation of $Q_n$ involving an integral on $\partial B_{2n}$ in the vein of (\ref{mix-D-f2}). This goal can be achieved in at least two different ways as shown by \lemref{dj} and \lemref{djj} below.

Let $A\in{\mathcal K_1}(\mathbb C^n)$ admit the origin as interior point so that $h_{A}(u)>0$, for every $u\in\mathbb C^n\setminus\{0\}$. For every $u\in \mathbb C^n\setminus\{0\}$, we associate to $A$ two $n$-by-$n$ complex matrices, namely \begin{equation}\label{matriceM}
M_{A}(u)=
\begin{pmatrix}
	u_1
	\\
	\vdots
	\\
	u_n
\end{pmatrix}
\begin{pmatrix}
	\dfrac{\partial h_A}{\partial z_1}(u)&\ldots&\dfrac{\partial h_A}{\partial z_n}(u)
\end{pmatrix}
=
\begin{pmatrix}
u_1\dfrac{\partial h_A}{\partial z_1}(u)&\ldots&u_1\dfrac{\partial h_A}{\partial z_n}(u)
\\
\vdots&\ddots&\vdots\\
u_n\dfrac{\partial h_A}{\partial z_1}(u)&\ldots&u_n\dfrac{\partial h_A}{\partial z_n}(u)
\end{pmatrix}
\end{equation}
and the hermitian matrix
\begin{equation}\label{matriceN}
N_A(u)=\overline M_A(u)+{}^t\!M_A(u)
=
\begin{pmatrix}
\bar u_1\dfrac{\partial h_A}{\partial \bar z_1}(u)
+
u_1\dfrac{\partial h_A}{\partial z_1}(u)
&\ldots&
\bar u_1\dfrac{\partial h_A}{\partial \bar z_n}(u)
+
u_n\dfrac{\partial h_A}{\partial z_1}(u)
\\
\vdots&\ddots&\vdots
\\
\bar u_n\dfrac{\partial h_A}{\partial \bar z_1}(u)
+
u_1\dfrac{\partial h_A}{\partial z_n}(u)
&\ldots&
\bar u_n\dfrac{\partial h_A}{\partial \bar z_n}(u)
+
u_n\dfrac{\partial h_A}{\partial z_n}(u)
\end{pmatrix}
\,.
\end{equation}
The matrix $M_A(u)$ has rank $1$ and so $N_A(u)$ has rank at most $2$ whence, for any $3\leqslant k\leqslant n$, all the principal $k$-by-$k$ minors vanish. 
The trace is
\begin{align}
\text{tr}\left(N_A(u)\right)
&=
\sum_{\ell=1}^n
\bar u_\ell
\dfrac{\partial h_A}{\partial \bar z_\ell}(u)
+
u_\ell
\dfrac{\partial h_A}{\partial z_\ell}(u)
\\
&=
2
\sum_{\ell=1}^n
\re
\left(
\bar u_\ell
\dfrac{\partial h_A}{\partial \bar z_\ell}(u)
\right)
\\
&=
\re\langle\nabla h_A(u),u\rangle
\\
&=
h_A(u)>0\,,
\end{align}
where the last equality comes from Euler's theorem on positively homogeneous functions. The principal $2$-by-$2$ minor obtained by selecting the $j$-th and the $k$-th rows and columns of $N_A(u)$, with $j<k$, is
\begin{align}
&\left(
\bar u_j\dfrac{\partial h_A}{\partial \bar z_j}(u)
+
u_j\dfrac{\partial h_A}{\partial z_j}(u)
\right)
\left(
\bar u_k\dfrac{\partial h_A}{\partial \bar z_k}(u)
+
u_k\dfrac{\partial h_A}{\partial z_k}(u)
\right)
-
\left\vert
\bar u_j\dfrac{\partial h_A}{\partial \bar z_k}(u)
+
u_k\dfrac{\partial h_A}{\partial z_j}(u)
\right\vert^2
\\
&=
\bar u_j\bar u_k 
\dfrac{\partial h_A}{\partial \bar z_j}(u)\dfrac{\partial h_A}{\partial \bar z_k}(u)
+
\bar u_j u_k 
\dfrac{\partial h_A}{\partial \bar z_j}(u)\dfrac{\partial h_A}{\partial z_k}(u)
+
u_j\bar u_k 
\dfrac{\partial h_A}{\partial z_j}(u)
\dfrac{\partial h_A}{\partial \bar z_k}(u)
+
u_j u_k
 \dfrac{\partial h_A}{\partial z_j}(u)
 \dfrac{\partial h_A}{\partial z_k}(u)
\\
&-
\left\vert \bar u_j\dfrac{\partial h_A}{\partial  z_k}(u)\right\vert^2 
-
\bar u_j \bar u_k 
\dfrac{\partial h_A}{\partial \bar z_k}(u)\dfrac{\partial h_A}{\partial \bar z_j}(u)
-
u_j u_k 
\dfrac{\partial h_A}{\partial z_j}(u)
\dfrac{\partial h_A}{\partial z_k}(u)
-
\left\vert u_k\dfrac{\partial h_A}{\partial\bar z_j}(u)\right\vert^2 
\\
&=
2\re
\left(
\bar u_j u_k
 \dfrac{\partial h_A}{\partial z_k}(u)
 \dfrac{\partial h_A}{\partial\bar z_j}(u)
 \right)
 -
\left\vert \bar u_j\dfrac{\partial h_A}{\partial  z_k}(u)\right\vert^2 
-
\left\vert u_k\dfrac{\partial h_A}{\partial  \bar z_j}(u)\right\vert^2 
\\
&=
2
\re\left(
\bar u_j\dfrac{\partial h_A}{\partial  z_k}(u)
\right)
\re\left(
u_k\dfrac{\partial h_A}{\partial \bar z_j}(u)
\right)
-
2
\im\left(
\bar u_j\dfrac{\partial h_A}{\partial  z_k}(u)
\right)
\im\left(
u_k\dfrac{\partial h_A}{\partial \bar  z_j}(u)
\right)
\\
&-
\left[
\re 
\left(\bar u_j\dfrac{\partial h_A}{\partial  z_k}(u)\right)
\right]^2
-
\left[
\im 
\left(\bar u_j\dfrac{\partial h_A}{\partial  z_k}(u)\right)
\right]^2
-
\left[
\re 
\left(u_k\dfrac{\partial h_A}{\partial \bar  z_j}(u)\right)
\right]^2
-
\left[
\im 
\left(u_k\dfrac{\partial h_A}{\partial \bar z_j}(u)\right)
\right]^2
\\
&=
-\left[
\re
\left(\bar u_j\dfrac{\partial h_A}{\partial  z_k}(u)\right)
-
\re
\left(u_k\dfrac{\partial h_A}{\partial \bar z_j}(u)\right)
\right]^2
-
\left[
\im
\left(\bar u_j\dfrac{\partial h_A}{\partial  z_k}(u)\right)
+
\im
\left(u_k\dfrac{\partial h_A}{\partial \bar z_j}(u)\right)
\right]^2
\\
&=
-\left[
\re
\left(
\bar u_j\dfrac{\partial h_A}{\partial  z_k}(u)
-
\bar u_k\dfrac{\partial h_A}{\partial z_j}(u)\right)
\right]^2
-
\left[
\im
\left(
\bar u_j\dfrac{\partial h_A}{\partial  z_k}(u)
-
\bar u_k\dfrac{\partial h_A}{\partial  z_j}(u)\right)
\right]^2
\\
&=
-\left\vert
\bar u_j\dfrac{\partial h_A}{\partial  z_k}(u)
-
\bar u_k\dfrac{\partial h_A}{\partial z_j}(u)
\right\vert^2
\\
&=
-\left\vert
u_j\dfrac{\partial h_A}{\partial \bar  z_k}(u)
-
u_k\dfrac{\partial h_A}{\partial\bar  z_j}(u)
\right\vert^2
\,,
\end{align}
it is a nonpositive real number and so is the sum of all these minors. 

The characteristic polynomial of $N_A(u)$ equals $(-\lambda)^{n-2}$ times the quadratic polynomial 
\begin{equation}
\lambda^2-h_A(u)\lambda-
\sum_{1\leqslant j<k\leqslant n}
\left\vert
u_j\dfrac{\partial h_A}{\partial\bar  z_k}(u)
-
u_k\dfrac{\partial h_A}{\partial\bar  z_j}(u)
\right\vert^2\,.
\end{equation}
By the well known Faddeev-Le Verrier algorithm, the preceding quadratic polynomial can be expressed in a more intrinsic way as
\begin{equation}
\lambda^2
-
h_A(u)\lambda
+
\dfrac{1}{2}
\left[
h_A(u)^2
-
\text{tr} \left(N_A(u)^2\right)
\right].
\end{equation}
The (ordinary) discriminant $\Delta(A,u)$ of this polynomial is positive, indeed
\begin{align}
\Delta(A,u)
&=
2\text{tr} \left(N_A(u)^2\right)-h_A(u)^2
\\
&=
h_A(u)^2+
4
\sum_{1\leqslant j<k\leqslant n}
\left\vert
u_j\dfrac{\partial h_{A_1}}{\partial\bar  z_k}(u)
-
u_k\dfrac{\partial h_{A_1}}{\partial\bar  z_j}(u)
\right\vert^2
\\
&\geqslant 
h_A(u)^2>0\,,
\end{align}
so the preceding polynomial has, for every $u\neq0$, the distinct roots 
\begin{equation}
\lambda(A,u)_1
=
\dfrac{h_A(u)}{2}+\dfrac{\sqrt{\Delta(A,u)}}{2}>0\,,
\end{equation}
\begin{equation}
\lambda(A,u)_2
=
\dfrac{h_A(u)}{2}-\dfrac{\sqrt{\Delta(A,u)}}{2}\leqslant 0\,.
\end{equation}
It follows that the hermitian matrix $N_A(u)$ is indefinite and it is positive semidefinite if and only if $\lambda(A,u)_2=0$. If $A=\varepsilon B_{2n}$, for some $\varepsilon>0$, a simple computation shows that $\lambda_1(\varepsilon B_{2n},u)=\varepsilon h_{B_{2n}}(u)$ whereas $\lambda_2(\varepsilon B_{2n},u)=0$. Indeed,
\begin{align}
u_j\dfrac{\partial h_{\varepsilon B_{2n}}}{\partial\bar  z_k}(u)
-
u_k\dfrac{\partial h_{\varepsilon B_{2n}}}{\partial\bar  z_j}(u)
=
\dfrac{u_j\varepsilon u_k}{2\Vert u\Vert}-\dfrac{u_k\varepsilon u_j}{2\Vert u\Vert}=0,
\end{align}
for every $1\leqslant j,k\leqslant n$, so $\Delta(\varepsilon B_{2n},u)=h_{\varepsilon B_{2n}}^2=\varepsilon^2h_{B_{2n}}^2$, $\lambda_1(\varepsilon B_{2n},u)=\varepsilon h_{B_{2n}}(u)$ and $\lambda_2(\varepsilon B_{2n},u)=0$. The negativity of $\lambda_2(A,u)$ provides a measures of how different is a generic $A\in{\mathcal K}_1(\mathbb C^n)$ with respect to a full-dimensional ball about the origin.

\begin{lem}\label{dj}
Let $n\geq 2$ and $A_1,A_2,\ldots,A_n\in{\mathcal K}_2(\mathbb C^n)$ with $0\in \relint A_1$, then 
\begin{equation*}
\iota^*_{\partial B_{2n}}(\dc h_{A_1}\wedge\ddc h_{A_2}\wedge\ldots\wedge \ddc h_{A_n})
=
4^{n-1}n!
\iota^*_{\partial B_{2n}}
D_n(N_{A_1},\Hess_\mathbb C h_{A_2},\ldots,\Hess_\mathbb C h_{A_n})\upsilon_{\partial B_{2n}}\,,
\end{equation*} 
In particular
\begin{equation}\label{Kmix-D-f1}
Q_n(A_1,A_2,\ldots,A_n)
=
\dfrac{4^{n-1}}{\varkappa_n}
\int_{\partial B_{2n}}
\iota^*_{\partial B_{2n}}
D_n(N_{A_1},\Hess_\mathbb C h_{A_2},\ldots,\Hess_\mathbb C h_{A_n})\upsilon_{\partial B_{2n}}\,.
\end{equation}
\end{lem}
\noindent
\pf
In the global coordinates of $\mathbb C^n$, for every $u\in\partial B_{2n}$, the standard volume form on $\partial B_{2n}$ at the point $u$ reads  
\begin{align*}
\upsilon_{\partial B_{2n}}
&=
\iota^*_{\partial B_{2n}}
\sum_{\ell=1}^n 
(\re u_\ell) \upsilon_{2n}[x_\ell]-(\im u_\ell)\upsilon_{2n}[y_\ell]
\\
&=
\iota^*_{\partial B_{2n}}
\sum_{\ell=1}^n 
u_\ell \upsilon_{2n}[z_\ell]-\bar u_\ell\upsilon_{2n}[\bar z_\ell]
\,.
\end{align*}
On each open subset of the family $\{z\in\partial B_{2n}\mid \re z_r\neq0\}$, (resp. $\{z\in\partial B_{2n}\mid \im z_r\neq0\}$), as $1\leq r\leq n$, both $\upsilon_{\partial B_{2n}}$ and $\iota^*_{\partial B_{2n}}(
\dc h_{A_1}\wedge\ddc h_{A_2}\wedge\ldots\wedge \ddc h_{A_n})$ get simpler expressions because they both become multiples of $\upsilon_{2n}[x_r]$, (resp. $\upsilon_{2n}[y_r]$). Indeed, for $u\in\partial B_{2n}$, the tangent space $T_u\partial B_{2n}$ is the zero-set of the equation
$$
\sum_{m=1}^n(\re u_m)x_m+(\im u_m)y_m=0\,,
$$ 
so if $1\leq r\leq n$ is such that $\re u_r\neq 0$ then, on $T_u\partial B_{2n}$, one has 
\begin{equation*}
\de x_r=
-
\sum_{\substack{m=1 \\ m\neq r}}^{n}\dfrac{\re u_m}{\re u_r}\de x_m
-
\sum_{m=1}^n \dfrac{\im u_m}{\re u_r}\de y_m\,,
\end{equation*}
\begin{align*}
\upsilon_{2n}[x_\ell]
&=
\de y_\ell
\wedge 
\left(-
\sum_{\substack{m=1 \\ m\neq r}}^{n}\dfrac{\re u_m}{\re u_r}\de x_m
-
\sum_{m=1}^n \dfrac{\im u_m}{\re u_r}\de y_m
\right)
\wedge
\de y_r
\wedge
\bigwedge_{\substack{j=1\\ r\neq j\neq\ell}}^n\de x_j\wedge \de y_j
\\
&=
\de y_\ell
\wedge 
\left(-
\dfrac{\re u_\ell}{\re u_r}\de x_\ell
\right)
\wedge
\de y_r
\wedge
\bigwedge_{\substack{j=1\\ r\neq j\neq\ell}}^n\de x_j\wedge \de y_j
\\
&=
\left(
\dfrac{\re u_\ell}{\re u_r}
\right)
\de y_r
\wedge
\bigwedge_{\substack{j=1\\ j\neq r}}^n\de x_j\wedge \de y_j
\\
&=
\left(
\dfrac{\re u_\ell}{\re u_r}
\right)
\upsilon_{2n}[x_r]\,,
\end{align*}
for $\ell\neq r$, and 
\begin{align*}
\upsilon_{2n}[y_\ell]
&=
\de x_\ell
\wedge 
\left(-
\sum_{\substack{m=1 \\ m\neq r}}^{n}\dfrac{\re u_m}{\re u_r}\de x_m
-
\sum_{m=1}^n \dfrac{\im u_m}{\re u_r}\de y_m
\right)
\wedge
\de y_r
\wedge
\bigwedge_{\substack{j=1\\ r\neq j\neq\ell}}^n\de x_j\wedge \de y_j
\\
&=
\de x_\ell
\wedge 
\left(-
\dfrac{\im u_\ell}{\re u_r}\de y_\ell
\right)
\wedge
\de y_r
\wedge
\bigwedge_{\substack{j=1\\ r\neq j\neq\ell}}^n\de x_j\wedge \de y_j
\\
&=
\left(
-\dfrac{\im u_\ell}{\re u_r}
\right)
\de y_r
\wedge
\bigwedge_{\substack{j=1\\ j\neq r}}^n\de x_j\wedge \de y_j
\\
&=
\left(-
\dfrac{\im u_\ell}{\re u_r}
\right)
\upsilon_{2n}[x_r]\,,
\end{align*}
for every $\ell$, whence
\begin{equation*}
\upsilon_{\partial B_{2n}}
=
\left[\re u_r
+
\dfrac{(\im u_r)^2}{\re u_r}
+
\sum_{\substack{\ell=1\\ \ell\neq r}}^n
\dfrac{(\re u_\ell)^2}{\re u_r}
+
\dfrac{(\im u_\ell)^2}{\re u_r}
\right]
\upsilon_{2n}[x_r]
=
\dfrac{1}{\re u_r}
\upsilon_{2n}[x_r]\,.
\end{equation*}
A similar computation on the open subset $\im z_r\neq 0$ yields
\begin{equation*}
\de y_r=
-
\sum_{m=1}^{n}\dfrac{\re u_m}{\im u_r}\de x_m
-
\sum_{\substack{m=1 \\ m\neq r}}^n \dfrac{\im u_m}{\im u_r}\de y_m
\,,
\end{equation*}
\begin{align*}
\upsilon_{2n}[x_\ell]
&=
\de y_\ell
\wedge 
\de x_r
\wedge
\left(-
\sum_{m=1}^{n}\dfrac{\re u_m}{\im u_r}\de x_m
-
\sum_{\substack{m=1 \\ m\neq r}}^n \dfrac{\im u_m}{\im u_r}\de y_m
\right)
\wedge
\bigwedge_{\substack{j=1\\ r\neq j\neq\ell}}^n\de x_j\wedge \de y_j
\\
&=
\de y_\ell
\wedge 
\de x_r
\wedge 
\left(-
\dfrac{\re u_\ell}{\im u_r}\de x_\ell
\right)
\wedge
\bigwedge_{\substack{j=1\\ r\neq j\neq\ell}}^n\de x_j\wedge \de y_j
\\
&=
\left(-
\dfrac{\re u_\ell}{\im u_r}
\right)
\de x_r
\wedge
\bigwedge_{\substack{j=1\\ j\neq r}}^n\de x_j\wedge \de y_j
\\
&=
\left(-
\dfrac{\re u_\ell}{\im u_r}
\right)
\upsilon_{2n}[y_r]\,,
\end{align*}
for every $\ell$, and
\begin{align*}
\upsilon_{2n}[y_\ell]
&=
\de x_\ell
\wedge 
\de x_r
\wedge
\left(-
\sum_{m=1}^{n}\dfrac{\re u_m}{\im u_r}\de x_m
-
\sum_{\substack{m=1 \\ m\neq r}}^n \dfrac{\im u_m}{\im u_r}\de y_m
\right)
\wedge
\bigwedge_{\substack{j=1\\ r\neq j\neq\ell}}^n\de x_j\wedge \de y_j
\\
&=
\de x_\ell
\wedge 
\de x_r
\wedge 
\left(-
\dfrac{\im u_\ell}{\im u_r}\de y_\ell
\right)
\wedge
\bigwedge_{\substack{j=1\\ r\neq j\neq\ell}}^n\de x_j\wedge \de y_j
\\
&=
\left(
\dfrac{\im u_\ell}{\im u_r}
\right)
\de x_r
\wedge
\bigwedge_{\substack{j=1\\ j\neq r}}^n\de x_j\wedge \de y_j
\\
&=
\left(
\dfrac{\im u_\ell}{\im u_r}
\right)
\upsilon_{2n}[y_r]\,,
\end{align*}
for $\ell\neq r$, whence 
\begin{equation*}
\upsilon_{\partial B_{2n}}
=
\left[
-\dfrac{(\re u_r)^2}{\im u_r}
-
\im u_r
-
\sum_{\substack{\ell=1\\ \ell\neq r}}^n
\dfrac{(\re u_\ell)^2}{\im u_r}
+
\dfrac{(\im u_\ell)^2}{\im u_r}
\right]
\upsilon_{2n}[y_r]
=
\dfrac{-1}{\im u_r}
\upsilon_{2n}[y_r]\,.
\end{equation*}
On the other hand, as any $(2n-1)$-form on $\mathbb C^n$, the value of the form on the point $u$
$$
\dc h_{A_1}\wedge\ddc h_{A_2}\wedge\ldots\wedge \ddc h_{A_n}
$$
can be expressed in the basis (\ref{basex}), (\ref{basey}). In fact, 
\begin{align*}
&
\dc h_{A_1}\wedge\ddc h_{A_2}\wedge\ldots\wedge \ddc h_{A_n}
\\
&=
i
\left[
\left(
\sum_{j_1=1}^n
\dfrac{\partial h_{A_1}}{\partial \bar z_{j_1}}\,\de\bar z_{j_1}
\right)
-
\left(
\sum_{\ell_1=1}^n
\dfrac{\partial h_{A_1}}{\partial z_{\ell_1}}\,\de z_{\ell_1}
\right)
\right]
\wedge
\left[
\bigwedge_{k=2}^n
2i\sum_{\ell_k,j_k=1}^n
\dfrac{\partial^2 h_{A_k}}{\partial z_{\ell_k}\partial\bar z_{j_k}}
\,\de z_{\ell_k}\wedge\de\bar z_{j_k}
\right]
\\
&=
2^{n-1}i^n
\left(
\sum_{j_1=1}^n
\dfrac{\partial h_{A_1}}{\partial \bar z_{j_1}}\,\de\bar z_{j_1}
\right)
\wedge
\sum_{\sigma,\tau\in\mathfrak{S}_n}
\left(
\prod_{k=2}^n
\dfrac{\partial^2 h_{A_k}}{\partial z_{\sigma(k)}\partial\bar z_{\tau(k)}}
\right)
\de z_{\sigma(2)}\wedge\de \bar z_{\tau(2)}
\wedge\ldots\wedge
\de z_{\sigma(n)}\wedge\de\bar z_{\tau(n)}
\\
&-
2^{n-1}i^n
\left(
\sum_{\ell_1=1}^n
\dfrac{\partial h_{A_1}}{\partial z_{\ell_1}}\,\de z_{\ell_1}
\right)
\wedge
\sum_{\sigma,\tau\in\mathfrak{S}_n}
\left(
\prod_{k=2}^n
\dfrac{\partial^2 h_{A_k}}{\partial z_{\sigma(k)}\partial\bar z_{\tau(k)}}
\right)
\de z_{\sigma(2)}\wedge\de \bar z_{\tau(2)}
\wedge\ldots\wedge
\de z_{\sigma(n)}\wedge\de\bar z_{\tau(n)}
\\
&=
2^{n-1}i^n
\sum_{\sigma,\tau\in\mathfrak{S}_n}
\left(
\dfrac{\partial h_{A_1}}{\partial \bar z_{\tau(1)}}
\prod_{k=2}^n
\dfrac{\partial^2 h_{A_k}}{\partial z_{\sigma(k)}\partial\bar z_{\tau(k)}}
\right)
\de\bar z_{\tau(1)}\wedge
\de z_{\sigma(2)}\wedge\de \bar z_{\tau(2)}
\wedge\ldots\wedge
\de z_{\sigma(n)}\wedge\de\bar z_{\tau(n)}
\\
&-
2^{n-1}i^n
\sum_{\sigma,\tau\in\mathfrak{S}_n}
\left(
\dfrac{\partial h_{A_1}}{\partial z_{\sigma(1)}}
\prod_{k=2}^n
\dfrac{\partial^2 h_{A_k}}{\partial z_{\sigma(k)}\partial\bar z_{\tau(k)}}
\right)
\de z_{\sigma(1)}\wedge
\de z_{\sigma(2)}\wedge\de \bar z_{\tau(2)}
\wedge\ldots\wedge
\de z_{\sigma(n)}\wedge\de\bar z_{\tau(n)}
\,,
\end{align*}
where all the derivatives are evaluated on $u$.
Now, by omitting the differentials in the brackets, one has 
\begin{align*}
&\de\bar z_{\tau(1)}\wedge\de z_{\sigma(2)}\wedge\de \bar z_{\tau(2)}
\wedge\ldots\wedge
\de z_{\sigma(n)}\wedge\de\bar z_{\tau(n)}
\\
&=
(-1)^{\frac{(n-1)(n-2)}{2}}
\de\bar z_{\tau(1)}\wedge
\de z_{\sigma(2)}\wedge\ldots\wedge\de z_{\sigma(n)}
\wedge
\de \bar z_{\tau(2)}\wedge\ldots\wedge\de\bar z_{\tau(n)}
\\
&=
(-1)^{\frac{(n-1)(n-2)}{2}}
(-1)^{(n-1)}
\text{sgn}(\tau)
\de z_{\sigma(2)}\wedge\ldots\wedge\de z_{\sigma(n)}
\wedge
\de\bar z_1\wedge\ldots\wedge\de\bar z_n
\\
&=
(-1)^{\frac{n(n-1)}{2}}
(-1)^{\sigma(1)-1}
\text{sgn}(\sigma\cdot\tau)
\de z_1\wedge\ldots\wedge\left[\de z_{\sigma(1)}\right]\wedge\ldots\wedge\de z_n\wedge\de\bar z_1\wedge\ldots\wedge\de\bar z_n
\\
&=
(-1)^{\frac{n(n-1)}{2}}
\text{sgn}(\sigma\cdot\tau)
\de\bar z_{\sigma(1)}\wedge
\de z_1\wedge\ldots\wedge\left[\de z_{\sigma(1)}\right]\wedge\ldots\wedge\de z_n
\wedge
\de\bar z_1\wedge\ldots\wedge\left[\de\bar z_{\sigma(1)}\right]\wedge\ldots\wedge\de\bar z_n
\\
&=
\text{sgn}(\sigma\cdot\tau)
\de\bar z_{\sigma(1)}\wedge
\bigwedge_{\substack{\ell=1\\ \ell\neq\sigma(1)}}^n
\de z_\ell\wedge\de\bar z_\ell\\
&=
(-2i)^n
\text{sgn}(\sigma\cdot\tau)
\upsilon_{2n}[z_{\sigma(1)}]
\\
&=
(-2i)^n
\text{sgn}(\sigma\cdot\tau)
\left(
\dfrac{\upsilon_{2n}[x_{\sigma(1)}]+i\upsilon_{2n}[y_{\sigma(1)}]}{2}
\right)
\end{align*}
and
\begin{align*}
&
\de z_{\sigma(1)}\wedge\de z_{\sigma(2)}\wedge\de \bar z_{\tau(2)}
\wedge\ldots\wedge
\de z_{\sigma(n)}\wedge\de\bar z_{\tau(n)}
\\
&=
(-1)^{\frac{n(n-1)}{2}}
\de z_{\sigma(1)}\wedge\de z_{\sigma(2)}\wedge\ldots\wedge\de z_{\sigma(n)}\wedge\de\bar z_{\tau(2)}\wedge\ldots\wedge\de\bar z_{\tau(n)}
\\
&=
(-1)^{\frac{n(n-1)}{2}}
\text{sgn}(\sigma)
\de z_1\wedge\ldots\wedge\de z_n\wedge\de\bar z_{\tau(2)}\wedge\ldots\wedge\de\bar z_{\tau(n)}
\\
&=
(-1)^{\frac{n(n-1)}{2}}
(-1)^{\tau(1)-1}
\text{sgn}(\sigma)
\de z_{\tau(1)}
\wedge\de z_1\ldots\wedge\left[\de z_{\tau(1)}\right]\wedge\ldots\wedge\de z_n
\wedge\de\bar z_{\tau(2)}\wedge\ldots\wedge\de\bar z_{\tau(n)}
\\
&=
(-1)^{\frac{n(n-1)}{2}}
\text{sgn}(\sigma\cdot\tau)
\de z_{\tau(1)}\wedge\de z_1\wedge\ldots\wedge\left[\de z_{\tau(1)}\right]\wedge\ldots\wedge\de z_n
\wedge
\de\bar z_1\wedge\ldots\left[\de\bar z_{\tau(1)}\right]\wedge\ldots\wedge \de\bar z_n
\\
&=
\text{sgn}(\sigma\cdot\tau)
\de z_{\tau(1)}
\wedge
\bigwedge_{\substack{\ell=1\\ \ell\neq\tau(1)}}^n
\de z_\ell\wedge\de\bar z_\ell
\\
&=
(-2i)^n
\text{sgn}(\sigma\cdot\tau)
\upsilon_{2n}[\bar z_{\tau(1)}]
\\
&=
(-2i)^n
\text{sgn}(\sigma\cdot\tau)
\left(
\dfrac{-\upsilon_{2n}[x_{\tau(1)}]+i\upsilon_{2n}[y_{\tau(1)}]}{2}
\right)
\,.
\end{align*}
If $\re u_r\neq 0$, on $T_u\partial B_{2n}$ one gets
\begin{equation*}
\upsilon_{2n}[x_{\sigma(1)}]+i\upsilon_{2n}[y_{\sigma(1)}]
=
\left(
\dfrac{\re u_{\sigma(1)}}{\re u_r}-i\dfrac{\im u_{\sigma(1)}}{\re u_r}
\right)\upsilon_{2n}[x_r]
=
\left(\dfrac{\bar u_{\sigma(1)}}{\re u_r}\right)\upsilon_{2n}[x_r]
=
\bar u_{\sigma(1)} \upsilon_{\partial B_{2n}}
\end{equation*}
and
\begin{equation*}
-\upsilon_{2n}[x_{\tau(1)}]+i\upsilon_{2n}[y_{\tau(1)}]
=
\left(
-\dfrac{\re u_{\tau(1)}}{\re u_r}-i\dfrac{\im u_{\tau(1)}}{\re u_r}
\right)\upsilon_{2n}[x_r]
=
\left(\dfrac{- u_{\tau(1)}}{\re u_r}\right) \upsilon_{2n}[x_r]
=
- u_{\tau(1)} \upsilon_{\partial B_{2n}}\,.
\end{equation*}
On the open subset $\im u_r\neq0$ one gets
\begin{equation*}
\upsilon_{2n}[x_{\sigma(1)}]+i\upsilon_{2n}[y_{\sigma(1)}]
=
\left(
-\dfrac{\re u_{\sigma(1)}}{\im u_r}+i\dfrac{\im u_{\sigma(1)}}{\im u_r}
\right)\upsilon_{2n}[y_r]
=
\left(\dfrac{-\bar u_{\sigma(1)}}{\im u_r}\right)\upsilon_{2n}[y_r]
=
\bar u_{\sigma(1)} \upsilon_{\partial B_{2n}}
\end{equation*}
and
\begin{equation*}
-\upsilon_{2n}[x_{\tau(1)}]+i\upsilon_{2n}[y_{\tau(1)}]
=
\left(
\dfrac{\re u_{\tau(1)}}{\im u_r}+i\dfrac{\im u_{\tau(1)}}{\im u_r}
\right)\upsilon_{2n}[y_r]
=
\left(\dfrac{u_{\tau(1)}}{\im u_r}\right) \upsilon_{2n}[y_r]
=
- u_{\tau(1)} \upsilon_{\partial B_{2n}}\,.
\end{equation*}
In both cases, it follows that
\begin{align*}
&
\iota^*_{\partial B_{2n}}(\dc h_{A_1}\wedge\ddc h_{A_2}\wedge\ldots\wedge \ddc h_{A_n})
\\
&=
4^{n-1}
\sum_{\sigma,\tau\in\mathfrak{S}_n}
\text{sgn}(\sigma\cdot\tau)
\left(
\bar u_{\sigma(1)}\dfrac{\partial h_{A_1}}{\partial \bar z_{\tau(1)}}
+
u_{\tau(1)}\dfrac{\partial h_{A_1}}{\partial z_{\sigma(1)}}
\right)
\left(
\prod_{k=2}^n
\dfrac{\partial^2 h_{A_k}}{\partial z_{\sigma(k)}\partial \bar z_{\tau(k)}}
\right)
\upsilon_{\partial B_{2n}}
\\
&=
4^{n-1}
n!
D_n(N_{A_1},\Hess_\mathbb C h_{A_2},\ldots,\Hess_\mathbb C h_{A_n})
\upsilon_{\partial B_{2n}}\,,
\end{align*}
where in the last equality we have used property 1 of $D_n$. 
The lemma follows by the arbitrary choice of $u\in\partial B_{2n}$. In particular, by the definition of $Q_n$, equality \eqref{Kmix-D-f1} follows at once.
\findim

\begin{rmq}
It should be remarked that, although the hermitian matrix $N_{A_1}$ figuring in~(\ref{Kmix-D-f1}) is not positive defined, the integrand $D_n(N_{A_1}(u),\Hess_\mathbb C h_{A_2}(u),\ldots,\Hess_\mathbb C h_{A_n}(u))$ is non negative at any point $u\in\partial B_{2n}$. Indeed, by virtue of~\cororef{Alex3}, there exists a positive function $a:\partial B_{2n}\to\mathbb R$ such that
\begin{equation}
D_n(N_{A_1}(u),\Hess_\mathbb C h_{A_2}(u),\ldots,\Hess_\mathbb C h_{A_n}(u))
\geqslant
a(u)\tr N_{A_1}(u)=a(u)h_{A_1}(u).
\end{equation}
The hypothesis $0\in\relint A_1$ implies the positivity of $h_{A_1}$ on  $\partial B_{2n}$.
\end{rmq}

For every $A\in{\mathcal K}_2(\mathbb C^n)$, let $Y_A(u)$ be the $n$-by-$n$ matrix whose $(\ell,j)$-entry $Y_A(u)_{\ell,j}$ is given by
\begin{equation}\label{matriceY}
\bar u_\ell\left[
\left(
\overline{\Hess_\mathbb C h_A (u)}
\right)_j \,u
\right]
+
u_j\left[
\left(
\Hess_\mathbb C h_A (u)
\right)_\ell \,\bar u
\right],
\end{equation}
where, for every $1\leqslant r\leqslant n$, by $\left(\Hess_\mathbb C h_{A_1}(u)\right)_r$ we denote the $r$-th row of the complex hessian matrix of $h_A$ at the point $u$. If 
\begin{equation}
W(u)=\bar u\; {}^tu=
\left(\begin{array}{ccc}
\bar u_1 u_1 & \ldots & \bar u_1 u_n 
\\
\vdots & \ddots & \vdots 
\\
\bar u_n u_1 & \ldots & \bar u_n u_n
\end{array}\right)
\end{equation}
it is easy to realize that $Y_A(u)={}^t[W(u)\Hess_\mathbb C h_A(u)]+\overline{W(u)\Hess_\mathbb C h_A(u)}$, so that 
\begin{align*}
{}^tY_A(u)
&=
W(u)\Hess_\mathbb C h_A(u)
+
{}^t\left[\overline{W(u)\Hess_\mathbb C h_A(u)}\right]
\\
&=
W(u)\Hess_\mathbb C h_A(u)
+
\Hess_\mathbb C h_A(u)W(u)
\\
&=
\overline{\overline{W(u)\Hess_\mathbb C h_A(u)}}
+
\overline{{}^t\left[W(u)\Hess_\mathbb C h_A(u)\right]}
\\
&=
\overline{Y_A(u)}
\end{align*}

i.e. $Y_A(u)$ is hermitian.

\begin{lem}\label{djj}
Let $n\geq 2$ and $A_1,A_2,\ldots,A_n\in{\mathcal K}^2_2(\mathbb C^n)$ with $0\in \relint A_1$, then 
\begin{equation}\label{Kmix-D-f22}
Q_n(A_1,A_2,\ldots,A_n)
=
\dfrac{4^{n-1}}{\varkappa_n}
\int_{\partial B_{2n}}
\iota^*_{\partial B_{2n}}
D_n(Y_{A_1},\Hess_\mathbb C h_{A_2},\ldots,\Hess_\mathbb C h_{A_n})\upsilon_{\partial B_{2n}}\,.
\end{equation}
\end{lem}
\noindent
\pf By \lemref{koush}, the form $\dc h_{A_1}$ is cohomologous to the form $(\nu_{\partial A_1}^{-1})^*\alpha_{\partial A_1}$, which by (\ref{koushnirenko}) equals
\begin{equation}\label{coomologa}
\dfrac{i}{2}
\sum_{j=1}^n
\left[
\sum_{\ell=1}^n
u_\ell\dfrac{\partial^2 h_{A_1}}{\partial \bar z_j\partial z_\ell}(u)
-
\bar u_\ell\dfrac{\partial^2 h_{A_1}}{\partial \bar z_j\partial \bar z_\ell}(u)
\right]
\de\bar z_j
+
\left[
\sum_{\ell=1}^n
u_\ell\dfrac{\partial^2 h_{A_1}}{\partial  z_j\partial z_\ell}(u)
-
\bar u_\ell\dfrac{\partial^2 h_{A_1}}{\partial  z_j\partial \bar z_\ell}(u)
\right]
\de z_j.
\end{equation}
Euler's relation for $h_{A_1}$ in complex coordinates simply reads
\begin{equation}
h_{A_1}(z)=\sum_{\ell=1}^n z_\ell \dfrac{\partial h_{A_1}}{\partial z_\ell}(z)+
\bar z_\ell \dfrac{\partial h_{A_1}}{\partial \bar z_\ell}(z).
\end{equation}
A differentiation with respect to $\bar z_j$ or $z_j$, $1\leqslant j\leqslant n$, followed by an evaluation on the point $u$ yields
\begin{equation}
\sum_{\ell=1}^n
u_\ell\dfrac{\partial^2 h_{A_1}}{\partial \bar z_j\partial z_\ell}(u)
+
\bar u_\ell\dfrac{\partial^2 h_{A_1}}{\partial \bar z_j\partial \bar z_\ell}(u)
=0
\end{equation}
or
\begin{equation}
\sum_{\ell=1}^n
u_\ell\dfrac{\partial^2 h_{A_1}}{\partial z_j\partial z_\ell}(u)
+
\bar u_\ell\dfrac{\partial^2 h_{A_1}}{\partial z_j\partial \bar z_\ell}(u)
=0,
\end{equation}
then (\ref{coomologa}) becomes
\begin{align}
&
\dfrac{i}{2}
\sum_{j=1}^n
\left[
2
\sum_{\ell=1}^n
u_\ell\dfrac{\partial^2 h_{A_1}}{\partial \bar z_j\partial z_\ell}(u)
\right]
\de \bar z_j
+
\left[
-2
\sum_{\ell=1}^n
\bar u_\ell\dfrac{\partial^2 h_{A_1}}{\partial  z_j\partial \bar z_\ell}(u)
\right]
\de z_j
\\
&=
i
\sum_{j=1}^n
\left[
\sum_{\ell=1}^n
u_\ell\dfrac{\partial^2 h_{A_1}}{\partial \bar z_j\partial z_\ell}(u)
\right]
\de \bar z_j
-
\left[
\sum_{\ell=1}^n
\bar u_\ell\dfrac{\partial^2 h_{A_1}}{\partial  z_j\partial \bar z_\ell}(u)
\right]
\de z_j
\\
&=
i
\sum_{j_1=1}^n
\left[
\left(\overline{\Hess_\mathbb C h_{A_1}(u)}\right)_{j_1}\,u
\right]
\de \bar z_{j_1}
-
i
\sum_{\ell_1=1}^n
\left[
\left(\Hess_\mathbb C h_{A_1}(u)\right)_{\ell_1}\,\bar u
\right]
\de z_{\ell_1}.
\end{align}
It follows that the form $\dc h_{A_1}\wedge \ddc h_{A_2}\wedge\ldots\wedge\ddc h_{A_n}$ is cohomologous to the form $(\nu_{\partial A_1}^{-1})^*\alpha_{\partial A_1}\wedge \ddc h_{A_2}\wedge\ldots\wedge\ddc h_{A_n}$ and the latter can be computed as the former in the proof of \lemref{dj}:
\begin{align}
&
i
\sum_{j_1=1}^n
\left[
\left(\overline{\Hess_\mathbb C h_{A_1}(u)}\right)_{j_1}\,u
\right]
\de \bar z_{j_1}
\wedge
\left[
\bigwedge_{k=2}^n
2i
\sum_{\ell_k,j_k=1}^n
\dfrac{\partial^2 h_{A_k}}{\partial z_{\ell_k}\partial\bar z_{j_k}}
\,\de z_{\ell_k}\wedge\de\bar z_{j_k}
\right]
\\
&
-
i
\sum_{\ell_1=1}^n
\left[
\left(\Hess_\mathbb C h_{A_1}(u)\right)_{\ell_1}\,\bar u
\right]
\de z_{\ell_1}
\wedge
\left[
\bigwedge_{k=2}^n
2i\sum_{\ell_k,j_k=1}^n
\dfrac{\partial^2 h_{A_k}}{\partial z_{\ell_k}\partial\bar z_{j_k}}
\,\de z_{\ell_k}\wedge\de\bar z_{j_k}
\right]
\\
&=
2^{n-1}i^n
\sum_{\sigma,\tau\in\mathfrak{S}_n}
\left[
\left[
\left(\overline{\Hess_\mathbb C h_{A_1}(u)}\right)_{\tau(1)}\,u
\right]
\prod_{k=2}^n
\dfrac{\partial^2 h_{A_k}}{\partial z_{\sigma(k)}\partial\bar z_{\tau(k)}}
\right]
\de\bar z_{\tau(1)}\wedge
\bigwedge_{k=2}^n\de z_{\sigma(k)}\wedge\de \bar z_{\tau(k)}
\\
&-
2^{n-1}i^n
\sum_{\sigma,\tau\in\mathfrak{S}_n}
\left[
\left[
\left(\Hess_\mathbb C h_{A_1}(u)\right)_{\sigma(1)}\,\bar u
\right]
\prod_{k=2}^n
\dfrac{\partial^2 h_{A_k}}{\partial z_{\sigma(k)}\partial\bar z_{\tau(k)}}
\right]
\de z_{\sigma(1)}\wedge
\bigwedge_{k=2}^n\de z_{\sigma(k)}\wedge\de \bar z_{\tau(k)}
\\
&=
4^{n-1}
\sum_{\sigma,\tau\in\mathfrak{S}_n}
\text{sgn}(\sigma\cdot\tau)
Y_{A_1}(u)_{\sigma(1),\tau(1)}
\left(
\prod_{k=2}^n
\dfrac{\partial^2 h_{A_k}}{\partial z_{\sigma(k)}\partial \bar z_{\tau(k)}}
\right)
\upsilon_{\partial B_{2n}}
\\
&=
4^{n-1}
n!
D_n\left(Y_{A_1}(u),\Hess_\mathbb C h_{A_2}(u),\ldots,\Hess_\mathbb C h_{A_n}(u)\right)
\upsilon_{\partial B_{2n}}.
\end{align}
As a consequence one obtains (\ref{Kmix-D-f22}).
\findim

\color{black}
Let us set, for $A_2,\ldots,A_n\in{\mathcal K}_2(\mathbb C^n)$, $u\in\mathbb C^n\setminus\{0\}$ and $1\leqslant j,k\leqslant n$,
\begin{equation}
Z_{j,k}(A_2,\ldots,A_n)(u)
=
\dfrac{(-1)^{j+k}}{n} D_{n-1}\left(\Hess_\mathbb C h_{A_2}^{[j,k]}(u),\ldots,\Hess_\mathbb C h_{A_n}^{[j,k]}(u)\right)
\end{equation}
and 
\begin{equation}
\Xi_{j,k}(A_2,\ldots,A_n)(u)
=
\begin{cases}
Z_{j,k}(A_2,\ldots,A_n)(u),& \text{if } j\neq k\,,
\\~\\
-
\sum\limits_{\substack{\ell=1\\ \ell\neq j}}^n 
Z_{\ell,\ell}(A_2,\ldots,A_n)(u),& \text{if } j= k\,.
\end{cases}
\end{equation}
The $n$-by-$n$ hermitian matrix $\Xi(A_2,\ldots,A_n)(u)$\label{matriceXi} whose $(j,k)$-entry is $\Xi_{j,k}(A_2,\ldots,A_n)(u)$ will be helpful in the sequel.
\begin{lem}\label{lemma-Xi}
Let $n\geq 2$, $A_1,A_2,\ldots,A_n\in{\mathcal K}_2(\mathbb C^n)$ and $u\in\mathbb C^n\setminus\{0\}$. Then
\begin{align}
&D_n(N_{A_1},\Hess_\mathbb C h_{A_2}(u),\ldots,\Hess_\mathbb C h_{A_n}(u))
\\
&=
 h_{A_1}(u) D_n\left( I_n,\Hess_\mathbb C h_{A_2}(u),
\ldots,
\Hess_\mathbb C h_{A_n}(u)
\right)
+\re\langle\nabla h_{A_1}(u),\Xi(A_2,\ldots,A_n)(u)\,u\rangle.
\end{align}
\end{lem}

\noindent
\pf 
As $\text{tr}\left(N_{A_1}(u)\right)=h_{A_1}(u)$, if all the derivatives are computed at $u$, one has
\begin{align}
&D_n(N_{A_1},\Hess_\mathbb C h_{A_2},\ldots,\Hess_\mathbb C h_{A_n})
\\
&=
D_n(h_{A_1}I_n+(N_{A_1}-h_{A_1}I_n),\Hess_\mathbb C h_{A_2},\ldots,\Hess_\mathbb C h_{A_n})
\\
&=
h_{A_1}
D_n(I_n,\Hess_\mathbb C h_{A_2},\ldots,\Hess_\mathbb C h_{A_n})
+
D_n(N_{A_1}-h_{A_1}I_n,\Hess_\mathbb C h_{A_2},\ldots,\Hess_\mathbb C h_{A_n}).
\end{align}
For the sake of notation, let $Z_{j,k}=Z_{j,k}(A_2,\ldots,A_n)(u)$, $\Xi_{j,k}=\Xi_{j,k}(A_2,\ldots,A_n)(u)$ and $\Xi=\Xi(A_2,\ldots,A_n)(u)$.
By (\ref{comoda})
\begin{align}
&D_n(N_{A_1}-h_{A_1}I_n,\Hess_\mathbb C h_{A_2},\ldots,\Hess_\mathbb C h_{A_n})
\\
&=
-
\sum_{j=1}^n
Z_{j,j}
\sum_{\substack{k=1\\ k\neq j}}^n
\left(
\bar u_k\dfrac{\partial h_{A_1}}{\partial \bar z_k}
+
u_k\dfrac{\partial h_{A_1}}{\partial z_k}
\right)
+
\sum_{\substack{1\leqslant j,k\leqslant n\\
j\neq k}}
Z_{j,k}
\left(
\bar u_j\dfrac{\partial h_{A_1}}{\partial \bar z_k}
+
u_k\dfrac{\partial h_{A_1}}{\partial z_j}
\right).\label{ultima}
\end{align}
Now observe that
\begin{align}
-\sum_{j=1}^n
Z_{j,j}
\sum_{\substack{k=1\\ k\neq j}}^n
\left(
\bar u_k\dfrac{\partial h_{A_1}}{\partial \bar z_k}
+
u_k\dfrac{\partial h_{A_1}}{\partial z_k}
\right)
&=
-\sum_{k=1}^n
\left(
\bar u_k\dfrac{\partial h_{A_1}}{\partial \bar z_k}
+
u_k\dfrac{\partial h_{A_1}}{\partial z_k}
\right)
\sum_{\substack{j=1\\ j\neq k}}^n
Z_{j,j}
\\
&=
2\re
\sum_{k=1}^n
\dfrac{\partial h_{A_1}}{\partial \bar z_k}
\left(
\sum_{\substack{j=1\\ j\neq k}}^n
-\overline{Z_{j,j}}
\right)
\bar u_k
\\
&=
2\re
\sum_{k=1}^n
\dfrac{\partial h_{A_1}}{\partial \bar z_k}
\overline{\Xi_{k,k}}
\bar u_k\label{termine1}
\end{align}
and
\begin{align}
\sum_{\substack{1\leqslant j,k\leqslant n\\
j\neq k}}
\left(
\bar u_j\dfrac{\partial h_{A_1}}{\partial \bar z_k}
+
u_k\dfrac{\partial h_{A_1}}{\partial z_j}
\right)
Z_{j,k}
&=
\sum_{k=1}^n
\dfrac{\partial h_{A_1}}{\partial \bar z_k}
\left(
\sum_{\substack{j=1\\
j\neq k}}^n
Z_{j,k}
\bar u_j
\right)
+
\sum_{j=1}^n
\dfrac{\partial h_{A_1}}{\partial z_j}
\left(
\sum_{\substack{k=1\\
k\neq j}}^n
Z_{j,k}
u_k
\right)
\\
&=
\sum_{k=1}^n
\dfrac{\partial h_{A_1}}{\partial \bar z_k}
\left(
\sum_{\substack{j=1\\
j\neq k}}^n
\overline{Z_{k,j}}
\bar u_j
\right)
+
\sum_{k=1}^n
\dfrac{\partial h_{A_1}}{\partial z_k}
\left(
\sum_{\substack{j=1\\
j\neq k}}^n
Z_{k,j}
u_j
\right)
\\
&=
2\re
\sum_{k=1}^n
\dfrac{\partial h_{A_1}}{\partial \bar z_k}
\left(
\sum_{\substack{j=1\\
j\neq k}}^n
\overline{Z_{k,j}}
\bar u_j
\right)
\\
&=
2\re
\sum_{k=1}^n
\dfrac{\partial h_{A_1}}{\partial \bar z_k}
\left(
\sum_{\substack{j=1\\
j\neq k}}^n
\overline{\Xi_{k,j}}
\bar u_j
\right)\label{termine2}
\,.
\end{align}
It follows that $
D_n(N_{A_1}-h_{A_1}I_n,\Hess_\mathbb C h_{A_2},\ldots,\Hess_\mathbb C h_{A_n})
=
\re\langle\nabla h_{A_1},\Xi u\rangle$.
\findim

\begin{coro}\label{corobuono}
Let $A_1,\ldots,A_n\in{\mathcal K}_2(\mathbb C^n)$ such that 
$$
\int_{\partial B_{2n}}\re\langle\nabla h_{A_1}(u),\Xi(A_2,\ldots, A_n)(u)\, u\rangle\upsilon_{\partial B_{2n}}=0\,,
$$ 
then
\begin{equation}\label{formula-buona}
Q_n(A_1,\ldots,A_n)
=
\dfrac{4^{n-1}}{\varkappa_n}
\int_{\partial B_{2n}}
h_{A_1}
D_n(I_n,\Hess_\mathbb C h_{A_2},\ldots,\Hess_\mathbb C h_{A_n})
\upsilon_{\partial B_{2n}}\,.
\end{equation}
\end{coro}

\noindent
\pf The statement is a direct consequence of \lemref{lemma-Xi}.
\findim


\begin{lem}\label{lemma-RXi}
Let $A_1,A_2,\ldots,A_n\in{\mathcal K}_2(\mathbb R^n)$ and $u\in\mathbb C^n\setminus\{0\}$. Then 
\begin{equation}
\Xi(A_2,\ldots,A_n)(u)\re u=\re\left[\Xi(A_2,\ldots,A_n) (u)u\right]=0\,,
\end{equation}
in particular $\re\langle\nabla h_{A_1},\Xi(A_2,\ldots,A_n)(u)\, u\rangle=0$.
\end{lem}

\noindent
\pf
For every $2\leqslant r\leqslant n$, since $A_r\subset\mathbb R^n$, $h_{A_r}$ does not depend on $y_1,\ldots,y_n$, so that $4\Hess_\mathbb C h_{A_r}=\Hess_\mathbb R h_{A_r}$, $\Xi(A_2,\ldots,A_n)$ is a real matrix and so $\Xi(A_2,\ldots,A_n)(u)\re u=\re[\Xi(A_2,\ldots,A_n) (u)u]$. If $\re u=0$ the lemma is trivial, so let us suppose $\re u\neq 0$. We now prove that, for every  $1\leqslant j\leqslant n$, $(\Xi(A_2,\ldots,A_n)(u)\re u)_j=0$. Let us first  consider the case $A_2=\ldots=A_n=A$. In this case  
\begin{align}
&
(\Xi(A,\ldots,A)(u)\re u)_j
\\
&=
\Xi_{j,j}(A,\ldots,A)\re u_j+\sum_{\substack{k=1\\ k\neq j}}^n \Xi_{j,k}(A,\ldots,A)\re u_k
\\
&=
\sum_{\substack{k=1\\ k\neq j}}^n 
Z_{j,k}(A,\ldots,A)\re u_k- Z_{k,k}(A,\ldots,A)\re u_j
\\
&=
\dfrac{1}{4^n n}
\sum_{\substack{k=1\\ k\neq j}}^n 
(-1)^{j+k}
\re u_k
\det\Hess_\mathbb R h_A^{[j,k]}
-
\re u_j
\det\Hess_\mathbb R h_A^{[k,k]}.
\end{align}
For every $k\neq j$, the matrix $\Hess_\mathbb R h_A^{[j,k]}$ admits the row 
\begin{equation}
R=\left(
\dfrac{\partial^2 h}{\partial x_k\partial x_1},
\ldots,
\dfrac{\partial^2 h}{\partial x_k\partial x_{k-1}}
,
\dfrac{\partial^2 h}{\partial x_k\partial x_{k+1}}
,\ldots,
\dfrac{\partial^2 h}{\partial x_k\partial x_n}
\right)
\end{equation}
so let $M_{j,k}$ be the matrix obtained from $\Hess_\mathbb R h^{[j,k]}$ by replacing $R$ with $R \re u_k$ and interchanging this new row with the first row of $\Hess_\mathbb R h^{[j,k]}$. 
In the same manner, the matrix $\Hess_\mathbb R h_A^{[k,k]}$ admits the row 
\begin{equation}
S=
\left(\dfrac{\partial^2 h}{\partial x_j\partial x_1},
\ldots,
\dfrac{\partial^2 h}{\partial x_j\partial x_{k-1}}
,
\dfrac{\partial^2 h}{\partial x_j\partial x_{k+1}}
,\ldots,
\dfrac{\partial^2 h}{\partial x_j\partial x_n}
\right)
\end{equation}
so let $M_{k,k}$ be the matrix obtained from $\Hess_\mathbb R h^{[k,k]}$ by replacing $S$ with $S\re u_j$ and interchanging this new row with the first row of $\Hess_\mathbb R h^{[k,k]}$. 
By construction the matrices $M_{j,k}$ and $M_{k,k}$ differ in their first rows only,
\begin{equation}
\det M_{j,k}=
\begin{cases}
(-1)^{k-1}\re u_k\det \Hess_\mathbb R h^{[j,k]}, &\text{if } k<j,
\\
(-1)^{k-2}\re u_k\det \Hess_\mathbb R h^{[j,k]}, &\text{if } k>j
\end{cases}
\end{equation}
whereas
\begin{equation}
\det M_{k,k}=
\begin{cases}
(-1)^{j-2}\re u_k\det \Hess_\mathbb R h^{[k,k]},&\text{if } k<j,
\\
(-1)^{j-1}\re u_k\det \Hess_\mathbb R h^{[k,k]},&\text{if } k>j.
\end{cases}
\end{equation}
For $k\neq j$, let $L_{j,k}$ be the matrix whose first row is the sum of the first rows of $M_{j,k}$ and $M_{k,k}$ and whose remaining rows coincide with those of $M_{j,k}$ (or $M_{k,k}$); it follows that
\begin{align}
(\Xi(A,\ldots,A)(u)\re u)_j
&=
(-1)^{j-1}
\sum_{k=1}^{j-1}
\det M_{j,k}+\det M_{k,k}
+
(-1)^{j}
\sum_{k=j+1}^n
\det M_{j,k}+\det M_{k,k}
\\
&=
(-1)^{j-1}
\sum_{k=1}^{j-1}
\det L_{j,k}
+
(-1)^{j}
\sum_{k=j+1}^n
\det L_{j,k}.
\end{align}
For every $k\neq j$, $\det L_{j,k}=0$ because the first row of $L_{jk}$ is a linear combination of the other ones. Indeed, by \lemref{Eulero}
\begin{equation}
\sum_{s=1}^n
\dfrac{\partial h_A}{\partial x_s\partial x_r}(u) \re u_s=0 
\end{equation} 
and it is enough to write the preceding equality as
\begin{equation}
\dfrac{\partial h_A}{\partial x_k\partial x_r}(u) \re u_k
+
\dfrac{\partial h_A}{\partial x_j\partial x_r}(u) \re u_j
=
-\sum_{\substack{s=1\\s\neq j\\ s\neq k}}^n
\dfrac{\partial h_A}{\partial x_s\partial x_r}(u) \re u_s
\end{equation}
to realize that the first row of $L_{j,k}$ is a linear combination of the other ones.
We can now pass to the general case. Let $t_2,\ldots,t_n$ be positive real numbers and set $A=\sum_{r=2}^n t_r A_r$. Then 
\begin{equation}
\Hess_\mathbb R h_A=\sum_{r=2}^n t_r \Hess_\mathbb R h_{A_r}
\end{equation}
whence
\begin{align}
\Xi_{j,k}(A_2,\ldots,A_n)
&=
\dfrac{(-1)^{j+k}}{n 4^n}
D_{n-1}
\left(
\Hess_\mathbb R h_{A_2}^{[j,k]},\ldots,\Hess_\mathbb R h_{A_n}^{[j,k]}
\right)
\\
&=
\dfrac{(-1)^{j+k}}{n! 4^n}
\dfrac{\partial}{\partial t_2\cdots\partial t_n}
\det\left(\sum_{r=2}^n t_r \Hess_\mathbb R h_{A_r}^{[j,k]}\right)
\\
&=
\dfrac{(-1)^{j+k}}{n! 4^n}
\dfrac{\partial}{\partial t_2\cdots\partial t_n}
\det\left(\Hess_\mathbb R h_A^{[j,k]}\right),
\end{align}
for every $k\neq j$, and
\begin{align}
\Xi_{j,j}(A_2,\ldots,A_n)
&=
-\dfrac{1}{n 4^n}
\sum_{\substack{\ell=1\\ \ell\neq j}}^n
D_{n-1}
\left(
\Hess_\mathbb R h_{A_2}^{[\ell,\ell]},\ldots,\Hess_\mathbb R h_{A_n}^{[\ell,\ell]}
\right)
\\
&=
-\dfrac{1}{n! 4^n}
\sum_{\substack{\ell=1\\ \ell\neq j}}^n
\dfrac{\partial}{\partial t_2\cdots\partial t_n}
\det\left(\sum_{r=2}^n t_r \Hess_\mathbb R h_{A_r}^{[\ell,\ell]}\right)
\\
&=
-\dfrac{1}{n! 4^n}
\dfrac{\partial}{\partial t_2\cdots\partial t_n}
\sum_{\substack{\ell=1\\ \ell\neq j}}^n
\det\left(\Hess_\mathbb R h_A^{[\ell,\ell]}\right),
\end{align}
so that 
\begin{align}
&
(\Xi(A_2,\ldots,A_n)(u)\re u)_j
\\
&=
\Xi_{j,j}(A_2,\ldots,A_n)\re u_j+\sum_{\substack{k=1\\ k\neq j}}^n \Xi_{j,k}(A_2,\ldots,A_n)\re u_k
\\
&=
\dfrac{1}{4^n n!}
\dfrac{\partial}{\partial t_2\cdots\partial t_n}
\sum_{\substack{k=1\\ k\neq j}}^n 
(-1)^{j+k}
\re u_k
\det\Hess_\mathbb R h_A^{[j,k]}
-
\re u_j
\det\Hess_\mathbb R h_A^{[k,k]}.
\end{align}
By the first part of the proof, each summand in the preceding expression is zero. The last statement is a consequence of the equality
\begin{equation}
\nabla h_{A_1}=\left(
\dfrac{\partial h_{A_1}}{\partial x_1},\ldots,\dfrac{\partial h_{A_1}}{\partial x_n}\right),
\end{equation}
indeed 
\begin{align}
\re\langle\nabla h_{A_1}(u),\Xi(A_2,\ldots,A_n)(u)\,u\rangle
&=
\sum_{j=1}^n\dfrac{\partial h_{A_1}}{\partial x_j}(u)(\Xi(A_2,\ldots,A_n)(u)\,\re u)_j=0
\end{align}
since $\Xi(A_2,\ldots,A_n)(u)\,\re u=0$.
\findim 

We end the section by proving one more formula for the mixed volume on ${\mathcal K}_2(\mathbb R^n)$ in the vein of \cororef{corobuono}.

\begin{lem}
For every $A_1\in{\mathcal K}_1(\mathbb R^n)$ and every $A_2,\ldots,A_n\in{\mathcal K}_2(\mathbb R^n)$,
\begin{align}\label{nuova}
V_n(A_1,\ldots,A_n)
=
\dfrac{1}{\varkappa_n}
\int_{\partial B_{2n}}
h_{A_1}
D_n(I_n,\Hess_\mathbb R h_{A_2},\ldots,\Hess_\mathbb R h_{A_n})\upsilon_{\partial B_{2n}}.
\end{align}
\end{lem}

\noindent
\pf 
By (\ref{mix-D-f}) it's enough to prove that the right-hand side of (\ref{nuova}) equals
\begin{equation}
\int_{\partial B_{n}}
h_{A_1}
D_n(I_n,\Hess_\mathbb R h_{A_2},\ldots,\Hess_\mathbb R h_{A_n})\upsilon_{\partial B_{n}},\label{intreale}
\end{equation}
but if $A_2,\ldots,A_n\in{\mathcal K}_3(\mathbb R^n)$ this is a consequence of the following \lemref{nonni}. The case with $A_2,\ldots,A_n\in{\mathcal K}_2(\mathbb R^n)$ follows by continuity. \findim

\begin{lem}\label{nonni}
Let $f\in{\mathcal C}^1(\mathbb R^n\setminus\{0\})$ a positively homogeneous function of degree $(2-n)$. Then
\begin{equation}
\dfrac{1}{\varkappa_n}
\int_{\partial B_{2n}}
f\upsilon_{\partial B_{2n}}
=
\int_{\partial B_n}
f\upsilon_{\partial B_n}.
\end{equation} 
\end{lem}

\noindent
\pf
Looking at $f$ as a function on $\mathbb R^n\setminus\{0\}+i\mathbb R^n$ depending on the real variables $x_1,\ldots,x_n$ only, the condition on the degree of $f$ imply $f\in L^1(B_n)$ and  $f\in L^1(B_{2n})$. By Stokes' theorem, for every $0<\varepsilon\leqslant 1$, 
\begin{equation}
\int_{\partial B_{2n}}
f\upsilon_{\partial B_{2n}}
=
\int_{B_{2n}\setminus\varepsilon B_{2n}}
\de(f\upsilon_{\partial B_{2n}})
+
\int_{\partial (\varepsilon B_{2n})}
f\upsilon_{\partial (\varepsilon B_{2n})}
\end{equation} 
and, as $\varepsilon\to0^+$, one gets
\begin{equation}
\int_{\partial B_{2n}}
f\upsilon_{\partial B_{2n}}
=
\int_{B_{2n}}
\de(f\upsilon_{\partial B_{2n}}).
\end{equation} 
Recall that, with respect to the basis (\ref{basex}) and (\ref{basey}), $\upsilon_{\partial B_{2n}}=\sum_{\ell=1}^n
x_\ell\upsilon_{2n}[x_\ell]-y_\ell\upsilon_{2n}[y_\ell]$, with $\de\upsilon_{2n}[x_\ell]=\de\upsilon_{2n}[y_\ell]=0$ for every $1\leqslant \ell\leqslant n$, so
\begin{align}
\de(f\upsilon_{\partial B_{2n}})
&=
\sum_{\ell=1}^n
\de (f x_\ell \upsilon_{2n}[x_\ell])-\de  (f y_\ell \upsilon_{2n}[y_\ell])
\\
&=
\sum_{\ell=1}^n
\de(fx_\ell)\wedge \upsilon_{2n}[x_\ell]
-
\de(fy_\ell)\wedge\upsilon_{2n}[y_\ell].
\end{align}
As $f$ does not depend on $y_1,\ldots,y_n$, 
\begin{align}
\sum_{\ell=1}^n
\de(fx_\ell)\wedge \upsilon_{2n}[x_\ell]
&=
\sum_{\ell=1}^n
\left[
\sum_{k=1}^n
\left(
\dfrac{\partial f}{\partial x_k}x_\ell+f\delta_{k,\ell}
\right)
\de x_k
\right]
\wedge\upsilon_{2n}[x_\ell]
\\
&=
\sum_{\ell=1}^n
\sum_{k=1}^n
\left(
\dfrac{\partial f}{\partial x_k}x_\ell+f\delta_{k,\ell}
\right)
\delta_{k,\ell}\,
\upsilon_{2n}
\\
&=
\left[
\sum_{\ell=1}^n
\dfrac{\partial f}{\partial x_\ell}x_\ell+f
\right]
\upsilon_{2n}
\\
&=
[(2-n)f+nf]\upsilon_{2n}\label{Eu1}
\\
&=
2f\upsilon_{2n}
\end{align}
where the equality (\ref{Eu1}) is a consequence of Euler's theorem on homogeneous functions, and
\begin{align}
\sum_{\ell=1}^n
\de(fy_\ell)\wedge \upsilon_{2n}[y_\ell]
&=
\sum_{\ell=1}^n
\left[
\sum_{k=1}^n
\dfrac{\partial f}{\partial x_k}y_\ell
\de x_k
\right]
\wedge\upsilon_{2n}[y_\ell]
+
\left[
\sum_{k=1}^n
\left(
\dfrac{\partial f}{\partial y_k}y_\ell+f\delta_{k,\ell}
\right)
\de y_k
\right]
\wedge\upsilon_{2n}[y_\ell]
\\
&=
\sum_{\ell=1}^n
\left[
\sum_{k=1}^n
f\delta_{k,\ell}
\right]
(-\delta_{k,\ell})\upsilon_{2n}
\\
&=
-nf\upsilon_{2n},
\end{align}
so that $\de(f\upsilon_{\partial B_{2n}})=(n+2)f\upsilon_{2n}$. We now pass to polar coordinates in $\mathbb R^n$ only by setting
\begin{equation}
x_\ell=\rho\cos\vartheta_\ell\prod_{j=1}^{\ell-1}\sin\vartheta_j
\quad
\text{for every}
\quad
1\leqslant \ell\leqslant n-1
\quad
\text{and}
\quad 
x_n=\rho\prod_{j=1}^{n-1}\sin\vartheta_j,
\end{equation}
with $\rho\geqslant 0$, $\vartheta_1,\ldots,\vartheta_{n-2}\in[0,\pi)$ and $\vartheta_{n-1}\in[0,2\pi)$. In these coordinates the unit ball $B_{2n}$ is given by $\rho^2+\Vert y\Vert^2\leqslant 1$ and the standard volume form $\upsilon_{2n}$ becomes
\begin{equation}
\upsilon_{2n}
=
\rho^{n-1}
J(\vartheta)\de\rho\wedge\de y_1\wedge\bigwedge_{\ell=1}^{n-1}\de\vartheta_\ell\wedge\de y_{\ell+1}
\end{equation}
with 
\begin{equation}
J(\vartheta)=
\prod_{j=1}^{n-2}
(\sin\vartheta_j)^{n-j-1}.
\end{equation}
If $\tilde f$ is the expression of $f$ in the new coordinates, remark that by homogeneity $\tilde f(\rho,\vartheta)=\rho^{2-n}\tilde f(1,\vartheta)$, remark also that $J(\vartheta)\de\vartheta$ is the standard volume form on an $(n-1)$-dimensional sphere. It follows that
\begin{align}
\dfrac{1}{\varkappa_n}
\int_{\partial B_{2n}}
f\upsilon_{\partial B_{2n}}
&=
\dfrac{n+2}{\varkappa_n}
\int_{B_{2n}}
\tilde f\upsilon_{2n}
\\
&=
\dfrac{n+2}{\varkappa_n}
\int_{B_{2n}}
\rho f(1,\vartheta)J(\vartheta)\de\rho\wedge\de y_1\wedge\bigwedge_{\ell=1}^{n-1}\de\vartheta_\ell\wedge\de y_{\ell+1}.
\end{align}
If $\de\vartheta=\de\vartheta_1\wedge\ldots\wedge\de\vartheta_{n-1}$ and $\de y=\de y_1\wedge\ldots\wedge\de y_n$, by Fubini's theorem
\begin{align}
\dfrac{1}{\varkappa_n}
\int_{\partial B_{2n}}
f\upsilon_{\partial B_{2n}}
&=
\dfrac{n+2}{\varkappa_n}
\int_0^1
\rho\left[
\int_{[0,\pi]^{n-2}\times[0,2\pi]}
\tilde f(1,\vartheta)J(\vartheta)
\left(
\int_{\{y\in\mathbb R^n\mid\Vert y\Vert^2\leqslant 1-\rho^2\}}
\de y
\right)
\de\vartheta
\right]
\de\rho
\\
&=
\dfrac{n+2}{\varkappa_n}
\int_0^1
\rho\left[
\int_{[0,\pi]^{n-2}\times[0,2\pi]}
\tilde f(1,\vartheta)J(\vartheta)
\left(\sqrt{1-\rho^2}\right)^n\varkappa_n
\de\vartheta
\right]
\de\rho
\\
&=
(n+2)
\left[\int_0^1
\rho
\left(\sqrt{1-\rho^2}\right)^n
\de\rho
\right]
\int_{[0,\pi]^{n-2}\times[0,2\pi]}
\tilde f(1,\vartheta)J(\vartheta)
\de\vartheta
\\
&=
(n+2)
\left[
\dfrac{1}{n+2}
\right]
\int_{[0,\pi]^{n-2}\times[0,2\pi]}
\tilde f(1,\vartheta)J(\vartheta)
\de\vartheta
\\
&=
\int_{\partial B_n}
f\upsilon_{\partial B_n}.
\end{align}
The lemma is thus proved.\findim

\section{Examples}\label{esempi}

This section provides the explicit computation of the pseudovolume for some convex bodies of $\mathbb C^n$. If $A\in{\mathcal K}(\mathbb C^n)$, its pseudovolume $P_n(A)$ does not merely depend on the shape of $A$ but it also depends on the position of $A$ with respect to the complex structure. Though it would be interesting to compute the extrema of $P_n(R(A))$, as $R$ runs in the group $\text{O}(2n)$ of orthogonal transformations of $\mathbb R^{2n}$, we will content ourselves to compute pseudovolumes in the simplest conditions. We just mention that for a a fulldimensional body $A\in{\mathcal K}(\mathbb C^n)$ the minimum of $P_n(R(A))$, for $R\in\text{O}(2n)$, can be strictly positive.

\begin{exe}\label{pseudoball}
The $n$-pseudovolume of the unit full-dimensional ball $B_{2n}$ of $\mathbb C^n$.

\noindent
{\rm
As $\rho_{B_{2n}}(z)=-1+\Vert z\Vert^2=-1+\sum_{\ell=1}^n z_\ell\bar z_\ell $, we have $\Vert\nabla \rho_{B_{2n}}(z)\Vert=2\Vert z\Vert=2$ for every $z\in\partial B_{2n}$, so 
\begin{align*}
P_n(B_{2n})
&=
\frac{i^n}{n!\varkappa_n}\int_{\partial B_{2n}}
\left(\frac{1}{2}\sum_{\ell=1}^n z_\ell\de\bar z_\ell-\bar z_\ell\de z_\ell\right)\wedge\left(\sum_{\ell=1}^n\de z_\ell\wedge\de\bar z_\ell\right)^{\wedge(n-1)}
\\
&=
\frac{i^n}{n!\varkappa_n}\int_{B_{2n}} \left(\sum_{\ell=1}^n\de z_\ell\wedge\de\bar z_\ell\right)^{\wedge n}
\\
&=
\frac{2^n}{\varkappa_n }\int_{B_{2n}}\upsilon_{2n}
\\
&=
\frac{2^n\varkappa_{2n}}{\varkappa_n}
\\
&=
\frac{2^n\varGamma(1/2)^n\varGamma(1+(n/2))}{n!}
\,.
\end{align*} 
Table~\ref{tab1} collects the values of $P_n(B_{2n})$ for $1\leq n\leq 10$.

\begin{table}[htb]
\begin{center}
\caption{Pseudovolume of full-dimensional unit balls of $\mathbb C^n$.}
\label{tab1}
{\small
\begin{tabular}{ccccccccccc}\toprule
$n$&1&2&3&4&5&6&7&8&9&10\\\midrule
$\displaystyle P_n(B_{2n})$&$\pi$&$2\pi$&$\pi^2$&$\dfrac{4\pi^2}{3}$&$\dfrac{\pi^3}{2}$&$\dfrac{8\pi^3}{15}$&$\dfrac{\pi^4}{6}$&$\dfrac{16\pi^4}{105}$&$\dfrac{\pi^5}{24}$&$\dfrac{32\pi^5}{945}$\\[2ex]
\bottomrule
\end{tabular}
}
\end{center}
\end{table}%

The computation through the Levi form is easier. Indeed ${\mathcal E}_{\partial B_{2n}}\equiv 1$ on $\partial B_{2n}$, so
\begin{equation*}
P_n(B_{2n})=\frac{2^{n-1}(n-1)!}{n!\varkappa_n}\,2n\varkappa_{2n}=\frac{2^n\varkappa_{2n}}{\varkappa_n}\,.
\end{equation*}
The computation of $P_n(B_{2n})$ as a convex body is possible too. In fact, $h_{B_{2n}}(z)=\Vert z\Vert$ so 
\begin{equation*}
(\ddc h_{B_{2n}})^{\wedge n}=4^n n!\det\left(
\frac{\partial^2 \Vert z\Vert}{\partial z_\ell\partial\bar z_k}\right) \upsilon_{2n}\,,
\end{equation*}
and \lemref{hess} yields
\begin{equation*}
\det\left(
\frac{\partial^2 \Vert z\Vert}{\partial z_\ell\partial\bar z_k}\right)
=
2^{-(n+1)}\Vert z\Vert^{-n }\,.
\end{equation*}
By using \eqref{wallis}, it follows that
\begin{align*}
\frac{1}{n!\varkappa_n}\int_{B_{2n}}(\ddc h_{B_{2n}})^{\wedge n}
&=
\frac{1}{\varkappa_n}\int_{B_{2n}}
 2^{n-1}\Vert z\Vert^{-n }\,\upsilon_{2n}
 \\
 &=
\frac{2^{n-1}}{\varkappa_n}\int_0^1
 \rho^{n-1}\de\rho\int_0^{2\pi}\de\vartheta_{2n-1}\prod_{\ell=1}^{2n-2}\int_0^\pi(\sin\vartheta_\ell)^{2n-\ell-1}\de\vartheta_\ell
 \\
 &=
  \frac{2^n\pi}{n\varkappa_n}\prod_{\ell=1}^{2n-2}\varGamma(1/2)\frac{\varGamma\big((2n-\ell)/2\big)}{\varGamma\big((2n-\ell+1)/2\big)}
  \\
  &=
   \frac{2^n\pi}{n\varkappa_n}\frac{ \pi^{n-1}}{(n-1)!}
   \\
   &=
   \frac{2^n\varkappa_{2n}}{\varkappa_n}\,,
\end{align*}
just as expected. As a last easiest computation, we mention the following one based on (\ref{strana2}):
\begin{align*}
P_n(B_{2n})
&=
\dfrac{1}{2^nn!\varkappa_n}
\int_{B_{2n}}
(\ddc h^2_{B_{2n}})^{\wedge n}
\\
&=
\dfrac{4^n}{2^nn!\varkappa_n}
\int_{B_{2n}}
n!\upsilon_{2n}
\\
&=
\dfrac{2^n\varkappa_{2n}}{\varkappa_{2n}}\,.
\end{align*}
Since $2/\varkappa_1=1$, it follows that the surface area $\varkappa_2$ of the unit disk in $\mathbb C$ equals its $1$-dimesional pseudovolume. However $2^n/\varkappa_n>1$ for $n>1$, so in higher dimension the $2n$-dimensional volume $\varkappa_{2n}$ of $B_{2n}$ is strictly smaller that its $n$-dimensional pseudovolume. Moreover $P_n(B_{2n})$ increases with $n$ for $n\leq 6$, then it decreases monotonically to $0$ as $n\to+\infty$. Observe that, for every $R\in \text{O}(2n)$, $R(B_{2n})=B_{2n}$ so the pseudovolume of a unit full-dimensional ball in $\mathbb C^n$ depends on its shape only, regardless of its position with respect to the complex structure.
}
\end{exe}

\begin{exe}\label{exebd} The $n$-pseudovolume of odd-dimensional unit balls $B_{2n-1}\subset\mathbb C^n$. 

\noindent
{\rm
If $E\subset\mathbb C^n$ is a real hyperplane, the intersection $E\cap B_{2n}$ is a $(2n-1)$-unit ball $B_{2n-1}$ in $E$. Up to a unitary transformation of $\mathbb C^n$, we may suppose that $E=\{z\in\mathbb C^n\mid \re z_1=0\}$. In this case 
\begin{equation*}
h_{B_{2n-1}}(z)=\sqrt{\left(\frac{z_1-\bar z_1}{2i}\right)^2+\sum_{\ell=2}^{n}\vert z_\ell\vert^2}\,.
\end{equation*}
If $n=1$, $B_1=[-i,i]\subset\mathbb C$ then, by \eqref{magnifica2}, $P_1(B_1)=2$. For $n>1$, \lemref{hess2} shows that
\begin{equation*}
\det\left(
\frac{\partial^2 h_{B_{2n-1}}}{\partial z_\ell\partial\bar z_k}\right)
=
\frac{\sum\nolimits_{\ell=2}^{n}\vert z_\ell\vert^2}{(2h_{B_{2n-1}})^{n+2}}
\,,
\end{equation*} 
which, by virtue of \eqref{polar}, in polar coordinates becomes
\begin{equation*}
\frac{\rho^2(\sin\vartheta_1)^2(\sin\vartheta_2)^2}{2^{n+2}\rho^{n+2}(\sin\vartheta_1)^{n+2}}
=
\frac{(\sin\vartheta_2)^2}{2^{n+2}\rho^n(\sin\vartheta_1)^n}
\,.
\end{equation*} 
It follows that $ P_n(B_{2n-1})$ equals
\begin{equation*}
\frac{1}{n!\varkappa_n}
\int_{B_{2n}}
4^n n!
\frac{(\sin\vartheta_2)^2}{2^{n+2}\rho^n(\sin\vartheta_1)^n}\,
\rho^{2n-1}
\prod_{\ell=1}^{2n-2}
(\sin\vartheta_\ell)^{2n-\ell-1}
\de\rho\wedge\de\vartheta
\end{equation*}
i.e.
\begin{equation*}
\frac{2^{n-2}}{\varkappa_n}
\int_{B_{2n}}
\rho^{n-1}(\sin\vartheta_1)^{n-2}(\sin\vartheta_2)^{2n-1}
\prod_{\ell=3}^{2n-2}
(\sin\vartheta_\ell)^{2n-\ell-1}
\de\rho\wedge\de\vartheta\,,
\end{equation*}
where, for the sake of notation, $\de\vartheta=\de\vartheta_1\wedge\ldots\wedge\de\vartheta_{2n-1}$.
By \eqref{wallis}
\begin{align*}
\int_0^\pi(\sin\vartheta_1)^{n-2}\de\vartheta_1
&=
\frac{\varGamma(1/2)\varGamma((n-1)/2)}{\varGamma(n/2)}\,,
\\
\int_0^\pi(\sin\vartheta_2)^{2n-1}\de\vartheta_2
&=
\frac{\varGamma(1/2)(n-1)!}{\varGamma(n+(1/2))}\,,
\\
\prod_{\ell=3}^{2n-2}
\int_0^\pi \sin^{2n-\ell-1}\vartheta_{\ell}\de\vartheta_\ell
&=
\frac{\varGamma(1/2)^{2n-4}}{(n-2)!}\,,
\end{align*}
so that
\begin{align*}
P_n(B_{2n-1})
&=
\frac{2^{n-2}2\pi\varGamma(1/2)^{2n-2}\varGamma((n-1)/2)(n-1)!}{\varkappa_n n\varGamma(n/2)\varGamma(n+(1/2))(n-2)!}
\\
&=
\frac{2^{n-2}\varGamma(1/2)^{n}\varGamma((n-1)/2)(n-1)}{\varGamma(n+(1/2))}
\\
&=
\dfrac{2^{n-1}\pi^{n/2}\varGamma((n+1)/2)}{\varGamma(n+(1/2))}\,.
\end{align*}
The table~\ref{tab2} collects the values of $P_n(B_{2n-1})$ and $P_n(B_{2n})$, for $1\leq n\leq 10$. A simple check shows that $P_n(B_{2n-1})<P_n(B_{2n})$ for every $n\in\mathbb N^*$. Moreover, $P_n(B_{2n-1})$ increases with $n$ for $n\leq 6$, then it decreases monotonically to $0$, as $n\to+\infty$.

\begin{table}[htb]
\caption{$n$-pseudovolume of $B_{2n-1}$ and $B_{2n}$.}
\label{tab2}
\begin{center}
{\small
\begin{tabular}{ccccccccccc}\toprule
$n$&1&2&3&4&5&6&7&8&9&10\\\midrule
$P_n(B_{2n-1})$&$2$&$\dfrac{4\pi}{3}$&$\dfrac{32\pi}{15}$&$\dfrac{32\pi^2}{35}$&$\dfrac{1024\pi^2}{945}$&$\dfrac{256\pi^3}{693}$&$\dfrac{16384\pi^3}{45045}$&$\dfrac{2048\pi^4}{19305}$&$\dfrac{1048576\pi^4}{11486475}$&$\dfrac{16384\pi^5}{692835}$\\\midrule
$P_n(B_{2n})$&$\pi$&$2\pi$&$\pi^2$&$\dfrac{4\pi^2}{3}$&$\dfrac{\pi^3}{2}$&$\dfrac{8\pi^3}{15}$&$\dfrac{\pi^4}{6}$&$\dfrac{16\pi^4}{105}$&$\dfrac{\pi^5}{24}$&$\dfrac{32\pi^5}{945}$
\\\bottomrule
\end{tabular}
}
\end{center}
\end{table}

}
\end{exe}

\begin{exe}
The $n$-pseudovolume of the standard $2n$-cube of $\mathbb C^n$.

\noindent
Consider the standard $2n$-cube $I_{2n}\subset\mathbb C^n$. $I_{2n}$ is the cartesian product of $n$ copies of the convex subset $\conv(\{1+i,-1+i,-1-i,1-i\})\subset\mathbb C$. As follows from \exeref{ncubo}, the $n$-th component of the face vector $f(I_{2n})$ is given by
\begin{equation*}
f_n(I_{2n})=2^n{2n\choose n}\,.
\end{equation*}
In $\mathbb R^{2n}$ there are ${2n\choose n}$ $n$-dimensional coordinate subspaces and the $n$-faces of $I_{2n}$ are organized into $2^n$ collections, each one consisting of faces parallel to a same real $n$-dimensional coordinate subspace of $\mathbb R^{2n}\simeq\mathbb C^n$. As shown in \exeref{realsim}, of such real subspaces only $2^n$ are equidimensional (and in fact real similar), so that the cardinal number of ${\mathcal B}_{\rm ed}(I_{2n},n)$ is $4^n$. For each $\Delta\in{\mathcal B}_{\rm ed}(I_{2n},n)$, one has $\varrho(\Delta)=1$, $\vol_n(\Delta)=2^n$ and $\psi_{I_{2n}}(\Delta)=2^{-n}$, whence $P_n(I_{2n})=4^n$. Notice that $P_n(I_{2n})=\vol_{2n}(I_{2n})$, quite an uncommon circumstance! Indeed, for every $\lambda>0$, $P_n(\lambda I_{2n})=\lambda^n 4^n$ whereas $\vol_{2n}(\lambda I_{2n})=\lambda^{2n}4^n$.
\end{exe} 

\begin{exe}
The $n$-pseudovolume of the standard $(2n-1)$-cube of $\mathbb C^n$.

\noindent
Consider a $(2n-1)$dimensional cube $\Gamma$ in $\mathbb C^n$. Up to a unitary transformation, we may suppose that $E_\Gamma=\{z\in\mathbb C^n\mid\im z_n=0\}$. If we take $\Gamma$ as the cartesian product of $n-1$ copies of the convex subset $\conv(\{1+i,-1+i,-1-i,1-i\})\subset\mathbb C$ with the interval $[-1,1]$. As follows from \exeref{ncubo}, the $n$-th component of the face vector $f(\Gamma)$ is given by
\begin{equation*}
f_n(\Gamma)=2^{n-1}{2n-1\choose n}\,.
\end{equation*}
In $\mathbb R^{2n-1}$ there are ${2n-1\choose n}$ $n$-dimensional coordinate subspaces and the $n$-faces of $\Gamma$ are organized into $2^{n-1}$ collections, each one consisting of faces parallel to a same real $n$-dimensional coordinate subspace of $\mathbb R^{2n-1}\subset\mathbb C^n$. Of such real subspaces only $2^{n-1}$ are equidimensional (and in fact real similar): they are obtained by adding the vector $e_n$ to each of the $2^{n-1}$ real similar $(n-1)$-dimensional coordinate subspaces of $\mathbb C^{n-1}$. So that the cardinal number of ${\mathcal B}_{\rm ed}(\Gamma,n)$ is $4^{n-1}$. For each $\Delta\in{\mathcal B}_{\rm ed}(\Gamma,n)$, one has $\varrho(\Delta)=1$, $\vol_n(\Delta)=2^n$ and $\psi_\Gamma(\Delta)=2^{-n+1}$, whence $P_n(\Gamma)=2^{2n-1}$ and, by homogeneity, we can deduce the $n$-pseudovolume of any $(2n-1)$-cube. 
\end{exe}

\begin{exe}
The $2$-pseudovolume of the standard unit $3$-crosspolytope of $\mathbb C^2$.

\noindent
Consider the standard unit $2n$-crosspolytope  $\Theta_{2n}\subset\mathbb C^n$. As shows \exeref{ntetra}, the $n$-th component of the face vector $f(\Theta_{2n})$ is given by
\begin{equation*}
f_n(\Theta_{2n})=2^{n+1}{2n\choose n+1}\,.
\end{equation*}
In the case $n=2$, the $4$-crosspolytope $\Theta_4$ has only equidimensional $2$-faces, $16$ of which are real similar whereas the remaining ones have a coefficient of area distortion equal to $2/3$. The outer angle of any $2$-face equals $1/6$ and the area of any $2$-face is $\sqrt3/2$, so that $P_2(\Theta_4)=(20/9)\sqrt3$. Observe by the way that the $3$-dimensional crosspolytope, (i.e. the octahedron) $\Theta_3\subset \mathbb C\times\re\mathbb C$ has $8$ equidimensional facets with $\varrho=2/3$, so that $P_2(\Theta_3)=4\sqrt3/3$.
\end{exe}

\begin{exe}
The $2$-pseudovolume of the standard unit $4$-crosspolytope of $\mathbb C^2$.

\noindent
Consider the standard unit $2n$-crosspolytope  $\Theta_{2n}\subset\mathbb C^n$. As shows \exeref{ntetra}, the $n$-th component of the face vector $f(\Theta_{2n})$ is given by
\begin{equation*}
f_n(\Theta_{2n})=2^{n+1}{2n\choose n+1}\,.
\end{equation*}
In the case $n=2$, the $4$-crosspolytope $\Theta_4$ has only equidimensional $2$-faces, $16$ of which are real similar whereas the remaining ones have a coefficient of area distortion equal to $2/3$. The outer angle of any $2$-face equals $1/6$ and the area of any $2$-face is $\sqrt3/2$, so that $P_2(\Theta_4)=(20/9)\sqrt3$. Observe by the way that the $3$-dimensional crosspolytope, (i.e. the octahedron) $\Theta_3\subset \mathbb C\times\re\mathbb C$ has $8$ equidimensional facets with $\varrho=2/3$, so that $P_2(\Theta_3)=4\sqrt3/3$.
\end{exe}

\begin{exe}\label{B34} The mixed $2$-pseudovolume of $B_3$ and $B_4$.

\noindent
The form $\ddc h_{B_3}\wedge\ddc h_{B_4}$ is locally integrable and it is differentiable almost everywhere on $B_4$. 
The computation of the derivatives of $h_{B_3}$ and $h_{B_4}$ carried on in the proofs of \lemref{hess} and \lemref{hess2} show that 
\begin{equation*}
Q_2(B_3,B_4)
=
\frac{1}{2\pi}
\int_{B_4}
\dfrac{2y_1^2+3x_2^2+3y_2^2}{h_{B_3}^3h_{B_4}}
\,\upsilon_4=\dfrac{16}{3}\,.
\end{equation*}
\end{exe}

\begin{exe}\label{B24} The mixed $2$-pseudovolume of a $2$-dimensional unit ball and $B_4$.

\noindent
Let $A\subset \mathbb C^2$ a $2$-dimensional unit ball. Since $Q_2(A,B_4)$ is translation invariant and orthogonally invariant, by \thmref{azzok1}, 
$Q_2(A,B_4)=(4/3)\textsl{v}_1(A)=4\pi/3$.
\end{exe}

\begin{exe}
If $n>1$, the Kazarnovski{\v\i} $n$-dimensional pseudovolume is not orthogonally invariant.

\noindent
{\rm
Let $I_n\subset\mathbb R^n\subset\mathbb C^n$ be the standard unit $n$-cube. We know that $P_n(I_n)=\vol_n(I_n)=2^n$. 
If $F:\mathbb R^{2n}\to\mathbb R^{2n}$ is the orthogonal mapping considered in \rmqref{qinv}, it is easy to realize that
the image of $F(I_n)$ via $F$ is a $n$-cube in $\mathbb C^n$ such that $E_{F(I_n)}^\mathbb C\neq\{0\}$, so that $P_n(F(I_n))=0$. In particular, the present example implies that the $1$-form $\alpha_{\partial A}$ cannot be orthogonally invariant.
}
\end{exe}

\begin{exe}
The $n$-dimensional Kazarnovski{\v\i} mixed psudovolume is not rational on lattice polytopes.  

\noindent
{\rm
Given $n$ polytopes $\Gamma_1,\ldots,\Gamma_n\subset\mathbb C^n$ whose vertices have coordinates in the ring of Gauss' integers, then $n!Q_n(\Gamma_1,\ldots,\Gamma_n)$ is, generally, not an integer and $(2n)!Q_n(\Gamma_1,\ldots,\Gamma_n)$ need not be either. If $n=1$, consider the polytope $\Gamma\subset\mathbb C$ with vertices $1,-1,i,-i$; then $1!P_1(\Gamma)=2\sqrt2\notin\mathbb Z$ and $2!P_1(\Gamma)=4\sqrt2\notin\mathbb Z$. The  example can be generalized to $\mathbb C^n$.
}
\end{exe}

\section{Open questions}\label{open}

In this section we present three open questions about $Q_n$ together with some comments and partial answers . 
\begin{enumerate}
\item Prove or disprove the non-vanishing condition of~\cororef{ndeg2} in the case of arbitrary convex bodies;
\item Find monotonicity conditions for $Q_n$;
\item Find weaker forms of Alexandroff-Fenchel inequality for $Q_n$;
\item Characterize the convex bodies for which Alexandroff-Fenchel inequality for $Q_n$ holds true.
\end{enumerate}

\subsection{Non vanishing condition}
In the polytopal case a non vanishing condition has been already stated in~\cite{Ka1} and proved in~\cite{Ka10}; our~\cororef{ndeg2} gives an alternative proof.
This non vanishing condition in the polytopal case is necessary and sufficient. However, in the general case, it is just sufficient and it is not known if it is also necessary.

\begin{lem}\label{buonino}
Let $A_1,\dots,A_n\in{\mathcal K}(\mathbb C^n)$ and $A=A_1+\ldots+A_n$ their Minkowski sum. If $\dim_\mathbb C A<n$ then $Q_n(A_1,\ldots,A_n)=0$.
\end{lem}
\noindent
\pf
For $1\leq\ell\leq n$, let $(\Gamma_\ell^{(m)})_m$ be a sequence of polytopes with vertices belonging to $\partial A_\ell$ that converges to $A_\ell$ for the Hausdorff metric. As $\Gamma_\ell^{(m)}\subseteq A_\ell$, it follows that, for every $m\in\mathbb N$, $\dim_\mathbb C (\Gamma_1^{(m)}+\ldots+\Gamma_n^{(m)})\leq\dim_\mathbb C A <n$ so, by \cororef{ndeg} and \cororef{azzok2}, $Q_n(\Gamma_1^m,\ldots,\Gamma_n^m)=0$. The continuity of $Q_n$ implies $Q_n(A_1,\ldots,A_n)=0$.\findim

The difficulties with the necessity of the condition come from the limiting procedure. Indeed it may happen that $A_1,\dots,A_n\in{\mathcal K}(\mathbb C^n)$ are approximated by $\Gamma_1^{(m)},\ldots,\Gamma_n^{(m)}\in{\mathcal P}(\mathbb C^n)$ such that $Q_n(\Gamma_1^{(m)},\ldots,\Gamma_n^{(m)})>0$ for every $m\in\mathbb N$, but 
$$
\lim_{m\to\infty} Q_n(\Gamma_1^{(m)},\ldots,\Gamma_n^{(m)})=Q_n(A_1,\ldots,A_n)=0^+.
$$

\subsection{Monotonicity}
The question of monotonicity of $Q_n$ (with respect to inclusion in each argument) is also very interesting and it seems related to non degeneracy and to the valuation property. Let us start with a simple general fact.

\begin{lem}\label{monotono}
Let $A_1,\ldots, A_m,K_2,\ldots,K_n\in{\mathcal K}(\mathbb C^n)$ and $A=A_1+\ldots+A_m$. 
For every $1\leq \ell\leq m$ and $1\leq k\leq n$, 
\begin{equation}\label{mon1}
Q_n(A_\ell[k],K_{k+1},\ldots,K_n)\leq \sum_{j=1}^m Q_n(A_j[k],K_{k+1},\ldots,K_n)\leq Q_n(A[k],K_{k+1},\ldots,K_n)\,. 
\end{equation}
\end{lem}
\noindent
\pf
The last statement is a monotonicity result and is a straightforward consequence of the~\eqref{mon1}, so it is enough to prove the latter inequalities. By direct computation:
\begin{align*}
Q_n(A_\ell[k],K_{k+1},\ldots,K_n)
&=
\dfrac{1}{n!\varkappa_n}\int_{B_n}(\ddc h_{A_\ell})^{\wedge k}\wedge \ddc h_{K_{k+1}}\wedge\ldots\wedge\ddc h_{K_n}
\\
&\leq
\dfrac{1}{n!\varkappa_n}\int_{B_n}\sum_{j=1}^m(\ddc h_{A_j})^{\wedge k}\wedge \ddc h_{K_{k+1}}\wedge\ldots\wedge\ddc h_{K_n}
\\
&=
\sum_{j=1}^m Q_n(A_j[k],K_{k+1},\ldots,K_n)
\end{align*} 
and 
\begin{align*}
&
\sum_{j=1}^m Q_n(A_j[k],K_{k+1},\ldots,K_n)
\\
&=
\dfrac{1}{n!\varkappa_n}\int_{B_n}\sum_{j=1}^m(\ddc h_{A_j})^{\wedge k}\wedge \ddc h_{K_{k+1}}\wedge\ldots\wedge\ddc h_{K_n}
\\
&\leq
\dfrac{1}{n!\varkappa_n}\int_{B_n}\sum_{\substack{0\leq k_1,\ldots, k_m\leq k \\ k_1+\ldots+k_m=k}}\dfrac{k!}{k_1!\cdots k_m!}\bigwedge_{\ell=1}^m(\ddc h_{A_{\ell}})^{\wedge k_\ell}\wedge \ddc h_{K_{k+1}}\wedge\ldots\wedge\ddc h_{K_n}
\\
&=
\dfrac{1}{n!\varkappa_n}
\int_{B_n}(\ddc h_{A_1}+\ldots +\ddc h_{A_m})^{\wedge k}\wedge \ddc h_{K_{k+1}}\wedge\ldots\wedge\ddc h_{K_n}
\\
&=
\dfrac{1}{n!\varkappa_n}
\int_{B_n}(\ddc h_{A})^{\wedge k}\wedge \ddc h_{K_{k+1}}\wedge\ldots\wedge\ddc h_{K_n}
\\
&=
Q_n(A[k],K_{k+1},\ldots,K_n)\,.
\end{align*} 
The lemma is thus proved.\findim

Observe that, by symmetry,~\lemref{monotono} is still valid in any of the arguments of $Q_n$. In particular if $A_\ell\subset A$ for some $\ell$, then \lemref{monotono} provides a monotonicity result for $Q_n$ on ${\mathbb K}(\mathbb C^n)$. \exeref{nonmon} suggests that, in such a generality, there is no other monotonicity result.    

\color{black}
\begin{exe}\label{nonmon}
The Kazarnovski{\v\i} $n$-dimensional pseudovolume is, generally, not monotonically incresing.

\noindent
{\rm Consider the the square $$\Delta=\conv\{2+2i,-2+2i,-2-2i,2-2i\}\subset\mathbb C=\mathbb C\times\{0\}\subset\mathbb C^2$$
and let, for every $\lambda>0$, $\Gamma^\lambda=\conv\{\Delta;(0,2\lambda)\}\subset\mathbb C^2$. The polytope $\Gamma^\lambda$ is a 3-dimensional pyramid on the base $\Delta$ with apex in the point $(0,2\lambda)$. Let also 
$$
K^\lambda=\conv\{(i,0),(- i,0),(i,\lambda),(-i,\lambda)\}\subset \mathbb C^2\,,
$$ 
then $K^\lambda\subset\Gamma^\lambda$, for every $\lambda>0$. Figure~\ref{fig5} depicts both $K^\lambda$ and $\Gamma^\lambda$ when $\lambda=1/4$.

\begin{figure}[ht]
\begin{center}
\begin{tikzpicture}[scale=2.5]
\draw[->] (0,0,-3) -- (0,0,3) node[below] {$\re z_1$};
\draw[->] (-2.5,0,0) -- (2.5,0,0) node[below] {$\im z_1$};
\draw[->] (0,0,0) -- (0,1.5,0) node[above] {$\re z_2$};
\filldraw[thick,fill=green!20!white,draw=blue,semitransparent,line cap=round, line join=round] (-2,0,2)--(2,0,2)--(2,0,-2)--(-2,0,-2)--cycle;
\draw[thick,fill=green!20!white,draw=blue,semitransparent,line cap=round, line join=round] (-2,0,-2)--(0,0.5,0)--(2,0,-2);
\filldraw[thick, fill=blue!20!white, draw=orange,very nearly opaque,line cap=round, line join=round] (-1,0,0)--(1,0,0)--(1,0.25,0)--(-1,0.25,0)--cycle;
\draw[thick,fill=green!20!white,draw=blue,nearly transparent,line cap=round, line join=round] (-2,0,2)--(0,0.5,0)--(2,0,2);
\draw[thick,fill=green!20!white,draw=blue,nearly transparent,line cap=round, line join=round] (-2,0,2)--(0,0.5,0)--(-2,0,-2);
\draw[thick,,draw=blue,fill=green!20!white,nearly transparent,line cap=round, line join=round] (2,0,2)--(0,0.5,0)--(2,0,-2);
\draw[dashed,blue!40!white,line cap=round, line join=round] (-2,0,0)--(0,0.5,0)--(2,0,0);
\end{tikzpicture}
\caption{The rectangle $K^\lambda$ and the pyramid $\Gamma^\lambda$ for $\lambda=1/4$.}
\label{fig5}
\end{center}
\end{figure}
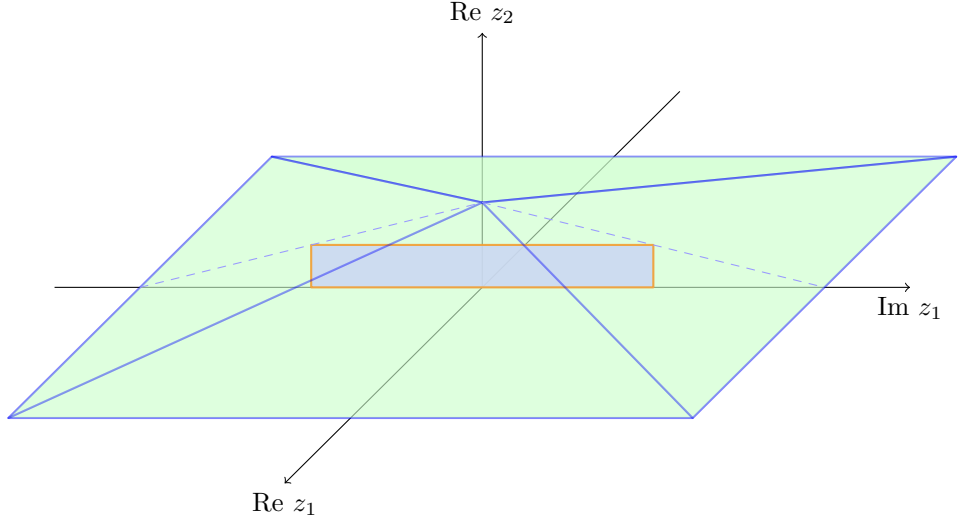

It is very easy to show that $P_2(K^\lambda)=2\lambda$. In order to compute $P_2(\Gamma^\lambda)$ notice that  $\Gamma^\lambda$ has five $2$-dimensional faces, namely 
\begin{align*}
\Delta_1&=
\conv\{(2+2 i,0),(-2+2 i,0),(0,2\lambda)\},
\\
\Delta_2&=
\conv\{(2-2 i,0),(-2-2 i,0),(0,2\lambda)\},
\\
\Delta_3&=
\conv\{(2+2 i,0),(2-2 i,0),(0,2\lambda)\},
\\
\Delta_4&=
\conv\{(-2+2 i,0),(-2-2 i,0),(0,2\lambda)\}.
\end{align*}
and $\Delta_5=\Delta$. The face $\Delta_5$ does not contribute to $P_2(\Gamma^\lambda)$ because $\varrho(\Delta_5)=0$, whereas
\begin{equation*}
\varrho(\Delta_\ell)=\frac{\lambda^2}{1+\lambda^2}\,,\quad
\vol_2(\Delta_\ell)=4\sqrt{1+\lambda^2}\,,\quad
\psi(\Delta_\ell)=\frac{1}{2}\,,
\end{equation*}
for every $\ell\in\{1,2,3,4\}$, 
so that 
\begin{equation*}
P_2(\Gamma^\lambda)
=
\frac{8\lambda^2}{\sqrt{1+\lambda^2}}\,.
\end{equation*}
This computation shows that $P_2(\Gamma^\lambda)<P_2(K^\lambda)$ as soon as $0<\lambda<1/\sqrt{15}$. The preceding example is a variation of one provided by B.~Ya.~Kazarnovski{\v\i} (in a private communication) and can be generalized to higher dimensions. Observe that $\dim_\mathbb R K^\lambda<\dim_\mathbb R \Gamma^\lambda$, i.e. $K^\lambda$ is dimensionally smaller than $\Gamma^\lambda$, however it can have a bigger pseudovolume than $\Gamma^\lambda$.\findim
}
\end{exe}

Under suitable assumptions, special monotonicy results can be obtained by using the valuation property.
The following lemma is an example of this kind.

\begin{lem}\label{luglio2019}
Let $A\in{\mathcal K}(\mathbb C^n)$ such that $\dim_\mathbb R A=2n-1$ and let $H\subset\mathbb C^n$ a real affine hyperplane such that $E_A\cap E_H=E_A^\mathbb C$. Then $P_n(A\cap H)=0$ and $P_n(A)=P_n(A\cap H^+)+P_n(A\cap H^-)$. In particular $P_n(A)\geq \max\{P_n(A\cap H^+),P_n(A\cap H^-)\}$.
\end{lem}
\noindent
\pf
The hyperplane $H$ cannot be included in $\aff_\mathbb R A$, since otherwise (for dimensional reasons) it should coincide with $\aff_\mathbb R A$, then $E_A=E_H$ and $E_H=E_A=E_A^\mathbb C$, which is impossible because $\dim_\mathbb R E_H=2n-1$ whereas $\dim_\mathbb R E_A^\mathbb C=2n-2$. The  convex body $A\cap H$ has zero pseudovolume. In fact, any sequence of polytopes with vertices on the relative boundary of $A\cap H$ and converging to $A\cap H$ has the property that each polytope of the sequence has a complex dimension that is smaller than $n$. By~\cororef{ndeg},~\cororef{azzok2} and~\rmqref{continuo2}, it follows that $P_n(A\cap H)=0$. Thanks to~\cororef{valutazione-completa}, $P_n(A)=P_n(A\cap H^+)+P_n(A\cap H^-)$, then $P_n(A)\geq \max\{P_n(A\cap H^+),P_n(A\cap H^-)\}$.  \findim
 
\begin{rmq}
By \lemref{luglio2019}, one gets $P_n(A)>0$ as soon as $P_n(A\cap H^+)>0$ or $P_n(A\cap H^-)>0$. Unless $A$ is a polytope, the latter inequalities are hardly checked, so \lemref{luglio2019} does not imply any non vanishing criterion for $P_n(A)$.
\end{rmq}

\subsection{Alexandroff-Fenchel type inequalities}\label{afk}

If $n\geq 2$ and $A_1,A_2,A_3,\ldots,A_n\in{\mathcal K}(\mathbb C^n)$, the general Alexandroff-Fenchel inequality for Kazarnovsk{\v\i}i  mixed pseudovolume would read as follows:
\begin{equation}\label{KAF}
Q_n(A_1,A_2,A_3,\ldots,A_n)^2\geq Q_n(A_1[2],A_3,\ldots,A_n)Q_n(A_2[2],A_3,\ldots,A_n)\,,
\end{equation}
however, as shown by \rmqref{nokaf}, the latter equality is generally false, so we look for weaker inequalities like
\begin{equation}\label{KAFw}
Q_n(A_1,A_2,A_3,\ldots,A_n)^2\geq 
c\,Q_n(A_1[2],A_3,\ldots,A_n)Q_n(A_2[2],A_3,\ldots,A_n)\,,
\end{equation}
for some constant $c\in [0,1]$.
The following instances of (\ref{KAFw}) recall Minkowski's quadratic inequality for mixed volumes. 

\begin{propos}
Let $1\leqslant k\leqslant n$ be an integer and let $\Gamma_1,\ldots,\Gamma_k\in{\mathcal P}(\mathbb C^n)$. If $\Gamma=\sum_{\ell=1}^k\Gamma_k$, then 
\begin{align}\label{bestkaf0}
Q_n(\Gamma_1,\ldots,\Gamma_k,B_{2n}[n-k])^2
\geqslant
\min_{\Delta\in{\mathcal B}(\Gamma,k)}
\varrho(\Delta)^2
\prod_{j=1}^2
Q_n(\Gamma_j[2],\Gamma_3,\ldots,\Gamma_k,B_{2n}[n-k])
\,.
\end{align}
In particular, if $A_1,\ldots,A_k\in{\mathcal K}(\mathbb R^n)$, then 
\begin{align}\label{bestkaf1}
Q_n(A_1,\ldots,A_k,B_{2n}[n-k])^2
\geqslant
\prod_{j=1}^2
Q_n(A_j[2],A_3,\ldots,A_k,B_{2n}[n-k])
\,.
\end{align}
\end{propos}

\noindent
\pf
By \cororef{azzok2}, (\ref{comb}), (\ref{AF}) and (\ref{diseg}),
\begin{align}
&Q_n(\Gamma_1,\ldots,\Gamma_k,B_{2n}[n-k])^2
\\
&=
\left[
\dfrac{2^{n-k}\varkappa_{2n-k}}{\varkappa_n\displaystyle{n\choose k}}
\right]^2
\left[
\sum_{\Delta\in{\mathcal B}(\Gamma,k)}
\varrho(\Delta)
V_k(\Delta_1,\ldots,\Delta_k)
\psi_\Gamma(\Delta)
\right]^2
\\
&\geqslant
\left[
\dfrac{2^{n-k}\varkappa_{2n-k}}{\varkappa_n\displaystyle{n\choose k}}
\right]^2
\min_{\Delta\in{\mathcal B}(\Gamma,k)}
\varrho(\Delta)^2
\left[
\sum_{\Delta\in{\mathcal B}(\Gamma,k)}
V_k(\Delta_1,\ldots,\Delta_k)
\psi_\Gamma(\Delta)
\right]^2
\\
&=
\left[
\dfrac{2^{n-k}\varkappa_{2n-k}}{\varkappa_n\displaystyle{n\choose k}}
\dfrac{\displaystyle{2n\choose k}}{\varkappa_{2n-k}}
\right]^2
\min_{\Delta\in{\mathcal B}(\Gamma,k)}
\varrho(\Delta)^2
\,
V_{2n}(\Gamma_1,\ldots,\Gamma_k,B_{2n}[2n-k])^2
\\
&\geqslant
\left[
\dfrac{2^{n-k}\displaystyle{2n\choose k}}{\varkappa_n\displaystyle{n\choose k}}
\right]^2
\min_{\Delta\in{\mathcal B}(\Gamma,k)}
\varrho(\Delta)^2
\prod_{j=1}^2
V_{2n}(\Gamma_j[2],\Gamma_3,\ldots\Gamma_k,B_{2n}[2n-k])
\\
&\geqslant
\left[
\dfrac{2^{n-k}\displaystyle{2n\choose k}}{\varkappa_n\displaystyle{n\choose k}}
\dfrac{\varkappa_{2n-k}}{\displaystyle{2n\choose k}}
\right]^2
\min_{\Delta\in{\mathcal B}(\Gamma,k)}
\varrho(\Delta)^2
\prod_{j=1}^2
V^\varrho_k(\Gamma_j[2],\Gamma_3,\ldots,\Gamma_k)
\\
&=
\left[
\dfrac{2^{n-k}\varkappa_{2n-k}}{\varkappa_n\displaystyle{n\choose k}}
\right]^2
\min_{\Delta\in{\mathcal B}(\Gamma,k)}
\varrho(\Delta)^2
\prod_{j=1}^2
V^\varrho_k(\Gamma_j[2],\Gamma_3,\ldots,\Gamma_k)
\\
&=
\min_{\Delta\in{\mathcal B}(\Gamma,k)}
\varrho(\Delta)^2
\prod_{j=1}^2
Q_n(\Gamma_j[2],\Gamma_3,\ldots,\Gamma_k,B_{2n}[n-k])\,.
\end{align}
If $\Gamma_1,\ldots,\Gamma_k\in{\mathcal P}(\mathbb R^n)$, then $\min_{\Delta\in{\mathcal B}(\Gamma,k)}\varrho(\Delta)^2=1$ and (\ref{bestkaf1}) follows by continuity.\findim

Observe that if every $k$-dimensional face of $\Gamma$ is equidimensional, then $\min_{\Delta\in{\mathcal B}(\Gamma,k)}\varrho(\Delta)>0$ and so (\ref{bestkaf0}) is non trivial, whereas (\ref{bestkaf1}) provides an instance of (\ref{KAF}) because, in fact, it deals with ordinary mixed volume in $\mathbb R^{2n}$. The following proposition shows that for $k=1$ a better inequality holds true.

\begin{propos}\label{bestkaf2}
For every $n\geqslant 2$ and $A\in{\mathcal K}(\mathbb C^n)$
\begin{equation}\label{wKAF}
Q_n(A,B_{2n}[n-1])^2\geqslant 
\left(\dfrac{2n-2}{2n-1}\right)Q_n(A[2],B_{2n}[n-2])P_n(B_{2n})\,.
\end{equation}
\end{propos}

\noindent
\pf 
If $\Gamma\in{\mathcal P}(\mathbb C^n)$, by (\ref{f-azzok}), (\ref{meta}), (\ref{AF}), (\ref{kth}), (\ref{ineq}) and \exeref{pseudoball}, one gets
\begin{align}
Q_n(\Gamma,B_{2n}[n-1])^2
&=
\left[
\dfrac{2^{n-1}\varkappa_{2n-1}}{n\varkappa_n}\textsl{v}_1^\varrho(\Gamma)\right]^2
\\
&=
\left[
\dfrac{2^{n-1}\varkappa_{2n-1}}{n\varkappa_n}\dfrac{2n}{\varkappa_{2n-1}}V_{2n}(\Gamma,B_{2n}[2n-1])\right]^2
\\
&\geqslant
\dfrac{2^{2n}\varkappa_{2n}}{\varkappa_n^2}
V_{2n}(\Gamma[2],B_{2n}[2n-2])
\\
&=
\dfrac{2^{2n}\varkappa_{2n}}{\varkappa_n^2}
\dfrac{\varkappa_{2n-2}}{n(2n-1)}
\textsl{v}_2(\Gamma)
\\
&=
\dfrac{2^n\varkappa_{2n-2}}{n(2n-1)\varkappa_n}P_n(B_{2n})\textsl{v}_2(\Gamma)
\\
&\geqslant
\dfrac{2^n\varkappa_{2n-2}}{n(2n-1)\varkappa_n}P_n(B_{2n})
\textsl{v}^\varrho_2(\Gamma)
\\
&=
\dfrac{2^n\varkappa_{2n-2}}{n(2n-1)\varkappa_n}P_n(B_{2n})
\dfrac{n(n-1)\varkappa_n}{2^{n-1}\varkappa_{2n-2}}
Q_n(\Gamma[2],B_{2n}[n-2])
\\
&=
\left(\dfrac{2n-2}{2n-1}\right)
Q_n(\Gamma[2],B_{2n}[n-2])
P_n(B_{2n}).
\end{align} 
By continuity, the inequality can be extended from the polytopal case to the general one.\findim

\begin{rmq}\label{nokaf}
The coefficient $(2n-2/2n-1)$ on the right-hand side of~(\ref{wKAF}) cannot be increased to $1$. Indeed, if $n=2$ and $\Gamma\subset\mathbb R^2$ is a regular polygon with $m$ edges, each of which has length $1$, it is easy to show that 
\begin{equation}
Q_2(\Gamma,B_4)=\dfrac{2\varkappa_3}{2\varkappa_2}\textsl{v}_1^\varrho(\Gamma)=\dfrac{4}{3}\textsl{v}_1(\Gamma)=\dfrac{2}{3}m
\end{equation}
and
\begin{equation}
P_2(\Gamma)P_2(B_4)=\vol_2(\Gamma)2\pi=\dfrac{\pi}{2}m\cot\left(\dfrac{\pi}{m}\right),
\end{equation}
so that~(\ref{KAF}) holds true for $m\in\{3,4,5\}$, but it doesn't as soon as $m\geqslant 6$. An even simpler example is provided by the $2$-dimensional unit ball $B_2\in{\mathcal K}(\mathbb R^2)$. Indeed, by \exeref{B24}, one gets 
$$
Q_2(B_2,B_4)^2
=
\dfrac{16}{9}\pi^2
<
2\pi^2
=
P_2(B_2)P_2(B_4).
$$
\end{rmq}

We now address the question of equality cases in (\ref{wKAF}) for $n=2$. If $A$ is not reduced to a single point, $Q_2(A,B_4)=(4/3)\textsl{v}_1(\Gamma)>0$ and the equality $$Q_2(A,B_4)^2=\dfrac{4\pi}{3} P_2(A)\,,$$ requires $P_2(A)>0$. Such an equality is then equivalent to 
$$
\dfrac{\textsl{v}_1(A)^2}{\textsl{v}^\varrho_2(A)}=\dfrac{3}{4}\pi\,.
$$ 
The latter condition requires non real bodies, in fact for real ones it is impossible because the smallest value of the so called \emph{isoperimetric ratio}, $\textsl{v}_1(A)^2/\textsl{v}_2(A)=\textsl{v}_1(A)^2/\vol_2(A)$, is attained when $A$ is a $2$-dimensional real ball for which it equals $\pi$. If $A$ is a generic $2$-dimensional ball about the origin of $\mathbb C^2$, then 
$$
\dfrac{\textsl{v}_1(A)^2}{\textsl{v}^\varrho_2(A)}
=
\dfrac{\textsl{v}_1(A)^2}{\varrho(A)\vol_2(A)}
\geqslant
\pi\,,
$$ 
because $\varrho(A)\leqslant 1$, so we have to look for a convex body whose dimension is $3$ or $4$. Nor $B_4$ neither $B_3$ make (\ref{wKAF}) become an equality, as $P_2(B_4)>0$ and, by \exeref{pseudoball}, \exeref{exebd} and \exeref{B34}, 
$$
Q_2(B_3,B_4)^2=\dfrac{256}{9}> \dfrac{16}{9}\pi^2=\dfrac{2}{3}P_2(B_3)P_2(B_4).
$$ 
Finally we observe that no equality in (\ref{wKAF}) can involve a polytope. Indeed, if $\Gamma\in{\mathcal P}(\mathbb C^2)$ is such that $\dim_\mathbb R\Gamma\geqslant 3$, then $\dim_\mathbb C\Gamma=2$, and by \thmref{iok}, $\max_{\Delta\in{\mathcal B}(\Gamma,2)}\varrho(\Delta)>0$ so
\begin{align*}
\dfrac{\textsl{v}_1(\Gamma)^2}{\textsl{v}_2^\varrho(\Gamma)}
&=
\dfrac{16\varkappa_3^{-2}V_4(\Gamma,B_4[3])^2}{\sum\limits_{\Delta\in{\mathcal B}(\Gamma,2)}\varrho(\Delta)\vol_2(\Delta)\psi_\Gamma(\Delta)}
\\
&\geqslant
\dfrac{16\varkappa_3^{-2}V_4(\Gamma,B_4[3])^2}{\left[\max\limits_{\Delta\in{\mathcal B}(\Gamma,2)}\varrho(\Delta)\right]\sum\limits_{\Delta\in{\mathcal B}(\Gamma,2)}\vol_2(\Delta)\psi_\Gamma(\Delta)}
\\
&\geqslant
\dfrac{16\varkappa_3^{-2}V_4(\Gamma,B_4[3])^2}{6\varkappa_2^{-1}\left[\max\limits_{\Delta\in{\mathcal B}(\Gamma,2)}\varrho(\Delta)\right]V_4(\Gamma[2],B_4[2])}
\\
&=
\dfrac{16\varkappa_3^{-2}\varkappa_4}{6\varkappa_2^{-1}\left[\max\limits_{\Delta\in{\mathcal B}(\Gamma,2)}\varrho(\Delta)\right]}
\cdot
\dfrac{V_4(\Gamma,B_4[3])^2}{\varkappa_4V_4(\Gamma[2],B_4[2])}
\\
&=
\dfrac{3}{4}\pi
\left[
\max\limits_{\Delta\in{\mathcal B}(\Gamma,2)}\varrho(\Delta)
\right]^{-1}
\dfrac{V_4(\Gamma,B_4[3])^2}{\vol_4(B_4)V_4(\Gamma[2],B_4[2])}\,.
\end{align*}
Now the term in brackets is greater or equal to $1$ and, even if it were equal to $1$, the last fraction involving mixed volumes is always strictly greater than $1$ because $V_4(\Gamma,B_4[3])^2\geqslant  \vol_4(B_4)V_4(\Gamma[2],B_4[2])$ by (\ref{AF}) and, as  proved by Bol~\cite{Bol}, an equality would absurdly require  $\Gamma$ and $B_4$ be homothetic. It follows that the equality case in (\ref{wKAF}) possibly requires bodies with a more  involved boundary structure.

Apart from the case of real bodies or that of $n$ coincident complex ones, despite \rmqref{nokaf} shows that (\ref{KAF}) is generally false on the whole ${\mathcal K}(\mathbb C^n)$, one can find non trivial examples for which (\ref{KAF}) in fact holds true.

An interesting investigation in this direction stems from  \cororef{corobuono}. Let $n>1$ and fix $(n-1)$ convex bodies $A_2,\ldots,A_n\in{\mathcal K}_2(\mathbb C^n)$. Observe that, by \lemref{dj} and \lemref{lemma-Xi}, for every convex body $A_1\in{\mathcal K}_2(\mathbb C^n)$,
\begin{align}
&Q_n(A_1,A_2,A_3,\ldots,A_n)
\\
&=
\dfrac{4^{n-1}}{\varkappa_n}
\int_{\partial B_{2n}}
D_n(N_{A_1},\Hess_\mathbb C h_{A_2},\Hess_\mathbb C h_{A_3},\ldots,\Hess_\mathbb C h_{A_n})\upsilon_{2n}
\\
&=
\dfrac{4^{n-1}}{\varkappa_n}
\int_{\partial B_{2n}}
h_{A_1}D_n(I_n,\Hess_\mathbb C h_{A_2},\Hess_\mathbb C h_{A_3},\ldots,\Hess_\mathbb C h_{A_n})\upsilon_{2n}
\\
&+
\dfrac{4^{n-1}}{\varkappa_n}
\int_{\partial B_{2n}}
\re\langle\nabla h_{A_1}(u),
\Xi(A_2,A_3,\ldots,A_n)(u)u\rangle\upsilon_{2n}\,.
\end{align}
The preceding computation suggests to consider the set ${\mathcal X}(A_2,\ldots,A_n)$\label{X} of convex bodies $A_1$ such that 
\begin{equation}\label{integrandX}
\int_{\partial B_{2n}}\re\langle\nabla h_{A_1}(u),\Xi(A_2,A_3\ldots,A_n)(u)\,u\rangle\upsilon_{\partial B_{2n}}=0\,.
\end{equation} 
Recall that $\Xi(A_2,A_3\ldots,A_n)(u)$ is a hermitian matrix and $\nabla h_{A_1}(u)=\nu^{-1}_{\partial A_1}(u)$, so 
$$
\re\langle\nabla h_{A_1}(u),\Xi(A_2,A_3\ldots,A_n)(u)\,u\rangle
=
\re\langle\Xi(A_2,A_3\ldots,A_n)(u)\,\nu^{-1}_{\partial A_1}(u),u\rangle
$$
i.e. the integral in (\ref{integrandX}) is proportional to the average on $\partial B_{2n}$ of the projection of the vector $\Xi(A_2,A_3\ldots,A_n)(u)\,\nu^{-1}_{\partial A_1}(u)$ on the unit vector $u$.

If $n=2$, it is worth noting that $B_4\notin {\mathcal X}(B_4)$, indeed on $\partial B_4$, one gets $\nabla h_{B_4}(u)=u$, 
\begin{equation}
\Xi(B_4)(u)\,u
=
\frac{1}{8}
\left(
\begin{array}{ccc}
\vert u_2\vert^2-2 && u_1\bar u_2
\\
&&
\\
\bar u_1 u_2 && \vert u_1\vert^2 -2
\end{array}
\right)
\left(
\begin{array}{c}
u_1
\\
~
\\
u_2
\end{array}
\right)
=
-\frac{1}{4}
\left(
\begin{array}{c}
\vert u_1\vert^2u_1
\\
~
\\
\vert u_2\vert^2u_2
\end{array}
\right)
\end{equation}
and $\re\langle\nabla h_{B_4}(u),\Xi(B_4)(u)\,u\rangle=\langle\nabla h_{B_4}(u),\Xi(B_4)(u)\,u\rangle=-(\vert u_1\vert^4+\vert u_2\vert^4)/4\neq 0$.

If ${\mathcal X}(\mathbb C^n)$ denotes the largest subspace of ${\mathcal K}_2(\mathbb C^n)^n$ on which, for almost all $u\in \partial B_{2n}$, the equality $\re\langle\nabla h_{A_1}(u),\Xi(A_2,\ldots,A_n)\rangle=0$ holds true,  it seems reasonable to think that the inclusion ${\mathcal K}_2(\mathbb R^n)^n\subseteq {\mathcal X}(\mathbb C^n)$ is strict. Thanks to \cororef{corobuono}, the space ${\mathcal X}(\mathbb C^n)$ could be a starting point towards a proof of the inequality~(\ref{KAF}) on ${\mathcal X}(\mathbb C^n)$. The reason why we think so relies on Alexandroff’s original proof of (\ref{AF}) contained in \cite{A2}. That proof is based on the representation (\ref{mix-D-f}) of mixed volume, formally analogous to (\ref{formula-buona}) on ${\mathcal X}(\mathbb C^n)$, and on the following remarkable inequality enjoyed by mixed discriminants.
\begin{thm}[Alexandroff]\label{Ada}
Given complex $n$-by-$n$ real symmetric matrices $M,N,M_3,\ldots, M_n$ with
$M, M_3, \ldots , M_n$ positive definite and $N$ arbitrary, then
\begin{align}\label{ada}
D_n(M,N,M_3,\ldots,M_n)^2 \geqslant D_n(M,M,M_3,\ldots,M_n)D_n(N,N,M_3,\ldots,M_n),
\end{align} 
with equality if and only if $N =\lambda M$ for some $\lambda\in\mathbb R$.
\end{thm}

We recall that in the second page, lines 20-23, of the original paper \cite{A2}, Alexandroff claims that an hermitian version \thmref{Ada} can be proved by just adapting the proof he gives for real symmetric ones. 
Proofs of the hermitian version of \thmref{Ada} which differ from that indicated by Alexandroff have been found by Shenfeld and 
Van Handel \cite{SvH} as well as Li \cite{Ping}, who also observed that, by continuity, the inequality (\ref{ada}) remains true if $M,M_3,\ldots,M_n$ are only positive semi-definite. In view of (\ref{formula-buona}) and the hermitian version of \thmref{Ada}, it seems natural to investigate the extent to which the space ${\mathcal X}(\mathbb C^n)$ is adequate for (\ref{KAF}) to hold.

\section{Useful constructions and computations}\label{u}

We collect here some constructions and computations which have been used in the preceding sections. 

\subsection{Polar coordinates} 
Consider the following coordinate transformation in $\mathbb R^m$
\begin{equation*}
\xi_1=\rho\cos\vartheta_1\,,\quad \xi_\ell=\rho(\cos\vartheta_\ell)\prod_{j=1}^{\ell-1}\sin\vartheta_j\,,\quad{\rm if}\quad 2\leq\ell\leq m-1,\quad{\rm and}\quad\xi_{m}=\rho\prod_{j=1}^{m-1}\sin\vartheta_j\,,
\end{equation*}
where $\rho\in(0,+\infty)$, $\vartheta_{m-1}\in(0,2\pi)$ and $\vartheta_\ell\in(0,\pi)$, for every $1\leq\ell\leq m-2$. The corresponding jacobian is given by
\begin{equation*}
\rho^{m-1}\prod_{\ell=1}^{m-2}(\sin\vartheta_\ell)^{m-\ell-1}\,.
\end{equation*}
Observe that, for every $1\leq k\leq m$,
\begin{equation}\label{polar}
\sum_{\ell=k}^{m}
\xi_\ell^2=\left(
\rho\prod_{j=1}^{k-1}\sin\vartheta_j
\right)^2\,.
\end{equation}
Recall also that, for every $n\in\mathbb N$, the $n$-th Wallis' integral can be computed via Euler's gamma function as follows 
\begin{equation}\label{wallis}
\int_0^\pi (\sin\vartheta)^{n}\de\vartheta=\frac{\varGamma(1/2)\varGamma((n+1)/2)}{\varGamma((n+2)/2)}\,.
\end{equation}

\subsection{Complex hessians}

\begin{lem}\label{hess}
For every $n\in\mathbb N^*$ and every $z\in \mathbb C^n\setminus\{0\}$, $\det \Hess_\mathbb C h_{B_{2n}}=2^{-(n+1)}\Vert z\Vert^{-n }$.
\end{lem}
\noindent\pf 
An easy computation shows that, for every $1\leq \ell,k\leq n$,
\begin{equation*}
\frac{\partial^2 h_{B_{2n}}}{\partial z_\ell\partial\bar z_k}=\frac{2\Vert z\Vert^2\delta_{\ell,k}-z_k\bar z_\ell}{4\Vert z\Vert^3}\,,
\end{equation*}
so that we are led to show that
\begin{equation}\label{det}
2^{2n}\Vert z\Vert ^{3n}\det\Hess_\mathbb C h_{B_{2n}}
=
2^{n-1}\Vert z\Vert^{2n}
\,.
\end{equation}
Consider the matrix $M=(2\Vert z\Vert^2\delta_{\ell,k}-z_k\bar z_\ell)_{\ell,k}$
and observe that it admits the eigenvalue $\lambda_1=2\Vert z\Vert^2$. Indeed, the matrix $M-\lambda_1I_n=(-z_k\bar z_\ell)_{\ell,k}$ has rank $1$ because, for every $1\leq k\leq n$, the $k$-th column equals $-z_k{}^t(\bar z_1,\ldots,\bar z_n)$. It follows that the dimension of the corresponding eigenspace is $n-1$. A direct computation shows that the other eighenvalue of $M$ is $\lambda_2=\Vert z\Vert^2$. Indeed, let $a_\ell$ denote the $\ell$-th row of $M-\lambda_2I_n$ and suppose $z_k\neq0$, then $a_k=\sum_{\ell\neq k}a_\ell(-z_\ell/z_k)$.
Since the columns $a_\ell$, with $\ell\neq k$, are linearly independent, it follows that $M-\lambda_2I_n$ has rank $n-1$. \findim

\begin{lem}\label{hess2}
For every $n\in\mathbb N\setminus\{0,1\}$ and every $z\in \{z\in\mathbb C^n\mid \re z\neq0\}$ 
\begin{equation*}
\det
\Hess_\mathbb C  h_{B_{2n-1}}
=
\frac{\displaystyle\sum\nolimits_{\ell=2}^{n}\vert z_\ell\vert^2}{(2h_{B_{2n-1}})^{n+2}}\,.
\end{equation*}
\end{lem}
\noindent\pf 
We already know that 
\begin{equation*}
h_{B_{2n-1}}(z)=\sqrt{\left(\frac{z_1-\bar z_1}{2i}\right)^2+\sum_{\ell=2}^{n}\vert z_\ell\vert^2}
\end{equation*}

and an easy computation shows that, for every $2\leq \ell,k\leq n$,
\begin{equation*}
\frac{\partial^2 h_{B_{2n-1}}}{\partial z_\ell\partial\bar z_k}=\frac{2h_{B_{2n-1}}^2\delta_{\ell,k}-z_k\bar z_\ell}{4 h_{B_{2n-1}}^3}\,,
\end{equation*}
\begin{equation*}
\frac{\partial^2 h_{B_{2n-1}}}{\partial z_\ell\partial\bar z_1}
=
\frac{(\bar z_1- z_1)\bar z_\ell}{8 h_{B_{2n-1}}^3}\,,
\qquad
\frac{\partial^2 h_{B_{2n-1}}}{\partial z_1\partial\bar z_k}
=
\frac{(z_1-\bar z_1) z_k}{8 h_{B_{2n-1}}^3}
\,,
\end{equation*}
\begin{equation*}
\frac{\partial^2 h_{B_{2n-1}}}{\partial z_1\partial\bar z_1}
=
\frac{\displaystyle\sum\nolimits_{\ell=2}^{n}\vert z_\ell\vert^2}{4h_{B_{2n-1}}^{3}}\,,
\end{equation*}
so it's enough to show that
\begin{equation*}
4^nh_{B_{2n-1}}^{3n}
\det
\Hess_\mathbb C  h_{B_{2n-1}}
=
2^{n-2}(h_{B_{2n-1}}^2)^{n-1}
\sum\nolimits_{\ell=2}^{n}\vert z_\ell\vert^2\,.
\end{equation*}
Let $M_n=4h_{B_{2n-1}}^{3} \Hess_\mathbb C  h_{B_{2n-1}}$.  
Observe that $\lambda_1=2h^2_{B_{2n-1}}$ is an eigenvalue of $M_n$, in fact, for every $2\leq\ell\leq n$, the $\ell$-th row of $M_n- \lambda_1 I_n$ is equal to 
\begin{equation*}
-\bar z_\ell\left(\frac{\bar z_1- z_1}{2},z_2,\ldots,z_n\right)\,,
\end{equation*}
i.e. a multiple of $v=((\bar z_1- z_1)/2,z_2,\ldots,z_n)$. Since none of the last $n-1$ rows is identically zero on $B_{2n-1}\setminus\{0\}$ and as the first row of $M_n-\lambda_1I_n$ is not a multiple of $v$, it follows that $M_n- \lambda_1 I_n$ has rank $2$, i.e. $M_n$ admits a diagonal form with $\lambda_1$ occurring on the last $n-2$ diagonal elements. The remaining eigenvalues are 
\begin{equation*}
\lambda_2=h^2_{B_{2n-1}}-h_{B_{2n-1}}\left(\frac{z_1-\bar z_1}{2i}\right)
\qquad
{\rm and}
\qquad
\lambda_3=h_{B_{2n-1}}^2+h_{B_{2n-1}}\left(\frac{z_1-\bar z_1}{2i}\right)\,.
\end{equation*}
Indeed, for $j=2,3$, the first row of the matrix $M_n-\lambda_j I_n$ is a linear combination of the remaining ones, (which are all non zero). For $j=2,3$ and $2\leq \ell\leq n$, the coefficient of the the $\ell$-th row in such a linear combination is
\begin{equation*}
i z_\ell\left[h_{B_{2n-1}}+\left(\frac{z_1-\bar z_1}{2i}\right)\right]^{-1}
\qquad
{\rm or}
\qquad
-i z_\ell\left[h_{B_{2n-1}}-\left(\frac{z_1-\bar z_1}{2i}\right)\right]^{-1}\,,
\end{equation*}
respectively.
It follows that 
\begin{align*}
\det M_n
&=
(2h_{B_{2n-1}}^2)^{n-2}
h_{B_{2n-1}}^2
\left[
h_{B_{2n-1}}^2-\left(\frac{z_1-\bar z_1}{2i}\right)^2
\right]
\\
&=
2^{n-2}(h_{B_{2n-1}}^2)^{n-1}
\sum\nolimits_{\ell=2}^{n}\vert z_\ell\vert^2
\,
\end{align*}
as required. \findim

\subsection{Linear groups}

For every $n\in\mathbb N^*$, let ${\rm M}(n,\mathbb C)$\label{matricinpernc}, (resp. ${\rm M}(2n,\mathbb R)$)\label{matrici2nper2nr}, denote the set of $n$-by-$n$ matrices (resp. $2n$-by-$2n$ matrices) with entries in $\mathbb C$ (resp. $\mathbb R$). As usual
\begin{align*}
{\rm GL}(n,\mathbb C)
&=
\{M\in M(n,\mathbb C)\mid \det M\neq0\}\,,
\\
{\rm GL}(2n,\mathbb R)
&=
\{M\in M(2n,\mathbb R)\mid \det M\neq0\}\,,
\\
{\rm U}(n)
&=
\{M\in M(n,\mathbb C)\mid M\,{}^t\!\overline{M}=I_n\}\,, 
\\
{\rm O}(2n)
&=
\{M\in M(2n,\mathbb R)\mid M\,{}^t\!M=I_n\}\,,
\\
{\rm SO}(2n)
&=
\{M\in{\rm O}(2n)\mid \det M=1\}\,.
\end{align*}

\begin{lem}\label{unitario1}
The mapping $\psi_1:{\rm M}(1,\mathbb C)=\mathbb C\to {\rm M}(2,\mathbb R)$\label{psi1} defined, for $m\in \mathbb C$, by
\begin{align*}
\psi_1(m)=
\left(\begin{array}{cc}
\re m & -\im m \\
\im m & \re m
\end{array}\right)\,,
\end{align*}
has the following properties for every $m,m^\prime\in\mathbb C$:
\begin{enumerate}
\item $\psi_1(m+m^\prime)=\psi_1(m)+\psi_1(m^\prime)$, 
\item $\psi_1(mm^\prime)=\psi_1(m)\psi_1(m^\prime)$, 
\item $\det\psi_1(m)=\vert m\vert^2$,
\item ${}^t\psi_1(m)=\psi_1(\overline{m})$.
\end{enumerate}
In particular, $\psi_1({\rm U}(1))= {\rm SO}(2,\mathbb R)$.
\end{lem}
\noindent\pf The proof is straightforward. For the last statement, observe that, if $m\in{\rm U}(1)$, then $$\psi_1(m)\,{}^t\psi_1(m)=\psi_1(m)\psi_1(\overline{m})=\psi_1(m\overline{m})=\psi_1(1)=I_2\,$$ so that $\psi_1(m)\in  {\rm SO}(2,\mathbb R)$, i.e.  $\psi_1({\rm U}(1))\subseteq {\rm SO}(2,\mathbb R)$. Conversely, an element from ${\rm SO}(2,\mathbb R)$ can always be written as
$$
\left(\begin{array}{cc}
\cos\theta & -\sin\theta\\
\sin\theta & \cos\theta
\end{array}\right)\,,
$$
for some $\theta\in\mathbb R$, and such a matrix is nothing but $\psi_1(e^{i\theta})$.
\findim

\begin{lem}\label{unitarion}
If $n\in\mathbb N\setminus\{0,1\}$, the mapping $\psi_n:{\rm M}(n,\mathbb C)\to {\rm M}(2n,\mathbb R)$ defined, for~$M=(m_{\ell j})\in {\rm M}(n,\mathbb C)$, by\label{psin}
\begin{align*}
\psi_n(M)
=
\left(\begin{array}{ccc}
\psi_1(m_{11}) & \ldots  & \psi_1(m_{1 n}) \\
\vdots & \ddots & \vdots\\
\psi_1(m_{n1}) & \ldots & \psi_1(m_{nn})
\end{array}\right)\,,
\end{align*}
has the following properties for every $M,M^\prime\in{\rm M}(n,\mathbb C)$:
\begin{enumerate}
\item $\psi_n(M+M^\prime)=\psi_n(M)+\psi_n(M^\prime)$, 
\item $\psi_n(MM^\prime)=\psi_n(M)\psi_n(M^\prime)$, 
\item $\det\psi_n(M)=\vert \det M\vert^2$,
\item ${}^t\psi_n(M)=\psi_n({}^t\!\overline{M})$. In particular:
\begin{enumerate}
\item if $M$ is hermitian, then $\psi_n(M)$ is symmetric;
\item $\psi_n({\rm U}(n))\subset {\rm SO}(2n,\mathbb R)$, but the inclusion is strict.
\end{enumerate}
\end{enumerate}
\end{lem}
\noindent\pf The equality $\psi_n(MM^\prime)=\psi_n(M)\psi_n(M^\prime)$ is straightforward. If $\re M=(\re m_{\ell,j})$ and $\im M=(\im m_{\ell,j})$, then $M=\re M+i\im M$ and
\begin{align*}
\det\psi_n(M)
&=
(-1)^{\frac{n(n-1)}{2}}
\det
\left(\begin{array}{ccccc}
\re m_{11} & -\im m_{11} & \ldots  & \re m_{1 n} & -\im m_{1 n} \\
\vdots & \vdots & \ddots & \vdots &  \vdots\\
\re m_{n1} & -\im m_{n1} & \ldots  & \re m_{n n} & -\im m_{nn} \\
\im m_{11} & \re m_{11} & \ldots & \im m_{1n} & \re m_{1n} \\
\vdots & \vdots & \ddots & \vdots &  \vdots\\
\im m_{n1} & \re m_{n1} & \ldots & \im m_{nn} & \re m_{nn}
\end{array}\right)
\\
&=
\det
\left(\begin{array}{cc}
\re M & -\im M \\
\im M & \re M 
\end{array}\right)\,,
\end{align*}
whence
\begin{align*}
\det
\left(\begin{array}{cc}
\re M & -\im M \\
\im M & \re M 
\end{array}\right)
&=
\det
\left(\begin{array}{cc}
\re M +i\im M& -\im M +i\re M\\
\im M & \re M 
\end{array}\right)
\\
&=
\det
\left(\begin{array}{cc}
\re M +i\im M& 0\\
\im M & \re M -i\im M
\end{array}\right)
\end{align*}
so that
\begin{align*}
\det
\left(\begin{array}{cc}
M& 0\\
\im M& \overline{M}
\end{array}\right)
&=
(\det M)(\det \overline{M})
\\
&=
(\det M)(\overline{\det M})
\\
&=
\vert\det M\vert^2\,.
\end{align*}
Also
\begin{align*}
{}^t\psi_n(M)
&=
\left(\begin{array}{ccc}
{}^t\psi_1(m_{11}) & \ldots  & {}^t\psi_1(m_{n1}) \\
\vdots & \ddots & \vdots\\
{}^t\psi_1(m_{1n}) & \ldots & {}^t\psi_1(m_{nn})
\end{array}\right)
\\
&=
\left(\begin{array}{ccc}
\psi_1(\overline{m_{11}}) & \ldots  & \psi_1(\overline{m_{n1}}) \\
\vdots & \ddots & \vdots\\
\psi_1(\overline{m_{1n}}) & \ldots & \psi_1(\overline{m_{nn}})
\end{array}\right)
\\
&=
\psi_n({}^t\overline{M})\,.
\end{align*}
so that if $M={}^t\!\overline{M}$ then $\psi_n(M)=\psi({}^t\!\overline{M})={}^t\psi_n(M)$. Moreover, 
$M\in{\rm U}(n)$ implies $I_{2n}=\psi_n(I_n)=\psi_n(M\,{}^t\!\overline{M})=\psi_n(M)\psi_n({}^t\!\overline{M})=\psi_n(M)\,{}^t\psi_n(M)$, i.e. $\psi_n(M)\in {\rm O}(2n,\mathbb R)$. Since $M$ is unitary, $\det \psi_n(M)=\vert\det M\vert^2=1$, then $\psi_n(M)\in {\rm SO}(2n,\mathbb R)$. If $A\in {\rm M}(2,\mathbb R)$ denotes the matrix
\begin{align*}
A=\left(\begin{array}{cc}
0 & 1\\
1 & 0
\end{array}\right)\,,
\end{align*}
when $n=2k$, the matrix with $2k$ diagonal blocks equal to $A$ provides an element of ${\rm SO}(4k,\mathbb R)\setminus\psi_{2k}({\rm U}(2k))$. If $n=2k+1$, the matrix with $2k$ diagonal blocks equal to $A$ and a block equal to $I_2$, yields an element of ${\rm SO}(4k+2,\mathbb R)\setminus\psi_{2k+1}({\rm U}(2k+1))$.\findim

Remark that, for every $n\in\mathbb N^*$, $\psi_n(iI_n)$ is the matrix of the operator of $\mathbb R^{2n}$ defining the usual complex structure.

\cleardoublepage
\addcontentsline{toc}{section}{List of notation}

\subsection*{List of notation}
\begin{table}[H]
\begin{tabular}{lp{10cm}l}
Notation& Meaning &Page\\
\toprule
$\chi_A$& Indicator function of the set $A$, $\chi(x)=1$ if $x\in A$ and $\chi(x)=0$ if $x\notin A$&\pageref{chi}\\
$(\,,)$& Standard scalar product in $\mathbb R^n$&\pageref{scalarr}\\
$\aff_\mathbb R A$& Affine subspace spanned by the subset $A$ of $\mathbb R^n$ or $\mathbb C^n$ over $\mathbb R$&\pageref{affr}\\ 
$E_A$&Linear subspace parallel to $\aff_\mathbb R A$ &\pageref{giacituraA}\\
$\relint A$& Relative interior of $A$ &\pageref{relintA}\\
$\relbd A$& Relative boundary of $A$&\pageref{relbdA}\\
${\mathcal G}_k(\mathbb R^n)$& Grassman manifold of $k$-dimensional linear subspaces of $\mathbb R^n$&\pageref{grassk}\\
${\mathcal G}(\mathbb R^n)$& Topological sum of ${\mathcal G}_k(\mathbb R^n)$ &\pageref{wgrass}\\
$E^\perp$& Orthocomplement of the linear subspace $E\subset \mathbb R^n$ with respect to the standard scalar product $(\,,)$&\pageref{perp}\\
${\mathcal K}(\mathbb R^n)$& Convex subsets of $\mathbb R^n$&\pageref{konvexr},\pageref{konvexcr}\\
${\mathcal C}(\mathbb R^n)$& Non empty compact subsets of $\mathbb R^n$&\pageref{compactr}\\
${\mathcal P}(\mathbb R^n)$& Convex polytopes of $\mathbb R^n$&\pageref{konvexr},\pageref{konvexcr}\\
$d_A$& Real dimension of $E_A$&\pageref{rdimA}\\
$[v,w]$& Segment with $v$ and $w$ as endpoints&\pageref{segmento}\\
$h_A$& Support function of the convex subset $A$&\pageref{supportf}\\
$B_n$& Full-dimensional unit ball in $\mathbb R^n$ about the origin&\pageref{fullnball}\\
$H_A(v)$& Supporting hyperplane of the convex subset $A$ in the direction $v$&\pageref{supportH}\\
$\preccurlyeq\,,\prec\,,\succcurlyeq\,, \succ $& Order relations for faces &\pageref{facce}\\
${\mathcal B}(\Gamma)$&Boundary complex of the polytope $\Gamma$&\pageref{bc}\\
${\mathcal B}(\Gamma,k)$& Set of $k$-dimensional faces of the polytope $\Gamma$&\pageref{bck}\\
$f(\Gamma)$&Face vector of the polytope $\Gamma$&\pageref{facevector}\\
$u_{\Delta_{k-1},\Delta_k}$&Outer unit normal vector to the face $\Delta_{k-1}$&\pageref{ounv}\\
$K_{\Delta,\Gamma}$&Dual cone to the face $\Delta$ of the polytope $\Gamma$&\pageref{dualcone}\\
$\Sigma_{k,\Gamma}$&$k$-star of the polytope $\Gamma$&\pageref{kstar}\\
$\psi_\Gamma(\Delta)$&Outer angle of the face $\Delta\preccurlyeq\Gamma$ with respect to $\Gamma$&\pageref{outerangle}\\
$\vol_\ell$&$\ell$-dimensional Lebesgue measure in $\mathbb R^\ell$&\pageref{vollell}\\
$\varkappa_\ell$&$\ell$-dimensional Lebesgue measure of $B_\ell$&\pageref{varkappaell}\\
$\Sigma(\Gamma)$&Normal fan of the polytope $\Gamma$&\pageref{nfan}\\
$(A)_\varepsilon$& $\varepsilon$-neighborhood of $A$&\pageref{eintorno}\\
${\mathcal P}_\circ(\mathbb R)$& Neighborhoods of polytopes in $\mathbb R^n$&\pageref{pinfinitor},\pageref{konvexcr}\\
$\Delta_n$&Standard $n$-simplex&\pageref{nsimplesso}\\
$I_n$&Standard $n$-cube&\pageref{ncubo}\\
$\Theta_n$&Standard $n$-crosspolytope&\pageref{ntetra}\\
$q_\Delta$&Orthogonal projection from $\mathbb R^n$ to $E_\Delta^\perp$&\pageref{q-Delta}\\
$\Subd h_A$& Subdifferential of $h_A$&\pageref{subd}\\
${\mathcal S}^2(\mathbb R^n)$& piecewise $2$-regular subsets of $\mathbb R^n$&\pageref{skreg}\\
$\rho_A$&Defining function of $A$&\pageref{deffunct}\\
${\mathcal K}^2(\mathbb R^n)$& Convex piecewise $2$-regular subsets of $\mathbb R^n$&\pageref{kreg2}\\
${\mathcal K}_1(\mathbb R^n)$& Strictly convex subsets of $\mathbb R^n$&\pageref{strict}\\
${\mathcal K}_\ell(\mathbb R^n)$& Convex bodies with ${\mathcal C}^\ell(\mathbb R^n\setminus\{0\})$ support function&\pageref{ell-strict},\pageref{ell-strict},\pageref{konvexcr}\\
${\mathcal K}^2_\ell(\mathbb R^n)$& piecewise $2$-regular convex bodies with ${\mathcal C}^\ell(\mathbb R^n\setminus\{0\})$ support function&\pageref{2-ell-strict},\pageref{konvexcr}\\
${\mathfrak h}$&Haussdorf metric&\pageref{hausm}\\
$R_\varepsilon A$&$\varepsilon$-regularization of $A$&\pageref{ereg}\\
$\varphi_\varepsilon$&Cut-off function&\pageref{ecutoff}\\
$\upsilon_n$&Standard volume $n$-form on $\mathbb R^n$&\pageref{vformn}\\
$A[k]$&Sequence of $k$ subsets equal to $A$&\pageref{aktimes}\\
$\textsl{v}_k$&$k$-th intrinsic volume&\pageref{vintrk}\\
$V_n$&$n$-dimensional (Minkowski) mixed volume&\pageref{nvmix}\\
$\upsilon_{\partial B_n}$ & Standard volume form on $\partial B_n$&\pageref{vformpB}\\
$\Hess_\mathbb R h$ & Real Hessian matrix of $h$&\pageref{real-hess}\\
$D$ &sum of principal mixed minors corresponding to the nonzero entries of $I_n$ in the Laplace expansion of the mixed discriminant $D_n(I_n,\Hess_\mathbb R h_{A_2},\ldots,\Hess_\mathbb R h_{A_n})$ along $I_n$.&\pageref{md}, \pageref{mix-D-f2}\\
\end{tabular}
\end{table}%

\vfill\newpage

\begin{table}[ht]
\begin{tabular}{lp{10cm}l}
Notation& Meaning &Page\\
\toprule
$\textsl{v}_k^\varphi$&$k$-th intrisic $\varphi$-volume&\pageref{vintrfk}\\
$V_k^\varphi$&$k$-dimensional mixed $\varphi$-volume&\pageref{vintrfk}\\
$\varrho(E_1, E_2)$&Coefficient of volume distortion under the projection of $E_1$ on $E_2$&
\pageref{vdist}\\
$\dim_\mathbb R E$& Dimension of the $\mathbb R$-linear subspace $E\subseteq\mathbb C^n$&\pageref{rdim}\\
$\dim_\mathbb C E$& Complex dimension of the $\mathbb C$-linear subspace $\lin_\mathbb C E\subseteq\mathbb C^n$&\pageref{cdim}\\
$E^\mathbb C$& Maximal complex subspace of the $\mathbb R$-linear subspace $E\subseteq \mathbb C^n$&\pageref{cmax}\\
$E^\perp$& Orthocomplement of the linear subspace $E\subseteq\mathbb C^n$ with respect to the standard scalar product $\re\langle\,,\rangle$&\pageref{rperp}\\
$E^{\perp_\mathbb C}$& Orthocomplement of the linear subspace $E\subseteq\mathbb C^n$ with respect to the standard hermitian product $\langle\,,\rangle$&\pageref{rperp}\\
$E^\prime$& $E^\perp\cap \lin_\mathbb C E$ &\pageref{eprime}\\
${\mathcal G}(\mathbb C^n,\mathbb R)$& Set of $\mathbb R$-linear subspaces of $\mathbb C^n$&\pageref{wcgrass}\\
$\varrho(E)$& Coefficient of volume distortion under the projection of $E$ on $iE^\prime$& \pageref{roe}\\
$\equiv_4$& Equivalence modulo $4$&\pageref{pd}\\
$\upsilon_E$& Volume form on the $\mathbb R$-linear subspace $E\subset \mathbb C^n$&\pageref{vformE}\\
${\mathcal K}(\mathbb C^n)$& Convex subsets of $\mathbb C^n$&\pageref{konvexc}\\
${\mathcal P}(\mathbb C^n)$& Convex polytopes of $\mathbb C^n$&\pageref{konvexc}\\
${\mathcal S}^2(\mathbb C^n)$& Smooth compact subsets of $\mathbb C^n$&\pageref{konvexc}\\
${\mathcal K}^2_\ell(\mathbb C^n)$& $k$-regular convex bodies with ${\mathcal C}^\ell(\mathbb C^n\setminus\{0\})$ support function&\pageref{2-ell-strict},\pageref{konvexc}\\
$\varrho(A)$& Coefficient of volume distortion under the projection of $E_A$ on $iE^\prime_A$&\pageref{roA}\\
${\mathcal B}_{\rm ed}(\Gamma,k)$& Set of $k$-dimensional equidimensional faces of the polytope $\Gamma$&\pageref{bcked}\\
$\upsilon_{2n}$ &Standard volume $2n$-form on $\mathbb C^n$&\pageref{vform2n}\\
$\dc$& $i(\bar\partial-\partial)$& \pageref{dc}\\
$\ddc$ &The Bott-Chern operator $2i\partial\bar\partial$ &\pageref{ddc}\\
$*$& Hodge $*$-operator&\pageref{*-op}\\
$\upsilon_E$ & Volume form on the oriented real linear subspace $E\subset\mathbb C^n$ &\pageref{vformE}\\
$\nu_{\partial A}$& Gauss map of $\partial A$&\pageref{Gaussmap}\\
$\upsilon_{\partial A}$ & Volume form on $\partial A$ & \pageref{vformpA}\\
$\alpha_{\partial A}$ &  Differential $1$-form defined on $\partial A$ &\pageref{alphaf}\\
$P_n$& $n$-dimensional Kazarnovski{\v \i} pseudovolume&\pageref{nkaza}\\
$\sec_{\partial A}$& Second quadratic form of $\partial A$ & \pageref{second}\\
$\sec_{\partial A}^{\mathbb C}$ & $\mathbb R$-bilinear symmetric form associated to $\sec_{\partial A}$ & \pageref{secondc}\\
${\mathcal L}_{\partial A}$ & Levi form of $\partial A$ & \pageref{levi}\\
${\mathcal E}_{\partial A}$ & Product of the eigenvalues of ${\mathcal L}_{\partial A}$& \pageref{kappalevi}\\
${\mathfrak S}_n$ & Symmetric group over the set $\{1,\ldots,n\}$ &\pageref{gruppo-simm}\\
$\text{mix}^\sigma(M_1,\ldots,M_n)$& mixed version of matrices&\pageref{mix}\\
$\text{mix}_\sigma(M_1,\ldots,M_n)$& mixed version of matrices&\pageref{mix}\\
$\Hess_\mathbb C h$ & Complex hessian matrix of $h$ &\pageref{c-hess}\\
$D_n$ & Mixed discriminant&\pageref{dm}\\
$M^{[j,k]}$ & sub-matrix of $M$ & \pageref{sottomatrice}\\
$T$ & Current &\pageref{current}\\
$\varsigma$ & Test form &\pageref{tform}\\
$\langle\!\langle T,\varsigma\rangle\!\rangle$& Pairing between the current $T$ and the testform $\varsigma$ &\pageref{pairing}\\
$\lambda_\Delta$ & $(k,k)$-current associated to the equidimensional $k$-face $\Delta$ of a polytope $\Gamma$ & \pageref{lambdak}\\
$\iota_\Delta$ & Inclusion mapping $K_\Delta\to\mathbb C^n$  & \pageref{lambdak}\\
$\Xi(A_2,\ldots,A_n)(u)$& $n$-by-$n$ matrix whose entries are (sums) of mixed discriminants of complex hessian matrices&\pageref{matriceXi}\\
${\mathcal X}(\mathbb C^n)$& Largest subspace of ${\mathcal K}_2(\mathbb C^n)^n$ on which, for almost all $u\in\partial B_{2n}$, the equality $\re\langle\nabla h_{A_1}(u),\Xi(A_2,\ldots,A_n)(u)\,u\rangle=0$ holds true&\pageref{X}\\
$\text{M}(n,\mathbb C)$& Set of complex $n$-by-$n$ matrices&\pageref{matricinpernc}\\
$\text{M}(2n,\mathbb R)$&Set of real $2n$-by-$2n$ matrices&\pageref{matrici2nper2nr}\\
$\psi_1$&Mapping from $\text{M}(1,\mathbb C)$ to $\text{M}(2,\mathbb R)$&\pageref{psi1}\\
$\psi_n$&Mapping from $\text{M}(n,\mathbb C)$ to $\text{M}(2n,\mathbb R)$&\pageref{psin}\\
\end{tabular}
\end{table}%
\vfill
\newpage

\cleardoublepage

\end{document}